\documentclass[letterpaper,11pt,oneside,reqno]{amsbook}
\usepackage{bm}
\usepackage{mathrsfs}
\usepackage{amsfonts,amsmath, amssymb,amsthm,amscd,stmaryrd}
\usepackage[height=9.6in,width=5.95in]{geometry}
\usepackage[mpexclude,DIV13]{typearea}
\usepackage{verbatim}
\usepackage{hyperref}
\usepackage{graphicx}
\usepackage[latin1]{inputenc}
\usepackage{latexsym}
\usepackage{lscape}
\usepackage{epsfig}
\include{bibtex}

\usepackage{setspace}

\usepackage{parskip}

\newcommand{\ovbar}[1]{\mkern 1.5mu\overline{\mkern-1.5mu#1\mkern-1.5mu}\mkern 1.5mu}

\begin{document}
\bibliographystyle{alpha}
\newcommand{\cn}[1]{\overline{#1}}
\newcommand{\e}[0]{\varepsilon}
\newcommand{\phimac}{\varphi}
\newcommand{\bbf}[0]{\mathbf}

\newcommand{\Pfree}[5]{\ensuremath{\mathbb{P}^{#1,#2,#3,#4,#5}}}
\newcommand{\PfreeShort}{\ensuremath{\mathbb{P}^{BB}}}

\newcommand{\WH}[8]{\ensuremath{\mathbb{W}^{#1,#2,#3,#4,#5,#6,#7}_{#8}}}
\newcommand{\Wfree}[5]{\ensuremath{\mathbb{W}^{#1,#2,#3,#4,#5}}}
\newcommand{\WHShort}[3]{\ensuremath{\mathbb{W}^{#1,#2}_{#3}}}
\newcommand{\WHShortCouple}[2]{\ensuremath{\mathbb{W}^{#1}_{#2}}}

\newcommand{\walk}[3]{\ensuremath{X^{#1,#2}_{#3}}}
\newcommand{\walkupdated}[3]{\ensuremath{\tilde{X}^{#1,#2}_{#3}}}
\newcommand{\walkfull}[2]{\ensuremath{X^{#1,#2}}}
\newcommand{\walkfullupdated}[2]{\ensuremath{\tilde{X}^{#1,#2}}}

\newcommand{\PH}[8]{\ensuremath{\mathbb{Q}^{#1,#2,#3,#4,#5,#6,#7}_{#8}}}
\newcommand{\PHShort}[1]{\ensuremath{\mathbb{Q}_{#1}}}
\newcommand{\PHExp}[8]{\ensuremath{\mathbb{F}^{#1,#2,#3,#4,#5,#6,#7}_{#8}}}

\newcommand{\D}[8]{\ensuremath{D^{#1,#2,#3,#4,#5,#6,#7}_{#8}}}
\newcommand{\DShort}[1]{\ensuremath{D_{#1}}}
\newcommand{\partfunc}[8]{\ensuremath{Z^{#1,#2,#3,#4,#5,#6,#7}_{#8}}}
\newcommand{\partfuncShort}[1]{\ensuremath{Z_{#1}}}
\newcommand{\bolt}[8]{\ensuremath{W^{#1,#2,#3,#4,#5,#6,#7}_{#8}}}
\newcommand{\boltShort}[1]{\ensuremath{W_{#1}}}
\newcommand{\boltNew}{\ensuremath{W}}
\newcommand{\QTLH}{\ensuremath{\mathfrak{H}}}
\newcommand{\QTLHgen}{\ensuremath{\mathfrak{L}}}

\newcommand{\whitenoise}{\ensuremath{\mathscr{\dot{W}}}}
\newcommand{\mf}{\mathfrak}

\newcommand{\EE}{\ensuremath{\mathbb{E}}}
\newcommand{\PP}{\ensuremath{\mathbb{P}}}
\newcommand{\var}{\textrm{var}}
\newcommand{\N}{\ensuremath{\mathbb{N}}}
\newcommand{\R}{\ensuremath{\mathbb{R}}}
\newcommand{\C}{\ensuremath{\mathbb{C}}}
\newcommand{\Z}{\ensuremath{\mathbb{Z}}}
\newcommand{\Q}{\ensuremath{\mathbb{Q}}}
\newcommand{\T}{\ensuremath{\mathbb{T}}}
\newcommand{\E}[0]{\mathbb{E}}
\newcommand{\OO}[0]{\Omega}
\newcommand{\F}[0]{\mathfrak{F}}
\def \Ai {{\rm Ai}}
\newcommand{\G}[0]{\mathfrak{G}}
\newcommand{\ta}[0]{\theta}
\newcommand{\w}[0]{\omega}
\newcommand{\ra}[0]{\rightarrow}
\newcommand{\vectoro}{\overline}
\newcommand{\crairy}{\mathcal{CA}}
\newcommand{\nc}{\mathsf{NoTouch}}
\newcommand{\ncf}{\mathsf{NoTouch}^f}
\newcommand{\wxy}{\mathcal{W}_{k;\bar{x},\bar{y}}}
\newcommand{\AP}{\mathfrak{a}}
\newcommand{\cm}{\mathfrak{c}}
\newtheorem{theorem}{Theorem}[chapter]
\numberwithin{figure}{chapter}
\newtheorem{partialtheorem}{Partial Theorem}[section]
\newtheorem{conj}[theorem]{Conjecture}
\newtheorem{lemma}[theorem]{Lemma}
\newtheorem{proposition}[theorem]{Proposition}
\newtheorem{corollary}[theorem]{Corollary}
\newtheorem{claim}[theorem]{Claim}
\newtheorem{experiment}[theorem]{Experimental Result}

\renewcommand{\thesection}{\thechapter.\arabic{section}}

\def\todo#1{\marginpar{\raggedright\footnotesize #1}}
\def\change#1{{\color{green}\todo{change}#1}}
\def\note#1{\textup{\textsf{\color{blue}(#1)}}}

\theoremstyle{definition}
\newtheorem{rem}[theorem]{Remark}

\theoremstyle{definition}
\newtheorem{com}[theorem]{Comment}

\theoremstyle{definition}
\newtheorem{definition}[theorem]{Definition}

\theoremstyle{definition}
\newtheorem{definitions}[theorem]{Definitions}

\theoremstyle{definition}
\newtheorem{conjecture}[theorem]{Conjecture}

\newcommand{\airysh}{\mathcal{A}}
\newcommand{\hfixed}{\mathcal{H}}
\newcommand{\afixed}{\mathcal{A}}
\newcommand{\canopynoarg}{\mathsf{C}}
\newcommand{\canopy}[3]{\ensuremath{\mathsf{C}_{#1,#2}^{#3}}}
\newcommand{\argmax}{x_{{\rm max}}}
\newcommand{\zmax}{z_{{\rm max}}}

\newcommand{\Rkle}{\ensuremath{\mathbb{R}^k_{>}}}
\newcommand{\Ronele}{\ensuremath{\mathbb{R}^k_{>}}}
\newcommand{\ewxy}{\mathcal{E}_{k;\bar{x},\bar{y}}}

\newcommand{\bxyf}{\mathcal{B}_{\bar{x},\bar{y},f}}
\newcommand{\bxyflr}{\mathcal{B}_{\bar{x},\bar{y},f}^{\ell,r}}

\newcommand{\bxyfone}{\mathcal{B}_{x_1,y_1,f}}

\newcommand{\ptac}{p}
\newcommand{\ptact}{v}
\newcommand{\nF}{H}

\newcommand{\fext}{\mathfrak{F}_{{\rm ext}}}
\newcommand{\gext}{\mathfrak{G}_{{\rm ext}}}
\newcommand{\xext}{{\rm xExt}(\mathfrak{c}_+)}

\newcommand{\bmotion}{X}

\newcommand{\dd}{\, {\rm d}}
\newcommand{\signc}{\Sigma}
\newcommand{\wxylr}{\mathcal{W}_{k;\bar{x},\bar{y}}^{\ell,r}}
\newcommand{\wxylrprime}{\mathcal{W}_{k;\bar{x}',\bar{y}'}^{\ell,r}}
\newcommand{\Rklezero}{\ensuremath{\mathbb{R}^k_{>0}}}
\newcommand{\XYfM}{\textrm{XY}^{f}_M}

\newcommand{\upright}{D}
\newcommand{\energy}{E}
\newcommand{\xmax}{{\rm max}_1}
\newcommand{\ymax}{{\rm max}_2}
\newcommand{\lppls}{\mathcal{L}}
\newcommand{\lpplsre}{\mathcal{L}^{{\rm re}}}
\newcommand{\lpplsarg}[1]{\mathcal{L}_{n}^{\fa \to #1}}
\newcommand{\larg}[3]{\mathcal{L}_{n}^{#1,#2;#3}}
\newcommand{\BP}{M}
\newcommand{\weight}{\mathsf{Wgt}}
\newcommand{\pairweight}{\mathsf{PairWgt}}
\newcommand{\sumweight}{\mathsf{SumWgt}}
\newcommand{\mpgood}{\mathcal{G}}
\newcommand{\mpg}{\mathsf{Fav}}
\newcommand{\mcgone}{\mathsf{Fav}_1}
\newcommand{\radnik}[2]{\mathsf{RN}_{#1,#2}}
\newcommand{\size}[2]{\mathsf{S}_{#1,#2}}
\newcommand{\pdr}{\mathsf{PolyDevReg}}
\newcommand{\pwr}{\mathsf{PolyWgtReg}}
\newcommand{\lwr}{\mathsf{LocWgtReg}}
\newcommand{\hwp}{\mathsf{HighWgtPoly}}
\newcommand{\maxswf}{\mathsf{MaxScSumWgtFl}}
\newcommand{\emaxswf}{\e \! - \! \maxswf}
\newcommand{\minswf}{\mathsf{MinScSumWgtFl}}
\newcommand{\eminswf}{\e \! - \! \minswf}
\newcommand{\surreg}{\mathcal{R}}
\newcommand{\scf}{\mathsf{FavSurgCond}}
\newcommand{\disjtpoly}{\mathsf{DisjtPoly}}
\newcommand{\intint}[1]{\llbracket 1,#1 \rrbracket}
\newcommand{\maxsym}{*}
\newcommand{\polynum}{\#\mathsf{Poly}}
\newcommand{\dlp}{\mathsf{DisjtLinePoly}}
\newcommand{\lowb}{\underline{B}}
\newcommand{\highb}{\overline{B}}
\newcommand{\tottt}{t_{1,2}^{2/3}}
\newcommand{\tot}{t_{1,2}}

\newcommand{\mc}{\mathcal}
\newcommand{\vect}{\mathbf}
\newcommand{\bt}{\mathbf{t}}
\newcommand{\scB}{\mathscr{B}}
\newcommand{\scBres}{\mathscr{B}^{\mathrm{re}}}
\newcommand{\rightshadow}{\mathrm{RS}Z}
\newcommand{\dbm}{L}
\newcommand{\dysonbm}{DBM}
\newcommand{\edgedbm}{\mc{L}^{\scal}}
\newcommand{\edgedysonbm}{D^{\rm edge}}
\newcommand{\gue}{\mathrm{GUE}}
\newcommand{\edgegue}{\mathrm{GUE}^{\mathrm{edge}}}
\newcommand{\eqdist}{\stackrel{(d)}{=}}
\newcommand{\geqdist}{\stackrel{(d)}{\succeq}}
\newcommand{\leqdist}{\stackrel{(d)}{\preceq}}
\newcommand{\scal}{{\rm sc}}
\newcommand{\fa}{x_0}
\newcommand{\hit}{H}
\newcommand{\scaledle}{\mathsf{Sc}\mc{L}}
\newcommand{\cenleup}{\mathscr{L}^{\uparrow}}
\newcommand{\cenledown}{\mathscr{L}^{\downarrow}}
\newcommand{\eln}{T}
\newcommand{\xmin}{{\rm Corner}^{\mfl,\mc{F}}}
\newcommand{\ymin}{{\rm Corner}^{\mfr,\mc{F}}}
\newcommand{\gxmin}{{\rm Corner}^{\gfl,\mc{F}}}
\newcommand{\gymin}{{\rm Corner}^{\gfr,\mc{F}}}
\newcommand{\barxmin}{{\rm \overline{C}}{\rm orner}^{\mfl,\mc{F}}}
\newcommand{\barymin}{{\rm \overline{C}}{\rm orner}^{\mfr,\mc{F}}}
\newcommand{\gbarxmin}{{\rm \overline{C}}{\rm orner}^{\gfl,\mc{F}}}
\newcommand{\gbarymin}{{\rm \overline{C}}{\rm orner}^{\gfr,\mc{F}}}
\newcommand{\qmin}{{\rm Corner}^{\mfl,\mc{F}^1}}
\newcommand{\barqmin}{{\rm \overline{C}}{\rm orner}^{\mfl,\mc{F}^1}}
\newcommand{\test}{T}
\newcommand{\mfl}{\mf{l}}
\newcommand{\mfr}{\mf{r}}
\newcommand{\gfl}{\ell}
\newcommand{\gfr}{r}
\newcommand{\jre}{J}
\newcommand{\highfl}{{\rm HFL}}
\newcommand{\flyleap}{\mathsf{FlyLeap}}
\newcommand{\touch}{\mathsf{Touch}}
\newcommand{\notouch}{\mathsf{NoTouch}}
\newcommand{\close}{\mathsf{Close}}
\newcommand{\boxclose}{\mathsf{BoxClose}}
\newcommand{\abovepar}{\mathsf{High}}
\newcommand{\vecint}{\bar{\iota}}
\newcommand{\cornthree}{{\rm Corner}^\mc{G}_{k,\mfl}}
\newcommand{\cornfour}{{\rm Corner}^\mc{H}_{k,\fa}}
\newcommand{\mpgg}{\mathsf{Fav}_{\mc{G}}}

\newcommand{\lefta}{{\rm SideLeft}}
\newcommand{\righta}{{\rm SideRight}}
\newcommand{\mida}{{\rm Middle}}
\newcommand{\xnmac}{z_n}
\newcommand{\cdor}{{\rm C}}
\newcommand{\alphapr}{\alpha'}
\newcommand{\betapr}{\beta'}
\newcommand{\barqpr}{\bar{q}'}
\newcommand{\pspr}{\pairsep_{\alphapr,\betapr,\barqpr}}

\newcommand{\wien}{W}
\newcommand{\pole}{P}
\newcommand{\pp}{p}

\newcommand{\maxpoly}{\mathrm{MaxDisjtPoly}}

\newcommand{\low}{\mathsf{Low}}
\newcommand{\nolow}{\mathsf{NoLow}}
\newcommand{\up}{\mathsf{Up}}
\newcommand{\neargeod}{\mathsf{NearGeod}}

\newcommand{\lshift}{\mc{L}^{\rm shift}}

\newcommand{\const}{D_k}
\newcommand{\numcone}{14}
\newcommand{\numctwo}{13}
\newcommand{\numcthree}{4}
\newcommand{\cone}{c_1}
\newcommand{\Cone}{C_1}
\newcommand{\rsC}{C}
\newcommand{\rsc}{c}
\newcommand{\ctemp}{d_0}
\newcommand{\smallc}{c_0}
\newcommand{\smallcprime}{c_1}
\newcommand{\smallcanother}{c_2}
\newcommand{\smallcnew}{c_3}
\newcommand{\rcon}{r_0}
\newcommand{\Cstrong}{E}
\newcommand{\formerE}{C}
\newcommand{\Chat}{\hat{C}}
\newcommand{\Cwb}{K}
\newcommand{\constnew}{\gamma_k}

\newcommand{\dist}{\vert\vert}
\newcommand{\fik}{\mc{F}_i^{[0,\ipdval]^c}}
\newcommand{\mcfa}{\mc{H}[\fa]}
\newcommand{\tent}{{\rm Tent}}
\newcommand{\goodk}{\mathsf{G}_{0,\ipdval}}
\newcommand{\pairsep}{{\rm PS}}
\newcommand{\mbf}{\mathsf{MBF}}
\newcommand{\nbd}{\mathsf{NoBigDrop}}
\newcommand{\bd}{\mathsf{BigDrop}}
\newcommand{\jleft}{j_{{\rm left}}}
\newcommand{\jright}{j_{{\rm right}}}
\newcommand{\smalljfluc}{\mathsf{SmallJFluc}}
\newcommand{\mfone}{M_{\mc{F}^1}}
\newcommand{\mfthree}{M_{\mc{G}}}
\newcommand{\fev}{\mathsf{F}}

\newcommand{\para}{Q}
\newcommand{\ipd}{d_{{ip}}}
\newcommand{\ipdval}{d}
\newcommand{\deltamac}{\zeta}
\newcommand{\deltapi}{\theta}
\newcommand{\emac}{\eta}
\newcommand{\grabell}{\rho}

\newcommand{\staircase}{SC}
\newcommand{\coninit}{\Psi}
\newcommand{\initcond}{\mathcal{I}}
\newcommand{\maxmin}{\pwr}

\newcommand{\nmac}{n}

\newcommand{\rmreg}{{\rm Reg}}

\newcommand{\boundgood}{\mathsf{G}}

\newcommand{\strongjump}{J_+}
\newcommand{\overstrongjump}{{\overline{J}}_+}

\makeatletter
\newsavebox\myboxA
\newsavebox\myboxB
\newlength\mylenA

\newcommand*\xoverline[2][0.75]{%
    \sbox{\myboxA}{$\m@th#2$}%
    \setbox\myboxB\null
    \ht\myboxB=\ht\myboxA%
    \dp\myboxB=\dp\myboxA%
    \wd\myboxB=#1\wd\myboxA
    \sbox\myboxB{$\m@th\overline{\copy\myboxB}$}
    \setlength\mylenA{\the\wd\myboxA}
    \addtolength\mylenA{-\the\wd\myboxB}%
    \ifdim\wd\myboxB<\wd\myboxA%
       \rlap{\hskip 0.5\mylenA\usebox\myboxB}{\usebox\myboxA}%
    \else
        \hskip -0.5\mylenA\rlap{\usebox\myboxA}{\hskip 0.5\mylenA\usebox\myboxB}%
    \fi}
\makeatother



\title[Brownian bridge regularity for the  Airy line ensemble]{Brownian regularity for the  Airy line ensemble, \\ and multi-polymer watermelons \\ in Brownian last passage percolation}



\author[A. Hammond]{Alan Hammond}
\address{A. Hammond\\
  Department of Mathematics and Statistics\\
 U.C. Berkeley \\
  899 Evans Hall \\
  Berkeley, CA, 94720-3840 \\
  U.S.A.}
  \email{alanmh@berkeley.edu}
  \thanks{The author is supported by NSF grant DMS-$1512908$. }

\begin{abstract} 
The Airy line ensemble is a positive-integer indexed system of random continuous curves whose finite dimensional distributions are given by the multi-line Airy process. It is a natural object in the KPZ universality class: for example, its highest curve, the Airy$_2$ process, describes after the subtraction of a parabola the limiting law of the scaled energy of a geodesic running from the origin to a variable point on an anti-diagonal line in such problems as Poissonian last passage percolation. The ensemble of curves resulting from the Airy line ensemble after the subtraction of the same parabola enjoys a simple and explicit spatial Markov property, the  {\em Brownian Gibbs} property. 

In this paper, we employ the Brownian Gibbs property to make a close comparison between the Airy line ensemble's curves after affine shift and Brownian bridge, proving the finiteness of a superpolynomially growing moment bound on Radon-Nikodym derivatives. 

We also determine the value of a natural exponent describing in Brownian last passage percolation the decay in probability for the existence of several near geodesics that are disjoint except for their common endpoints, where the notion of `near' refers to a small deficit in scaled geodesic energy, with the 
parameter specifying this nearness tending to zero. 

To prove both results, we introduce a technique that may be useful elsewhere for finding upper bounds on probabilities of events concerning random systems of curves enjoying the Brownian Gibbs property. 

Several results in this article play a fundamental role in a further study of Brownian last passage percolation in three companion papers, ~\cite{ModCon},~\cite{NonIntPoly} and~\cite{Patch},  in which geodesic coalescence and geodesic energy profiles are investigated in scaled coordinates. 
\end{abstract}

\maketitle

\newpage

\tableofcontents

\chapter{Introduction}

\section{Kardar-Parisi-Zhang universality}

The topic of Kardar-Parisi-Zhang [KPZ] universality is concerned with random growth processes, in which a one-dimensional interface evolves randomly in time, and universal structures that describe,  independently of the microscopic details that specify the random growth,   the large-scale behaviour of interfaces that emerges at advanced time. The class of processes that are expected to lie in the KPZ universality class -- that is, for which these universal structures offer an accurate limiting description -- is very broad. The basic characteristics that are expected to place a random growth model in the KPZ class are that 
growth occurs in a direction normal to the present surface, alongside two forces: a
surface tension whose effect is smoothening, and a
local random force whose effect is to roughen the surface.

Many last passage percolation models are expected to lie in the KPZ class. An LPP model is a totally asymmetric random growth model that comes equipped with a planar random environment, which is independent between disjoint regions. Directed planar paths, which are restricted say to move only in directions in the first quadrant, are each assigned an energy via the environment:
a path's energy is equal to (or perhaps some variant of) the integral of the environment along the path. For a given pair of planar points, the second displaced from the first upwards and to the right,
the path that attains the maximum energy among those directed paths with this pair of endpoints is called a geodesic. 
 If the first endpoint is held fixed, and the second is varied, say horizontally, then the geodesic's energy is a random function that may play the role of the random interface that we have mentioned. 
By physicists, this energy profile is called the narrow wedge, a term that applies when the first endpoint is held fixed. 

Scaling exponents play a basic role in KPZ universality. 
In  pursuing a discussion of LPP models,   it is useful  to describe these exponents. 
For such a model in the KPZ class, the random interface (for the narrow wedge case) is a function of the variable horizontal coordinate of the non-fixed, higher, endpoint.  
The vertical displacement between the two endpoints, which is being held fixed, is a parameter that plays the role of microscopic time in this random interface model. We label this parameter~$n > 0$.
The dominant mechanism of growth is linear, the 
value of the interface evaluated on
the diagonal $y=x$
growing at order $n$. The typical deviation of this random value from its mean scales as~$n^{1/3}$, while the horizontal displacement required in order to see a non-trivial correlation 
of the interface's value with its value on the diagonal is of the order of $n^{2/3}$. 
 These assertions have been rigorously demonstrated for only a few LPP models, each of which enjoys an integrable structure: the seminal work of Baik, Deift and Johansson~\cite{BDJ1999} 
 rigorously established the one-third exponent, and moreover obtained the GUE Tracy-Widom distributional limit, for the case of Poissonian last passage percolation, while the two-thirds power law for transversal fluctuation was derived for this model by Johansson~\cite{Johansson2000}.

Given these scaling exponents of one-third and two-thirds, the aim of investigating KPZ universality is advanced by introducing scaled coordinates. For example, the geodesic from the origin to $(n,n)$
may be expected to fluctuate from the diagonal by an order of $n^{2/3}$, while its energy fluctuates from its mean value by an order of $n^{1/3}$. The use of scaled coordinates in space involves affinely bringing $(n,n)$ to the location $(0,1)$ and then further squeezing the $x$-coordinate by a factor of $n^{2/3}$. Viewed in the new coordinates, the geodesic travels a unit distance upwards from the origin, and fluctuates horizontally on a unit-order, independently of large choices of $n$. 

\begin{figure}[ht]
\begin{center}
\includegraphics[height=5cm]{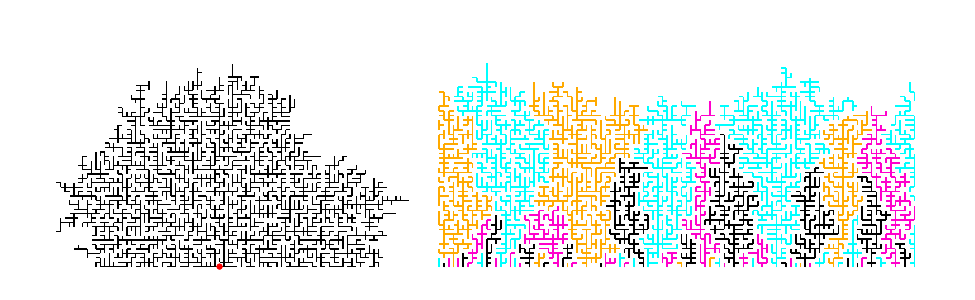}
\caption{
Consider a discrete model of random growth in which the initially healthy integer lattice sites in the upper half-plane become infected. At time zero, a certain subset of such sites on the $x$-axis are infected. At any given positive integer time, one uninfected site in the upper half-plane that is a nearest neighbour of a presently infected site is selected, uniformly at random. The site then becomes infected, with a `transmission' edge being added from the newly infected site into the already infected set, this edge selected uniformly from the available possibilities.
In this way, the infected region grows, one site at a time. 
 At any given moment, this region  is the collection of vertices abutting the present set of transmission edges. The transmission edge-set is partitioned into a collection of trees, each rooted at one of the initially infected sites.
This infection model, a variant of the Eden model introduced in~\cite{Eden}, is expected to lie in the KPZ class.
In the left sketch, it begins from a single infected site, the narrow wedge case. In the right, each site on the $x$-axis is initially infected, so that growth begins from a flat initial condition. A forest of competing trees forms, with the interface at a given time partitioned into the canopies of the trees surviving at that time.
By working with Brownian LPP, a model that, unlike the infection model, enjoys known integrable characteristics,
our study will demonstrate that the competition between trees respects KPZ scaling in such cases as the right sketch. For example, when the interface has height $n$, the canopy of a typical surviving tree has width $n^{2/3}$. Interface behaviour from a general initial condition will be understood by noting the similarities between the interface on the canopy of any given surviving tree and this interface for the single tree in the narrow wedge case seen in the left sketch.}
\label{f.eden}
\end{center}
\end{figure}

Brownian last passage percolation is an LPP model which has attractive integrable and probabilistic properties.
The present article forms with the companion papers~\cite{ModCon},~\cite{NonIntPoly} and~\cite{Patch} an inquiry into 
Brownian LPP viewed in these scaled coordinates. The scaled geodesic energy profile, begun from the narrow wedge, converges to one of the 
the canonical objects in the KPZ universality class, the Airy$_2$ process. A refined understanding of the locally Brownian nature of this process is achieved in the present article.
The scaled energy profile may also be studied when it is initiated by a far more general initial condition, and the articles prove a rather strong (and uniform) assertion of locally Brownian structure for
this class of Airy-like processes. Deeply implicated in the study is the field of scaled geodesics,  and uniqueness and coalescence properties for them, which lead to an understanding that this scaled geodesic field resembles a forest of coalescing geodesic branches whose geometry respects in a powerful way the scaled coordinates used to depict it. Figure~\ref{f.eden} illustrates some of these themes.

 A fundamental aspect of the approach to KPZ universality pursued by our four-paper study is to form an alliance between integrable results and probabilistic techniques. The geodesic energy profile in an integrable LPP model, such as Brownian LPP, may be embedded, via the Robinson-Schensted-Knuth correspondence, as the uppermost curve in an ordered system, or ensemble, of random curves, that collectively may be viewed as a collection of mutually avoiding random walks (or Brownian motions, in the Brownian LPP case). 
This observation, which is integrable in origin, becomes a critical probabilistic tool in our study: such ensembles may be resampled, and thus analysed, using a Brownian Gibbs resampling property.
Broadly, the technical contribution of this article is to analyse, via the Brownian Gibbs property, natural ensembles arising in Brownian LPP; doing so leads to 
 significant new inferences regarding, for example, the Airy$_2$ process. These Brownian LPP deductions will be exploited in the three companion papers where the study of this model in scaled coordinates is pursued.

In the next section of this introduction, we will offer a conceptual overview of the study undertaken in this work and the three companion papers. We will make informal statements of the principal conclusions and indicate crudely the ideas that drive the proofs of the results. The four papers have been written so that they may be read independently of one another, but the reader of any given one of them who wants to understand how the pieces fit together is encouraged to consult the next section. In the third and fourth sections of the introduction, we will then introduce and state some of the principal conclusions in the present article.

\subsection{Acknowledgments} 
I am very grateful to Riddhipratim Basu and Jeremy Quastel for extensive and valuable discussions regarding this paper and related ideas. I would like to thank Ivan Corwin, Shirshendu Ganguly, Milind Hegde and Jim Pitman for helpful 
discussions; two referees for thorough and perceptive comments; and Xuan Wu for her close reading and valuable corrections. 
I thank Judit Z{\'a}dor for the simulations shown in Figure~\ref{f.eden}.

\section{A conceptual overview of the scaled Brownian last passage percolation study}\label{s.concept}

Our discussion is of course heuristic, but, to give it some precision, we begin by offering the definition of Brownian last passage percolation; a geometric view of the model; and a specification of the scaled coordinate description of geodesics and their energies.

\subsection{Brownian last passage percolation}\label{s.brlpp}

\subsubsection{The model's definition.}
On a probability space carrying a law labelled~$\PP$, we let $B:\Z \times \R \to \R$ denote an ensemble of independent  two-sided standard Brownian motions $B(k,\cdot):\R\to \R$, $k \in \Z$.

Let $i,j \in \Z$ with $i \leq j$. 
We denote the integer interval $\{ i,\cdots  , j \}$ by $\llbracket i,j \rrbracket$.
Further let $x,y \in \R$ with $x \leq y$.
With these parameters given, we consider the collection of  non-decreasing lists 
 $\big\{ z_k: k \in \llbracket i+1,j \rrbracket \big\}$ of values $z_k \in [x,y]$. 
With the convention that $z_i = x$ and $z_{j+1} = y$,
we associate an energy to any such list, namely   $\sum_{k=i}^j \big( B ( k,z_{k+1} ) - B( k,z_k ) \big)$.
We may then define  the maximum energy, 
$M^1_{(x,i) \to (y,j)}$,  to be the supremum of the energies of all such lists. 
(The superscript `$1$' anticipates a generalization that we will later specify  under which the maximum over several lists is instead considered.)

 The process $M^1_{(0,1) \to (\cdot,n)}: [0,\infty) \to \R$ may have been considered first by~\cite{GlynnWhitt}; it was studied  further in~\cite{O'ConnellYor}.

\subsubsection{Staircases.}
Each such list is in bijection with a certain subset of $[x,y] \times [i,j] \subset \R^2$ that we call a {\em staircase}.
Staircases offer a geometric perspective on Brownian LPP and perhaps help in visualizing the problems in question.


 The staircase associated to the non-decreasing list $\big\{ z_k: k \in \llbracket i+1,j \rrbracket \big\}$ is specified as the union of certain horizontal planar line segments, and certain vertical ones.
The horizontal segments take the form $[ z_k,z_{k+1} ] \times \{ k \}$ for $k \in \llbracket i , j \rrbracket$.
Here, the convention that $z_i = x$ and  $z_{j+1} = y$ is again adopted. 
The right and left endpoints of each consecutive pair of horizontal segments are interpolated by a vertical planar line segment of unit length. It is this collection of vertical line segments that form
the vertical segments of the staircase.

The resulting staircase may be depicted as the range of an alternately rightward and upward moving path from starting point $(x,i)$ to ending point $(y,j)$. 
The set of staircases with these starting and ending points will be denoted by $\staircase_{(x,i) \to (y,j)}$.
Any staircase $\phi \in \staircase_{(x,i) \to (y,j)}$
is assigned an energy $E(\phi) = \sum_{k=i}^j \big( B ( k,z_{k+1} ) - B( k,z_k ) \big)$ via the non-decreasing $z$-list with which it is in bijection.


A staircase  $\phi \in \staircase_{(x,i) \to (y,j)}$ whose energy  attains the maximum value $M^1_{(x,i) \to (y,j)}$ is called a geodesic from $(x,i)$ to~$(y,j)$.

\subsubsection{Moving to scaled coordinates}

We want to depict geodesics in scaled coordinates and to measure their energy in a suitable scale. To understand what the right transformations will be, it is useful to mention how
the one-third and two-thirds KPZ scaling considerations that have been outlined  are manifest in Brownian LPP. When the ending height $j$ exceeds the starting height $i$ by a large quantity $n \in \N$, and the location $y$ exceeds $x$ also by $n$, then the maximum energy $M^1_{(x,i) \to (y,j)}$ grows linearly, at rate $2n$,
and has a fluctuation about this mean of order $n^{1/3}$. Moreover, if $y$ is permitted to vary from this location, then it is changes of $n^{2/3}$ in its value that result in a non-trivial correlation of the maximum energy from its original value.

On the basis of these facts, we want to introduce scaled coordinates suitable for the description of staircases whose height is of order~$n$.
The  relevant scaling transformation 
will be  $R_n:\R^2 \to \R^2$ given by
$$
 R_n \big(v_1,v_2 \big) = \Big( 2^{-1} n^{-2/3}( v_1 - v_2) \, , \,   v_2/n \Big) \, .
$$ 
In these scaled units, a vertical displacement of one unit corresponds to a displacement of $(n,n)$ in the original description. Horizontal displacement is magnified by a factor of twice $n^{2/3}$
in the original units.

The image of any staircase under $R_n$
will be called an $n$-zigzag. An $n$-zigzag has starting and ending points inherited from the preimage staircase.
A staircase moves alternately rightwards and upwards, while an $n$-zigzag correspondingly moves rightwards, directly, and northwesterly,
along sloping lines at the low gradient $-2n^{-1/3}$.



The image under $R_n$ of a geodesic is called an $n$-polymer, or often simply a polymer. 
Scaled geodesics play a central role in our study, so the succinct name `polymer' seems apt. This usage is not standard, however, since the term `polymer' is often used to refer to typical realizations of the path measure in LPP models at positive temperature.

We also want to describe a geodesic's energy in scaled units. 
The polymer whose endpoints (in scaled units) are $(x,0)$ and $(y,1)$ is for given $x,y \in \R$ almost surely unique, and will be denoted by~$\rho_{n;(x,0)}^{(y,1)}$.
(This uniqueness is proved as part of our study of Brownian LPP:  see~\cite[Lemma~$4.6(1)$]{Patch}.)
This polymer is the image under $R_n$ of the geodesic between $(2n^{2/3}x,0)$ to $(n  + 2n^{2/3}y,n)$.
The polymer will be ascribed a weight, which will be the geodesic's energy, measured in scaled units. The weight will be called   $\weight_{n;(x,0)}^{(y,1)}$. It is specified by the formula
\begin{equation}\label{e.weightmzeroone} 
  \weight_{n;(x,0)}^{(y,1)} \,     =  \,   2^{-1/2} n^{-1/3} \Big(  M^1_{(2n^{2/3}x,0) \to (n  + 2n^{2/3}y,n)} - 2n  -  2n^{2/3}(y-x) \Big) \, .
\end{equation}

The quantity  $\weight_{n;(x,0)}^{(y,1)}$ may be expected to be, for given real choices of $x$ and $y$, a unit-order random quantity, whose law is tight in the scaling parameter $n \in \N$. 

In fact, any $n$-zigzag $\phi$ receives a weight~$\weight(\phi)$, given by replacing the $M^1$ term in the formula~(\ref{e.weightmzeroone}) by the energy of the staircase that is mapped to $\phi$ under~$R_n$.

\subsection{Principal directions of inquiry}

Some of the leading players in our study have been introduced: the use of scaled coordinates, polymers with given endpoints such as $\rho_{n;(x,t_1)}^{(y,t_2)}$, and the polymer weight $\weight_{n;(x,t_1)}^{(y,t_2)}$. Having made the acquaintance of these characters,
we are ready to begin informally stating the main conclusions of this four-paper study.
We structure this presentation by next introducing four themes: the narrow wedge polymer weight profile; three KPZ scaling exponents viewed via Brownian LPP;  the general initial condition polymer weight profile; and the low probability that several disjoint polymers coexist in a unit-order region. 
In the ensuing subsections, we explain each theme and  present, albeit imprecisely, the principal associated conclusion reached by our study.  
Although results are labelled as theorems or corollaries,  all these assertions are merely informal versions intended to convey the rough meaning of the result in question. When the result is proved in one of the three companion papers, the paper and result name appear in square brackets at the result's statement.

It will not be apparent from this presentation what connections there are between the results. These connections are very substantial, as we will endeavour to communicate after presenting the four themes.

\subsubsection{The narrow wedge polymer weight profile}

 The random function $y \to   \weight_{n;(0,0)}^{(y,1)}$
may be viewed as the weight profile obtained by scaled maximizing paths that travel from the origin at  vertical coordinate zero to the variable location $y$ at vertical coordinate  one.
Since the first polymer endpoint is fixed, this is the narrow wedge case. This random function locally resembles Brownian motion, of diffusion constant one, 
but globally, it hews to the shape of the parabola $-2^{-1/2}y^2$, at least for values of $\vert y \vert = o(n^{1/9})$ that are not extremely large. 

Our study presents a very strong conclusion regarding the locally Brownian nature of the narrow wedge profile.

\begin{theorem}\label{t.bridge.informal}
Suppose that the weight profile  $y \to   \weight_{n;(0,0)}^{(y,1)}$ is viewed as a function of $y \in [-1,1]$
and is adjusted, by the addition of an affine function, in order that it vanishes at the endpoints $y = -1$ and $y = 1$.
The adjusted process takes values in the same space of continuous functions as does standard Brownian bridge $B:[-1,1] \to \R$ with $B(-1) = B(1) = 0$.
The adjusted weight profile is extremely similar to the bridge: for an event whose Brownian bridge probability equals $a \in (0,1)$, where $a$ is supposed to satisfy $a \geq \exp \big\{ - g n^{1/12} \big\}$,
the probability of the event for the adjusted profile is at most $a \exp \big\{  G    (\log a^{-1})^{5/6} \big\}$.
The positive constants $G$ and $g$ do not depend on $n$, provided that $n$ is large enough.  
\end{theorem}
This conclusion is reached in the present article, in Theorem~\ref{t.airytail.ln} and Proposition~\ref{p.lereg}.

\subsubsection{Three scaling principles: the powers one-third, two-thirds and one-half}

\begin{figure}[ht]
\begin{center}
\includegraphics[height=7cm]{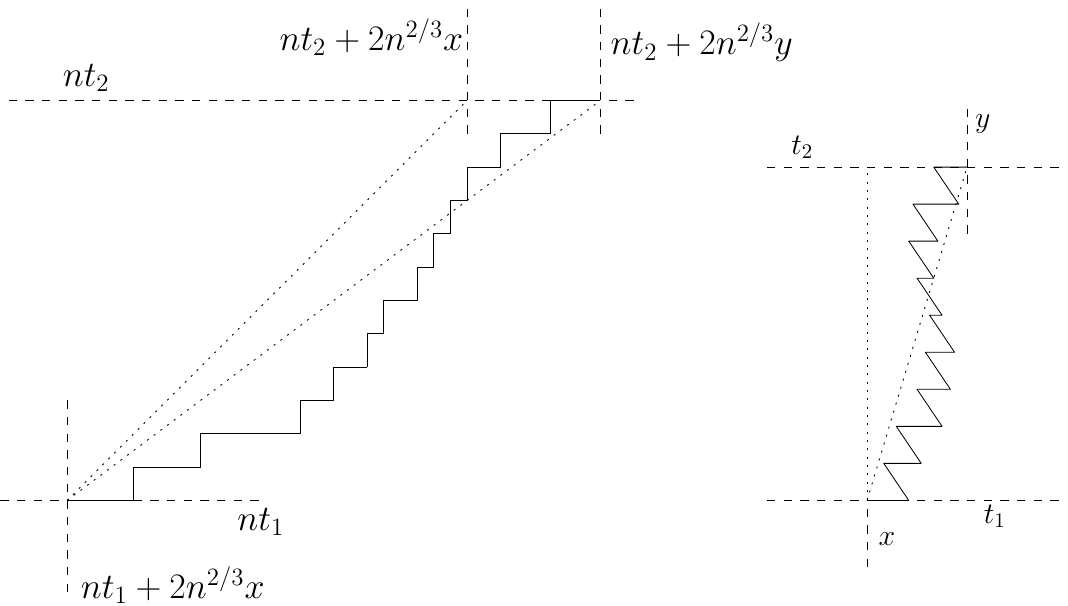}
\caption{Let $n \in \N$, $t_1,t_2 \in \R$, $t_1 < t_2$, and $x, y \in \R$. The endpoints of the geodesic in the left sketch have been selected so that, when the scaling map~$R_n$ is applied to produce the right sketch,
the $n$-polymer $\rho_{n;(x,t_1)}^{(y,t_2)}$ results.}
\label{f.scaling}
\end{center}
\end{figure}

Here we express three basic KPZ scaling principles. Each principle is manifest in scaled coordinate Brownian LPP and, in each case, a result to this effect has been proved in the four-paper study.
These results are certainly not the first expressions of these principles, which are fundamental to the whole arena of KPZ, but they offer strongly on-scale articulations of the principles and are valuable in our analysis of scaled Brownian LPP.
Before stating the results informally, we mention that the polymer $\rho_{n;(x,0)}^{(y,1)}$
may be viewed as having lifetime $[0,1]$, moving during its life from $x \in \R$ to $y \in \R$. It is equally possible to choose two general times $t_1,t_2 \in \R$ with $t_1 < t_2$, and speak of the polymer from 
$(x,t_1)$ to $(y,t_2)$. This lifetime-$[t_1,t_2]$ polymer will be denoted   $\rho_{n;(x,t_1)}^{(y,t_2)}$: see Figure~\ref{f.scaling}.
We will write $\tot = t_2 - t_1$ for the polymer's duration. The quantity $n \tot$ appears in some upcoming statements: this is the polymer lifetime in the original frame.

 \noindent{\em A power of two-thirds dictates polymer geometry.}  A polymer  whose lifetime is $[t_1,t_2]$, and thus has duration $\tot$, 
 fluctuates laterally, away from the planar line segment that interpolates its endpoints, by an order of $\tot^{2/3}$.

\begin{theorem}\label{t.polyfluc.informal}\cite[Theorem~$1.5$]{NonIntPoly} 
The polymer $\rho_{n;(x,t_1)}^{(y,t_2)}$ journeys from $(x,t_1)$ to $(y,t_2)$. Let $\ell$ denote the planar line segment that interpolates this pair of endpoints.
Let $a$ be a small positive quantity, and consider the polymer at time $t$, where $t$ is either $t_1 + a$ or $t_2 - a$. We may speak of the polymer's location at this time: this means its horizontal coordinate at height $t$.
The probability that the polymer location  differs from the time-$t$ location of $\ell$ by more than $a^{2/3}r$ is at most $\exp \big\{ - d r^\alpha \big\}$. Here, $\alpha = 3/4$.
The parameter $r$ is measuring polymer fluctuation on the scale suitable given the duration~$a$ from the endpoint time. It is supposed to verify the condition $r = o\big( (n \tot)^{1/36} \big)$.
The constant $d$ is universal.
\end{theorem}

\cite[Theorem~$2$]{BSS17} offers a similar assertion for exponential last passage percolation, with $\alpha = 1$. 

 \noindent{\em  A power of one-third dictates polymer weight.} A polymer whose lifetime is $[t_1,t_2]$ has a weight of order $\tot^{1/3}$. 
 Actually, this is only true if the polymer makes no significant lateral movement. For example, $\rho_{n;(x,t_1)}^{(y,t_2)}$ may be expected to have a weight of order $\tot^{1/3}$
 provided that the endpoints verify $\vert y - x \vert = O(\tot^{2/3})$. 
 Beyond this scale, parabolic curvature dictates polymer weight.

\begin{corollary}\label{c.maxminweight.informal}\cite[Corollary~$2.1$]{ModCon}
Consider the parabolically adjusted  polymer weight 
$$
\tot^{-1/3} \weight_{n;(x,t_1)}^{(y,t_2)} +  2^{-1/2} \tot^{-4/3}  \big(y - x  \big)^2 \, .
$$
Allow $x$ and $y$ to vary over a unit interval of values
in such a way that the normalized difference $\vert y - \vert \tot^{-2/3}$
 is at most order~$(n \tot)^{1/18}$.
The probability that the maximum absolute value of the parabolically adjusted weight exceeds $r$ is at most $\exp \big\{ - O(1) r^{3/2} \big\}$, uniformly in high choices of $n \tot$.
The parameter $r$ is supposed to be at most of order $(n \tot)^{1/18}$.
\end{corollary}

See Propositions~$10.1$ and~$10.5$ in~\cite{SlowBond}
 for comparable bounds concerning exponential and Poissonian LPP.
 
 \noindent{\em A power of one-half dictates polymer weight differences.} Consider two polymers of unit duration both of whose endpoints differ by a small quantity $\e$. 
 Then the polymers' weights typically differ by an order of $\e^{1/2}$. 

\begin{theorem}\label{t.differenceweight.informal}\cite[Theorem~$1.1$]{ModCon} 
For $x,y \in \R$ given, consider the maximum difference in weight between any pair of $n$-polymers travelling between $[x,x+\e] \times \{ 0 \}$ and $[y,y+\e] \times \{ 1 \}$.
The probability that this maximum exceeds $\e^{1/2} R$ is at most $\exp \big\{ - O(1) R^{3/2} \big\}$, uniformly in high choices of $n$.
We suppose that  $\vert y - x \vert$ and $R$ are at most of order $n^{1/18}$.
\end{theorem}

\subsubsection{The polymer weight profile initiated from general initial data}

The narrow wedge weight profile  $y \to   \weight_{n;(0,0)}^{(y,1)}$ concerns polymers all of which are forced to begin at the origin at time zero. We may make a much more general definition of a polymer weight profile.  This will be the  $f$-rewarded line-to-point polymer weight  $\weight_{n;(*:f,0)}^{(y,1)}$. Here, $f: \R \to \R \cup \{ -\infty \}$ is an initial condition. 
Paths (or rather $n$-zigzags) may begin anywhere on the real line at time zero; they travel to $y \in \R$ at time one. (Because they are free at the beginning and fixed at the end, we refer to these paths as `line-to-point'.) They begin with a reward given by evaluating $f$ at the starting location, and then gain the weight associated to the journey they make.  The value $\weight_{n;(*:f,0)}^{(y,1)}$, which we will define momentarily, denotes the maximum $f$-rewarded weight of all such paths.   In the notation $\weight_{n;(*:f,0)}^{(y,1)}$, we again use subscript and superscript expressions to refer to space-time pairs of starting and ending locations. The starting spatial location is being denoted $*:f$. The star is intended to refer to the free time-zero endpoint, which may be varied, and the $:f$ to the reward offered according to where this endpoint is placed.

Let $\mc{I}$ denote the space of $f$ that grow at most linearly. This is in essence  the broadest class of $f$ suitable for a study of the weight profiles $y \to \weight_{n;(*:f,0)}^{(y,1)}$ for all sufficiently high $n \in \N$.

For $f \in \mc{I}$, we define the $f$-rewarded line-to-point polymer weight  $\weight_{n;(*:f,0)}^{(y,1)}$ to be 
$$
    \sup_{x \in (-\infty,2^{-1}n^{1/3} + y]} \big(  \weight_{n;(x,0)}^{(y,1)}    + f(x) \big)  \, .
$$

The polymer that attains this weight is, naturally enough, denoted by  $\rho_{n;(*:f,0)}^{(y,1)}$. The $f$-rewarded line-to-point polymer uniqueness implicitly asserted by the use of this notation is proved in
~\cite[Lemma~$4.6(2)$]{Patch}.

Note that the narrow wedge case, in which all polymers depart from the origin, 
corresponds to the choice that $f:\R \to \R \cup \{-\infty\}$ equals zero at the origin and is otherwise minus infinity.

There are two main conclusions that the four-paper study reaches in regard to these general weight profiles.

First, in \cite[Theorem~$1.3$]{ModCon}, 
 the modulus of continuity of the profiles is studied. The maximum deviation of a profile over any interval of length $x$ lying inside a given bounded interval
is shown to be at most of the order of $x^{1/2} \big( \log x^{-1} \big)^{2/3}$. This result holds for $n$ high enough, and in essence it holds uniformly over the class of initial data $f$ lying in the function space~$\mc{I}$.

The second conclusion is one of the most conceptually central results in the study. It asserts a strong similarity to Brownian motion for general weight profiles.

\begin{theorem}\label{t.unifpatchcompare.informal}\cite[Theorem~$1.2$]{Patch}  
Let $f \in \mc{I}$. Consider the weight profile $y \to \weight_{n;(*:f,0)}^{(y,1)}$ as a function of $y \in [-1,1]$.
The domain $[-1,1]$ may be broken into a  random but finite number of subintervals called {\em patches}.
The weight profile when restricted to any given patch coincides on the patch with a random process --  call it the fabric process -- defined on $[-1,1]$ that withstands a very demanding Brownian comparison  on this interval. The comparision is made after the fabric process is affinely adjusted so that it vanishes at the interval endpoints $-1$ and $1$. To wit,
 suppose that $a \geq \exp \big\{ - O(1) n^{1/12}\big\}$ and that an event has probability $a > 0$ for {\em Brownian bridge} on $[-1,1]$; then its probability for the {\em affinely adjusted fabric process} is at most $a^{2/3 + o(1)}$.
This lower bound hypothesis on $a$ is mild, rapidly evaporating as $n$ rises; neglecting it, we may reexpress the last assertion: the Radon-Nikodym derivative of the latter process with respect to the former lies in any $L^p$ space for $p \in (1,3)$.
 The number of patches needed is random, but its probability of exceeding $\ell \in \N$ is at most $\ell^{-2 + o(1)}$.
There are some constants hidden in this description: in the $L^p$-norms for Brownian comparison, in the $+o(1)$ term in the exponent of the patch number decay rate, and in the lower bound on $n$ that is in fact needed.
All these constants can be chosen in essence uniformly in $f \in \mc{I}$.
\end{theorem}

As such, the weight profile $[-1,1] \to \R: y \to \weight_{n;(*:f,0)}^{(y,1)}$ may be called a {\em Brownian patchwork quilt} -- it is woven together from a small set of pieces of fabric, each vividly Brownian.

The ideas that drive the proof of this result form a conceptual highway for many parts of the four-paper study. We will overview them after discussing the final arena of main results.

\subsubsection{The rarity of coexistence of several disjoint polymers}

Consider a polymer $\rho_{n;(x,0)}^{(y,1)}$ with given endpoints $x,y \in \R$.
Considerations of scaling and polymer uniqueness mean that, if we perturb the time-zero endpoint by a short distance, we may expect that the perturbed polymer will merge after a short while with the original one, and that the two will then run together for the rest of their lives. Indeed, the two-thirds power law for polymer geometry suggests that the polymers  $\rho_{n;(x,0)}^{(y,1)}$ and $\rho_{n;(x + \e,0)}^{(y,1)}$ will typically merge by time of order $\e^{3/2}$.

A natural related question concerns $\maxpoly_{n;([x-\e,x+\e,0])}^{([y-\e,y+\e],1)}$, the maximum  possible cardinality of a pairwise disjoint set of $n$-polymers all of which begin in $[x-\e,x+\e]$ at time zero and end in $[y-\e,y+\e]$  at  time one. How quickly does the probability that this random variable exceeds a given integer decay as a function of $\e$ in the limit of low~$\e$? Can a bound be found that holds uniformly in high $n$?
 Our study provides a rather strong estimate.

\begin{theorem}\label{t.disjtpoly.pop.informal}\cite[Theorem~$1.1$]{NonIntPoly}
Set $I = [x- \e,x+ \e]$ and $J = [y- \e,y+  \e]$, we have that  
$$
\PP \Big( \maxpoly_{n;(I,0)}^{(J,1)}  \geq k \Big) 
  \leq  
 \e^{(k^2 - 1)/2  + o(1) } \, .
$$
The result holds for $k \in \N$ given, for $n$ large enough, and for $\vert x - y \vert \leq \e^{-1/2}$.
\end{theorem}

This theorem leads rather directly to an interesting consequence.
\begin{theorem}\label{t.maxpoly.pop.informal}\cite[Theorem~$1.4$]{NonIntPoly}  
Let $I$ and $J$ be unit-length intervals. Then the probability that  $\maxpoly_{n;(I,0)}^{(J,1)}  \geq k$ decays at least as fast as $k^{-O(1) (\log \log k )^2}$, uniformly in high~$n$.
\end{theorem}

As we will later indicate, Theorem~\ref{t.disjtpoly.pop.informal}  plays a key role in deriving the patch description of the general weight profile, Theorem~\ref{t.unifpatchcompare.informal}.  
Theorem~\ref{t.disjtpoly.pop.informal}  is a consequence of the solution of another natural Brownian LPP problem, which solution is provided in the present article. 
The problem concerns a topic that may called the {\em rarity of high-weight multi-polymer watermelons}.
For $k \in \N$, a $k$-tuple multi-polymer watermelon consists of $k$ $n$-zigzags, which share their pair of endpoints, but which are otherwise disjoint. The watermelon's weight is the sum of the weight of the constituent zigzags. 
 If the shared endpoints are $(0,0)$ and $(x,1)$, then the watermelon weight may be at most  $k \cdot \weight_{n;(0,0)}^{(x,1)}$. After all, each member zigzag   travels between $0$ at time zero and $x$ at time one. 
 We now define  the event $\neargeod_{n,k;(0,0)}^{(x,1)}(\eta)$ that there exists a $k$-tuple multi-polymer watermelon with endpoints $(0,0)$ and $(x,1)$ whose weight falls short of the maximum $k \cdot \weight_{n;(0,0)}^{(x,1)}$
 by at most~$\eta$. The problem of the rarity of the high-weight watermelon asks: what is the decay rate of the probability of this event as $\eta \searrow 0$? Here is the solution.

\begin{theorem}\label{t.neargeod.informal}
We have that 
$$
 \PP \Big(  
\neargeod_{n,k;(0,0)}^{(x,1)} \big(   \eta \big) \Big) =  \eta^{k^2 - 1 + o(1)} \, .
$$
For given $k \in \N$,
the result holds uniformly for  $\eta$ small,  $n$ high and $\vert x \vert = o(n^{1/9})$. 
\end{theorem}
The official version of this result will be presented in this article: see Theorem~\ref{t.neargeod}. 

\begin{figure}[ht]
\begin{center}
\includegraphics[height=9cm]{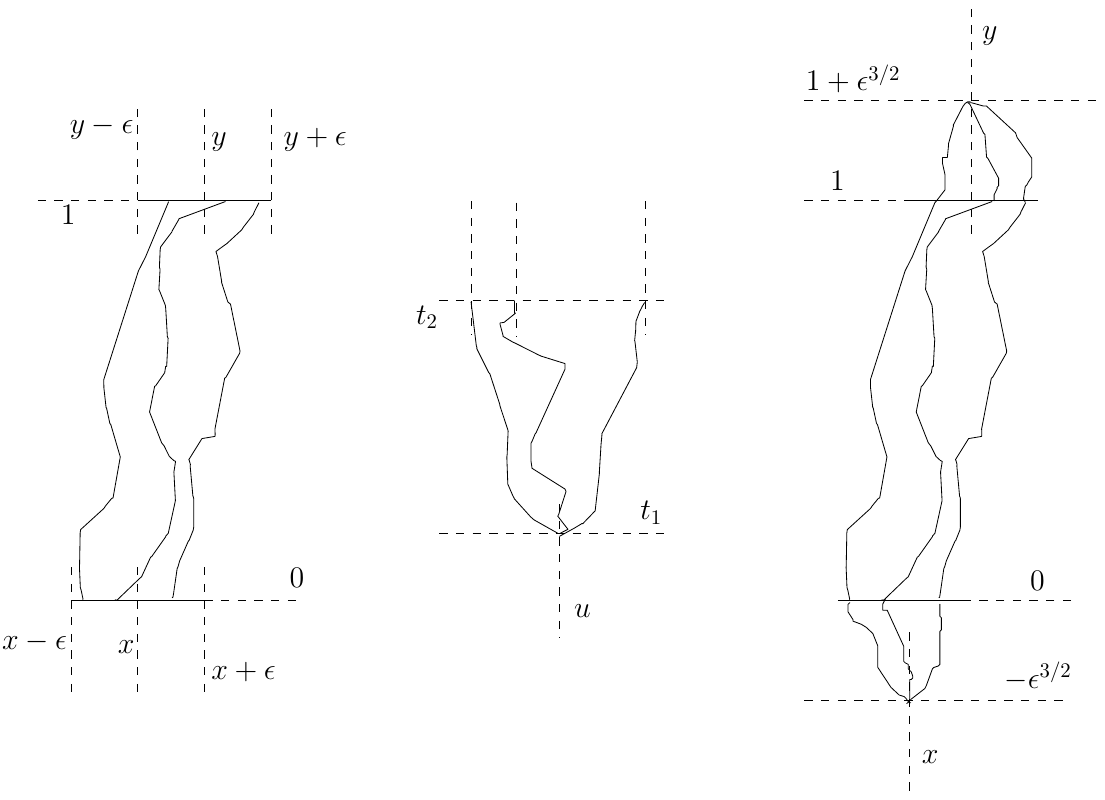}
\caption{{\em Left:}  The three polymers concerned in the event $\maxpoly_{n;([x-\e,x+\e],0)}^{([y-\e,y+\e],1)} \geq 3$. {\em Middle:} When surgery is undertaken with the aim of bringing together the endpoints of the three polymers, a bouquet of three paths with a shared endpoint will be used. {\em Right:} Surgery is performed with the use of two bouquets, facing in opposite directions, which extend the original three polymers into longer paths, each beginning and ending at $(x,-\e^{3/2})$ and $(y,1 +\e^{3/2})$. The outcome is a multi-polymer watermelon.}
\label{f.triple}
\end{center}
\end{figure}  

The relevance of the new result for the proof of  Theorem~\ref{t.disjtpoly.pop.informal} lies in the similarity in the two problems addressed by these theorems.
The latter concerns the presence of several disjoint polymers with nearby endpoints, whereas the former addresses the coexistence of several {\em near} polymers that actually share their endpoints but are otherwise disjoint. 
Indeed, the idea of the proof of Theorem~\ref{t.disjtpoly.pop.informal} is to relate the result to Theorem~\ref{t.neargeod.informal} by setting the latter's parameter $\eta$ to be of the order of the former's $\e^{1/2}$. A system of $k$ pairwise disjoint polymers with $\e$-near endpoints 
may typically be surgically altered to manufacture a system of near polymers with endpoints that are shared, but which otherwise maintain disjointness; in this way,  the probability bound in Theorem~\ref{t.neargeod.informal} may be harnessed.
The technique of surgery, which is the backbone of the article~\cite{NonIntPoly}, is depicted in Figure~\ref{f.triple}.
The three results concerning the KPZ scaling exponents all have a role to play in ensuring that surgery may be performed at a suitable, on-scale, cost. Because endpoints need to be moved by an order of $\e$ in the horizontal coordinate, the two-thirds power law for polymer geometry dictates that changes should be made in a region of height $\e^{3/2}$ near the endpoint locations. The one-third power law for polymer weight implies that the concerned changes involve short paths of typical weight $\big( \e^{3/2} \big)^{1/3} = \e^{1/2}$. These considerations determine that the choice $\eta = \Theta( \e^{1/2} )$ is suitable and thus why a factor of one-half appears in the exponent $(k^2 - 1)/2$ in Theorem~\ref{t.disjtpoly.pop.informal}.

\subsection{General weight profiles and polymer forests: a snapshot}

We have offered informal assertions of many of the principal conclusions of the four-paper study. We do not propose here to attempt a detailed overview of their proofs. In keeping with the aim of this section as a global snapshot of our study, we will now explain a certain critical concept that underlies the proof in~\cite{Patch} of the Brownian comparison result for general polymer weight profiles, reviewed here as Theorem~\ref{t.unifpatchcompare.informal}. 
Although the basic idea that we will explain has required significant refinement in order to be implemented rigorously, it remains conceptually crucial to the overall approach in the study.
 As such, we hope that presenting the idea now will serve to explain at least roughly how the main components of the study interact. This article and the three companion papers contain further expository material; this material may be read alongside the snapsnot we now present for a more detailed overview. (A few words about such further reading follow, in the next subsection.)

Consider again the polymer weight profile $y \to \weight_{n;(*:f,0)}^{(y,1)}$ for an initial condition~$f$ in the function space~$\mc{I}$.
To obtain Theorem~\ref{t.unifpatchcompare.informal}, we want to break the domain $[-1,1]$ into a random but controlled number of patches and show that the patch-restricted profiles withstand an exacting comparison to Brownian motion.

We are aided in this task by the first major input of our study. In the special, narrow wedge, case, where $f$ equals zero at a given location and is otherwise minus infinity,
Theorem~\ref{t.bridge.informal} shows that the affinely adjusted weight profile, when compared to Brownian bridge, has a Radon-Nikodym derivative that lies in every $L^p$-space for $p \in (1,\infty)$, uniformly in high~$n$.

\begin{figure}[ht]
\begin{center}
\includegraphics[height=12cm]{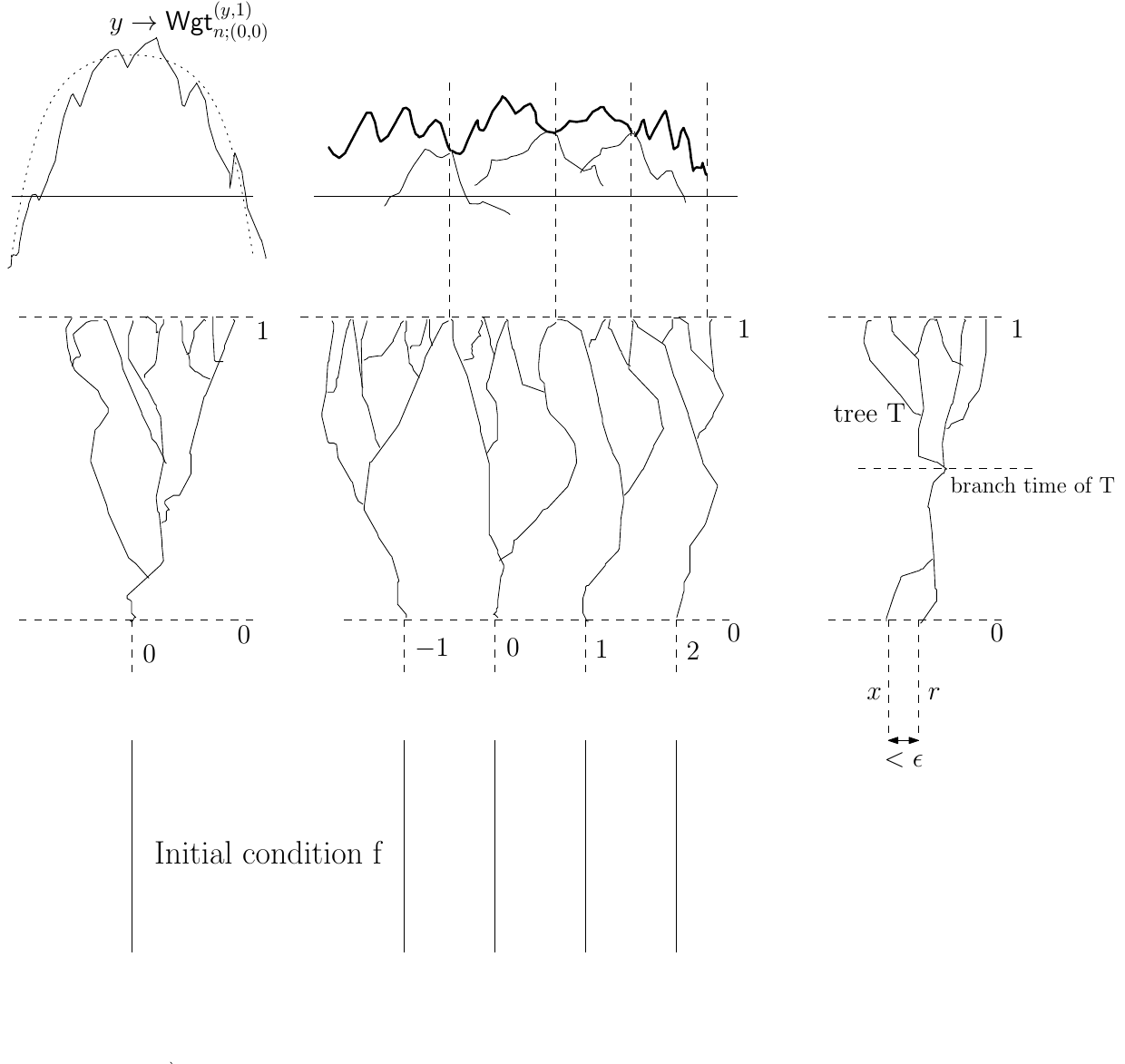}
\caption{In the left and middle sketches, the weight profile $y \to \weight_{n;(*:f,0)}^{(y,1)}$, the $f$-rewarded polymer forest, and the function $f$ are depicted.
 On the left, $f$ is zero at zero and otherwise minus infinity, so that the profile is the narrow wedge case $y \to \weight_{n;(0,0)}^{(y,1)}$.
In the middle, $f$ is instead zero at integer points. The  profile $y \to \weight_{n;(*:f,0)}^{(y,1)}$ is depicted in bold. The dashed vertical lines that contact the line at height one indicate canopy boundary points.
In a given canopy with root $r$, the emboldened profile takes the form $y \to \weight_{n;(r,0)}^{(y,1)}$. The right sketch depicts one tree, rooted at $r$, in an $f$-rewarded polymer forest. The proposed rerooting of this tree to a nearby discrete mesh element $x \in \e \Z$ is illustrated.}
\label{f.manytrees}
\end{center}
\end{figure}

Figure~\ref{f.manytrees}
illustrates several key ideas in the conceptual highway that leads from Brownian comparison for narrow wedge profiles to the much more general case.
To begin on that highway, we point out a geometric view of the narrow wedge weight profile:  
see the left sketch of Figure~\ref{f.manytrees}. Note that $\weight_{n;(0,0)}^{(y,1)}$
 is the weight of the polymer $\rho_{n;(0,0)}^{(y,1)}$. These polymers, indexed by $y \in \R$, all stream out of the origin at time zero, to arrive at their various ending locations $y$ at time one. Almost sure uniqueness of polymers with given endpoints suggests that, once separated, polymers will not meet again. Thus, the system of polymers should be viewed as a tree, with a root at $(0,0)$, and a canopy $\R \times \{ 1 \}$ of ending locations.

 Return now to a general initial condition $f: \R \to \R \cup \{ - \infty \}$.
 The $f$-rewarded line-to-point polymers may be traced backwards in time from locations $(y,1)$ with $y \in \R$. They arrive at time zero at a variety of locations, in contrast to the narrow wedge case. 
 They share something with that case, however: as time decreases from one to zero, two polymers that meet will stay together; this follows from $f$-rewarded line-to-point polymer uniqueness.

 This fact has the implication that we may view the collection of polymers  $\rho_{n;(*:f,0)}^{(y,1)}$ indexed by $y \in \R$ as a forest, which we may call the $f$-rewarded polymer forest.
 Each constituent tree has a root lying on the $x$-axis, and a canopy that consists of an interval lying in the line $y = 1$. 
 Indeed, we may partition the $y=1$ copy of the real line into this set of canopies. 
 The polymer weight profile $y \to \weight_{n;(*:f,0)}^{(y,1)}$, when restricted to any given canopy, 
 would seem to have much in common with the narrow wedge weight profile. After all, all the concerned polymers, ending at locations $(y,1)$, for points $y$ in the given canopy, share their starting location, namely $(r,0)$, where $r \in \R$ is the root of the canopy in question: see the middle sketch of  Figure~\ref{f.manytrees}.

 This then is the basis of our description of the general weight profile 
 $y \to \weight_{n;(*:f,0)}^{(y,1)}$ as a patchwork quilt of pieces, on each of which the profile withstands a tough comparison to Brownian motion. The patches of this quilt should be the canopies in the $f$-rewarded polymer forest. 
 On any given patch, the weight profile takes the form $y \to \weight_{n;(r,0)}^{(y,1)}$ where $r$ is the root of the tree in the polymer forest whose canopy is the patch. Since $r$ is fixed, the weight profile should withstand a demanding comparison to Brownian motion, because the narrow wedge profile does.
 
 This plan is  a crude outline of how Theorem~\ref{t.unifpatchcompare.informal} will be proved. There are serious obstacles along the way, whose resolution requires that we draw on several of the main conclusions that we summarised in the preceding section. We do not explain the obstacles or their resolution in any detail in this overview. Or rather, we do so only in the merest outline,  by voicing two basic objections to the plan, and then indicating something of how our apparatus will answer them.
 
\begin{itemize}
\item The decomposition into patches in the general weight profile description that we seek calls for control on the number of patches per unit length. Certainly, we want to argue that this number is tight in $n \in \N$, so that the description can be said to be valid in scaled coordinates. We must understand, then, what prevents there being an asymptotically unlimited (in the high~$n$ limit) number of trees per unit length  in any $f$-rewarded polymer forest.
\item It is important for the plan to work that the weight profile  $y \to \weight_{n;(r,0)}^{(y,1)}$, where $r$ is the root of a given tree in the $f$-rewarded polymer forest,   is indeed statistically similar to the narrow wedge case, where say $r$ is zero.  Trouble looms here because the location $r$ is random, and it may be exceptional. As such, properties of the narrow profile    $y \to \weight_{n;(x,0)}^{(y,1)}$, such as locally Brownian structure, that obtain for typical values of $x \in \R$, may disappear for the possibly exceptional choices of root $r$ that we are considering.
\end{itemize}

We complete our tour of the conceptual highway of the four-paper study by indicating, very briefly, how the two objections may be answered; in so doing, we will see at least in a vague way how some of the principal conclusions of the study will play a useful role. 

In regard to the first difficulty,   consider for a further time the $f$-rewarded polymer forest. Pick a location in the interior of each canopy, so that a point process on $\R \times \{ 1 \}$ results.
 From each point, draw the $f$-rewarded line-to-point polymer back to time zero. This system of polymers is pairwise disjoint, because the time-one endpoints lie in different canopies.
 If the tightness alluded to in the first problem were to fail, we would be drawing, for high $n$, many disjoint polymers, all ending on a given unit-order interval. 
 It is here that the bound from~\cite{Patch}, quoted here as Theorem~\ref{t.maxpoly.pop.informal},
 asserting that it is a superpolynomial rarity in $k$ that $k$ disjoint polymers coexist in a unit-order region, serves by providing the desired tightness in $n$.

Regarding the second difficulty, that root locations $r$ may be exceptional, we may try to solve the problem by a rerooting procedure: see the right sketch in Figure~\ref{f.manytrees}. Close to the location $r$ is an element $x$ in the discrete $\e$-mesh $\e \Z$, where $\e > 0$ is small but fixed. The polymer weight profile that we are trying to describe,  $y \to \weight_{n;(r,0)}^{(y,1)}$, may be  compared to the narrow wedge profile  $y \to \weight_{n;(x,0)}^{(y,1)}$. Because $x$ lies in a discrete mesh, it may be viewed as a typical location, so that the objection raised in the second bullet point does not arise.
But why should the nearby profile  $y \to \weight_{n;(x,0)}^{(y,1)}$ be expected to offer an accurate description of  $y \to \weight_{n;(r,0)}^{(y,1)}$? 
We will now present an heuristic argument that, since $r$ and $x$ are close, these two functions should, with a high probability determined by this closeness, differ by a random constant as $y$ varies over the canopy of the tree of which $r$ is the root; Figure~\ref{f.manytrees}'s right sketch again  illustrates the heuristic. The tree in question has a {\em branch} time, with the constituent polymers following a shared course (the tree trunk, if you like) until that time, after which they may go their separate ways to the various locations in the canopy. 
Now, when the tree is rerooted the short distance from $(r,0)$ to $(x,0)$, we may expect that its structure changes by a modification in the form of the polymers only close to time zero.  Provided this modification has finished by the branch time, it will affect only the form of the tree trunk, and will be shared by all concerned polymers, no matter at  which canopy point $(y,1)$ they end. 
For this reason, the difference $y \to \weight_{n;(x,0)}^{(y,1)} - \weight_{n;(r,0)}^{(y,1)}$ may be expected to equal the discrepancy in weight between the new and the old tree trunks with high probability,
and thus typically be independent of the canopy location~$y$. 

In fact, when this potential resolution of the second objection is pursued, it becomes necessary to have powerful estimates on the probability that polymers from $(r,0)$ and the rerooted location $(x,0)$ merge quickly enough, not least because union bounds over the approximate root location are necessary. When the pair of polymers fail to so merge, two disjoint polymers with nearby endpoints coexist. Theorem~\ref{t.disjtpoly.pop.informal} plays a vital role in providing a suitably strong estimate, which ultimately enables the $L^{3-}$-space comparison seen in Theorem~\ref{t.unifpatchcompare.informal}. 
In fact, for reasons we will not attempt to explain here, this key application of  Theorem~\ref{t.disjtpoly.pop.informal} is made with $k=3$, so that a triple of polymers are concerned, rather than a pair.

\begin{figure}[ht]
\begin{center}
\includegraphics[height=16cm]{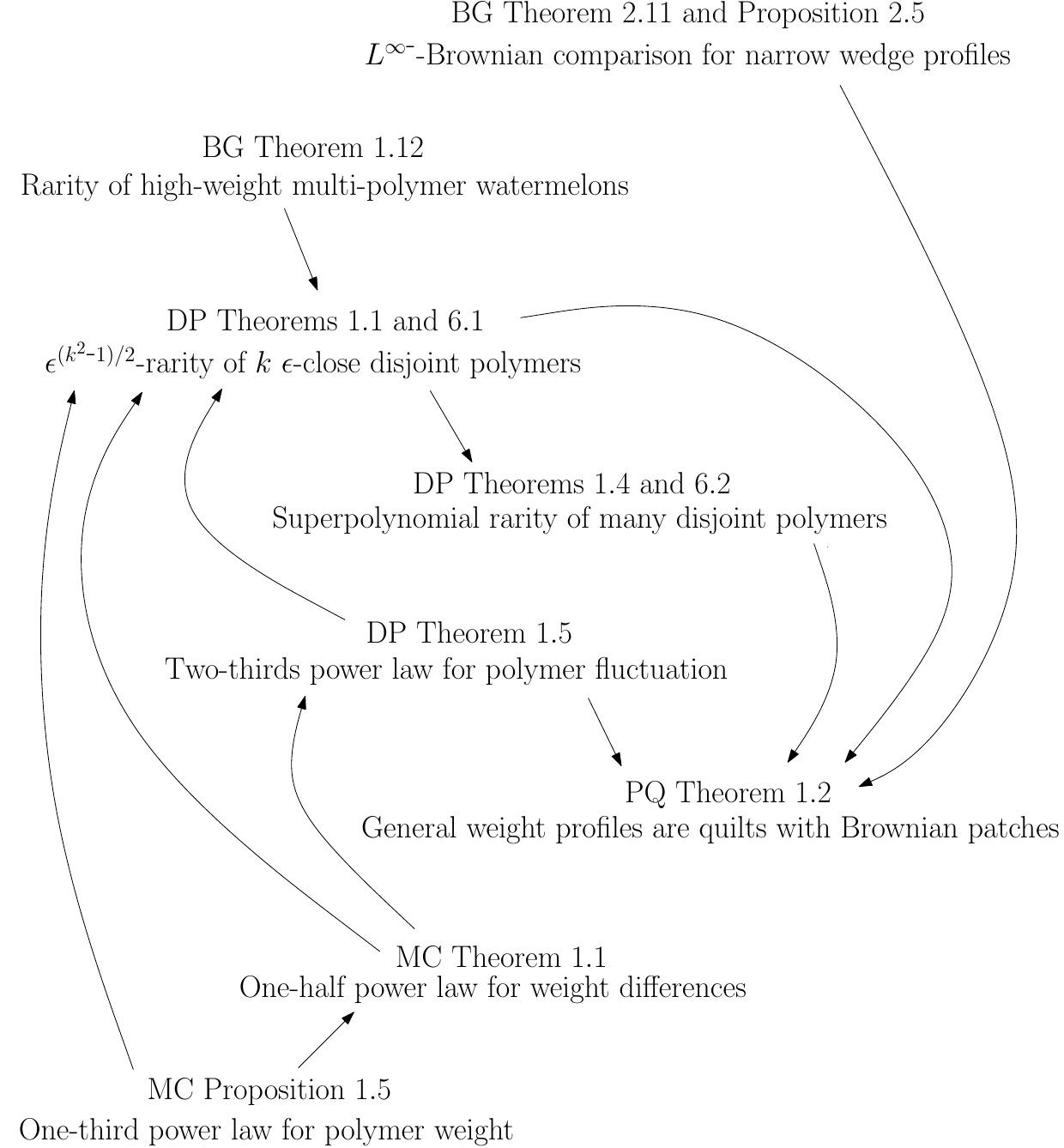}
\caption{The most significant conclusions in the four-paper study regarding scaled Brownian LPP are recorded. Arrows indicate that one result is used in the proof of another's.
The acronyms BG, MC, DP and PQ stand for Brownian Gibbs, Modulus of Continuity, Disjoint Polymers and Patchwork Quilt, and denote the present article,~\cite{ModCon},~\cite{NonIntPoly} and~\cite{Patch}. The sketch indicates only the more consequential results and implications; several significant technical results are not represented. It is apparent from the figure that~\cite[Theorem~$1.2$]{Patch} is a focal point of our study, and that the above order of acronym list would be a natural order for a reader who wishes to read the entire four-paper study.}
\label{f.map}
\end{center}
\end{figure}

Figure~\ref{f.map} may be consulted for an alternative summary of the four-paper study. 

\subsection{A brief discussion of further related topics}

\subsubsection{Directions to a more detailed overview}
We now make suggestions to the reader who has perused the preceding for a conceptual overview and wishes for a more detailed heuristic explanation regarding the Brownian LPP applications in the three companion papers. Leaving aside the contribution of the present work, the two most substantial concepts in the four-paper study concern the use of surgery in~\cite{NonIntPoly} to derive Theorem~\ref{t.disjtpoly.pop.informal} from Theorem~\ref{t.neargeod.informal}; and the polymer forest ideas for proving the patch description that we have just outlined.
The reader who wants more heuristic description of the first of these topics should read the `road map' in Section~$2.2$ in~\cite{NonIntPoly}. Regarding the second, we mention that our conceptual overview is similar to the `rough guide' offered in Section~$3$ of~\cite{Patch}. The heuristic story of how the patch description is built then continues in Section~$5$ of that article, where the rerooting procedure that is supposed to resolve the second of the two objections bullet-pointed above is explained more carefully.

\subsubsection{The limiting picture in high $n$}
The use of scaled coordinates to describe polymers and their weights is compatible with the expectation that it should be possible to take a high~$n$ limit of these structures and find a unique distributional limit.
The limiting field of polymer weights is in essence the space-time Airy sheet and the limiting system of polymers, the polymer fixed point. (The language of `fixed point' refers to the viewpoint that this structure is a fixed point of a renormalization operator.)
The space-time Airy sheet and the polymer fixed point are conjectural objects whose properties have been non-rigorously discussed in~\cite{CQR2015}.
A recent advance~\cite{MQR17} by Matetski, Quastel and Remenik rigorously constructs an important part of the space of information expected to be encoded in the space-time Airy sheet. These authors utilize a biorthogonal ensemble representation found by~\cite{Sas05,BFPS07} associated to the totally asymmetric exclusion process in order to  find Fredholm determinant formulas for the multi-point distribution of the height function of this growth process begun from an arbitrary initial condition. Using these formulas to take the KPZ scaling limit, the authors construct a scale invariant Markov process that lies at the heart of the KPZ universality class. The time-one evolution of this Markov process may be applied to very general initial data, and it is tempting to form a connection to our discussion by saying that this Markov process is the map
that sends $f \in \mc{I}$ to the random function $y \to \weight_{n;(*:f,0)}^{(y,1)}$, where $n$ assumes the value infinity. 
However,~\cite{MQR17}  is working with a totally asymmetric exclusion prelimit, rather than with Brownian LPP. If this new technique can be modified to apply to Brownian LPP, we will be able to speak of the limiting weight profile 
 $y \to \weight_{\infty;(*:f,0)}^{(y,1)}$. It may then be reasonable to attempt to bring the various conclusions of our Brownian LPP study to bear on the limiting profile. This view is reasonable because our assertions hold uniformly in high (but finite)~$n$, which is not to say that the existence of the limiting object would render the translation task from finite~$n$ to $n = \infty$ to be trivial in all cases. 

\subsubsection{How self-contained is the four-paper study?}
Are the proofs of the conclusions that are reached in the  three companion papers contained in those papers and the present article, without  inputs quoted from elsewhere?
The answer is no, but the outside inputs are fairly limited. Certain facts are needed about Brownian bridge,
and mutually avoiding systems of Brownian bridges, such as the Karlin-McGregor formula; these aside, the inputs are four assertions that will be recalled at the end of Section~\ref{s.neargeod} and in Section~\ref{s.deviation}:  O'Connell-Yor's observation that Brownian LPP has a Dyson Brownian motion ensemble structure; two bounds for the uppermost {\rm GUE} eigenvalue, due to Aubrun and Ledoux; and an observation recorded by Grabiner that relates these GUE eigenvalues to Dyson Brownian motion. 


\section{Non-intersecting line ensembles and their integrable and probabilistic analysis}\label{s.nonint}

We now leave the discussion of the four-paper study of Brownian LPP and turn to explaining the main results and ideas of the present article. In the next section, we will state the principal conclusions. In this one, we explain some of the overall themes and methods that are developed in the article.

A useful point of departure for this discussion is the scaled narrow wedge polymer weight profile in Brownian LPP, namely $\big[ - 2^{-1} n^{1/3} , \infty\big) \to \R: y \to \weight_{n;(0,0)}^{(y,1)}$.
A universal KPZ object may be expected to arise in the high $n$ limit. The limiting process will be defined on the entire real line. Calling it $\weight_{\infty;(0,0)}^{(y,1)}$, we would expect to have  $\weight_{\infty;(0,0)}^{(y,1)} = 2^{-1/2} \big( {\rm Airy}_2(x) - x^2  \big)$. The process ${\rm Airy}_2 : \R \to \R$ is a continuous and stationary process whose finite dimensional distributions may be expressed via Fredholm determinants.
Convergence to this process for the scaled narrow wedge polymer weight profile was proved for Poissonian last passage percolation by~\cite{PrahoferSpohn}  and~\cite{BDJ1999}, who showed convergence in finite dimensional distributions. Johansson~\cite{Johansson2003} derived a functional limit theorem, applying to LPP with an environment of geometric random variables, and proved the continuity of the  ${\rm Airy}_2$ process.

These results indicate how progress regarding the narrow wedge profile has been in large part due to the use of integrable techniques. Another vital perspective on the narrow wedge, which is also integrable in nature, involves writing the profile  $\weight_{n;(0,0)}^{(y,1)}$ in the form $\mc{L}_{n+1}(1,y)$, and embedding this process as the lowest indexed curve in a random ensemble $\mc{L}_{n+1}: \intint{n+1} \times  \big[ - 2^{-1} n^{1/3} , \infty\big) \to \R$ of continuous curves (of which there are $n+1$ in total). The ensemble is ordered, with $\mc{L}_{n+1}(1,\cdot)$ uppermost. We will specify this embedding precisely later, but it is defined via the Robinson-Schensted-Knuth correspondence, with the higher indexed curve values being specified in terms of the maximum weight of {\em multi-polymer watermelons}, namely  systems of point-to-point paths that are disjoint except for their shared endpoints.

The reason that embedding the scaled narrow wedge profile into this ensemble is valuable is that any given ensemble $\mc{L}_n$, for $n \in \N$, is in essence a system of Brownian motions conditioned on mutual avoidance.
Something similar is true in several LPP models, before or after the objects are scaled: in Poissonian LPP, for example, the unscaled narrow wedge geodesic energy profile is, as~\cite{PrahoferSpohn} demonstrated, the uppermost curve in a system of continuous time simple random walks conditioned on mutual avoidance, with a suitable boundary condition.

 Ensembles of this type, namely collections of one-dimensional Markov processes, such as random walks or Brownian motion, conditioned on mutual avoidance, are 
discussed the surveys \cite{Ferrari2010,FerrariSpohn2011} and  \cite{Johansson2006}. They form an important class of random models which arise in the study of random matrix theory, last passage percolation (or directed polymers), determinantal point processes, random tiling problems and asymmetric growth processes. 

An integrable technique is central to the analysis of such systems. The Karlin-McGregor formula, or, in different guises, the Lindstr{\"o}m-Gessel-Viennot formula and the method of free fermions in physics, expresses joint distributions of such systems as determinants whose entries are transition probabilities for the underlying Markov processes.
For example,  asymptotic analysis may be applied to these determinantal formulas, with exact expressions emerging for certain universal scaling limits; convergence to the limit is in the sense of finite dimensional distributions. This was the basis of the derivation of the limiting Airy$_2$ process from Poissonian LPP found by ~\cite{PrahoferSpohn}  and~\cite{BDJ1999}.
And indeed, such analysis may be performed not only for the polymer weight profile, which is the lowest indexed ensemble curve, but for the entire ensemble. In~\cite{PrahoferSpohn}, the high-$n$ finite distributional convergence of the Poissonian counterpart to  the ensemble that we have labelled $\mc{L}_{n+1}$ was proved; the limit is a  consistent family of finite dimensional distributions known as the multi-line Airy process (after the addition of a parabola).

Thus integrable techniques have been crucial to the analysis of these non-intersecting random walk or Brownian motion ensembles. The vital governing theme of the present article is that they are
also very well suited to analysis by probabilistic techniques. 
For Brownian LPP, the ensembles are mutually avoiding systems of Brownian motions (with certain boundary conditions).
Such ensembles enjoy an explicit resampling property, concerning say the conditional distribution of the top curve on a given interval given this curve's values at the interval's endpoints and the status of all the other curves.
(The rule is simple: the conditional distribution is Brownian bridge between the endpoints conditioned to remain above the given second curve.)
 In~\cite{AiryLE}, this {\em Brownian Gibbs} property was exploited in order to reexamine and develop the passage to the limit from mutually avoiding Brownian bridge systems at their edge to the multi-line Airy process. By establishing a certain uniform regularity at the edge of such systems, which holds after the KPZ scaling is taken, it was argued that convergence occurs, after a parabolic shift, to the Airy line ensemble, a positive integer-indexed ordered collection of random continuous  curves on the real line, stationary under horizontal shifts, whose curves are locally Brownian, and whose finite dimensional distributions are given by the multi-line Airy process.        The uppermost curve in the Airy line ensemble, namely the Airy$_2$ process, is thus a locally Brownian object.
 This statement may be interpreted either by taking a {\em local} limit, in which, for given $x \in \R$, the Gaussianity of $\e^{-1/2} \big( {\rm Airy}_2 (x + \e) - {\rm Airy}_2 (x) \big)$ is investigated after the low $\e$ limit is taken. The emergence of Brownian motion in this local limit of ${\rm Airy}_2$ has been proved by~\cite{Hagg} and~\cite{CatorPimentel}. An in essence stronger assertion of Brownian structure concerns a {\em unit-order} scale: for example, the absolute continuity of $[x,x+K]: y \to  {\rm Airy}_2(y) -  {\rm Airy}_2(x)$ with respect to Brownian motion (whose rate is two in view of the convention of definition for ${\rm Airy}_2$) on any compact interval~$[x,x+K]$. This last assertion was proved in~\cite{AiryLE} via the Brownian Gibbs technique.    

In the preceding section, we informally stated some of this article's principal conclusions, in the form of Theorem~\ref{t.bridge.informal} and~\ref{t.neargeod.informal}, with a view to their role in a broader study of scaled Brownian LPP. Now that we focus on the article's main conclusions apart from this role, we may summarise them thus:
\begin{itemize}
\item   In Theorem~\ref{t.rnbound}, an assertion of the close resemblance between  Brownian bridge and  affinely shifted curves in the Airy line ensemble, expressed as the finiteness of a super polynomial moment of the Radon-Nikodym derivative of the ensemble's curves. The same techniques yield  Theorem~\ref{t.bridge.informal}. These powerful, unit-order scale, results are achieved by pursuing the Brownian Gibbs technique introduced in~\cite{AiryLE}. 
\item In Theorem~\ref{t.airymodcon}(1), we prove that the {\rm Airy}$_2$ process has modulus of continuity of at most the order of $x^{1/2} \big( \log x^{-1} \big)^{1/2}$.
The polylogarithmic power law may be expected to be sharp.
 In Theorem~\ref{t.airymodcon}(2), we see that the same process varies by more than $\e^{1/2} K$
 on a given interval of length~$\e$ with probability at most $\exp \big\{ - O(1) K^{3/2} \big\}$ for any large enough $K$.
\item In Theorem~\ref{t.neargeod} and Corollary~\ref{c.neargeod}, we present the solution of a natural {\em high-weight multi-polymer watermelon rarity} problem in Brownian LPP.  This solution was informally stated in Theorem~\ref{t.neargeod.informal}.
\item In Theorem~\ref{t.othereigen},  control on the tails of edge eigenvalues in the Gaussian unitary ensemble that extends results presently known only for the top eigenvalue.
\end{itemize}

Leaving aside this article's role as a supplier of important tools for the broader Brownian LPP study, the guiding concept of the article is
 the use of the Brownian Gibbs property as a valuable probabilistic tool in the study of
 non-intersecting line ensembles.
  Beyond the concrete consequences that we have mentioned, the article introduces a general technique, the {\em jump ensemble method}, for proving upper bounds on the probabilities of events associated to ensembles of random curves that enjoy the Brownian Gibbs property. The method achieves outcomes, such as the membership in all $L^p$-spaces of the Radon-Nikydom derivative of the affinely adjusted Airy$_2$ process with respect to Brownian bridge, that have hitherto been out of reach; the article raises the prospect of reaching further such conclusions, and the jump ensemble method has been presented in a guise that we hope will permit its further application.
 
 We conclude this section by discussing how the Brownian Gibbs and other resampling properties have been used to date.
 First however, in an aside, we explain why the property is tied specifically  to {\em Brownian} last passage percolation.
 The RSK correspondence associates to various LPP models a mutually avoiding system of random curves. 
 All these presumably enjoy the Brownian Gibbs property after the scaling limit of $n \to \infty$ is taken. But Brownian LPP has the virtue of actually enjoying the Brownian Gibbs property in the prelimit;
 other LPP models would only satisfy an approximate version for high~$n$. The symmetries and invariance properties of Brownian motion make the Brownian Gibbs property a more attractive one for the purpose of study. This consideration explains why our four papers investigate {\em Brownian} LPP.
 



When a certain parabola is subtracted from each curve in the Airy line ensemble, the resulting non-intersecting ensemble enjoys the Brownian Gibbs property. The property 
 is shared~\cite{AiryLE} by perturbations of the parabolically shifted ensemble in which the highest curves are lifted away from the parabola far from the origin to become `Airy wanderers'~\cite{BBP2005}. It may be expected that it is shared by $\Z$-indexed tacnode processes~\cite{Johansson2013,DKZ2011,BorodinDuits}, which may be visualized crudely by inverting one parabolically shifted Airy line ensemble and placing it above another so that the curves of the two are forced into a deformation in and around a bounded region due to their mutual avoidance; by the Pearcey process~\cite{TracyWidom2006}, another $\Z$-indexed ensemble in which mutually avoiding Brownian motions split at a `cusp' into two packets in a neighbourhood of the origin, with the higher group surging upwards and the lower group falling away. 
The Airy line ensemble emerges from
the tacnode~\cite{Girotti2014} at a generic edge location, and from the
 Pearcey process far from the cusp; see~\cite{BC2013} for examples of such relationships. Dyson Brownian motion~\cite{Dyson1962}, a $\Z$-indexed collection of mutually avoiding Brownian motions that may be expected to arise as a bulk scaling limit of the Airy, Pearcey and tacnode line ensembles, is a closely related example. It would be interesting to construct rigorously the Dyson, Pearcey and tacnode line ensembles, as the Airy line ensemble was in~\cite{AiryLE}, and derive ensemble relationships between them that extend their known correlation kernel and finite dimensional distributional relations using a notion such as weak convergence of line ensembles in the upcoming Definition~\ref{maindef}.  

The Brownian Gibbs technique has found application 
in scaling limits that lead to softenings of the hard-core avoidance constraints among curves.  
The KPZ equation has been predicted since~\cite{KPZ} to model in a universal way surface growth with local randomness, smoothness and slope dependent growth speed, as the survey~\cite{IvanSurvey} discusses.
The equation  
has a narrow wedge solution, which models growth initiated at a point, that may be centred at late time $t$ and scaled by the characteristic factors of $(t^{2/3},t^{1/3})$. The resulting scaled solution for time $t > 0$ was shown  in~\cite{KPZLE} to be embedded as the lowest indexed curve in a KPZ$_t$ line ensemble that enjoys a softening of the Brownian Gibbs property determined by a Hamiltonian~${\bf H}_t$ that energetically penalizes, but does not forbid, curve crossing. Strong inferences result regarding the Brownian regularity of the scaled solution.
The counterpart to the prelimiting mutually avoiding Brownian bridge systems seen in the zero-temperature, Airy, theory is played by pairwise repulsive diffusions associated to the quantum Toda lattice Hamiltonian. These diffusions and their connections to directed polymers were discovered by O'Connell~\cite{OConnellToda}.

 The narrow wedge KPZ solution is the logarithm of the solution of the stochastic heat equation with multiplicative white noise with Dirac delta initial condition. In~\cite{OCW2016}, a multi-layer extension of this solution was constructed, with a conjectural Markovian evolution in time that was recently proved in~\cite{LunWarren}.  The logarithm of the multi-layer continuum system is, after scaling, given by the KPZ line ensemble: this has recently been shown by~\cite{Nica2016}, building on techniques of~\cite{CorwinNica}, as a consequence of this multi-layer continuum system being exhibited as the limit of the above mentioned multi-layer directed polymer model associated to the quantum Toda lattice.
In this limit, the number of ensemble curves is sent to infinity in a regime of intermediate disorder, already discussed in~\cite{AKQ2014}, 
 in which
 the temperature  scales to zero as a function of the system size.

The authors of~\cite{OsadaTanemura}
have constructed 
an infinite-dimensional stochastic differential equation modelling 
Dyson's Brownian motion and proved 
existence and pathwise uniqueness of solutions. They have also shown~\cite{OsadaTanemura2016} that the multi-line Airy process satisfies these stochastic dynamics, so that the SDE description and the multi-line Airy process are the same.

Systems of non-intersecting random walks are known as vicious walkers and were introduced by de Gennes~\cite{deGennes1968} as a soluble model of long thin fibres at thermal equilibrium subject to unidirectional stretching. These walks furnish examples of  line ensembles with discrete variants of the Brownian Gibbs property; 
the survey~\cite{Spohn2005} discusses examples of growth processes associated to such line ensembles.
Non-intersecting systems of discrete and continuous Markov processes are connected  to the theory of symmetric functions~\cite{Gessel1990}, including Schur functions, and to Young tableaux~\cite{AvM2005}, and have interpretations in terms of two-dimensional Yang-Mills theory~\cite{FMS2011}.

In the case of the asymmetric simple exclusion process and the stochastic six-vertex model, begun from step initial conditions, a discrete Gibbs property has been harnessed in~\cite{CorwinDimitrov} to verify the KPZ prediction of a $T^{2/3}$-scaling in space at advanced time $T$. Indeed, these models correspond to Hall-Littlewood processes, which are measures on plane partitions; and the levels lines of these plane partitions form ensembles of non-intersecting paths that verify a discrete Gibbs property. That Gibbs property controls the local behaviour of curves according to the {\em derivative} of the distance between adjacent curves -- and, in this regard, it differs from both the Brownian Gibbs property and its softened, ${\bf H}_t$, version; certain monotonicites that we will review in Lemmas~\ref{l.monotoneone} and~\ref{l.monotoneone} are unavailable, but a weaker form of near monotonicity is demonstrated in~\cite{CorwinDimitrov} as a tool for proving the $T^{2/3}$-scaling.


\section{The article's main results}

In the next three subsections, we present these results. We will also  take the opportunity to present the principal inputs needed for  their proofs, since these inputs serve to explain some of the connections between our results.

\subsection{Brownian bridge regularity for curves in the Airy line ensemble}\label{s.airyle}

The Airy line ensemble was constructed in \cite[Theorem 3.1]{AiryLE}. It is a random collection $\mc{A}: \N \times \R \to \R$  of continuous curves $\mc{A}(j,\cdot)$, indexed by the {\em positive} integers $j \in \N$, defined under a probability measure that we label~$\PP$.  For any finite set $I \subset \R$, we may define the random variable $\mc{A}[I]$ under~$\PP$ to be the point process  on $I \times \R$ given by $\big\{ \big( s,\mc{A}(j,s) \big) : j \in \N  \, , \, s \in I \big\}$. The law of the ensemble~$\mc{A}$ is the unique distribution supported on such collections of continuous curves such that, for each finite $I = \{ t_1,\cdots,t_m \}$, the process $\mc{A}[I]$ is a determinantal point process whose kernel is the extended Airy$_2$ kernel~$K^{{\rm ext}}_2$, specified by
$$
K^{{\rm ext}}_2 \big( s_1,x_1;s_2,x_2  \big)  = \begin{cases}
  \int_0^\infty e^{-\lambda(s_1 - s_2)}   {\rm Ai}\big(x_1+\lambda\big){\rm Ai}\big(x_2+\lambda\big) \dd \lambda \, \, & \textrm{if $s_1 \geq s_2$}
 \, , \\
 - \int_{-\infty}^0 e^{-\lambda(s_1 - s_2)}   {\rm Ai}\big(x_1+\lambda\big){\rm Ai}\big(x_2+\lambda\big) \dd \lambda \, \, & \textrm{if $s_1 < s_2$} \, , 
\end{cases}
$$ 
where ${\rm Ai}:\R \to \R$ is the Airy function. The Airy line ensemble's curves are ordered, with $\mc{A}(1,\cdot)$ uppermost.

Under the law~$\PP$, we further define 
\begin{equation}\label{e.lairy}
\mc{L}: \N \times \R \to \R  \,  , \,\,\, \, \mc{L}(i,x) = 2^{-1/2} \big( \mc{A}(i,x) - x^2 \big) \, \, \, \textrm{for $(i,x) \in \N \times \R$} \, .
\end{equation}
The ensemble $\mc{A}$ is stationary, and indeed ergodic~\cite{CorwinSun}, under horizontal shifts. This symmetry is of course lost for $\mc{L}$. However, $\mc{L}$ is a natural object: for example, its top curve $\mc{L}(1,\cdot)$ has the limiting law of a suitably scaled polymer weight profile in such models as geometric LPP~\cite[Theorem 1.2]{Johansson2003}. (The factor of $2^{-1/2}$ used in specifying $\mc{L}$ is employed in order that the curves of $\mc{L}$ locally resemble Brownian motion with a diffusion parameter equal to {\em one}.) 

Inherently related to the construction of the ensemble $\mc{A}$
in~\cite{AiryLE} is the assertion that the ensemble~$\mc{L}$ satisfies the Brownian Gibbs property. From this assertion, it readily follows that, if we
affinely shift the curves in $\mc{L}$ in order to compare them with Brownian bridge, the Radon-Nikodym derivative associated to this comparsion is almost surely finite.
Our first result presents a conclusion about how close this comparison is in terms of a moment bound on this Radon-Nikodym derivative.  
\begin{theorem}\label{t.rnbound}
Letting $K \in \R$ and $\ipdval \geq 1$, we write 
 $\mc{C} =  \mc{C}[K,K+\ipdval]$ for the space of real-valued continuous functions on $[K,K+\ipdval]$ whose endpoint values vanish, endowed with the topology of uniform convergence.  For $k \in \N$, we define 
 $\mc{L}^{[K,K+ \ipdval ]} \big( k, \cdot  \big):[K,K+\ipdval] \to \R$,
$$
\mc{L}^{[K,K+\ipdval]}(k,x) = \mc{L}(k,x) -  (K+\ipdval-x) \ipdval^{-1}  \mc{L}(k,K) - (x-K) \ipdval^{-1} \mc{L}(k,K+\ipdval) \, ,
$$
this being the affine translation of $\mc{L}(k,\cdot)$ that lies in $\mc{C}$. We further write $\mc{B} = \mc{B}^{[K,K+\ipdval]}$ for the law of Brownian bridge $B:[K,K+\ipdval] \to \R$ with $B(K) = B(K+\ipdval) = 0$;
that is, the law of $W(\cdot - K):[K,K+\ipdval] \to \R$ where $W:[0,\ipdval] \to \R$ is  standard Brownian motion conditioned to vanish at time~$\ipdval$.

The distribution of  $\mc{L}^{[K,K+\ipdval]} \big( k, \cdot  \big)$ under~$\PP$, and the law $\mc{B}$, are both supported on the space~$\mc{C}$. Let $f_k = f_{k,K,\ipdval}:\mc{C} \to [0,\infty)$ denote the Radon-Nikodym derivative of the former distribution with respect to the latter. There exists a sequence 
$\big\{ \alpha_k : k \in \N \big\}$, dependent on~$\ipdval$ but not on~$K$,  whose terms lie in $(0,1]$, and which satisfies $\inf \alpha_k^{1/k} > 0$, such that
 the quantity
$$
 j_k =  \int_{\mc{C}} \exp \Big\{ \alpha_k  \big( \log f_k \big)^{6/5} \Big\} \dd \mc{B}  \, ,
$$
which is independent of $K \in \R$, is finite for each $k \in \N$. Specifically, for each $k \in \N$,  $f_k$ has finite $L^p\big( \mc{C},  \mc{B} \big)$-norm for all $p \in [1,\infty)$ (with this $p$-norm being independent of $K \in \R$).
\end{theorem}
(The property that $\inf \alpha_k^{1/k} > 0$ has been stated as a convenient summary of the deduction made in proving Theorem~\ref{t.rnbound}. However, the main interest in the  theorem is probably that it makes an assertion for any given value of $k \in \N$; efforts to take $k$ to infinity using this result should not be expected to  yield sharp conclusions.)

In the next result, the first statement is in essence a restatement of Theorem~\ref{t.rnbound} that seeks to explain its meaning. The second is an example of an inference that can be drawn from this assertion and should be compared to the Brownian bridge probability 
 $\mc{B}^{[0,1]} \big( \sup_{x \in [0,1]}  \vert B ( x ) \vert \geq s  \big)$ lying in the interval $\big[ e^{-2s^2},2e^{-2s^2} \big]$
for all $s > 0$, (a fact that follows from the upcoming Lemma~\ref{l.maxfluc}).
\begin{theorem}\label{t.airytail}
Let $\ipdval \geq 1$.
There exist sequences  $\big\{ \beta_k : k \in \N \big\}$ and  $\big\{ \gamma_k : k \in \N \big\}$, without dependence on $\ipdval$ and satisfying $\limsup \beta_k^{1/k} < \infty$ and  $\liminf \gamma_k^{1/k} > 0$, 
such that the following hold. 
\begin{enumerate}
\item  For each $k \in \N$ and any measurable $A \subset \mc{C}[K,K+\ipdval]$, denote $a = \mc{B}^{[K,K+\ipdval]}(A)$; then the condition 
$a < \gamma_k \wedge \exp \big\{ -\ipdval^3 \big\}$ implies that
$$
\PP \Big( \mc{L}^{[K,K + \ipdval]} \big( k, \cdot  \big) \in A  \Big)  \leq a \cdot
  2103 \, \ipdval^{1/2}    \exp \Big\{ \beta_k \ipdval^2 \big( \log a^{-1} \big)^{5/6} \Big\}   \, .
$$ 
Specifically, this probability is $a  \, \cdot \, \exp  \big\{ (\log a^{-1})^{5/6} O_k(1) \big\}$, where $O_k(1)$ denotes a $k$-dependent term that is independent of $a$.
\item
For each $k \in \N$ and  $s \geq 2^{-1/2} \big( (\log \gamma_k^{-1}) \vee \ipdval^3 \,  + \,  \log 2 \big)^{1/2}$,
$$
 \PP \Big( \sup_{x \in [K,K+\ipdval]}  \big\vert \mc{L}^{[K,K+\ipdval]} ( k, x ) \big\vert \geq s \, \ipdval^{1/2} \Big) \leq 4206  \exp \Big\{ - 2s^2 \big( 1 - 2^{-1/6} \beta_k s^{-1/3} \big) \Big\}
  \, .
$$ 
\end{enumerate}
\end{theorem}

Our final main result concerning KPZ universal objects has two parts: an upper bound of the order of $x^{1/2} \big( \log x^{-1} \big)^{1/2}$ on the modulus of continuity of the  {\rm Airy}$_2$ process on a compact interval, and a local version with a pure $x^{1/2}$ power law.
Since this modulus of continuity is actually enjoyed by Brownian motion, it may be expected that the result in the first part is sharp: for example, that the polylogarithmic power of one-half cannot be improved.
The resampling techniques in this article might be address this question of a lower bound, but we do not pursue this.

\begin{theorem}\label{t.airymodcon}
Let $\mc{A}: \R \to \R$ denote the {\rm Airy}$_2$ process.
\begin{enumerate}
\item Let $I \subset \R$ be a bounded interval.
Let~$S$ denote the supremum of the ratio of 
 $\big\vert \mc{A}(y) - \mc{A}(x) \big\vert$ and $(y-x)^{1/2} \big( \log (y-x)^{-1} \big)^{1/2}$ as the variables $x$ and $y$ vary over~$I$ subject to $x < y$.
 Then $S$ is almost surely finite.
 \item 
There exist positive constants $M$ and $m$
such that, for any
 $x \in \R$, $\e \in (0,1]$  and $K \geq M$,
$$
 \PP \Big(  \sup_{h \in (0,\e)} \big\vert \mc{A}(x + h) - \mc{A}(x) \big\vert  \geq K\e^{1/2} \Big) \leq  M \exp \big\{ - m    K^{3/2}  \big\} \, . 
$$
 \end{enumerate}
\end{theorem}

The theorem's first part may be compared to the  $x^{1/2} \big( \log x^{-1} \big)^{2/3}$ modulus of continuity result~\cite[Theorem~$1.3$]{ModCon}   
which
applies to polymer weight profiles begun from general initial conditions in scaled Brownian LPP.

\subsection{Near geodesics with common endpoints in Brownian last passage percolation}\label{s.neargeod}


We now specify some notation for Brownian LPP, further to that in Section~\ref{s.brlpp}, 
state our principal inference Theorem~\ref{t.neargeod} regarding the scaled behaviour of this model, and 
introduce the associated unscaled and scaled ensembles $L_n$ and $\mc{L}_n^{\rm sc}$.


Recall that Brownian LPP is specified in terms of an ensemble $B:\Z \times \R \to \R$ of independent  two-sided standard Brownian motions $B(k,\cdot):\R\to \R$, $k \in \Z$. 
Let $i,j \in \N$ with $i \leq j$, and let $x > 0$. We have defined the set of staircases  $\staircase_{(0,i) \to (x,j)}$  with endpoints $(0,i)$ and $(x,j)$.
(Negatively indexed $B$-curves will not be used in this paper.)
To each member staircase~$\phi$, we have associated an energy~$E(\phi)$,
by summing the increments of the Brownian motions that run over $\phi$'s horizontal planar segments.

We now introduce notation for collections of mutually avoiding staircases and for the energy associated to them. Let  $n \in \N$, $\ell \in \intint{n}$ and~$x > 0$. 
(The integer interval notation $\intint{n}$ was introduced when Brownian LPP was.)
An $\ell$-tuple of staircases on $[0,x] \times \intint{n}$ is a vector  $\big( \phi_1,\cdots,\phi_\ell \big)$ where
\begin{itemize}
\item $\phi_j \in  \staircase_{(0,j) \to (x,n - \ell + j)}$ for $j \in \intint{\ell}$;
\item and the union of the horizontal planar segments of the $\phi_j$ are pairwise disjoint as $j$ varies. 
\end{itemize}
Let $\staircase^\ell_{(0,1) \to (x,n)}$ denote the set of such $\ell$-tuples. 
Each of the $\ell$ elements of any $\ell$-tuple in $\staircase^\ell_{(0,1) \to (x,n)}$ has an already specified energy. Define the energy $E\big( \phi \big)$ of any $\phi = \big( \phi_1,\cdots,\phi_\ell \big) \in \staircase^\ell_{(0,1) \to (x,n)}$ to be $\sum_{j=1}^\ell E(\phi_j)$.

Define the maximum $\ell$-tuple energy
\begin{equation}\label{e.mell}
 M_n^\ell(x) = \sup \Big\{ E(\phi): \phi \in  \staircase^\ell_{(0,1) \to (x,n)} \Big\} \, .
\end{equation}
When $\ell = 1$, a staircase attaining the maximum is called a geodesic.

\begin{figure}[ht]
\begin{center}
\includegraphics[height=12cm]{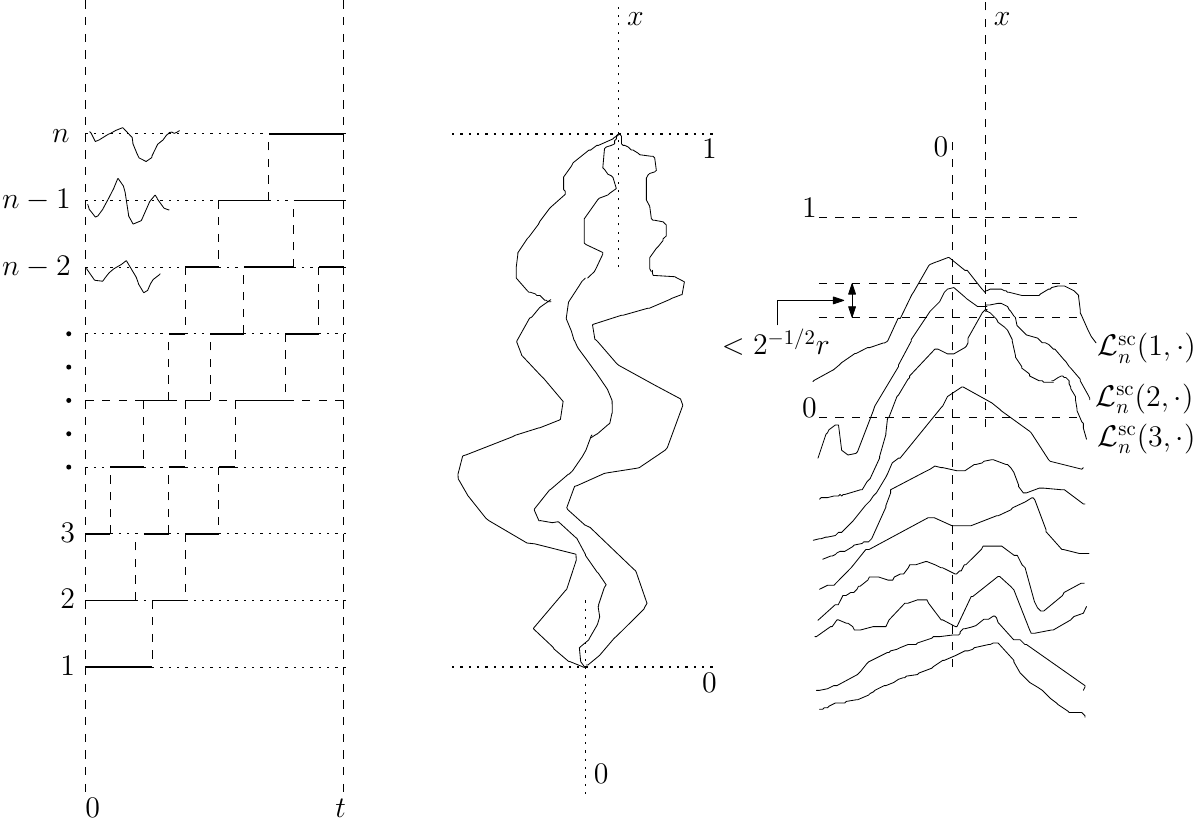}
\caption{{\em Left:} In the formation of the Brownian last passage percolation line ensemble, the maximum triple energy  $M^3_n(t) =\sum_{i=1}^3 L_n(i,t)$, for $t > 0$ given, is formed by considering the sum of the increments on the intervals indicated by horizontal solid black lines of the depicted independent Brownian motions and finding the maximum possible such value. {\em Middle:} Taking $n$ large and setting $t = n + 2n^{2/3}x$ for a given $x \in \R$, we may consider the maximizing triple and depict it after the change of coordinates $(x_1,x_2) \to \big( \tfrac{1}{2} n^{-2/3}(x_1-x_2), x_2 n^{-1} \big)$. If $n$ is high enough, the semi-discrete structure will be indiscernible in the new sketch, and the triple of paths --  a multi-polymer watermelon -- will appear to share the endpoints $(0,0)$ and $(1,x)$. 
{\em Right:} If in the scenario depicted in the middle sketch, the event $\neargeod_{n,3}\big( x,r \big)$ occurs for a given small $r > 0$, the elements in the path triple will have very similar energies, with a collective deficit of $r n^{1/3}$ over the total available in principle. Measuring the deficit in units of $2^{1/2} n^{1/3}$,  a $2^{-1/2}r$-near touch will arise between the top three curves in the scaled ensemble $\mc{L}^\scal_n$ over location~$x$.}
\label{f.threesketches}
\end{center}
\end{figure}

For $n,k \in \N$ with $n \geq k$, as well as $x \in \R$ with $x \geq - n^{1/3}/2$, and $r > 0$, we define
$$
\neargeod_{n,k}\big( x,r \big) = \Big\{  M_n^k \big( n + 2n^{2/3}x \big) \geq k \cdot M_n^1 \big( n + 2n^{2/3}x \big) \, - \, r n^{1/3} \Big\} \, .
$$
As the left and middle sketch in Figure~\ref{f.threesketches} depict, this event entails the existence of a $k$-tuple of disjoint paths whose collective energy is atypically close to the maximum possible, the shortfall being $2^{-1/2 }r$ when expressed in scaled units. The next theorem and corollary determine the first order decay for the probability of this event as $r \searrow 0$. Theorem~\ref{t.neargeod} 
is the official counterpart to the informal Theorem~\ref{t.neargeod.informal}: its proof will depend on the probabilistic apparatus of the jump ensemble to analyse the near-touch ensemble probability depicted in the right sketch of Figure~\ref{f.threesketches}, and the result then plays a foundational role in our broader Brownian LPP study. 
\begin{theorem}\label{t.neargeod}
There exist positive constants $K_0$, $K_1$, $a_0$ and $r_0$ and a positive sequence $\{ \beta_k: k \in \N \}$ with $\limsup \beta_k^{1/k} < \infty$
such that,
for $n,k \in \N$,  $x \in \R$ and $r \in \big(0,(r_0)^{k^2}\big)$ satisfying $k \geq 2$, $n \geq k \vee \, (K_0)^{k^2} \big( \log r^{-1} \big)^{K_0}$ and   $\vert x \vert \leq a_0 n^{1/9}$,
$$
 r^{k^2 - 1} \cdot \exp \big\{ - e^{K_1 k} \big\} \, \leq \, \PP \Big(  
\neargeod_{n,k}\big( x,r \big) \Big) \, \leq \, r^{k^2 - 1} \cdot 
\exp \Big\{ \beta_k \big( \log r^{-1} \big)^{5/6} \Big\} \, .
$$
\end{theorem}
\begin{corollary}\label{c.neargeod}
Let $k \geq 2$ and $x \in \R$. Consider the limit supremum and limit infimum of 
$$
 \big( \log r \big)^{-1} \cdot  \log  \PP \Big(  
\neargeod_{n,k}\big( x,r \big) \Big)
$$
as the limits $n \to \infty$ followed by $r \searrow 0$ are taken. The two limits exist and equal $k^2 - 1$.
Moreover, 
there exist positive constants $K_0$ and $a_0$ such that, for given $k \geq 2$, the limiting value is approached from above and below uniformly in the $r \searrow 0$ limit as the parameters $(n,x)$ vary over the set $\llbracket k  \vee \, (K_0)^{k^2} \big( \log r^{-1} \big)^{K_0} , \infty) \times  [-a_0 n^{1/9},a_0 n^{1/9}]$.  
\end{corollary}
\noindent{\bf Proof.} The result is immediate from Theorem~\ref{t.neargeod}. \qed

\medskip

In order to prove Theorem~\ref{t.neargeod}, we specify associated ensembles of curves. The $n$-indexed Brownian LPP line ensemble 
$L_n: \intint{n} \times \big[ 0 , \infty \big) \to \R$ is defined by insisting that, for $\ell \in \intint{n}$,
\begin{equation}\label{e.mlsum}
 M_n^\ell(x) = \sum_{i=1}^\ell L_n(i,x) \, .
\end{equation}

The {\em scaled} Brownian last passage percolation line ensemble
$$
 \mc{L}^{\rm sc}_{n}: \intint{n} \times \big[- \tfrac{1}{2} n^{1/3} , \infty \big) \to \R 
$$
is then specified
for $(i,x) \in \intint{n} \times \big[- \tfrac{1}{2} n^{1/3} , \infty \big)$ by setting
\begin{equation}\label{e.scl}
 \mc{L}^{\rm sc}_{n} \big( i,x \big)
  = 2^{-1/2} n^{-1/3} \Big( L_n \big( i, n + 2n^{2/3} x \big) - 2 n - 2 n^{2/3}x \Big) \, .
\end{equation}
The reader may wish to glance ahead to  Figure~\ref{f.lnscaledunscaled} to see a depiction of the ensembles $L_n$ and $\mc{L}^\scal_n$.

The right sketch in Figure~\ref{f.threesketches} shows how we will characterise the event~$\neargeod_{n,k}\big( x,r \big)$
in terms of the behaviour of the ensemble~$\mc{L}^\scal_n$.
To gauge this behaviour, we will certainly need to understand basic tightness and parabolic curvature properties of the sequence $\big\{ \mc{L}^\scal_n : n \in \N \big\}$.

In order to obtain these properties, we will  rely on the identity in law between $L_n$ and Dyson Brownian motion which we now recall.
For $n \in \N$, the $n$-indexed Dyson Brownian motion line ensemble $\dysonbm_n : \intint{n} \times [0,\infty) \to \R$ may formally be regarded as a system of $n$ mutually avoiding Brownian motions, each begun at the origin at time zero. 
This initial condition creates a singular conditioning, which may be interpreted using the theory of the Doob-$h$ transform (which is discussed in the text~\cite{RogersWilliamsVolOne}). Indeed, the function $h(\bar{x}) = \prod_{1 \leq i < j \leq n} (x_i- x_j)$, for $\bar{x} = (x_1,\cdots,x_n) \in \R^n$, is a strictly positive harmonic function for $n$-dimensional Brownian motion $B:\intint{n} \times [0,\infty) \to \R \cup \{ c \}$ that departs to a cemetery state $c$ on exit from the Weyl chamber $\big\{ x \in \R^n: x_1>\cdots >  x_n \big\}$. We may define $\dysonbm_n$ to be the Doob $h$-transform of $B$ with entrance point the origin in $\R^n$.

The next result is \cite[Theorem 7]{O'ConnellYor}.
\begin{proposition}\label{p.brlppdbm}
For $n \in \N$, the Brownian last passage percolation and Dyson Brownian motion line ensembles $L_n$ and $\dysonbm_n$, each of which maps $\intint{n} \times [0,\infty)$ to $\R$, are equal in law.
\end{proposition}
The equality in law between the top curves $L_n(1,\cdot)$ and $\dysonbm_n(1,\cdot)$
was earlier derived by~\cite{Baryshnikov} and~\cite{GTW}.

\subsection{Deviation inequalities for GUE eigenvalues at or near the edge}\label{s.deviation}

For $n \geq 1$ and $\sigma^2 \in (0,\infty)$, the Gaussian unitary ensemble~${\rm GUE}_n(\sigma^2)$ with entry variance $\sigma^2$ is the law on random $n \times n$ Hermitian matrices whose upper triangular entries are complex normal random variables $X_{i,j} \sim N(0,\sigma^2/2) + {\rm{i}\mkern1mu}  N(0,\sigma^2/2)$ 
and whose diagonal entries are real  normal random variables $X_{i,i} \sim N(0,\sigma^2)$. 

Under a probability measure that we denote by $\PP$, Hermitian Brownian motion ${\rm HBM}_n$ is the random process defined on $[0,\infty)$ and valued in $n \times n$ Hermitian matrices, whose  upper triangular entries ${\rm HBM}_n(t)_{ij}$ equal $B_{i,j;1}(t/2) + {\rm{i}\mkern1mu} B_{i,j;2}(t/2)$,
and whose diagonal entries ${\rm HBM}_n(t)_{ii}$  equal $B_{i,i}(t)$, where the $B$-processes are independent real-valued standard Brownian motions. Note that the law of~${\rm HBM}_n(t)$ equals 
${\rm GUE}_n(t)$ for any $t \in (0,\infty)$.

Let $\lambda_n:\intint{n} \times [0,\infty) \to \R$
be such that  $\big\{ \lambda_n(k,t): 1 \leq k \leq n \big\}$ is  a decreasing list of the eigenvalues of  the random matrix ${\rm HBM}_n(t)$.

The next result \cite[Theorem 3]{Grabiner} indicates the relevance of the Gaussian unitary ensemble for our study.

\begin{proposition}\label{p.hbmgrabiner}
For any $n \in \N$, the Hermitian Brownian motion eigenvalue process $\lambda_n$ and Dyson Brownian motion $\dysonbm_n$ are equal in law.
\end{proposition}


Upper bounds on the upper and lower tail of the scaled top ${\rm GUE}$  eigenvalue are known. 
For the lower tail, \cite[(5.16)]{Ledoux} states that there exist constants $C',c' > 0$ such that, for $n \geq 1$ and $\epsilon \in [n^{-2/3},1]$,
\begin{equation}\label{e.ledoux}
 \PP \big( \lambda_n\big( 1 , (4n)^{-1} \big) \leq 1 - \epsilon \big) \leq C' \exp \big\{ - c' n \epsilon^{3/2}  \big\} \, .
\end{equation}

For the upper tail,  there exist by  \cite[Proposition 1]{Aubrun} constants $\hat{C}$ and $\hat{c}$ such that, for $n \in \N$ and $t \geq 0$,
\begin{equation}\label{e.aubrun}
 \PP \Big( \lambda_n\big( 1, (4n)^{-1} \big) \geq 1 + t \Big) \leq \hat{C} \exp \big\{ - \hat{c} n t^{3/2}  \big\} \, .
\end{equation}
(The variance choice of $(4n)^{-1}$ is made in order that the edge eignvalues be close to unity. Regarding~(\ref{e.ledoux}), the lower bound in the condition $\epsilon \in [n^{-2/3},1]$ is insignificant, since an increase in the value of $C'$ can replace this condition by $\epsilon \in [0,1]$.)

The ordering of {\rm GUE} eigenvalues implies that~(\ref{e.ledoux}) holds for  $\lambda_n\big( i, (4n)^{-1} \big)$ for all eigenvalue indices $i \in \intint{n}$. 
Our next result extends the companion bound~(\ref{e.aubrun}) to eigenvalues other than the first near the top of the spectrum.
\begin{theorem}\label{t.othereigen}
There exists $n_0 \in \N$ such that, if  $(n,k) \in \N^2$ satisfies $n \geq k \vee n_0$, and $t \in [0,2^{1/2}n^{-11/18}]$, then
$$
\PP \Big( \lambda_n(k,\big(4n)^{-1}\big) \leq 1 - t \Big) \leq H_k \exp \big\{ - h_k n t^{3/2} \big\} \, ,
$$
where these $k$-indexed positive constants satisfy $\limsup H^{1/k^2}_k < \infty$ and $\liminf h^{1/k}_k > 0$. 
\end{theorem}

Expressed in terms of the scaled Brownian LPP ensemble $\mc{L}_n^\scal$, our point of departure (\ref{e.ledoux}) and~(\ref{e.aubrun}) is a quantified form of one-point tightness for the top curve $\mc{L}_n^\scal(1,\cdot)$. Theorem~\ref{t.othereigen} makes this inference for higher index curves. Its proof  is a variant of this paper's central theme: the analysis of such ensembles as  $\mc{L}_n^\scal$ by means of their Brownian Gibbs resampling property; the particular Brownian Gibbs argument used in this proof is a variant of one used to prove one of the key technical propositions in the construction~\cite{KPZLE} of the KPZ line ensemble.  

\chapter{Brownian Gibbs ensembles: definition and statements}

In this chapter, Brownian Gibbs ensembles are defined and an important regularity property that they may enjoy is identified. This regularity property has been designed to capture the behaviour of the scaled Brownian LPP ensembles  $\mc{L}_n^\scal$. Our principal results will be recast in terms of regular Brownian Gibbs ensembles. 
Our work will then be reduced to proving these reformulations in the later chapters.

Why do we reformulate results using regular ensembles? We have chosen to present several of our main results in terms of KPZ universal structures, such as the Airy line ensemble, because such statements make for what we hope is interesting headline reading. In the broader study of scaled Brownian LPP made in the three companion papers, however,
we will always retreat to the finite~$n$ prelimit. This prelimit is in essence described by many copies of the ensemble  $\mc{L}_n^\scal$. For this reason, it is the reformulated assertions, rather than the original ones, that will be quoted in these accompanying papers. 

The chapter has three sections. 
The first formulates the Brownian Gibbs property and the notion of a regular Brownian Gibbs ensemble and contains the assertion that the ensemble $\mc{L}_n^\scal$ is indeed regular. The second reformulates the article's main conclusions in these terms, and also presents some further tools using the regular ensemble framework. The third presents some general tools that are useful when working with the Brownian Gibbs property.
This third section ends with an overview of the structure of the remainder of the paper.

\section{Preliminaries: bridge ensembles and the Brownian Gibbs property}

\begin{figure}[ht]
\begin{center}
\includegraphics[height=12cm]{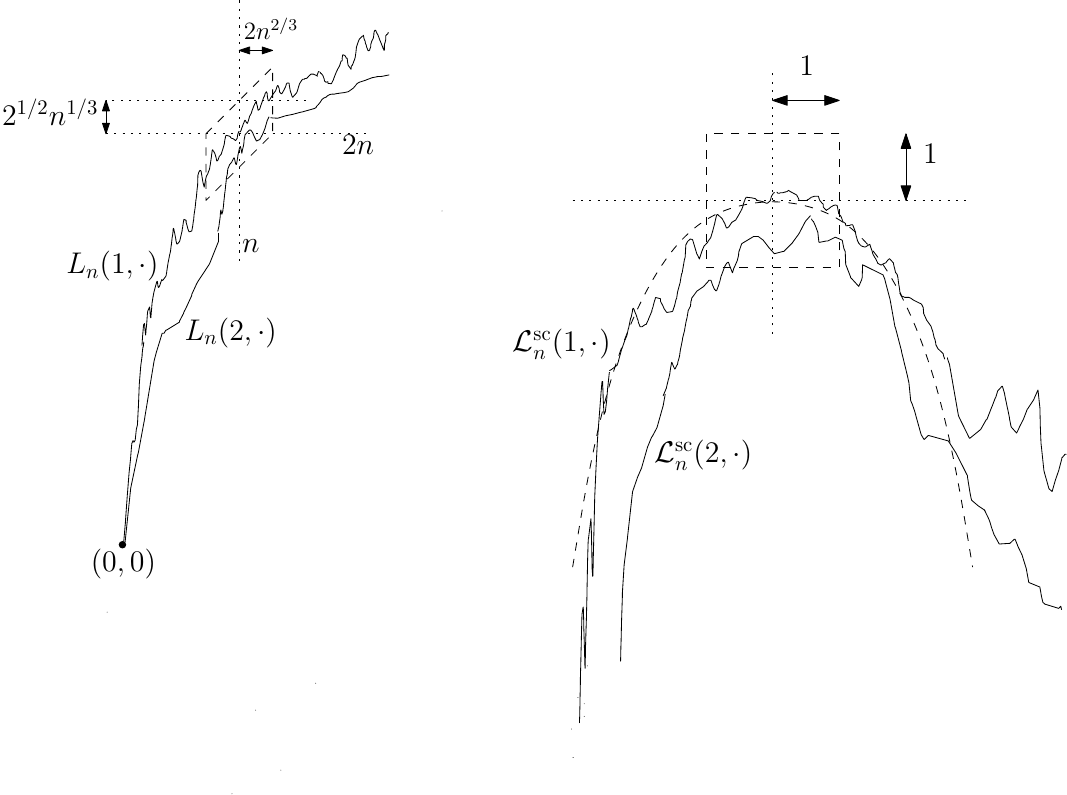}
\caption{A schematic depiction of the highest two curves in the unscaled and scaled Brownian last passage percolation line ensembles for a high value of~$n \in \N$. The dashed parallelogram on the left transforms into the dashed square on the right under the affine change of coordinates in~(\ref{e.scl}) by which $\mc{L}_n^{\scal}$ is formed from $L_n$. The dashed curve on the right equals $-Q(x) = - 2^{-1/2} x^2$. When $n$ is large, the highest curves in $\mc{L}^\scal_n$ surge upwards until, far to the left of the origin, they join a bounded channel about this parabola, which they then typically inhabit until far beyond the origin on the right, when the parabola drops away beneath them.}
\label{f.lnscaledunscaled}
\end{center}
\end{figure}

In this section, we introduce some basic notation and then discuss the Brownian Gibbs property enjoyed by such ensembles of random curves as the unscaled and scaled Brownian last passage percolation ensembles $L_n$ and $\mc{L}_n^\scal$ (for any given $n \in \N$). 
The scaled ensembles $\mc{L}^\scal_n$ form a sequence in~$n \in \N$ whose elements verify the Brownian Gibbs property and have curves that are locally Brownian but globally  parabolic, with the highest curves (which are those of lowest index) typically at unit-order height near the origin: see Figure~\ref{f.lnscaledunscaled}.  
Also in this section, we present a definition (of {\em regular} Brownian Gibbs ensembles) obtaining to such ensembles as the elements of the sequence $\big\{ \mc{L}^\scal_n : n \in \N \big\}$ that captures such behaviour; the definition has been introduced to provide a format for the statements of our results that is flexible enough to be convenient for the applications in this paper and in our study of scaled Brownian LPP.

\subsection{General notation}

We write $\N = \big\{ 1,2,\cdots \}$. For $i,j \in \N$ with $i \leq j$, the integer interval $\big\{ k \in \N: i \leq k \leq j \big\}$ will be (and has already been) denoted by $\llbracket i,j \rrbracket$.

For $x,y \in \R$, we write $x \wedge y = \min \{ x,y \}$ and $x \vee y = \max \{ x,y \}$. Division will take precedence over $\wedge$ (and $\vee$) so that $x \wedge y/2 = x \wedge (y/2)$.

Let $k \in \N$. We use an overhead bar notation, as in $\bar{x} \in \R^k$, to indicate a $k$-vector. We write $\bar{0} = (0,\cdots,0) \in \R^k$ and  $\vecint = (k-1,k-2,\cdots,1,0)$.
 A $k$-vector $\bar{x}=(x_{1},\ldots,x_{k} \big) \in \R^k$ is called a $k$-decreasing list if $x_i > x_{i+1}$ for $1 \leq i \leq k-1$. We write $\Rkle \subseteq \R^k$ for the set of $k$-decreasing lists. When $I \subseteq \R$, we write $I^k_>$ for the set of such lists each of whose components lies in $I$. We also write $I^k_\geq$ in the case that equality between consecutive elements is permitted.

For $\bar{x} \in \R^k$, $s > 0$ and $A \subseteq \R^k$, $\bar{x} + A = \big\{ \bar{x} + \bar{a}: \bar{a} \in A \big\}$ and $s \cdot A = \big\{ s\bar{a}: \bar{a} \in A \big\}$. 

With $\mathbb{Q}$ a probability measure and $A$ and $B$ events, we will write $\mathbb{Q}(A,B)$ for $\mathbb{Q}(A \cap B)$ and $\mathbb{Q}(\cdot \vert A,B)$ for  $\mathbb{Q}(\cdot \vert A \cap B)$. The event complementary to $A$ will be denoted by $\neg \, A$ or $A^c$.

\subsection{Mutually avoiding Brownian bridges: some definitions}

$\empty$

Throughout, Brownian motion and bridge have diffusion parameter one.

\begin{definition}\label{WBdef}
Let $k\in \N$, $a,b \in \R$ with $a < b$, and $\bar{x},\bar{y}\in \Rkle$. Write $\mc{B}_{k;\bar{x},\bar{y}}^{[a,b]}$ for the law of the
ensemble $B:\intint{k} \times [a,b] \to \R$ whose constituent curves $B(i,\cdot):[a,b] \to \R$, $i \in \intint{k}$, are 
independent Brownian bridges that satisfy $B(i,a) = x_i$ and $B(i,b) = y_i$. 

Let  $f:[a,b] \to \R \cup \{-\infty\}$ be a measurable function such that $x_k>f(a)$ and $y_k>f(b)$.
Define the non-touching event on an interval $A\subset [a,b]$ with lower boundary data $f$ by
\begin{equation*}
\notouch_f^A  =  \Big\{ \, \textrm{for all } x \in A \, ,  \,  B(i,x) > B(j,x) \textrm{ whenever } 1\leq i<j\leq k \, \, , \, \, \textrm{ and  $B(k,x) > f(x)$} \Big\} \, . 
\end{equation*}
We omit the subscript $f$ in the case that it equals $-\infty$ throughout $[a,b]$ (and thus plays no role). We omit the superscript $A$ in the case that $A = [a,b]$. With this convention, the event $\notouch$ always imposes {\em internal} curve avoidance, but only imposes {\em external} avoidance of the lower boundary condition when this is indicated in the subscript.

The conditional measure $\mc{B}^{[a,b]}_{k;\bar{x},\bar{y}} \big( \cdot \big\vert \notouch_f \big)$  is the  {\it mutually avoiding Brownian bridge ensemble on the interval $[a,b]$ with entrance data  $\vectoro{x}$, exit data $\vectoro{y}$ and lower boundary condition $f$}. 

We will occasionally refer to the acceptance probability, which is defined to be  $\mc{B}_{k;\bar{x},\bar{y}}^{[a,b]}\big(\notouch_f  \big)$.

\end{definition}

\subsection{Line ensembles and the Brownian Gibbs property}

The law $\mc{B}^{[a,b]}_{k;\bar{x},\bar{y}} \big( \cdot \big\vert \notouch \big)$ is a prototypical example of a line ensemble that enjoys the Brownian Gibbs property that we now define: this ensemble verifies the next definition with $\Sigma = \intint{k}$ and $\Lambda = [a,b]$.

\begin{definition}\label{maindef}
Let $\Sigma$ be an interval of $\Z$, and let $\Lambda$ be an interval of $\R$. Note that $\Sigma$ may be infinite and $\Lambda$ may have infinite length.
Consider the set $X$ of continuous functions $f:\Sigma\times \Lambda \rightarrow \R$ endowed with the topology of uniform convergence on compact subsets of $\Sigma\times\Lambda$. Let $\mathscr{C}$ denote the $\sigma$-algebra  generated by Borel sets in $X$.

A {\it $\Sigma$-indexed line ensemble} $\mathcal{L}$ is a random variable defined on a probability space $(\Omega,\mathscr{B},\PP)$, taking values in $X$ such that $\mathcal{L}$ is a $(\mathscr{B},\mathscr{C})$-measurable function. We view $\mathcal{L}$ as a collection of random continuous curves (despite using the word `line' to refer to them), indexed by $\Sigma$, each of which maps $\Lambda$ into $\R$. We will slightly abuse notation and write $\mathcal{L}:\Sigma\times \Lambda \rightarrow \R$, even though it is not $\mathcal{L}$ which is such a function, but rather $\mathcal{L}(\omega)$ for each $\omega \in \Omega$. 
Given a $\Sigma$-indexed line ensemble $\mathcal{L}$, and a sequence of such ensembles $\big\{ \mathcal{L}_n: n \in \N \big\}$,
we say that the sequence converges {\em weakly as a line ensemble}
to $\mathcal{L}$ if the measure on $(X,\mathscr{C})$ induced by $\mathcal{L}_n$ weak-* converges as $n \to \infty$ to the measure induced by $\mathcal{L}$.  This means that, for all bounded continuous functionals $f$, $\int {\rm d} \PP(\omega) f\big( \mathcal{L}_n(\omega) \big) \to \int {\rm d} \PP(\omega) f\big(\mathcal{L}(\omega) \big)$ as $n \to \infty$. A line ensemble is {\em ordered} if, for all $i,j \in \Sigma$, $i<j$, $\mathcal{L}(i,x)>\mathcal{L}(j,x)$ for all $x \in \Lambda$.
Note for example that $\mc{B}^{[a,b]}_{k;\bar{x},\bar{y}} \big( \cdot \big\vert \notouch \big)$ is ordered. Naturally,
statements such as this are understood as being asserted almost surely with respect to $\PP$.

We also mention that we will sometimes omit the term `line', so that `ensemble' is a synomym of `line ensemble'.
\end{definition}

We turn now to formulating the Brownian Gibbs property. 
\begin{definition}\label{maindefBGP}
For $n \in \N$ and an interval $\Lambda \subseteq \R$,  let $k \in \intint{n}$ and $a,b \in \Lambda$, with $a <b$. Set $f=\mathcal{L}_{k+1}$ unless $k=n$ when $f \equiv -\infty$. Write $D_{k;a,b} = \intint{k} \times (a,b)$ and $D_{K;a,b}^c = (\intint{n} \times \Lambda) \setminus  D_{K;a,b}$. Suppose that an ordered line ensemble $\mc{L}:\intint{n} \times \Lambda \to \R$ has the property that, for all such choices of $k$, $a$ and $b$,
\begin{equation*}
\textrm{Law}\Big(\mathcal{L} \big\vert_{D_{k;a,b}} \textrm{conditional on } \mathcal{L} \big\vert_{D_{K;a,b}^c}\Big) = \mc{B}^{[a,b]}_{k;\bar{x},\bar{y}} \big( \cdot \big\vert \notouch_f \big) \, ,
\end{equation*}
where on the right-hand side the entrance data $\bar{x}$ is taken equal to $\big(\mathcal{L}(1,a),\cdots,\mathcal{L}(k,a)\big)$, the exit data $\bar{y}$ to $\big(\mathcal{L}(1,b),\cdots,\mathcal{L}(k,b) \big)$, and where it is understood that the restriction of $f$ to $[a,b]$ is considered. Then the ensemble is said to have the Brownian Gibbs  property, or, more simply, to be Brownian Gibbs.

\end{definition}

\noindent{\em Remark.} 
For any $k$-curve ensemble $E: \intint{k} \times [a,b] \to \R$,
we use the bar $k$-vector notation in the form $\ovbar{E}(x) = \big( E(1,x),\cdots,E(k,x) \big)$ for $x \in [a,b]$.

\subsection{Regular Brownian Gibbs line ensembles}

Let $I$ be a closed interval in the real line and let $n \in \N$. In an extension of terminology, a line ensemble $\mc{L}_n: \intint{n} \times I \to \R$ that is ordered when restricted to $\intint{n} \times {\rm int}(I)$, where ${\rm int}(I)$ denotes the interior of $I$, is said to have the   Brownian Gibbs property if the ordered line ensemble given by the above restriction does so.  With this usage, the curves in a Brownian Gibbs line ensemble may all be equal, to zero for example, at say the left-hand endpoint of their common domain of definition.

\begin{definition}\label{d.regularsequence} 
%
Consider a Brownian Gibbs ensemble that has the form
$$
\mc{L}: \intint{\nmac} \times \big[ - \xnmac , \infty \big) \to \R   \, ,
$$
and which is defined on a probability space under the law~$\PP$.
The number $\nmac = \nmac(\mathcal{L})$ of ensemble curves and the absolute value $\xnmac$ of the finite endpoint may take any values in $\N$ and $[0,\infty)$.
(In fact, we may also take $\xnmac = \infty$, except that we would then take the domain of definition of $\mc{L}_n$ to be $\intint{n} \times \R$.)

Consider a given three-component vector $\bar\phimac \in (0,\infty)^3$, and 
let $\rsC$ and $\rsc$ be two positive constants. The ensemble $\mc{L}$
is said to be $(\bar\phimac,\rsc,\rsC)$-regular if the following conditions are satisfied.
\begin{enumerate}
\item {\bf Endpoint escape.} $\xnmac \geq  \rsc \nmac^{\phimac_1}$.
\item {\bf One-point lower tail.} If $z \geq -\xnmac$ satisfies $\vert z \vert \leq \rsc \nmac^{\phimac_2}$, then
$$
\PP \Big( \mc{L} \big( 1,z\big) + 2^{-1/2}  z^2 \leq - s \Big) \leq \rsC \exp \big\{ - \rsc s^{3/2} \big\}
$$
for all $s \in \big[1, \nmac^{\phimac_3} \big]$.
\item {\bf One-point upper tail.}  If $z \geq -\xnmac$ satisfies $\vert z \vert \leq \rsc \nmac^{\phimac_2}$, then
$$
\PP \Big( \mc{L} \big( 1,z\big) +  2^{-1/2} z^2 \geq  s \Big) \leq \rsC \exp \big\{ - \rsc s^{3/2} \big\}
$$
for all $s \in [1, \infty)$.
\end{enumerate}
\end{definition}

We will refer to these regular ensemble conditions in the form $\rmreg(i)$, for $i \in \intint{3}$.
Throughout the paper, we will reserve the symbols $c$ and $C$ for usage in denoting the constant parameters $(c,C)$ in a regular ensemble.

In alluding to an ensemble $\mc{L}$ that verifies these conditions, we will reserve the symbol $n$ for the number of its curves, and write $n$ in the subscript. When the notation $\mc{L}_n$ is used, reference is always being made to an $n$-curve ensemble that is  $(\bar\phimac,\rsc,\rsC)$-regular -- about the values of these five parameters, more momentarily.

The right sketch of Figure~\ref{f.threesketches}
illustrated how the near geodesic event 
$\neargeod_{n,k}( x,r )$ will be characterized by  near touching at distance of order~$r$ of the~$k$ highest curves in $\mc{L}^\scal_n$ at $x$.
As such,
the next proposition shows the relevance of the regular ensemble definition for the purpose of proving the near geodesic Theorem~\ref{t.neargeod}. 
\begin{proposition}\label{p.lereg}
 Let $\bar\phimac = \big(1/3,1/9,1/3 \big)$.
There exist choices of the positive constants $\rsc$ and $\rsC$ such that each of the scaled Brownian LPP line ensembles
$\mc{L}_n^{\scal}: \intint{n} \times \big[- \tfrac{1}{2} n^{1/3}  , \infty \big) \to \R$, $n \in \N$,  
is $\big(\bar\phimac,\rsc,\rsC\big)$-regular.
 \end{proposition}

The values of the three-vector  $\bar\phimac = \big(1/3,1/9,1/3 \big)$ and the constants $c$ and $C$ are set in all applications according to Proposition~\ref{p.lereg}. The definition of  $(\bar\phimac,\rsc,\rsC)$-regular has been recorded to permit other uses; this choice is however arguable, because it is not clear such uses will arise -- for this reason, we emphasise that $\bar\phimac$ equals $\big(1/3,1/9,1/3 \big)$ consistently in this work, and thus the reader is encouraged to focus on this value in interpreting statements, and to bear in mind that the value $\phimac_2 = 1/9$ dictates the value of certain quantities that appear on several occasions, such as $\phimac_1/4 \, \wedge \, \phimac_2/4 \, \wedge \,  \phimac_3/2$, which equals $1/36$ in applications. 

\subsection{The Airy line ensemble after parabolic shift enjoys the Brownian Gibbs property}

Note that the Airy line ensemble  $\mathcal{A}:\N\times \R\rightarrow \R$  recalled in Section~\ref{s.airyle} is a continuous ordered $\N$-indexed line ensemble.
The next proposition is asserted by the principal Theorem~$3.1$ of \cite{AiryLE}.

\begin{proposition}\label{p.airy}
The ensemble $\mathcal{L}:\N\times\R\rightarrow \R$ associated to the Airy line ensemble via~(\ref{e.lairy})  has the Brownian Gibbs property.
\end{proposition}

The ensemble $\mc{L} : \N \times \R \to \R$  is also regular,  as we shall demonstrate  in the proof of Theorem~\ref{t.rnbound} in the next section.

\section{Statements of principal results concerning regular ensembles}\label{s.ensemblesequences}

We now gather together our main conclusions concerning regular  ensembles. At the same time, we give the proofs of our main theorems, which follow directly from the new statements.

\subsection{The one-point lower tail estimate extends to curves of higher index}

Here we reformulate the eigenvalue bounds Theorem~\ref{t.othereigen}
as an assertion that extends  the one-point lower tail condition \rmreg(2) concerning the lowest indexed curve in a regular Brownian Gibbs ensemble to a counterpart result for other ensemble curves.

To present our result to this effect, Proposition~\ref{p.othercurves},
we associate to any $(\bar\phimac,\rsc,\rsC)$-regular ensemble two sequences $\big\{ C_k : k \in \N \big\}$ and $\big\{ c_k: k \in \N \big\}$. The dependence of the sequences on the regular ensemble  is communicated through the constant parameters $(c,C)$. The sequences are specified by setting, for each $k \geq 2$, 
  \begin{equation}\label{e.formere}
 \formerE_k = \max \Big\{  10 \cdot 20^{k-1} 5^{k/2} \Big( \tfrac{10}{3 - 2^{3/2}} \Big)^{k(k-1)/2} C \, , \, e^{c/2} \Big\} 
 \end{equation}
as well as $C_1 = 140 C$; and 
\begin{equation}\label{e.littlec}
 c_k =   \big( (3 - 2^{3/2})^{3/2} 2^{-1} 5^{-3/2} \big)^{k-1} c_1 \, ,
\end{equation}
with  $c_1 = 2^{-5/2} c \wedge 1/8$. 

Note that the sequences satisfy
 $\limsup \formerE_k^{1/k^2} < \infty$ and $\liminf c_k^{1/k} > 0$.

This usage of $C_k$ and $c_k$ will be made consistently in the paper.
  
\begin{proposition}\label{p.othercurves}
Let $\bar\phimac \in (0,\infty)^3$. 
Set $\delta = \phimac_1/2 \wedge  \phimac_2/2 \wedge \phimac_3$.  For any $n$-curve $(\bar\phimac,\rsc,\rsC)$-regular  ensemble $\mc{L}_n$, and  for each $k \in \N$,  the estimate
$$
\PP \Big( \, \mc{L}_n\big( k,x\big) + 2^{-1/2} x^2 \leq - s \, \Big) \leq \formerE_k \exp \big\{ - c_k s^{3/2} \big\}
$$
is valid provided that  $n \geq k 
  \vee  (c/3)^{-2(\phimac_1 \wedge \phimac_2)^{-1}} \vee  6^{2/\delta}$,
 $\vert x \vert \leq \rsc/2 \cdot n^{\delta}$ and $s \in \big[0, 2 n^{\delta} \big]$. 
\end{proposition}
We will prove the proposition by establishing the stronger Proposition~\ref{p.strongothercurves}, which asserts  a similar estimate for the lowest value adopted by $\mc{L}_n(k,\cdot)$ on a compact interval. 

Theorem~\ref{t.othereigen} is a consequence of Proposition~\ref{p.othercurves}. Indeed, the next result is a restatement of the theorem with a more explicit description of its constant parameters. 
\begin{corollary}\label{c.othereigen}
Suppose that   $(n,k) \in \N^2$ satisfies 
 $n \geq k 
  \vee  (c/3)^{-2\phimac_2^{-1}} \vee  6^{2/\delta}$, where $\phimac_2 = 1/9$ and $\delta = 1/{18}$. If $s \in [0,2^{1/2}n^{1/18-2/3}]$, then
$$
\PP \Big( \lambda_n\big( k,(4n)^{-1} \big) \leq 1 - s \Big) \leq \formerE_k \exp \big\{ - 2^{3/2} c_k n s^{3/2} \big\} \, ,
$$
where the parameters $(\rsc,\rsC)$ used to determine the sequences~(\ref{e.formere}) and~(\ref{e.littlec}) are
supplied by Proposition~\ref{p.lereg}.
\end{corollary}
\noindent{\bf Proof.} 
By Propositions~\ref{p.brlppdbm} and~\ref{p.hbmgrabiner} as well as Brownian scaling, $2n \lambda_n\big(k,(4n)^{-1}\big)$ has the law of~$L_n(k,n)$. Thus,~(\ref{e.scl}) implies that
$\lambda_k^n - 1$ has the law of $2^{-1/2} n^{-2/3} \mc{L}^\scal_n(k,0)$. Since $\mc{L}^\scal_n$ is a
 $\big(\bar\phimac,\rsc,\rsC\big)$-regular ensemble, with  $\bar\phimac = (1/3,1/9,1/3)$,
 by Proposition~\ref{p.lereg}, we may apply Proposition~\ref{p.othercurves} to find that, for $r \in [0,2n^{1/18}]$,
$$
\PP \big( \mc{L}^\scal_n(k,0) \leq - r \big) \, \leq \,  \formerE_k \exp \big\{ - c_k r^{3/2} \big\} \, .
$$
Setting $r = 2^{1/2} n^{2/3} s$, so that the left-hand side equals $\PP \big( \lambda_k^n - 1 \leq - s \big)$, we obtain the corollary. \qed

\subsection{The close encounter of several curves at the edge}
Now we reformulate Theorem~\ref{t.neargeod} in the language of regular Brownian Gibbs ensembles.

Let $m,k \in \N$ satisfy $m \geq k$, and let $a,b \in \R$ satisfy $b > a$. If $E:\intint{m} \times [a,b] \to \R$ is an ordered ensemble, $x \in [a,b]$ and $\phi > 0$, we let $\close\big( k; E , x , \phi \big)$ denote the event that 
$$
 E\big(1,x \big) \leq E\big(k,x\big) + \phi \, ,
 $$
 also adopting this definition when $[a,b]$ is replaced by $[a,\infty)$.

First we reformulate the upper bound in Theorem~\ref{t.neargeod}.
\begin{theorem}\label{t.airynt}
For $\bar\phimac \in (0,\infty)^3$ and $C,c > 0$, there exists a sequence $\big\{ \const = \const ( c ) : k \geq 2 \big\}$
satisfying $\sup_{k \geq 2} \const^{1/k} < \infty$
 such that the following holds. Let 
$$
\mc{L}_n:\intint{n} \times \big[-\xnmac,\infty\big) \to \R  
$$ 
be a $(\bar\phimac,\rsc,\rsC)$-regular ensemble defined under the law~$\PP$. Let $k \in \N$, $k \geq 2$, and $\e > 0$ satisfy
$\e <  2^{-1/2} (17)^{-1/k} C_k^{-1/k} \const^{-1}$. For $n \in \N$ satisfying $n \geq k 
  \vee  (c/6)^{-2(\phimac_1 \wedge \phimac_2)^{-1}} \vee  6^{2/\delta}$ (where $\delta = \phimac_1/2 \wedge \phimac_2/2 \wedge \phimac_3$),
and  
    $n^{\phimac_1/4 \, \wedge \, \phimac_2/4 \, \wedge \, \phimac_3/2}   \geq    \big( \rsc/4 \wedge 2^{1/2} \big)^{-1} 2^{1/3} \const \big( \log \e^{-1} \big)^{1/3}$,
\begin{enumerate}
\item the bound
$$
  \PP \Big( \close\big( k ; \mc{L}_n , \fa , \e \big)  \Big) 
  \leq     10^6 \exp \Big\{  8842 k^{7/2}  \const^{5/2} \big( \log \e^{-1} \big)^{5/6} \Big\}  \, \e^{k^2-1} 
$$ 
holds 
for any given $\fa \in c/2 \cdot n^{\phimac_1 \wedge \phimac_2} \cdot [-1,1]$; 
\item and the bound
\begin{eqnarray*}
 & & \PP \Big(  \exists \, x \in \R \, , \, \vert x - y \vert \leq \tfrac{1}{4} \const \big( \log \e^{-1} \big)^{1/3}:  \close\big( k ; \mc{L}_n , x , \e \big) \Big) \\
 & \leq &
\e^{k^2 - 3} \cdot 10^{46} \,  2^{6k^2} \const^{18} \exp \Big\{  8844 k^{7/2} \const^{5/2} \big( \log \e^{-1} \big)^{5/6} \Big\}  
\end{eqnarray*}
also holds whenever $\vert y \vert \leq c/2 \cdot n^{\phimac_1 \wedge \phimac_2}$.
\end{enumerate}
\end{theorem}
It should be noted that only the first part of Theorem~\ref{t.airynt} is needed for the purpose of proving Theorem~\ref{t.neargeod}.
The second part leads directly to the inference that  the probability that there exists a value of $x$ in a given compact interval 
such that  
$\neargeod_{n,k;(0,0)}^{(x,1)} \big(   \eta \big)$ occurs  has probability   $\eta^{k^2 - 3 + o(1)}$; this holds
for given $k \in \N$, uniformly for  $\eta$ small,  $n$ high and $\vert x \vert = o(n^{1/9})$.  Theorem~\ref{t.airynt}(2) is included because it may prove to be a useful tool which anyway follows rather directly from our results and techniques.
 
The lower bound in Theorem~\ref{t.neargeod} is reformulated next.
\begin{theorem}\label{t.airynt.lb}
 Let
$\bar\phimac \in (0,\infty)^3$ and $C,c > 0$, and let  
$$
\mc{L}_n:\intint{n} \times \big[-\xnmac,\infty\big) \to \R  
$$ 
be a $(\bar\phimac,\rsc,\rsC)$-regular ensemble  defined under the law~$\PP$. Set $\delta = \phimac_1/2 \wedge \phimac_2/2 \wedge \phimac_3$, and specify $s_k = \big( 8 \cdot 2 c_k^{-1} \log (28C_k) \big)^{2/3} \vee 2^{5/2}$ in terms of the sequences~(\ref{e.formere}) and~(\ref{e.littlec}) for $k \geq 2$. Then, for $\e \in \big( 0, k^{-2} s_k^{-1}/4 \big)$,
$$
  \PP \Big( \close\big( k ; \mc{L}_n , x , \e \big) \Big) 
  \geq    e^{- 52 s_k^2 k^3} \, \e^{k^2-1} 
$$ 
whenever $\vert x \vert \leq c/2 \cdot n^{\phimac_1 \wedge \phimac_2}$, $k \geq 2$ and  $n \geq k \vee (\rsc/6)^{-2(\phimac_1 \wedge \phimac_2)^{-1}} \vee 6^{2/\delta} \vee (s_k/2)^{1/\delta}$. 
\end{theorem}

\noindent{\bf Proof of Theorem~\ref{t.neargeod}.}
We consider $k \geq 2$, $n \geq k$ and $x \geq - n^{1/3}/2$. Note first that, by~(\ref{e.mlsum}),
$$
 k M_n^1 \big( n + 2 n^{2/3} x \big) - M_n^k  \big( n + 2 n^{2/3} x \big)  = \sum_{i=2}^k \Big( L_n\big(1,  n + 2 n^{2/3} x \big) -  L_n\big(i,  n + 2 n^{2/3} x \big)  \Big) \, .
$$
For any $r \geq 0$, the event $\neargeod_{n,k}(x,r)$
that the displayed quantity is at most $r n^{1/3}$ may by~(\ref{e.scl}) be expressed in the form
$$
 \Big\{  \sum_{i=2}^k \big( \mc{L}^\scal_n (1,x) - \mc{L}^\scal_n (i,x) \big) \leq 2^{-1/2} r  \Big\} \, .
$$
Since the ensemble $\mc{L}^\scal_n$ is ordered, we find that
$$
  \mc{L}^\scal_n (1,x) - \mc{L}^\scal_n (k,x) \leq   \sum_{i=2}^k \big( \mc{L}^\scal_n (1,x) - \mc{L}^\scal_n (i,x) \big) \leq ( k-1 ) \big(  \mc{L}^\scal_n (1,x) - \mc{L}^\scal_n (k,x)  \big) \, .
$$
Thus,
\begin{eqnarray*}
 & &  \PP \Big( \mc{L}^\scal_n \big(1,x\big) - \mc{L}^\scal_n \big(k,x\big) \leq (k-1)^{-1} 2^{-1/2}r  \Big) \\
  & \leq & \PP \Big( \neargeod_{n,k}\big(x,r\big) \Big) \leq \PP \Big( \mc{L}^\scal_n \big(1,x\big) - \mc{L}^\scal_n \big(k,x\big) \leq 2^{-1/2}r  \Big)  \, .
\end{eqnarray*}
By Proposition~\ref{p.lereg}, we are able to apply Theorem~\ref{t.airynt}(1) and Theorem~\ref{t.airynt.lb} to $\mc{L}^\scal_n$ to bound the right and left-hand terms here above and below, whenever $\vert x \vert \leq c/2 \cdot n^{1/9}$, where the constant parameter $\rsc > 0$ is provided by this regular ensemble. In bounding the left-hand term below, we use $\limsup s_k^{1/k} < \infty$ in order to obtain the doubly-exponential-in-$k$ correction factor on the left-hand side of the conclusion of Theorem~\ref{t.neargeod}.
In summary, the bounds on $\PP \big( \neargeod_{n,k}(x,r) \big)$ stated in Theorem~\ref{t.neargeod} are obtained by this means. \qed

\medskip

\noindent{\em Remark.} It is in fact technically incorrect to state that  $\mc{L}^\scal_n$ is ordered, because its curves are equal to zero at the left endpoint $-n^{1/3}/2$, so the ordering is not strict. This point is irrelevant for the application just made, as it will be subsequently, and we will call similar ensembles ordered later in the paper.

\subsection{Brownian bridge regularity for affinely translated ensemble curves}

Here we reformulate 
Theorem~\ref{t.airytail}(1).
\begin{theorem}\label{t.airytail.ln}
Let $\bar\phimac \in (0,\infty)^3$, $(\rsc,\rsC) \in (0,\infty)^2$ and $n \in \N$.
Suppose that $\mc{L}_n$ is an $n$-curve $\big(\bar\phimac,\rsc,\rsC\big)$-regular ensemble. Let $\ipdval \geq 1$ denote a parameter. Setting $\deltapi = \phimac_1 \wedge \phimac_2$, 
let $K \in \R$ satisfy $[K,K+ \ipdval] \subset \rsc/2 \cdot [-n^\deltapi,n^\deltapi]$,  and let $k \in \N$.
Let the parameter $\const = \const(\rsc)$ be specified by Theorem~\ref{t.airynt}, and further set $D_1 = D_2$.
Suppose that 
$n \geq k 
  \vee  (c/3)^{-2(\phimac_1 \wedge \phimac_2)^{-1}} \vee  6^{2/\delta}$,
 where  $\delta = \phimac_1/2 \wedge  \phimac_2/2 \wedge \phimac_3$.

For any measurable $A \subset \mc{C}[K,K+\ipdval]$, write  $a = \mc{B}^{[K,K+\ipdval]}(A)$. Suppose that $a$ satisfies the $k$-dependent upper bound 
$a < e^{-1} \wedge (17)^{-1/k} C_k^{-1/k} \const^{-1} \wedge \exp \big\{ -2^5 \ipdval^3 \const^{-3} \big\}$;
 as well as the $n$-dependent lower bound
\begin{equation}\label{e.alowerbound}
 a \geq \exp \Big\{   - \big( \rsc/2 \wedge 2^{1/2} \big) \const^{-1} n^{3\phimac_1/4 \, \wedge \,  3\phimac_2/4 \, \wedge \,  3\phimac_3/2} \Big\} \, .    
\end{equation} 
Then 
$$
\PP \Big( \mc{L}_n^{[K,K + \ipdval]} \big( k, \cdot  \big) \in A  \Big)  \leq a \cdot
  2103 \, \ipdval^{1/2}    \exp \Big\{ 6395 \, \ipdval^2 k^{7/2}  \const^2 \big( \log a^{-1} \big)^{5/6} \Big\}   \, .
$$ 
Specifically, this probability is $a  \cdot \exp \big\{ (\log a^{-1})^{5/6} O_k(1) \big\}$, where $O_k(1)$ denotes a $k$-dependent term that is independent of $a$.


\end{theorem}

The condition~(\ref{e.alowerbound}) may appear to be a nuisance. However, taking the form $a \geq \exp \big\{ - O(1) n^{1/12} \big\}$ in applications in view of Proposition~\ref{p.lereg}, it is seen to be rather mild; moreover, the condition is vacuous in the important case where $n=\infty$ and $\mc{L}_n$
is set equal to the parabolically shifted Airy line ensemble given by~(\ref{e.lairy}).

Theorems~\ref{t.airytail} and~\ref{t.airytail.ln}'s proofs appear 
in Section~\ref{s.brownianregularity}.

\subsection{Modulus of continuity for ensemble curves}

Here we reformulate the two parts of Theorem~\ref{t.airymodcon}
in the regular ensemble language. The reformulation of the first part, Theorem~\ref{t.weakbound}, will also later permit us to verify
Theorem~\ref{t.airynt}(2) from Theorem~\ref{t.airynt}(1).
\begin{definition}\label{d.klinemod}
For $k \in \N$, $a,b \in \R$ with $a < b$, and an ensemble $E:\intint{k} \times [a,b] \to \R$, define the ensemble's modulus of continuity 
$$
\omega_{k,[a,b]}(E,\delta) = \max_{i \in \intint{k}} \, \sup \Big\{ \big\vert E(i,x+s) - E(i,x) \big\vert :  (x,s) \in [a,b-\delta] \times [0,\delta] \Big\} \, .
$$
\end{definition}

\begin{theorem}\label{t.weakbound}
For $\bar\phimac \in (0,\infty)^3$ and $C,c > 0$, there exists a sequence $\big\{ \const = \const ( c ) : k \geq 2 \big\}$
satisfying $\sup_{k \geq 2} \const^{1/k} < \infty$
 such that the following holds. Let 
$\mc{L}_n:\intint{n} \times \big[-\xnmac,\infty\big) \to \R$ 
be a  $\big(\bar\phimac,\rsc,\rsC\big)$-regular ensemble defined under the law~$\PP$.
For $k \in \N$ with $k \geq 2$, set 
$$
g_k = 64 (k+2)k^{3/4} \const^{3/2} \vee 2^{-5/2} c_k^{1/2} \const^{3/2} \vee 50 k 2^k \, .
$$
Suppose that $k \geq 2$, $\Cwb \geq 96 \vee 118 g_k$ and 
$\e  < e^{-1} \wedge (17)^{-1/k} C_k^{-1/k} (381)^{-1} (k+2)^{2/3}k^{1/2} K^{-2/3}$; 
and further that 
 $n \in \N$ satisfies $n \geq k 
  \vee  (c/3)^{-2(\phimac_1 \wedge \phimac_2)^{-1}} \vee  6^{2/\delta}$ (where $\delta = \phimac_1/2 \wedge \phimac_2/2 \wedge \phimac_3$),
as well as  
\begin{equation}\label{e.epsilonandn}
n^{\phimac_1/4 \, \wedge \, \phimac_2/4 \, \wedge \, \phimac_3/2}   \geq    \big( \rsc/2 \wedge 2^{1/2} \big)^{-1} \tfrac{1}{380} (k+2)^{-2/3}k^{-1/2} K^{2/3} \big( \log \e^{-1} \big)^{1/3} \, .
\end{equation}
 Then 
$$
 \PP \Big( \omega_{k, I}\big(\mc{L}_n, \e \big) > 
 \Cwb \, \e^{1/2} \big( \log \e^{-1} \big)^{1/2}  
   \Big) \, \leq \,  25 k 2^k \, 
\e^{ 10^{-10} c_k (k+2)^{-2}k^{-3/2} K^{2} }  
\, , 
$$
where the interval $I$ is given by $I = \tfrac{1}{4} \const (\log \e^{-1})^{1/3} \cdot [-1,1]$.
\end{theorem}
\noindent{\bf Proof of Theorem~\ref{t.airymodcon}(1).}
To verify this assertion, recall that the Airy$_2$ process equals the uppermost curve $\mc{A}(1,\cdot)$ in the Airy line ensemble and that, as we saw in the proof of Theorem~\ref{t.rnbound}, the parabolically adjusted ensemble $\mc{L}:\N \times \R \to \R$
from~(\ref{e.lairy}) is regular.

Moreover, given the form of Theorem~\ref{t.airymodcon}(1) and the stationarity of $\mc{A}$, it suffices to demonstrate this result's statement for  the process $\mc{L}(1,\cdot): \R \to \R$ and with the interval $I$ chosen to equal $[-1,1]$. 
To do this, let $K > 0$ be large enough to permit the application of Theorem~\ref{t.weakbound} with $k=1$. 
We will first argue that 
 there exists $\e_0 > 0$ of the form $O(1) K^{-2/3}$ such that, for any dyadic scale $2^{-i_0} \leq \e_0$, 
 the supremum of the ratio of 
 $\big\vert \mc{L}(1,y) - \mc{L}(1,x) \big\vert$ and $(y-x)^{1/2} \big( \log (y-x)^{-1} \big)^{1/2}$ as the variables $x$ and $y$ vary over $[-1,1]$ subject to $x < y < x + 2^{-i_0}$
has probability at most $(\e_0)^{g K^2}$ of exceeding $K$.   
This claim is verified by applying Theorem~\ref{t.weakbound} with $n = \infty$ and $k=1$ to the regular ensemble~$\mc{L}$, with $\e$ being set equal to $2^{-i}$ for  all integer values of $i$ exceeding the specified $i_0 \in \N$.
The parameter $g > 0$ is  obtained by a suitable decrease from the value
 $10^{-10} c_k 3^{-2}$, the decrease made to cope with the union bound over dyadic scales. 
 
 From the claim just established, we see that the corresponding supremum, where $x$ and $y$ are instead permitted to vary over $[-1,1]$,
 may exceed $O(1)  K^{5/3}$ with probability at most $O(1)  K^{2/3} (\e_0)^{g K^2}$.
 Since $\e_0 = O(1) K^{-2/3}$, this supremum is certainly seen to be almost surely finite, so that the proof is complete. \qed

Theorem~\ref{t.weakbound} is a significantly less delicate result than Theorem~\ref{t.airytail.ln}, because it treats only fluctuations whose probability one may expect to be of  fast polynomial decay in $\e$. It is not implied by the latter result, however, because of the affine shift used in Theorem~\ref{t.airytail.ln}; the absence of this affine shift gives 
Theorem~\ref{t.weakbound} a practical usefulness that the more refined Theorem~\ref{t.airytail.ln} (in isolation from other bounds) lacks. 
We also draw attention to the $\e$-determined lower bound on $n$ in~(\ref{e.epsilonandn}): the index $n$ must rise in order that short scale fluctuation be understood. 
Theorem~\ref{t.weakbound}  offers an all-scale limiting description in the high~$n$ limit, 
something we took advantage of in the proof of Theorem~\ref{t.airymodcon}(1), but this description is not offered  in a uniform way for varying~$n$.

We now turn to maximal fluctuation near a given point and the reformulation of Theorem~\ref{t.airymodcon}(2), namely 
Theorem~\ref{t.aestimate}. 
In fact, 
Theorem~\ref{t.aestimate} concerns several ensemble curves, and is expressed in terms of the quantity 
$\omega_{k,[x,x+\delta]}(E,\delta)$, which we note is equal to  
$$
\max_{i \in \intint{k}} \, \sup \Big\{ \big\vert E(i,x+s) - E(i,x) \big\vert :  s \in  [0,\delta] \Big\}
$$
under the circumstances specified in Definition~\ref{d.klinemod}.

Theorem~\ref{t.aestimate} does not prove the presumably optimal $\exp \big\{ -c K^2 \big\}$ tail behaviour that we see in Theorem~\ref{t.weakbound}; rather we find a bound of the form~$\exp \big\{ -c K^{3/2} \big\}$. However, the new upper bound, being of the form $\exp \big\{-c K^{3/2} \big\}$ rather than $\e^{c K^2}$, does not require fast polynomial decay in $\e \searrow 0$ for the probability of the fluctuation event in question. 
\begin{theorem}\label{t.aestimate}
For $\bar\phimac \in (0,\infty)^3$, $C,c > 0$ and $n \in \N$, let 
$$
\mc{L}_n:\intint{n} \times \big[-\xnmac,\infty\big) \to \R  
$$ 
be a    $\big(\bar\phimac,\rsc,\rsC\big)$-regular ensemble defined under the law~$\PP$.
If $k \geq 1$, $x \in [-1,1]$, $\e \in (0,1/2)$ and $K \geq 1 \vee 2^{-3} c_k^2   k^{-9/2}  3^{-3} \vee 2^{19/2} k^{1/2}(k+2)$, then 
$$
 \PP \Big( \omega_{k,[x,x+\e]}\big(\mc{L}_n,\e \big) \geq K\e^{1/2} \Big) \leq  \big( 2^{3k/2} \pi^{k}  k \cdot 60     +   14 C_k \big) \exp \big\{ - c_k   k^{-3} 2^{-16}  K^{3/2}  \big\} \, . 
$$
whenever  $n \geq k \vee (\rsc/3)^{-2(\phimac_1 \wedge \phimac_2)^{-1}} \vee 6^{2/\delta} \vee \big( 2^{-8} K (k+2)^{-1} k^{-1/2} \big)^{1/\delta}$ (with $\delta = \phimac_1/2 \wedge \phimac_2/2 \wedge \phimac_3$).
\end{theorem}
\noindent{\bf Proof of Theorem~\ref{t.airymodcon}(2).}
Since the {\rm Airy}$_2$ process may, after the subtraction of a parabola, be embedded as the lowest indexed curve in a regular Brownian Gibbs ensemble with $n = \infty$,
this result follows from Theorem~\ref{t.aestimate}. \qed

Theorem~\ref{t.aestimate}  dispenses with the  $\e$-determined lower bound~(\ref{e.epsilonandn}) on the index $n$ 
that encumbers Theorem~\ref{t.weakbound}. The uniformity of the assertion remains compromised by a lower bound on $n$ that is $K$-determined, however.  
We present one further result, which dispenses with this shortcoming. The new result is not needed to prove the theorems in this article but we include it since it is a potentially useful result that  emerges directly {\em en route} to Theorem~\ref{t.aestimate}.

For $n \in \N$, $x \geq -\xnmac + 2$ and $t > 0$, define 
$$
\boundgood_t(x) \, = \, \bigcap_{x-2 \leq y \leq x+2} \Big\{ \, \mc{L}_n(1,y) + 2^{-1/2} y^2 \leq t \, , \, \mc{L}_n(k+1,y)  + 2^{-1/2} y^2 \geq - t \, \Big\} \, .
$$
\begin{proposition}\label{p.aestimateinference}
For $\bar\phimac \in (0,\infty)^3$, $C,c > 0$ and $n \in \N$, let 
$$
\mc{L}_n:\intint{n} \times \big[-\xnmac,\infty\big) \to \R  
$$ 
be a    $\big(\bar\phimac,\rsc,\rsC\big)$-regular ensemble defined under the law~$\PP$. 
If $k \geq 1$, $\vert x \vert \leq \rsc/2 \cdot n^{\phimac_1 \wedge \phimac_2}$, $\e \in (0,1/2)$ and $K \geq  2^{19/2} k^{1/2}(k+2)$, then, setting $t = 2^{-8} K (k+2)^{-1} k^{-1/2}$, 
$$
 \PP \Big( \omega_{k,[x,x+\e]}\big(\mc{L}_n,\e \big) \geq K\e^{1/2} \, , \, \boundgood_t(x) \Big) \leq 
 2^{3k/2} \pi^{k}  k \cdot 60 K^{-1}  \exp \big\{ - 2^{-12}  K^2 \big\}   
$$
whenever  $n \geq k + 1$.
\end{proposition}
In the bird's eye view tour of scaled Brownian LPP results presented in Section~\ref{s.concept}, we reviewed~\cite[Theorem~$1.1$]{ModCon},  a strong assertion of the one-half power law for polymer weight differences,
in the guise of the informal Theorem~\ref{t.differenceweight.informal}. Proposition~\ref{p.aestimateinference}, and its circumvention of the shortcoming we have mentioned, plays an important role in a proof of~\cite[Theorem~$1.1$]{ModCon}, though in fact this result is derived~\cite{ModCon} using a separate and simple Brownian Gibbs argument.

\section{Some generalities: notation and basic properties of Brownian Gibbs ensembles}\label{s.generalities}

In this  section, some useful general tools concerning Brownian Gibbs ensembles are presented.

\subsection{Some helpful lemmas and basic notation}\label{somehelp}

\subsubsection{Strong Gibbs property}
In order to explain this property of Brownian Gibbs line ensembles, we introduce the concept of a stopping domain.

\begin{definition}\label{defstopdom}
Consider a line ensemble $\mc{L}: \intint{n} \times [a,b]\to \R$. 
For $a<\ell<r<b$, and $k \in \intint{n}$, denote the $\sigma$-algebra generated by $\mathcal{L}$ outside $\intint{k} \times [\ell,r]$ by
\begin{equation*}
\mathcal{F}_{{\rm ext}}\big(k;\ell,r\big) \, = \, \sigma\Big\{ \, \mc{L} \textrm{ on } \intint{k} \times \big( [a,\ell] \cup [r,b]  \big) \, \, , \textrm{ and } \mathcal{L} \textrm{ on }  \llbracket k+1,n \rrbracket \times [a,b] \, \Big\} \, .
\end{equation*}

The random subset $\intint{k} \times \big( \mathfrak{l},\mathfrak{r} \big)$ of the domain $\intint{n} \times [a,b]$ of $\mc{L}$ is called a {\it stopping domain} if, for all $\ell<r$,
\begin{equation*}
\{\mathfrak{l} \leq \ell , \mathfrak{r}\geq r\} \in \mathcal{F}_{{\rm ext}}\big(k;\ell,r\big) \, .
\end{equation*}
In other words, the domain is determined by the information outside of it.
\end{definition}

We will make use of a version of the strong Markov property where the concept of stopping domain introduced in Definition \ref{defstopdom} plays the role of stopping time.

Let $C^k(\ell,r)$ denote the set of functions $f = (f_1,\ldots, f_k) : \intint{k} \times [\ell,r] \to \R$ with each $f_i:[\ell,r]\to \R$ continuous. Define
\begin{equation*}
C^k = \left\{ (\ell,r,f): \ell<r \textrm{ and } f \in C^k(\ell,r)\right\}.
\end{equation*}
Let $bC^k$ denote the set of Borel measurable functions from $C^k\to \R$.

\begin{lemma}\label{l.stronggibbslemma}
Consider a line ensemble $\mc{L}: \intint{n} \times [a,b]\to \R$ with the Brownian Gibbs property.
 Write $\PP$ and $\EE [\cdot]$ for the probability measure and expectation associated to $\mathcal{L}$.
 Fix $k \in \intint{n}$. For all stopping domains $\intint{k} \times \big( L,R \big)$, the following {\it strong Brownian Gibbs property} holds: for all $F\in bC^k$,  $\PP$ almost surely,
\begin{equation*}\label{strongeqn}
\EE\bigg[ \, F \Big( L,R,\mathcal{L}\big\vert_{\intint{k} \times (L,R)}  \Big) \, \bigg\vert \, \mathcal{F}_{ext}(k;L,R) \, \bigg] \, = \, \mc{B}_{k;\bar{x},\bar{y}}^{[L,R]} \Big[ \, F(L,R,B) \, \Big\vert \, \notouch_f \, \Big] \, ,
\end{equation*}
where $\bar{x} = \{\mathcal{L}(i,L)\}_{i=1}^{k}$, $\bar{y} = \{\mathcal{L}(i,R)\}_{i=1}^{k}$, $f(\cdot)=\mathcal{L}_{k+1}(\cdot)$ (or $-\infty$ if $k=n$).
On the right-hand side, a notational abuse is adopted under which 
$ \mc{B}_{k;\bar{x},\bar{y}}^{[L,R]} \big[ \cdot \, \big\vert \, \notouch_f \big]$ denotes conditional expectation with respect to this conditional measure.
\end{lemma}
This lemma is \cite[Lemma 2.5]{AiryLE}.
 The lemma's message is that the conditional law of a line ensemble inside a stopping domain is dictated by the domain's boundary data via the mutually avoiding Brownian bridge measure specified by this data.

\subsubsection{Monotonicity results}
The next two  lemmas, \cite[Lemmas 2.6 and 2.7]{AiryLE}, state two simple but important monotonicities exhibited by mutually avoiding Brownian bridge ensembles.

\begin{lemma}\label{l.monotoneone}
Fix $k\in \N$, $a<b$ and two measurable functions $f,g:[a,b]\rightarrow \R\cup\{-\infty\}$ such that for all $s\in [a,b]$, $f(s)\leq g(s)$. Let $\bar{x},\bar{y}\in \Rkle$ be two k-decreasing lists such that $x_k \geq g(a)$ and $y_k \geq g(b)$.
Recalling Definition~\ref{WBdef}, set $\PP_{k;f} = \mc{B}^{[a,b]} \big( \cdot \big\vert \notouch_{f} \big)$, and likewise define $\PP_{k;g}$. Then there exists a coupling of $\PP_{k;f}$ and $\PP_{k;g}$ such that almost surely $B_{k;f}(i,x) \leq B_{k;g}(i,x)$ for all $(i,x) \in \intint{k} \times [a,b]$.
\end{lemma}

\begin{lemma}\label{l.monotonetwo}
Fix $k\in \N$, $a<b$, a measurable function $f:[a,b]\rightarrow \R\cup\{-\infty\}$ and a measurable set $A\subseteq [a,b]$. Consider two pairs of $k$-decreasing lists $\bar{x},\bar{y}$ and $\bar{x}',\bar{y}'$ such that $x_k \wedge x'_k\geq f(a)$, $y_k \wedge y'_k \geq f(b)$ and $x_i' \geq x_i$ and $y_i' \geq y_i$ for each $i \in \intint{k}$.
Then the laws $\mc{B}_{k;\bar{x},\bar{y}}^{[a,b]} \big( \cdot \big\vert \notouch_f^A \big)$ and $\mathcal{B}^{[a,b]}_{k;\bar{x}',\bar{y}'} \big( \cdot \big\vert \notouch_f^A \big)$ may be coupled so that, denoting  by $B$ and $B'$ the ensembles defined under the respective measures, $B'(i,x) \geq B(i,x)$ for all $(i,x) \in \intint{k} \times [a,b]$.
\end{lemma}

\subsubsection{Gaussian random variables: notation and tail bounds}

For $k \in \N$, let $\bar{m} = (m_1,\cdots,m_k) \in \R^k$ be a $k$-vector, and let $\sigma^2 \in [0,\infty)$.
We will write $\nu^k_{\bar{m},\sigma^2}$ for the law of a $k$-vector $\ovbar{N} = (N_1,\cdots,N_k)$ of  independent Gaussian random variables, where $N_i$ has mean $m_i$ and variance $\sigma^2$  for $i \in \intint{k}$.

When $A \subseteq B \subseteq \R^k$, we denote by
$\nu^k_{\bar{m},\sigma^2}\big( A \, \big\vert \, B \big)$
the conditional probability under $\nu^k_{\bar{m},\sigma^2}$ given  $\ovbar{N} \in B$ that $\ovbar{N} \in A$. 

In the special case where $k=1$, and where we now take $m \in \R$,
 consider the choice $A = [y_1,y_2]$ and $B = [x,\infty)$ with $x \leq y_1 \leq y_2$. In this case, we will use the further shorthand $\nu^1_{m,\sigma^2}\big( y_1, y_2 \big\vert x,\infty \big)$ to denote $\nu^1_{m,\sigma^2} (A \vert B)$. We also write 
$\nu^1_{m,\sigma^2}(y_1,y_2)$ in place of  $\nu^1_{m,\sigma^2}([y_1,y_2])$. Furthermore, we will often omit the superscript in the notation $\nu^k_{\bar{m},\sigma^2}$ when $k=1$.

When $k=1$, and we take $m \in \R$, we will write 
$g_{m,\sigma^2}:\R \to [0,\infty)$, 
\begin{equation}\label{e.gaussiandensity}
g_{m,\sigma^2}(x) \, = \, (2\pi)^{-1/2} \sigma^{-1} \exp \Big\{ - \tfrac{(x-m)^2}{2\sigma^2} \Big\}
\end{equation} 
for the density of $\nu^1_{m,\sigma^2}$ at $x \in \R$.

Whenever one-dimensional Gaussian tail bounds are needed, the next lemma will be used. We will sometimes omit to mention that the lemma has been invoked.
\begin{lemma}\label{l.gaussiantail}
Let $m \in \R$ and $\sigma^2 \in [0,\infty)$ and let $x \in \R$ with $x \geq m$. Setting $t = (x-m)\sigma^{-1}$, we have
$$
 \nu^1_{m,\sigma^2} \big( x, \infty \big) \leq (2\pi)^{-1/2} \cdot  t^{-1} \exp \big\{ - t^2/2 \big\}   \, ;
$$
and if $t \geq 1$, then 
$$
\nu^1_{m,\sigma^2} \big( x, \infty \big) \geq
 (2\pi)^{-1/2} \cdot (2t)^{-1} \exp \big\{ - t^2/2 \big\} \, . 
$$
\end{lemma}
\noindent{\bf Proof.} The standard bounds 
$$
  (2\pi)^{-1/2}  \tfrac{t}{t^2 + 1} \exp \big\{ - t^2/2 \big\}   \leq  \nu^1_{0,1}(t,\infty) \leq  (2\pi)^{-1/2}  t^{-1} \exp \big\{ - t^2/2 \big\}
$$
for $t \geq 0$ may be found in~\cite[Section 14.8]{Williams}.
Note that $\tfrac{t}{t^2 + 1} \geq (2t)^{-1}$ for $t \geq 1$. \qed

On one occasion, we will also use the next fact.
\begin{lemma}\label{l.monotonenormal}
Let $a \in \R$ and $\sigma^2 > 0$.
For any $r > 0$, the map $\R \to \R: s \to  \frac{\nu_{a,\sigma^2}\big( s + r , \infty \big)}{\nu_{a,\sigma^2}\big( s , \infty \big)}$ is strictly decreasing. That is, for a normal random variable $N$ of any given mean and positive variance, the conditional probability $\PP \big( N \geq s+r \, \big\vert \, N \geq s \big)$ is strictly decreasing in $s \in \R$, for any given $r > 0$. 
\end{lemma}
\noindent{\bf Proof.} Let $N$ denote a  normal random variable of mean $a$ and variance $\sigma^2$. Note that
$$
\log \PP \big( N \geq s+r \, \big\vert \, N \geq s \big)
= \log \int_{s + r}^\infty \exp \big\{ - \tfrac{(x-a)^2}{2\sigma^2} \big\} \dd x \, - \,
 \log \int_{s}^\infty \exp \big\{ - \tfrac{(x-a)^2}{2\sigma^2} \big\} \dd x 
$$
has derivative in $s$ given by 
$$
 \frac{ \int_{s}^\infty \exp \big\{ - \tfrac{(s-a)^2 - (x-a)^2 - r^2}{2\sigma^2} \big\} \Big( - \exp \big\{ -  \tfrac{-(s-a)r}{\sigma^2} \big\}  +  \exp \big\{ - \tfrac{-r(x-a)}{\sigma^2} \big\}  \Big) \dd x}{\int_{s + r}^\infty \exp \big\{ - \tfrac{(s-a)^2}{2\sigma^2} \big\} \dd x  \int_{s}^\infty \exp \big\{ - \tfrac{(s-a)^2}{2\sigma^2} \big\} \dd x } \, .
$$
The integrand is strictly negative for all $x > s$. \qed


\medskip

We typically denote means and variances of Gaussian random variables by $m$ and $\sigma^2$ using subscripts for vector component indices. Each use of this notation is made locally with, we hope, little prospect of confusion between them.

\subsubsection{Standard bridges, general or Brownian}

Let $a,b \in \R$, with $a \leq b$. A continuous function $f:[a,b] \to \R$ with $f(a) = f(b)  = 0$ will be called a standard bridge. The space of such functions will be denoted $\mc{C}_{0,0} \big( [a,b], \R \big)$ (a usage that will supersede the earlier use of $\mc{C}[K,K+\ipdval]$ when $[a,b] = [K,K+\ipdval]$ with $K \in \R$ and $\ipdval \geq 1$).

When $f:[a,b] \to \R$ is general, the standard bridge obtained from $f$ by affine translation will be denoted $f^{[a,b]}$, so that
$$
 f^{[a,b]}(x) = f(x) - \tfrac{b-x}{b-a}f(a) -  \tfrac{x-a}{b-a}f(b) \, \, \, \, \textrm{for $x \in [a,b]$} \, .
$$
This notation extends to line ensembles: for example, if $\mc{L}_n: \intint{n} \times [-\xnmac,\infty) \to \R$ and $[a,b] \subset [-\xnmac,\infty)$, then  $\mc{L}_n^{[a,b]}:\intint{n} \times [a,b] \to \R$ is the line ensemble whose value  $\mc{L}_n^{[a,b]}(i,x)$ is specified by the displayed formula with $f = \mc{L}_n(i,\cdot)$.

When $k=1$ is taken in $\mc{B}_{k;\bar{x},\bar{y}}^{[a,b]}$,
we will write $B:[a,b] \to \R$ in place of $B(1,\cdot)$.
Note then that the process $B^{[a,b]}$ has law $\mc{B}_{1;0,0}^{[a,b]}$. It will be called {\em standard} Brownian bridge.

\subsubsection{Brownian bridge basics}

\begin{lemma}\label{l.brbr}
Let $a,b \in \R$ with $a < b$. For $l \in \N$, let $x_1,\cdots,x_l$ be an increasing sequence of elements of $(a,b)$. Also let $x,y \in \R$. Then the joint distribution under $\mc{B}_{1;x,y}^{[a,b]}$ of $\big( B(x_1),\cdots,B(x_l) \big) \in \R^l$ has density at $(z_1,\cdots,z_l) \in \R^l$ given by
$$
 Z^{-1} \prod_{i=0}^l g_{0,x_{i+1} - x_i}\big( z_{i+1} - z_i \big) \, ,
$$
where we take $x_0 = a$ and $x_{l+1} = b$, as well as $z_0 = x$ and $z_{l+1} = y$; the quantity $Z$ equals $g_{0,b-a}(y-x)$.
\end{lemma}
\noindent{\bf Proof.}
Let $W:[a,b] \to \R$ denote Brownian motion with $W(a) = x$.
The law $\mc{B}_{1;x,y}^{[a,b]}$ is the weak limit as $\phi \searrow 0$ of the law of $W$ given that $W(b) \in (y - \phi,y+\phi)$. 
The claimed formula for the density arises from Bayes' theorem after this limit is taken. \qed

\medskip

The next result is a trivial but sometimes useful consequence.
\begin{corollary}\label{c.brbr}
For choices of parameters given in Lemma~\ref{l.brbr}, let $A_1,\cdots,A_l \subseteq \R$ be a collection of intervals. Also set $A_0 = \{ x \}$ and $A_{l+1} = \{ y \}$. Writing $\mu$ for Lebesgue measure, we have that
$$
 \mc{B}_{1;x,y}^{[a,b]} \Big( \bigcap_{i=1}^k \big\{ B(x_i) \in A_i \big\} \Big) \, \geq \, \frac{1}{g_{0,b-a}(y-x)} \cdot  \prod_{i=1}^l \mu(A_i) \cdot \prod_{i=0}^{l} g_{0,x_{i+1} - x_i}(s_i) \, ,
$$
where $s_i = \sup \big\{ \vert z - z' \vert : z \in A_i , z' \in A_{i+1}  \big\}$.  
\end{corollary}
When applying the bound, also note that $g_{0,b-a}(y-x)^{-1} \geq \big(2\pi(b-a) \big)^{1/2}$.

\subsection{Bounds on maximum fluctuation of bridge ensembles}

\subsubsection{Maximum fluctuation of standard Brownian bridge}
\begin{lemma}\label{l.maxfluc}
Let $a,b \in \R$ with $b > a$. For any $h \in \R$ and $r > 0$,
$$
\mc{B}_{1;h,h}^{[a,b]} \bigg( \sup_{x \in [a,b]} B(x) \geq h + r \bigg) \, = \, \exp \Big\{ - 2 \tfrac{r^2}{b-a} \Big\} \, .
$$
Equality also holds when the condition $\inf_{x \in [a,b]} B(x)  \leq h - r$ is  considered instead.
\end{lemma}
\noindent{\bf Proof.}
By Brownian scaling and symmetry, the second statement reduces to the first and the first to the case where $a = 0$, $b=1$ and $h=0$.  The result then follows from equation~(3.40) in 
\cite[Chapter 4]{KaratzasShreve}. \qed
\subsubsection{Maximum fluctuation of mutually avoiding Brownian bridges}
\begin{lemma}\label{l.notouchfluc}
Let $k \in \N$, $\bar{x},\bar{y} \in \R^k_>$ and $[a,b] \subset \R$. Then, for any $r > 0$,
$$
 \mc{B}_{k;\bar{x},\bar{y}}^{[a,b]} \Big( \inf_{x \in [a,b]} B(k,x) < x_k \wedge y_k \, - \, \sqrt{2} (b-a)^{1/2} \big( k - 1 + r \big) \, \Big\vert \, \notouch^{[a,b]} \Big)
 \leq \big( 1 - 2 e^{-1} \big)^{-k}  e^{-4r^2} \, .
$$
\end{lemma}
\noindent{\bf Proof.}
Specify $\bar{z} \in \R^k_>$ so that $z_i = x_i \wedge y_i - \sqrt{2} (b-a)^{1/2}(i-1)$. Note that $z_i$ is at most both $x_i$ and $y_i$, and also that $z_i - z_{i-1} \leq - \sqrt{2}(b-a)^{1/2}$.
By Lemma~\ref{l.monotonetwo}, the probability in question may only increase if $\bar{x}$ and $\bar{y}$ are replaced by $\bar{z}$. Under the law~$\mc{B}_{k;\bar{z},\bar{z}}^{[a,b]}$, each curve is a vertical displacement of standard Brownian bridge, and there is a clear route to the mutual avoidance event $\notouch^{[a,b]}$: it is achieved if each of the standard bridges has a supremum in absolute value of at most $2^{-1/2} (b-a)^{1/2}$. Each bridge
has $\mc{B}_{k;\bar{z},\bar{z}}^{[a,b]}$-probability at least $1 - 2e^{-1}$ of being constrained in this way by Lemma~\ref{l.maxfluc}. Bayes' theorem, and another use of Lemma~\ref{l.maxfluc}, yields the result. \qed

\subsection{Some basic properties of regular ensembles}

\subsubsection{Basic parabolic symmetry of regular ensembles}
Let $\para:\R \to \R$ denote the parabola $\para(x) = 2^{-1/2} x^2$, and
let $l:\R^2 \to \R$ be given by $l(x,y) = - 2^{-1/2}y^2 - 2^{1/2}y(x-y)$. Note that $x \to l(x,y)$
is the tangent line of the parabola $x \to - \para(x)$ at the point $\big(y,-\para(y)\big)$. Note also that, for any $x,y \in \R$,
\begin{equation}\label{e.plp}
\para(x) = - l(x,y) + \para(x-y) \, .
\end{equation}

For $\xnmac \geq 0$, consider a regular ensemble $\mc{L}_n:\intint{n} \times [-\xnmac,\infty) \to \R$. 
For any $y_n > - \xnmac$, define $\lshift_{n,y_n}:\intint{n} \times [-\xnmac - y_n,\infty) \to \R$ to be the shifted ensemble given by 
$$
\lshift_{n,y_n}(i,x) = \mc{L}_{n}(i,x + y_n) - l(x+y_n,y_n)  \, .
$$
By~(\ref{e.plp}), $\lshift_{n,y_n}(i,x)$ equals
$\big( \mc{L}_n(i,x + y_n) + \para(x+y_n) \big) - \para(x)$. In the case of a regular ensemble~$\mc{L}_n$,  the last expression consists of the bracketed  process, for which the influence of curvature has been cancelled, at least when $x + y_n = n^{o(1)}$, to which is added the basic parabolic decay term.

\begin{lemma}\label{l.parabolicinv}
Let $\bar\phimac \in (0,\infty)^3$, $\rsc, \rsC > 0$ and $n \in \N$.
Suppose that  $\mc{L}_n:\intint{n} \times [-\xnmac,\infty) \to \R$
is  a    $\big(\bar\phimac,\rsc,\rsC\big)$-regular ensemble.
Set $\deltapi =  \phimac_1 \wedge \phimac_2$.
Whenever $y_n  \in \R$
satisfies $\vert y_n \vert \leq \rsc/2 \cdot n^\deltapi$, the ensemble $\lshift_{n,y_n}$ is    $\big(\bar\phimac,\rsc/2,\rsC\big)$-regular.
\end{lemma}
\noindent{\bf Proof.} Note that the domain of definition of $\lshift_{n,y_n}$ has left-hand endpoint $-\xnmac - y_n$ at most $- \rsc n^{\phimac_1} + \rsc/2 \cdot n^\deltapi \leq - \rsc/2 \cdot n^{\phimac_1}$. Thus~$\rmreg(1)$ holds.

Note that the parabolically adjusted random variable
$\lshift_{n,y_n}\big(1,z\big) + \para(z)$
 with which $\rmreg(2)$ and $(3)$ are concerned equals 
 $\mc{L}_n\big(1,z+y_n\big) - l(z+y_n,y_n) + \para(z) = \mc{L}_n\big(1,z+y_n\big) - Q(z+y_n)$ in view of (\ref{e.plp}).
Thus  $\rmreg(2)$ and $(3)$ for $\lshift_{n,y_n}$ follow from their counterparts for $\mc{L}_n$, for values of $\vert z \vert$ at most $\rsc n^{\phimac_2} - y_n$ and thus when $\vert z \vert \leq \rsc/2 \cdot n^{\phimac_2}$.  \qed

\subsubsection{The upper tail of the maximum of a Brownian Gibbs ensemble top curve}

We now present a result similar to the `no big max' Lemma~$5.1$ of \cite{AiryLE}. Our proposition roughly asserts that the upper tail estimate $\exp \{ - O(t^{3/2}) \}$ for the one-point law in $\rmreg(3)$ remains valid when we instead consider the maximum top curve value. 

\begin{proposition}\label{p.nobigmax}
Let $\bar\phimac \in (0,\infty)^3$, $\rsc, \rsC > 0$ and $n \in \N$.
Suppose that  $\mc{L}_n:\intint{n} \times [-\xnmac,\infty) \to \R$
is  a    $\big(\bar\phimac,\rsc,\rsC\big)$-regular ensemble. For $r \in \big[0,\rsc/2 \cdot n^{\phimac_1 \wedge \phimac_2}\big]$, $t \in \big[ 2^{7/2} , 2 n^{\phimac_3} \big]$ and $n \geq (2c)^{-2(\phimac_1 \wedge \phimac_2)^{-1}}$,
$$
\PP \Big( \sup_{x \in [-r,r]} \big( \mc{L}_n ( 1,x ) + 2^{-1/2}x^2 \big) \geq t \Big) \leq  (r + 1) \cdot  6  \rsC \exp \big\{ - 2^{-9/2} \rsc  t^{3/2} \big\} \, . 
$$
\end{proposition}
In fact, we now take the opportunity to record a version of this proposition which is slightly more flexible, in which the interval $[-r,r]$ may be centred elsewhere than at the origin. This more general version will not be applied in the present article but is valuable in some of the accompanying Brownian LPP papers, and its proof is immediate given the apparatus we have developed.

\begin{proposition}\label{p.nobigmax.gen}
Let $\bar\phimac \in (0,\infty)^3$, $\rsc, \rsC > 0$ and $n \in \N$.
Suppose that  $\mc{L}_n:\intint{n} \times [-\xnmac,\infty) \to \R$
is  a    $\big(\bar\phimac,\rsc,\rsC\big)$-regular ensemble. 
For  $\vert y \vert \leq c/2 \cdot n^{\phimac_1 \wedge \phimac_2}$, $r \in \big[0,\rsc/4 \cdot n^{\phimac_1 \wedge \phimac_2}\big]$, $t \in \big[ 2^{7/2} , 2 n^{\phimac_3} \big]$ and $n \geq c^{-2(\phimac_1 \wedge \phimac_2)^{-1}}$,
$$
\PP \Big( \sup_{x \in [y-r,y+r]} \big( \mc{L}_n ( 1,x ) + 2^{-1/2}x^2 \big) \geq t \Big) \leq  (r + 1) \cdot  6  \rsC \exp \big\{ - 2^{-11/2} \rsc  t^{3/2} \big\} \, . 
$$
\end{proposition}
\noindent{\bf Proof.} A consequence of Proposition~\ref{p.nobigmax} and the parabolic invariance Lemma~\ref{l.parabolicinv}. \qed

As we will shortly demonstrate, Proposition~\ref{p.nobigmax} reduces, also  via Lemma~\ref{l.parabolicinv}, to the next result, which concerns the top curve supremum over a bounded interval. The proof of Lemma~\ref{l.nobigmax} is an early illustration of the utility of the Brownian Gibbs property, showing how when a Brownian Gibbs ensemble attains a very high value at a possibly exceptional time, it in fact typically adopts rather high values on a neighbourhood of nearby times; thus, such behaviour may occur only with probability dictated by the one-point upper tail axiom~$\rmreg(3)$. 
\begin{lemma}\label{l.nobigmax}
Let $\bar\phimac \in (0,\infty)^3$, $\rsc, \rsC > 0$ and $n \in \N$.
Suppose that  $\mc{L}_n:\intint{n} \times [-\xnmac,\infty) \to \R$
is  a    $\big(\bar\phimac,\rsc,\rsC\big)$-regular ensemble.
For $t \in \big[ 2^{7/2} , 2 n^{\phimac_3} \big]$ and $n \geq \big(2c^{-1} \big)^{(\phimac_1 \wedge\phimac_2)^{-1}}$,
$$
\PP \Big( \sup_{x \in [-2,2]} \mc{L}_n \big( 1,x \big) \geq t \Big) \leq   6  \rsC \exp \big\{ - \rsc  t^{3/2}/8 \big\} \, . 
$$
\end{lemma}
\noindent{\bf Proof of Proposition~\ref{p.nobigmax}.} 
Note that
$$
 \sup_{x \in [-r,r]} \big( \mc{L}_n ( 1,x ) + 2^{-1/2}x^2 \big)  \leq  \max \sup_{x \in [m-1,m+1]} \big( \mc{L}_n ( 1,x ) + 2^{-1/2}x^2 \big)   \, ,
$$
where the maximum is taken over even integers $m$ with $\vert m \vert \leq r$. The supremum on the right-hand side equals
$$
\sup_{x \in [-1,1]} \big( \lshift_{n,m} ( 1,x ) + 2^{-1/2}x^2 \big) \, .
$$
As such, we will now be able to use Lemma~\ref{l.parabolicinv} to reduce Proposition~\ref{p.nobigmax} to Lemma~\ref{l.nobigmax}, 
using a union bound to find that the probability in the proposition is at most an  $(r + 1)\textsuperscript{st}$ multiple of 
$$
\PP \Big( \sup_{x \in [-1,1]} \big( \lshift_{n,m} ( 1,x ) + 2^{-1/2}x^2 \big) \geq t  \Big) \leq 
\PP \Big( \sup_{x \in [-1,1]} \lshift_{n,m} ( 1,x )  \geq t/2  \Big)  \, ,
$$
where we used 
$t \geq 2^{5/2}$. We now apply Lemma~\ref{l.nobigmax} to bound this quantity, noting from Lemma~\ref{l.parabolicinv} that the shifted ensembles $\lshift_{n,m}$ over the range $\vert m \vert \leq r$
are $\bar\phimac$-regular, because $r \leq \rsc/2 \cdot n^{\phimac_1 \wedge\phimac_2}$.
This completes the proof of Proposition~\ref{p.nobigmax}. \qed

\medskip

\noindent{\bf Proof of Lemma~\ref{l.nobigmax}.} 
We will prove the lemma in the guise
\begin{equation}\label{e.zeroelnsup}
 \PP \Big( \sup_{x \in [-2,0]} \mc{L}_n \big( 1,x \big) \geq t  \Big) \vee \PP \Big( \sup_{x \in [0,2]} \mc{L}_n \big( 1,x \big) \geq t  \Big) \leq  3  \rsC \exp \big\{ - \rsc   t^{3/2}/8 \big\}
\end{equation}
for $t \geq 2^{7/2}$.

We prove~(\ref{e.zeroelnsup}) only for the latter probability, the two proofs being the same.  Let $\chi$ denote the supremum of values $x \in [0,2]$ for which $\mc{L}_n \big( 1,x \big) \geq t$, with $\chi = -\infty$ if no such point exists. It is our aim then to bound $\PP \big( \chi \in [0,2] \big)$. For $r > 0$,
let $\up_{1}^{-2}(r)$ denote the event that $\mc{L}_n(1,-2) \geq - 2^{3/2}  - r$, so that the one-point lower tail bound $\rmreg(2)$ implies that 
\begin{equation}\label{e.negupbound}
 \PP \Big(  \neg \,  \up_{1}^{-2}\big( r \Big) \Big) \leq C \exp \big\{ - c r^{3/2} \big\}
\end{equation}
provided that $r \in [1,n^{\phimac_3}]$,
since $c n^{\phimac_2} \geq 2$.
Note that $\{ 1 \} \times [-2, \chi \vee -2]$ is a stopping domain. We now consider the choice $r = t/2 - 2^{3/2}$, this value lying in the interval $[1,n^{\phimac_3}]$ due to our assumption on~$t$. 
By the strong Gibbs Lemma~\ref{l.stronggibbslemma} and the monotonicity Lemma~\ref{l.monotoneone},
the conditional distribution of $\mc{L}_n\big( 1,\cdot \big):[-2,\chi] \to \R$ under $\PP$ given $\big\{ \chi \in [0,2] \big\} \cap \up_{1}^{-2}\big( t/2 - 2^{3/2} \big)$ stochastically dominates $\mc{B}_{1; -t/2, t}^{[-2,\chi]}$. Since zero is closer to $\chi$ than it is to $-2$, the probability that the latter bridge at zero exceeds $\tfrac{1}{2} \big( t - t/2 \big) = t/4$ is at least one-half. Thus,
$$
 \PP \Big( \mc{L}_n\big( 1,0 \big) \geq t/4 \,  \Big\vert \, \chi \in [0,2] \, , \, \up_{1}^{-2}\big(  t/2 - 2^{3/2} \big) \Big) \geq \tfrac{1}{2} \, ,
$$
from which we find that
$$
 \PP \big( \chi \in [0,2]  \big) \leq 2 \, \PP \Big( \mc{L}_n\big( 1,0 \big) \geq t/4 \Big) \, + \, \PP \Big(  \neg \, \up_{1}^{-2}\big(  t/2 - 2^{3/2} \big) \Big) \, .
$$
Applying \rmreg(3) using $2 \geq c n^{1/9}$ and (\ref{e.negupbound}), 
$$
 \PP \big( \chi \in [0,2]  \big) \leq 2  \rsC \exp \big\{ - \rsc  2^{-3} t^{3/2} \big\} \, + \, \rsC \exp \big\{ - \rsc (t/2 - 2^{3/2})^{3/2} \big\}  \, .  
$$ 
Since $t \geq 2^{7/2}$, we obtain~(\ref{e.zeroelnsup}). \qed

 \subsubsection{The collapse far from the origin  of the top curve in a regular ensemble}
 
 The three $\rmreg$ axioms make no comment on the behaviour of a regular ensemble's curves at very high distances from the origin. However, the conditions imply a rapid decay for the curves' value in this region, though not one dictated by the parabolic curvature seen at shorter range.
 
\begin{proposition}\label{p.collapsenearinfinity}
Let $\bar\phimac \in (0,\infty)^3$, $\rsc, \rsC > 0$ and $n \in \N$.
Suppose that  $\mc{L}_n:\intint{n} \times [-\xnmac,\infty) \to \R$
is  a    $\big(\bar\phimac,\rsc,\rsC\big)$-regular ensemble.

 For $\eta \in (0,\rsc]$, let
$\ell = \ell_\eta:\R \to \R$ denote the even function which is affine on $[0,\infty)$ and has gradient $ - 5 \cdot 2^{-3/2} \eta \nmac^{1/9}$ on this interval, and which satisfies
$\ell(\eta \nmac^{1/9}) = \big( - 2^{-1/2} + 2^{-5/2} \big) \eta^2 \nmac^{2/9}$. If  $\nmac \geq 2^{45/4} \rsc^{-9}$, then  
\begin{eqnarray*}
 & & \PP \Big( \mc{L} \big(1,z\big) > \ell(z) \, \, \textrm{for some} \, \,  z \in \big[ - \xnmac , \infty \big) \setminus \big[ - \eta  \nmac^{1/9} , \eta \nmac^{1/9} \big] \Big) \\
  & \leq &  
6C \exp \Big\{ - c \eta^3  2^{-15/4}   \nmac^{1/3} \Big\} \, .
\end{eqnarray*}
\end{proposition}

In fact, this proposition will not be applied in this article. It will, however, play a useful role in the companion paper scaled Brownian LPP applications, because control on curves in all of space is needed when global properties are considered.
For example, the proposition may be applied to find that the maximizer of  $x \to \mc{L}(1,x)$
is tight uniformly over the space of regular ensembles.

\subsection{A rough guide to the structure of the remainder of the paper}\label{s.roughguide}

Before surveying this structure, it is useful to remark how central we have made the apparatus of regular Brownian Gibbs ensembles.
Indeed, we have reduced Theorem~\ref{t.neargeod} to two results, Theorems~\ref{t.airynt}(1) and~\ref{t.airynt.lb},
that are expressed in terms of regular ensembles, while another main theorem, Theorem~\ref{t.rnbound}, which concerns Brownian bridge regularity of the Airy line ensemble, is a very close cousin of Theorem~\ref{t.airytail.ln}, a result phrased in terms of regular ensembles which, in place of using the language of Radon-Nikodym moment bounds, records  a bound that controls the deformation in probability of unlikely events.

The crucial technical results that we must prove, such as  the $k$-curve closeness Theorem~\ref{t.airynt} and the just mentioned Theorem~\ref{t.airytail.ln}, are upper bounds on regular Brownian Gibbs ensemble probabilities. 
They will be proved with the aid of a general tool for such upper bounds that we will call the {\em jump ensemble method}. 
This method is crucial for deriving the more subtle of the results in the article, and its development and application is the guiding theme of the proofs that we undertake.
The method is built on a more basic apparatus which we call  {\em missing closed middle reconstruction}. This more basic method is enough to prove the more straightforward $k$-curve closeness lower bound Theorem~\ref{t.airynt.lb} and the local maximal fluctuation Theorem~\ref{t.aestimate}.

These considerations dictate the overall structure of the remainder of the article.
In the next chapter, we will present the missing closed middle construction, and a related technique called the Wiener candidate approach, and prove Theorems~\ref{t.airynt.lb} and~\ref{t.aestimate}. 
The jump ensemble method refines the Wiener candidate approach. The foundations of the new method are presented in Chapter~\ref{c.jumpensemblefoundations}.
It is plausible that this technique may have numerous applications in understanding behaviour of Brownian Gibbs ensembles, and thus in KPZ universality, and we hope that the presentation of the method in a general form in this chapter may facilitate such applications. We provide two applications of the method in this article. These are the proofs of 
 Theorems~\ref{t.airynt} and Theorem~\ref{t.airytail.ln}, which appear in Chapter~\ref{c.jumpensembleapplications}.

The article ends with an appendix in which certain important properties of regular Brownian Gibbs ensembles are proved.
The concerned results include  Proposition~\ref{p.lereg}, the result that makes the necessary connection from  regular ensembles to Brownian last passage percolation. The higher curve index one-point lower tail bound Proposition~\ref{p.othercurves} is necessary to set up the jump ensemble method, and its proof also appears in the appendix. 
The just presented `collapse near infinity'
Proposition~\ref{p.collapsenearinfinity} has a  Brownian Gibbs proof that is easily expressed here given the notation and tools we anyway develop, and the appendix ends with its derivation.

\noindent{\bf Remark on notation.}
We end this overview with a notational comment. Regular Brownian Gibbs ensembles being critical for our purpose,
 {\em henceforth our standing and implicit assumption is} that $\mc{L}_n:\intint{n} \times [-\xnmac,\infty) \to \R$ is
a   $\big(\bar\phimac,\rsc,\rsC\big)$-regular ensemble, for a given admissible choice of these parameters.
(It is Proposition~\ref{p.lereg} that furnishes such an admissible choice in applications to Brownian last passage percolation.)


\chapter{Missing closed middle reconstruction and the Wiener candidate} 

This chapter is devoted to the proofs of two  bounds on Brownian Gibbs ensemble properties, Theorems~\ref{t.airynt.lb} and~\ref{t.aestimate}, that are available with a less refined approach than will be needed for our other principal conclusions, 
 Theorems~\ref{t.airynt} and Theorem~\ref{t.airytail.ln}.
 In fact, the resampling method by which we will prove these first two results will be a fundamental underpinning of the more subtle jump ensemble method that we will later develop.
A key role of the present chapter will be to introduce the methods in the chapter's title. 

It is convenient to treat the one-point $k$-curve closeness problem as an example that motivates the development of our resampling method. We begin the chapter by considering this problem:
  in Section~\ref{s.mutual}, we turn   to the beginnings of an explanation for the $k^2-1$ exponent that we have seen in this context,
   deriving it simply for $k$-curve systems of mutually avoiding Brownian bridges. In considering how to prove the lower bound Theorem~\ref{t.airynt.lb}
   in light of the solution of this toy problem, we are led to ask how it is that 
    that $k$-curve closeness in the $k$-curve system may be inherited by Brownian Gibbs ensembles with far more curves. 
    We answer this question by setting up a technique,  missing closed middle reconstruction, for 
    transferring ensemble behaviour from a small system to a large one.
    Section~\ref{s.missingclosedmiddle} presents the general apparatus of this technique and, in Section~\ref{s.wienercandidateapplications}, we use the technique to prove Theorem~\ref{t.airynt.lb} as well as Theorem~\ref{t.aestimate}.

\section{Close encounter between finitely many 
non-intersecting
Brownian bridges}\label{s.mutual}

We now study the $k$-curve closeness probability in $k$-curve systems, finding that the $k^2-1$ exponent arises for them. The principal conclusion of this section will be Proposition~\ref{p.expsimple}. 

First we need the next proposition, a basic input from integrable probability.
Recall that $\vecint$ denotes the $k$-decreasing list $(k-1,k-2,\cdots,1,0) \in \R^k_>$.

\begin{proposition}\label{p.grabiner}
Let $k \in \N$ and $\grabell,K > 0$, and  let  $\bar{y} \in [-K,K]^k_>$ and $\emac \in \big( 0, \grabell k^{-2} K^{-1}  \big)$.
\begin{enumerate}
\item
We have that
\begin{equation}\label{e.grabiner}
 \mc{B}_{k;\emac \vecint,\bar{y}}^{[0,\grabell]} \Big( \nc^{[0,\grabell]} \Big) \, = \, \emac^{k(k-1)/2} \cdot   \grabell^{-k(k-1)/2}  \cdot \prod_{1 \leq i < j \leq k} \big( y_i - y_j \big) \cdot \big( 1 + E \big) \, ,
\end{equation}
where the error term $E$ satisfies 
\begin{equation}\label{e.grabinererror}
 - \, 2  \emac \, \grabell^{-1} k^2  K  \, \, \leq
E \, \leq \, \, (e^2-1) \emac \,  \grabell^{-1} k^2 K   \, .
\end{equation}
\item
Suppose that the left endpoint vector $\emac \vecint$ in the left-hand side of~(\ref{e.grabiner}) is replaced by any vector $\bar{x} \in \R^k_>$ for which $x_1 < x_k + \emac$. Then the upper bound on the error term $E$ stated in~(\ref{e.grabinererror}) remains valid.
\item Suppose instead that this same vector $\emac \vecint$  is replaced by any vector $\bar{x} \in \emac \vecint + \R^k_\geq$: that is, by any vector with the property that each of its consecutive component differences is at least~$\emac$. Then the lower bound on $E$ stated in~(\ref{e.grabinererror}) is still valid.
\end{enumerate}
\end{proposition}
\noindent{\bf Proof: (1).}  Set $\bar{x} = \emac \vecint$. By the Karlin-McGregor formula,
$$
\mc{B}_{k;\emac \vecint,\bar{y}}^{[0,\grabell]} \big( \nc^{[0,\grabell]} \big)  = \Big( \prod_{i=1}^k h(x_i,y_i) \Big)^{-1}  {\rm det} \big( h(x_i,y_j)  \big)_{1\leq i,j \leq k} \, ,
$$
where $h : \R^2 \to [0,\infty)$ equals $h(x,y) = g_{0,\grabell}(y-x)  = (2\pi \grabell)^{-1/2} \exp \big\{ - \tfrac{(x - y)^2}{2\grabell} \big\}$. 
(This form of the Karlin-McGregor formula follows from equation~$(4)$ in their paper~\cite{KarlinMcGregor} in the case of Brownian motions,  by taking $n = k$ and $E_i = (y_i - \epsilon,y_i + \epsilon)$ and considering the limit $\epsilon \searrow 0$. In fact, the authors demand that the processes involved be stationary. The method of proof does not require this hypothesis in our case, however; if we insist on applying the result directly, we might do so by making use of Brownian motions on a circle and letting the circle's radius tend to infinity.)

Thus,
$$
\mc{B}_{k;\emac \vecint,\bar{y}}^{[0,\grabell]} \big( \nc^{[0,\grabell]} \big)  =   \exp \Big\{ - \grabell^{-1} \sum_{i=1}^k  x_i  y_i \Big\} \cdot
{\rm det} \big( e^{x_i y_j \grabell^{-1}}  \big)_{1\leq i,j \leq k} \, .
$$
The determinant is Vandermonde, and thus may be factored: 
$$
{\rm det} \big( e^{x_i y_j \grabell^{-1}}  \big)_{1\leq i,j \leq k} = 
{\rm det} \Big(  \big( e^{\emac y_j \grabell^{-1}} \big)^{i-1}  \Big)_{1\leq i,j \leq k} = \prod_{1 \leq i < j \leq k} \big( e^{\emac y_j \grabell^{-1}} - e^{\emac y_i \grabell^{-1}} \big) \, .
$$
Writing the term in the product in the form $\emac (y_j - y_i)\grabell^{-1} e^{\emac t_{i,j} \grabell^{-1}}$ for certain $t_{i,j} \in [y_i,y_j]$ using the mean value theorem, we find that $
\mc{B}_{k;\emac \vecint,\bar{y}}^{[0,\grabell]} \big( \nc^{[0,\grabell]} \big)$ equals  
$$  
 \big( \emac \grabell^{-1} \big)^{k(k-1)/2} \cdot  \bigg( \prod_{1 \leq i < j \leq k} ( y_i - y_j  ) \bigg) \cdot \big( 1 + E  \big)  \, ,
$$
where 
$$
 1 + E =  \exp \Big\{ - \grabell^{-1} \sum_{i=1}^k x_i y_i \Big\} \cdot \exp \big\{ \emac \grabell^{-1} t \big\}
$$
with
$t = \sum_{1 \leq i < j \leq k} t_{i,j}$.
The two right-hand factors will now be bounded by the same upper and lower bounds: to do so, we begin by seeing that the terms $\sum_{i=1}^k x_i y_i$
and $\emac t$ are in absolute value at most $\emac k^2 K$. In the latter case, this follows from $\vert  t \vert \leq k^2 K$, a fact which is due to $\bar{y} \in [-K,K]^k$. In the former, it is a consequence of Cauchy-Schwarz.
Since $\emac \grabell^{-1} k^2 K$ is at most one by hypothesis, the inequality  $e^x \leq 1 + (e-1)x$, $x \in [0,1]$, yields an upper bound on both of the last displayed factors; a lower bound on each is furnished by $1+x \leq e^{x}$, $x \in \R$. What we have learnt is that
$$
 \Big( 1 - \emac \grabell^{-1} k^2 K  \Big)^2 \, \leq \,   1 + E \, \leq \, 
 \Big( 1 + (e-1) \emac \grabell^{-1} k^2 K \Big)^2 \, . 
$$
Another use of $\emac \grabell^{-1 }k^2 K \leq 1$ yields   
Proposition~\ref{p.grabiner}(1). 

\noindent{\em Remark.} The preceding proof is a variant of arguments in \cite{Grabiner}, which derives asymptotic avoidance probabilities for Brownian ensembles with various constraints. 

\medskip

\noindent{\bf (2).} The differences between successive components of the vector $\bar{x}$ are all less than $\emac$. 
The law $\mc{B}_{k;\emac \vecint,\bar{y}}^{[0,\grabell]}$ may be obtained from  $\mc{B}_{k;\bar{x},\bar{y}}^{[0,\grabell]}$   by an affine translation of each curve.  Since the differences between successive components of the vector $\bar{x}$ are all less than $\emac$, the occurrence of  $\nc^{[0,\grabell]}$ is maintained by this procedure.

\medskip

\noindent{\bf (3).} The role-reversed observation applies. \qed

\begin{definition}\label{d.closeenc}
Let $k,m \in \N$ satisfy $m \geq k$. For $[a,b] \subseteq \R$, $x \in [a,b]$, $\phi > 0$,
and an ensemble $E:\intint{m} \times [a,b] \to \R$, we define the event $\close \big( k; E,x,\phi \big)$ that 
$$
 E(i,x) - E(k,x) \in (0,\phi) \, \, \, \, \textrm{for each $i \in \intint{k-1}$} \, .
$$
When the ensemble $E$ is ordered, this event is specified by the condition that $E(1,x) < E(k,x) + \phi$, consistently with the usage of this notation in Theorem~\ref{t.airynt}.

The parameter $k$ will be consistently used when the $\close$ event is studied, and we will use the shorthand $\close \big(  E,x,\phi \big)$. 
\end{definition}

Theorems~\ref{t.airynt}(1) and~\ref{t.airynt.lb} assert that this event's probability behaves as $\e^{k^2 -1 + o(1)}$ in a limit of $\e \searrow 0$. We now present a simple result that gives a first indication as to why this behaviour may be expected. A system of $k$ mutually avoiding Brownian bridges defined on an interval of unit-order length has such a probability of near-touch at the midpoint time, provided that the entrance and exit data are both well-spaced; moreover, this remains true in the presence of a lower boundary condition  that consistently remains a respectful distance below the endpoint locations.


\begin{proposition}\label{p.expsimple}
$\empty$
\begin{enumerate}
\item Suppose that $\bar{x},\bar{y} \in 2\vecint + [0,\infty)_\geq^k$ and $f:[-1,1] \to [-\infty,-1]$ is measurable. Moreover, we assume that $x_1 \vee y_1 \leq K$ for some given constant $K \geq 1$. If $\phi \in (0, k^{-2}K^{-1})$, then 
$$
\mc{B}_{k;\bar{x},\bar{y}}^{[-1,1]} \Big( \, \close \big( B,0,\phi \big) \, \, \Big\vert \, \, \nc_f^{[-1,1]}  \, \Big) \leq
 \phi^{k^2-1} \cdot \big( 1 - 2 e^{-1} \big)^{-k} \pi^{-(k-1)/2} K^{k^2} e^4  k^4 \, . 
$$
\item 
Suppose that $\bar{x},\bar{y} \in \vecint + [0,\infty)^k_\geq$; that  $\bar{x},\bar{y} \in [0,K]^k$ for some given $K \geq 1$; and also that  $f:[-1,1] \to \R$ is measurable and satisfies $f \leq -4 \sqrt{2}k$. If $\phi \in \big( 0 ,(2kK)^{-1} \big)$, then
\begin{eqnarray*}
 & & \mc{B}_{k;\bar{x},\bar{y}}^{[-1,1]} \Big( \, \close \big( B,0,\phi \big) \, \, \Big\vert \, \, \nc_f^{[-1,1]}  \, \Big) \\ & \geq &
  \phi^{k^2-1} \cdot \big( 1 -  \pi^{-1/2}e^{-1} \big) \, \tfrac{1}{16}  (2k)^{-(k^2-1)}  \pi^{-(k-1)}  \exp \big\{ - 2 (k-1) (2K+1)^2 \big\}  \, . 
\end{eqnarray*}
\end{enumerate}
\end{proposition}
Note that the case $f = -\infty$
of an absent lower boundary condition is included.

\medskip

\noindent{\bf Proof of Proposition~\ref{p.expsimple}.}
By Bayes' theorem, the conditional probability in question equals
\begin{equation}
  \mc{B}_{k;\bar{x},\bar{y}}^{[-1,1]} \Big( \,   \nc_f^{[-1,1]}  \, \Big)^{-1} \cdot \, \mc{B}_{k;\bar{x},\bar{y}}^{[-1,1]} \Big( \, \close \big( B,0,\phi \big)  \, , \,  \nc_f^{[-1,1]}  \, \Big) \, . \label{e.bayesconseq}
\end{equation}
The first factor is the reciprocal of the acceptance probability associated to the law $\mc{B}_{k;\bar{x},\bar{y}}^{[-1,1]}$ and the function $f$. 
This probability is bounded below for our parameter choices, as we now show.
\begin{lemma}\label{l.accproblb} 
When $\bar{x},\bar{y} \in 2\vecint + [0,\infty)_\geq^k$ and $f:[-1,1] \to \R$ satisfies $f \leq -1$,
$$
\mc{B}_{k;\bar{x},\bar{y}}^{[-1,1]} \Big( \,   \nc_f^{[-1,1]}  \, \Big) \geq \big(1 - 2 e^{-1} \big)^k \, .
$$
\end{lemma}
\noindent{\bf Proof.}
Consider the planar line segments $\ell_i \subset [-1,1] \times \R$ that interpolate the boundary data $(-1,x_i)$ and $(1,y_i)$ for $i \in \intint{k}$.
The corridor $\cdor_i \subset [a,b] \times \R$
consists of points whose vertical displacement from $\ell_i$ is less than one.
Under our assumption, the $k$ corridors are disjoint, with the lowest one, $\cdor_k$, lying above the graph of $f$.
Thus, for a sample of 
$\mc{B}_{k;\bar{x},\bar{y}}^{[-1,1]}$ to realize the event~$ \nc_f^{[-1,1]}$,   it is enough that each of the $k$ curves remain in the corridor that shares its index.
We may express this eventuality in terms of the standard bridges associated to the curves: in these terms, this event is
$$
 \bigcap_{i = 1}^k \Big\{ \big\vert B^{[-1,1]}\big( i, x \big) \big\vert < 1 \, \, \, \forall \, x \in [-1,1] \, \Big\} \, .
$$
The standard bridges are independent under 
$\mc{B}_{k;\bar{x},\bar{y}}^{[-1,1]}$. Thus, Lemma~\ref{l.maxfluc}
implies that the $\mc{B}_{k;\bar{x},\bar{y}}^{[-1,1]}$-probability of the displayed event is at least $\big(1 - 2 e^{-1} \big)^k$. \qed 

\medskip

We may write the second factor in~(\ref{e.bayesconseq}) in the form of a product of three terms
\begin{equation}\label{e.triplea}
 \mc{B}_{k;\bar{x},\bar{y}}^{[-1,1]} \Big( \, \close \big( B,0,\phi \big)  \, , \,  \nc_f^{[-1,1]}  \, \Big) =  A_1 \cdot A_2 \cdot A_3 \, ,
\end{equation}
where the first term is the probability of closeness
$$
A_1 = 
 \mc{B}_{k;\bar{x},\bar{y}}^{[-1,1]} \Big(   \, \close \big( B,0,\phi \big)  \, \Big) \, ;
$$
the second is the conditional probability of avoidance on the left interval $[-1,0]$ given closeness
$$
A_2 = 
 \mc{B}_{k;\bar{x},\bar{y}}^{[-1,1]} \Big(   \,  \nc_f^{[-1,0]}  \, \Big\vert \, \close \big( B,0,\phi \big) \Big) \, ,
$$
and the third is the conditional probability of avoidance on the right interval $[0,1]$,
$$
A_3 = 
 \mc{B}_{k;\bar{x},\bar{y}}^{[-1,1]} \Big(   \,  \nc_f^{[0,1]}  \, \Big\vert \, \close \big( B,0,\phi \big) \, , \,  \nc_f^{[-1,0]}  \Big) \, .
$$
As we will now show, these events have probabilities in the low $\phi$ limit described by the dominant terms $\phi^{k-1}$, $\phi^{k(k-1)/2}$ and  $\phi^{k(k-1)/2}$ .
Indeed, this triple product probability and the formula
$$
 \phi^{k^2 - 1} \, = \, \phi^{k-1} \cdot \phi^{k(k-1)/2} \cdot \phi^{k(k-1)/2} \, , 
$$
is a useful overview of the reason why we may at least begin to expect to see the exponent~$k^2 - 1$ in the one-point $k$-curve closeness estimate in Theorems~\ref{t.airynt}(1) and~\ref{t.airynt.lb}. 

We now substantiate these claims about dominant behaviour in the $\phi \searrow 0$ limit.
The next lemma treats upper bounds and in view of~(\ref{e.triplea}) completes the proof of Proposition~\ref{p.expsimple}(1).

\begin{lemma}\label{l.aub}
Let $f: [-1,1] \to \R \cup \{ - \infty \}$ be measurable.
\begin{enumerate}
\item
For any vectors $\bar{x},\bar{y} \in \R^k$, 
$$
A_1 \leq \pi^{-(k-1)/2} \phi^{k-1} \, .
$$
\item Suppose that $\phi < k^{-2}K^{-1}$. 
If $\bar{x},\bar{y} \in [0,K]^k$ for some given $K \geq 1$,
then $A_2$ and $A_3$ are both at most
$\phi^{k(k-1)/2} \cdot K^{k^2/2} e^2  k^2$.
\end{enumerate}
\end{lemma}
\noindent{\bf Proof: (1).}
We may condition on the value of $B(k,0)$; if it is taken equal to $x$, then the closeness event takes the form $\cap_{i=1}^{k-1} \big\{ B(i,0) \in [x,x+\phi] \big\}$. 
The random variable $B(i,0)$ is normally distributed under  $\mc{B}_{k;\bar{x},\bar{y}}^{[-1,1]}$ with mean $(x_i + y_i)/2$ and variance $1/2$; from~(\ref{e.gaussiandensity}), we thus see that it has probability at most $\phi \pi^{-1/2}$
to belong to any given interval of length $\phi$. The upper bound in the lemma's first part results from multiplying this bound over the $k-1$ indices. 

\noindent{\bf (2).} The quantity $A_2$
is at most the supremum of 
$$
 \mc{B}_{k;\bar{x},\bar{z}}^{[-1,0]} \Big(   \,  \nc_f^{[-1,0]}  \,  \Big) 
$$
as $\bar{x}$ varies over $[0,K]^k$
and $\bar{z}$ over vectors in some displacement of the set $[0,\phi]^k$. We may eliminate the lower boundary condition by setting $f = -\infty$, since it is trivial that doing so only increases the probability in question. We then apply Proposition~\ref{p.grabiner}(2) with $\emac = \phi$ and  $\grabell = 1$ to find that 
$$
A_2 \leq \phi^{k(k-1)/2} \cdot K^{k(k-1)/2} \big( 1 + (e^2 - 1)\phi k^2 K \big) \, . 
$$ 
Using $\phi \leq 1 \leq K$, we
confirm the stated upper bound on~$A_2$.

The quantity $A_3$
is at most the supremum of 
$$
 \mc{B}_{k;\bar{z},\bar{y}}^{[0,1]} \Big(   \,  \nc_f^{[0,1]}  \,  \Big) 
$$
as $\bar{y}$ varies over $[0,K]^k$
and $\bar{z}$ over vectors in some displacement of the set $[0,\phi]^k$. 
The preceding argument proves the upper bound on~$A_3$. \qed

\medskip 
 In order to prove the lower bound on $\mc{B}_{k;\bar{x},\bar{y}}^{[-1,1]} \big( \, \close ( B,0,\phi )  \, , \,  \nc_f^{[-1,1]}  \, \big)$ in Proposition~\ref{p.expsimple}(2), we first specify a convenient sufficient condition for the occurrence of $\close ( B,0,\phi )$. Let $D \subset \R^{k-1}$
 denote the $(k-1)$-dimensional box
 $$
 D \, = \, \big[ 2k - 3 , 2k - 2 \big] \times \cdots \times \big[3,4\big] \times \big[1,2\big] \, ;
 $$
specify the event 
$\boxclose \big(B,0,\phi \big)$ to be equal to 
$$
\Big\{ \Big( B(1,0) - B(k,0) , B(2,0) - B(k,0) , \cdots, B(k-1,0) - B(k,0) \Big) \in  \tfrac{1}{2k} \phi \cdot D \, , \, B(k,0) \geq -1 \Big\} \, ;
$$
and  
note that $  \boxclose \big(B,0,\phi \big)  \subset  \close \big( B,0,\phi \big)$.
In place of~(\ref{e.triplea}), we  write
$$
 \mc{B}_{k;\bar{x},\bar{y}}^{[-1,1]} \Big( \,  \boxclose \big( B,0,\phi \big)  \, , \,  \nc_f^{[-1,1]}  \, \Big) =  \nF_1 \cdot \nF_2 \cdot \nF_3 \, ,
$$
where
$$
\nF_1 = 
 \mc{B}_{k;\bar{x},\bar{y}}^{[-1,1]} \Big(   \,  \boxclose \big( B,0,\phi \big) \, \Big) \, ,
$$
$$
\nF_2 = 
 \mc{B}_{k;\bar{x},\bar{y}}^{[-1,1]} \Big(   \,  \nc_f^{[-1,0]}  \, \Big\vert \,  \boxclose \big( B,0,\phi \big) \, \Big) \, ,
$$
and
$$
\nF_3 = 
 \mc{B}_{k;\bar{x},\bar{y}}^{[-1,1]} \Big(   \,  \nc_f^{[0,1]}  \, \Big\vert \,  \boxclose \big( B,0,\phi \big) \, , \,  \nc_f^{[-1,0]}  \Big) \, .
$$

Proposition~\ref{p.expsimple}(2) thus follows from the next lemma. \qed
\begin{lemma}\label{l.alb}
$\empty$
\begin{enumerate}
\item
If $\phi \leq 1$ and  $\bar{x},\bar{y} \in [0,K]^k$ for some given $K \geq 1$,
then 
$$
\nF_1 \geq \phi^{k-1} \cdot \big( 1 -  \pi^{-1/2}e^{-1} \big) (2k)^{-(k-1)}  \pi^{-(k-1)/2}  \exp \big\{ - (k-1) (2K+1)^2 \big\} \, .
$$
\item If $\phi \leq (2kK)^{-1}$, $\vec{x},\vec{y} \in \vecint + [0,\infty)^k_\geq$ and $f:[-1,1] \to \R$ is measurable and satisfies $f \leq -4 \sqrt{2} k$, then both $\nF_2$ and $\nF_3$ are at least
$$
 \phi^{k(k-1)/2} \cdot  \tfrac{1}{4}  (2k)^{-k(k-1)/2} \, . 
$$ 
\end{enumerate}
\end{lemma}
\noindent{\bf Proof: (1).} Recall that $B(i,0)$ has law $\nu_{(x_i+y_i)/2,1/2}$ under 
 $\mc{B}_{k;\bar{x},\bar{y}}^{[-1,1]}$. Note that each of these means $(x_i + y_i)/2$ lies in $[0,K]$.
The probability that $B(k,0)$ adopts a value in $[-1,K+1]$ is at least $1 - 2 \nu_{0,1/2}(1,\infty) \geq 1 -  \pi^{-1/2}e^{-1}$. Supposing that it adopts such a value $x \in [-1,K+1]$, it is sufficient for $\boxclose(B,0,\phi)$ to occur that, for each
 $i \in \intint{k-1}$,  $B(i,0)$ lies in a certain interval of length $\tfrac{1}{2k}\phi$ inside $[x,x+\phi]$; for given $i$ and any such interval $I$,  this circumstance happens with probability at least 
$\int_I g_{0,1/2}(x) \dd x$ where note that the interval $I$ is of length $\tfrac{1}{2k} \phi$ and is contained in $[-2K,2K+\phi]$. Since $\phi \leq 1$, this expression is at least $\tfrac{1}{2k} \phi \pi^{-1/2} \exp \big\{ - (2K+1)^2 \big\}$. 

\noindent{\bf (2).} 
Note that $\nF_2$ is at least the infimum of 
$$
 \mc{B}_{k;\bar{x},\bar{z}}^{[-1,0]} \Big(   \,  \nc_f^{[-1,0]}  \,  \Big) 
$$
as $\bar{x}$ varies over elements of $\vecint + [0,\infty)^k_\geq$ and $\bar{z} \in \R^k_>$ over elements 
of $\R^k_>$ such that
$\big( z_1 - z_k, z_2 - z_k, \cdots, z_{k-1} - z_k \big) \in \tfrac{\phi}{2k}D$ and $z_k \geq - 1$. 
The displayed quantity may be represented 
$$
 \mc{B}_{k;\bar{x},\bar{z}}^{[-1,0]} \Big(   \,  \nc^{[-1,0]}  \,  \Big) \cdot 
 \mc{B}_{k;\bar{x},\bar{z}}^{[-1,0]} \Big( \,  \nc_f^{[-1,0]}  \,  \Big\vert \,  \nc^{[-1,0]}  \,  \Big) \, .
$$
A lower bound on the first term in the product is provided by Proposition~\ref{p.grabiner}(3) with the choice $\emac = \tfrac{1}{2k} \phi$ (and $\grabell = 1$). The assumption that $\phi \leq (2kK)^{-1}$ implies that the term $1+E$ is at least $1/2$.
The second term is bounded below by means of Lemma~\ref{l.notouchfluc}, taking $[a,b]$ equal to $[-1,0]$ and $\bar{x}$ equal to its present value, as well as $\bar{y} = \bar{z}$ and $r=k$. The 
lemma bounds above the probability of the complementary event, finding this probability to be at most $(1 - 2e^{-1})^{-k} e^{-4r^2}$, which is at most $1/2$. Thus, the second term is at least $1/2$. 

The lower bound on $\nF_3$ follows similarly. \qed

\medskip

We end this section by noting a related result that will be used later.

\begin{lemma}\label{l.browniannotouch}
For any $k \geq 2$, $\phi \in \big(0,k^{-3}/6\big)$, and $\bar{x},\bar{y} \in \R^k_>$,
$$
\mc{B}_{k;\bar{x},\bar{y}}^{[-1,1]} \Big( \notouch^{[-1,1]} \, , \, \close \big( B ,0, \phi \big) \Big)
\leq 4258 \cdot (36)^{k^2} (k^2-1)^{k^2}  \phi^{k^2-1} \big(\log \phi^{-1}\big)^{k^2/2} \, . 
$$
\end{lemma}
\noindent{\bf Proof.}
Suppose in the first instance that $(x_1 - x_k) \vee (y_1 - y_k) \leq \Chat (\log \e^{-1})^{1/2}$, where $\Chat = 6(k^2-1)^{1/2}$. By the affine symmetry of Brownian bridge, we may suppose that $\bar{x},\bar{y} \in [0,K]^k$ with $K = \Chat \big( \log \phi^{-1} \big)^{1/2} \geq 1$.
With this choice of $K$, we may apply Lemma~\ref{l.aub} to the formula~(\ref{e.triplea})
to learn that
$$
\mc{B}_{k;\bar{x},\bar{y}}^{[-1,1]} \Big( \notouch^{[-1,1]} \, , \, \close \big( B ,0, \phi \big) \Big) \leq  \phi^{k^2-1} \big(\log \phi^{-1}\big)^{k^2/2} \Chat^{k^2} \cdot \pi^{-(k-1)/2} e^4 k^4 \, .
$$
(The hypothesis $\phi < k^{-2} K^{-1}$ of Lemma~\ref{l.aub}(2) takes the form $6(k^2 - 1)^{1/2} k^2 \big( \log \phi^{-1} \big)^{1/2} \phi < 1$. The logarithmic term may be omitted here because $\phi < e^{-1}$. The present hypothesis $\phi \leq k^{-3}/6$ is then seen to imply the resulting bound.)

The post-$\cdot$ right-hand term in this display is at most $4258$, 
(because $\sup_{k \geq 2} \pi^{-k/2} k^4$ is attained by $k=7$
and is thus at most $44$).

Suppose instead that  $(x_1 - x_k) \vee (y_1 - y_k) > \Chat (\log \phi^{-1})^{1/2}$. In this case, note that
\begin{eqnarray*}
& &
\mc{B}_{k;\bar{x},\bar{y}}^{[-1,1]} \Big( \close\big(B,0,\phi\big)  \Big) 
 \leq  \sum_{i=1}^k
\mc{B}_{k;\bar{0},\bar{0}}^{[-1,1]} \Big( \big\vert B(i,0) \big\vert \geq    \tfrac{1}{6} \Chat (\log \phi^{-1})^{1/2}  \Big) \\  
& \leq &  k \cdot 2(2\pi)^{-1/2}  \Big( \tfrac{1}{3} \Chat (\log \phi^{-1})^{1/2}  \Big)^{-1} \exp \big\{ - \tfrac{1}{36} \Chat^2 \log \phi^{-1} \big\} \leq \phi^{\Chat^2/36} \, ,
\end{eqnarray*}
where in the first inequality, we used $\phi < \tfrac{1}{6} \Chat (\log \phi^{-1})^{1/2}$ (a bound due to $\phi \leq e^{-1}$ and $\Chat \geq 6$); and
then $\Chat (\log \phi^{-1})^{1/2} \geq 3 \sqrt{2} \pi^{-1/2} k$. Since $\Chat^2/36 = k^2 - 1$, we find that
$$
\mc{B}_{k;\bar{x},\bar{y}}^{[-1,1]} \Big( \close\big(B,0,\phi\big)  \Big) \leq \phi^{k^2 - 1} \, .
$$
This completes the proof of Lemma~\ref{l.browniannotouch}. \qed


\section{The reconstruction of the missing closed middle}\label{s.missingclosedmiddle}

Recall from the remark concerning notation at the end of Section~\ref{s.roughguide} that $\mc{L}_n:\intint{n} \times [-\xnmac,\infty) \to \R$ denotes 
a   $\big(\bar\phimac,\rsc,\rsC\big)$-regular ensemble.


As we prepare to introduce the reconstruction that, in the first instance, we will use to  prove the lower bound  Theorem~\ref{t.airynt.lb},  it may be useful to give a sense of the guiding ideas of this proof. The essence of these ideas is quite straightforward.  We may consider the eventuality that the highest~$k$ curves of the ensemble $\mc{L}_n$ are separated from one another to unit order above locations~$-1$ and~$1$, and are comfortably above the curve $\mc{L}_n(k+1,\cdot)$ throughout the interval $[-1,1]$. If we can argue that this circumstance occurs with a probability that is uniform in the parameter $n$, then we may note that the conditional distribution of $\mc{L}_n:\intint{k} \times [-1,1] \to \R$ verifies the hypotheses of Proposition~\ref{p.expsimple}; the lower bound in Theorem~\ref{t.airynt.lb} will then result.

(We have just made the first use  of a notational abuse that we will sometimes employ. The line ensemble $\mc{L}_n$ has domain of definition $\intint{n} \times [-\xnmac,\infty)$. When we are interested in an ensemble's behaviour on a subdomain, we will refer to such objects as $\mc{L}_n:\intint{k} \times [-1,1] \to \R$, even though this is technically incorrect.)

We will verify that this separation circumstance has uniformly positive probability by introducing a procedure by which the law of the entire ensemble $\mc{L}_n:\intint{n} \times [-\xnmac,\infty) \to \R$ may be sampled. The procedure may be viewed as a reconstruction: the law is realized, with the information that specifies it being represented in a certain convenient form. Some pieces of information in the representation are retained, and others are forgotten. The lost data is then reconstructed according to its conditional distribution given the retained data, so that a new, reconstructed, copy of the ensemble results.

Since the procedure and its jump ensemble elaboration will play such a key role in this paper, as we now present the procedure, we will discuss in some detail how it may be interpreted probabilistically.

\subsection{Specifying the missing closed middle reconstruction procedure}
Beyond the total curve number index $n$ and the fixed index $k \in \intint{n}$, the {\em missing closed middle} reconstruction procedure has three parameters: $\eln > 0$ and the {\em left} and {\em right} parameters $\gfl$ and $\gfr$ that satisfy 
$\gfl \in [-\eln,0]$
and $\gfr \in [0,\eln]$. The {\em middle interval} is $[\gfl,\gfr]$; it is straddled by the  {\em side intervals} $[-2\eln,\gfl]$ and $[\gfr,2\eln]$.

Let $n \in \N$ and  $a,b \in \R$ satisfy $-z_n \leq a \leq b$. Recall the standard bridge ensemble
$\mc{L}_{n}^{[a,b]} : \intint{n} \times [a,b] \to \R$ induced on~$[a,b]$ by $\mc{L}_n$, namely 
$$
\mc{L}_{n}^{[a,b]}\big( i, x \big) = 
\mc{L}_{n}\big( i, x \big) \, - \,
\ell_{n}^{[a,b]}\big( i, x \big) \, \, \, \textrm{for $(i,x) \in \intint{n} \times [a,b]$} \, , 
$$
where 
$\ell_{n}^{[a,b]}\big( i, \cdot \big)$ denotes the affine function whose values at $a$ and $b$ are $\mc{L}_n\big(i,a\big)$ and $\mc{L}_n\big(i,b\big)$.

Let $\mc{F}$ denote the {\em missing closed middle} $\sigma$-algebra, generated by the following collection of random variables:
\begin{itemize}
\item  the curves $\mc{L}_n:\llbracket k+1,n \rrbracket \times [-\xnmac,\infty) \to \R$ of index at least $k+1$;
\item the highest $k$ curves $\mc{L}_n:\intint{k}  \times \big( [-\xnmac,-2\eln] \cup [2\eln, \infty) \big) \to \R$ outside $(-2\eln,2\eln)$;
\item and the $2k$ standard bridges $\mc{L}^{[-2\eln,\gfl]}_{n}\big(i,\cdot \big)$ and $\mc{L}^{[\gfr,2\eln]}_{n}\big(i,\cdot \big)$, where $i \in \intint{k}$.
\end{itemize}
This $\sigma$-algebra clearly depends on the index~$k$, though we omit display of this dependence in our notation.

Let $\PP_{\mc{F}}$ denote the conditional probability given $\mc{F}$. That is, $\PP_{\mc{F}}(A) = \EE \big( {\bf 1}_A \, \big\vert \, \mc{F} \big)$ for any measurable event $A$.
The law $\PP_\mc{F}$ represents the information available to, and statistical uncertainty of, the 
observer who is informed by an experimenter who samples the law $\PP$ only of the data constituting~$\mc{F}$. To this observer, whom we may call the witness of $\mc{F}$, 
certain aspects of the behaviour of 
the line ensemble $\mc{L}_n:\intint{n} \times [-\xnmac,\infty) \to \R$ remain unknown, and random. 

The use of the law $\PP_\mc{F}$ is central to the proofs of our main results. It is thus valuable to carefully consider the nature of this object, a task which amounts to understanding the perspective of the witness of~$\mc{F}$ as he considers the conditional law of the entire line ensemble given the available data. In this regard, we should ask: what is known to the witness, and what is random? Does the conditional law of what is unknown have a convenient and explicit representation?

Clearly picturing the $\mc{F}$-witness's perspective is also useful, because we will later present arguments based on resamplings of the line ensemble $\mc{L}_n$ under which the data that is frozen is an augmentation of that specifying~$\mc{F}$. The concerned data in such arguments will be specified by a larger $\sigma$-algebra than $\mc{F}$, with the associated probability experiment having its own, more knowledgeable, witness. We will again evoke the perspective of these witnesses when we present these later arguments; the groundwork is set by clearly understanding the nature of $\PP_\mc{F}$.  

We should address then the posed questions about the perspective of the witness of~$\mc{F}$. 
What is known to this observer? Clearly, the three item list of data that specifies~$\mc{F}$.
What data is random, and how may we conveniently depict this randomness?
The $\mc{F}$-random data consists of 
\begin{itemize}
\item $\mc{L}_n:\intint{k} \times [\gfl,\gfr] \to \R$.
\end{itemize}
Alternatively, this data may be represented in the form:
\begin{itemize}
\item and the endpoint $k$-vectors $\big( \mc{L}_n(i,\gfl): i \in \intint{k} \big)$ and  $\big( \mc{L}_n(i,\gfr): i \in \intint{k} \big)$;
\item and the standard bridges $\mc{L}^{[\gfl,\gfr]}_{n}\big(i,\cdot \big)$ for $i \in \intint{k}$.
\end{itemize}

Both of these ways of writing the $\mc{F}$-random data will be used. We will call them the {\em one-piece} and {\em two-piece list} presentations. We speak of the missing {\em closed} middle because, as the two-piece list shows, the witness of~$\mc{F}$ is unaware of middle interval endpoint data as well as standard bridge data.

Consider for now the two-piece list presentation of the data.
 We first present a specification of a space of possible outcomes for the witness of~$\mc{F}$ of the form of such data.

Consider data 
$(\bar{x},\bar{y}) \in \R^k_> \times \R^k_>$ and a $k$-vector $\bar{b}$ whose components are standard bridges $b_i$ belonging to $\mc{C}_{0,0}\big( [\gfl,\gfr] , \R \big)$ for $i \in \intint{k}$.

The witness of $\mc{F}$ may consider the possibility that the $\mc{F}$-random data in two-piece list form adopts the values $(\bar{x},\bar{y})$ in the first item and $\bar{b}$ in the second. 
Indeed, the witness may reconstruct the top $k$ curves  in this eventuality:
for each $i \in \intint{k}$, the~$i\textsuperscript{th}$ curve so reconstructed may be denoted by 
$$
\larg{x_i}{y_i}{b_i} \big( i, \cdot\big):[-\xnmac,\infty) \to \R \, .
$$ 
It is specified by
\begin{equation}\label{e.reconcases}
 \larg{x_i}{y_i}{b_i} \big( i, s \big) = \begin{cases}
  \, \mc{L}_n\big(i,s \big) \, & s \in [-\xnmac,-2\eln) \, , \\
 \, \mc{L}_n^{[-2\eln,\gfl]}\big(i,s\big) + 
 \tfrac{\gfl - s}{\gfl + 2\eln}  \lppls_n\big(i, - 2\eln \big)
 + \tfrac{s + 2\eln}{\gfl + 2\eln}  x_i   \, & s\in \big[-2\eln,\gfl \big] \, , \\
 \, \tfrac{\gfr - s}{\gfr - \gfl} x_i + \tfrac{s - \gfl}{\gfr - \gfl} y_i \, + \, b_i(s)  \, & s \in \big[\gfl ,\gfr \big] \, , \\
 \, \mc{L}_n^{[\gfr,2\eln]}\big(i,s\big) + 
 \tfrac{s - \gfr}{2\eln - \gfr}  \lppls_n\big(i,2\eln\big)
 + \tfrac{2\eln - s}{2\eln - \gfr}  y_i   \, & s\in \big[\gfr,2\eln \big] \, ,  \\
   \, \mc{L}_n\big(i,s \big) \, & s \in (2\eln,\infty) \, .
\end{cases}
\end{equation}
Note that indeed the curve
$\larg{x_i}{y_i}{b_i}\big( i,\cdot \big)$ depends on its parameters only via the $i$-indexed variables $(x_i,y_i,b_i)$. The curve is  compatible with the data in $\mc{F}$ in the sense that
the curve is specified by its form outside $(-2\eln,2\eln)$,   standard bridges on $[-2\eln,\gfl]$ and on $[\gfr,2\eln]$, and its form on $[\gfl,\gfr]$; and of these three pieces of data, 
the first two coincide with the data specified by~$\mc{F}$.

As such, the set of $k$-vectors of reconstructed curves given by  
$$
 \Big( \larg{x_i}{y_i}{b_i} \big( i, \cdot\big):[-2\eln,2\eln]\to \R \, , \, i \in \intint{k} \Big)
$$ 
as the triple $\big(\bar{x},\bar{y},\bar{b}\big)$ varies over  $\R^k \times \R^k \times \mc{C}_{0,0}\big([\gfl,\gfr], \R \big)^k$  constitutes for the witness of~$\mc{F}$ a space of possible outcomes for the form of $\mc{L}_n$ on $\intint{k} \times [-2\eln,2\eln]$.  Note that here we have begun to neglect the region outside $[-2\eln,2\eln]$ where $\mc{F}$-data dictates the outcome.

What law on this set of triples gives the conditional distribution for the witness of~$\mc{F}$ of the outcome $\mc{L}_n:\intint{k} \times [-2\eln,2\eln] \to \R$? 

To answer this question, we begin by embuing the outcome set $\R^k \times \R^k \times \mc{C}_{0,0}\big([\gfl,\gfr], \R \big)^k$ with a {\em reference measure}, namely
the product measure of $k$-dimensional Lebesgue measure in the first two coordinates, and the standard $k$-dimensional Brownian bridge law~$\mc{B}^{[\gfl,\gfr]}_{k;\bar{0},\bar{0}}$ in the third. 

We may then find a form for the conditional distribution in question by finding its Radon-Nikodym distribution with respect to the reference measure. Indeed, we
let $h: \R^k \times \R^k \times \mc{C}_{0,0}\big([\gfl,\gfr], \R \big)^k \to [0,\infty)$ denote the density  with respect to the reference measure of the conditional law under $\PP_\mc{F}$ of the triple of  $k$-vectors 
$$
\bigg( \Big( \lppls_n \big(i,\gfl \big) : i \in \intint{k} \Big) \, , \,
\Big( \lppls_n \big(i,\gfr \big) : i \in \intint{k} \Big) \, ,
\, \Big( [\gfl,\gfr] \to \R: x \to \mc{L}^{[\gfl,\gfr]}_{n}\big(i,x \big) : i \in \intint{k}  \Big) \bigg) \, .  
$$
The Brownian bridge basic Lemma~\ref{l.brbr} may be used to explicitly compute $h$: $h\big(\bar{x},\bar{y}, \bar{b} \big)$ equals
\begin{eqnarray}
   & & Z^{-1}  \prod_{i=1}^k \, \exp\bigg\{-\tfrac{1}{2(\gfl + 2\eln)}\Big(\lppls_n\big(i,-2\eln \big) - x_i \Big)^2  - \,  \tfrac{1}{2(\gfr - \gfl)} \big( x_i - y_i \big)^2 \label{e.gdens.first.new} \\
 & &  \qquad \qquad \qquad \qquad \qquad \qquad  - \, \,  \tfrac{1}{2(2\eln - \gfr)} \Big( \lppls_n\big(i,2\eln  \big) - y_i \Big)^2 \bigg\} \,\cdot\,\BP_{1,k+1} \big( \bar{x},\bar{y}, \bar{b} \big) \, , \nonumber
\end{eqnarray}
where 
$\BP_{1,k+1}\big( \bar{x},\bar{y}, \bar{b} \big)$ denotes the indicator function of the event
\begin{eqnarray}
 & & \Big\{ \, \larg{x_i}{y_i}{b_i} \big(i,s\big) >  \larg{x_{i+1}}{y_{i+1}}{b_{i+1}} \big( i+1,s\big) \, \, \, \, \forall \, \, \big(s,i\big) \in [-2\eln,2\eln] \times \intint{k-1} \, \Big\} \label{e.intextavoid} \\
 & \cap &  \Big\{ \,  \larg{x_k}{y_k}{b_k} \big(k,s \big) >  \mc{L}_n \big( k+1,s \big)\, \, \, \, \forall \, \, s  \in [-2\eln,2\eln]   \,  \Big\}  \, ; \nonumber
\end{eqnarray}
the role of the $\mc{F}$-measurable quantity $Z \in (0,\infty)$ in~(\ref{e.gdens.first.new}) to normalize $h$ so that it is the density function of a probability measure. This means that $Z$ is specified by
\begin{eqnarray}
 Z & = & \int \prod_{i=1}^k \, \exp\bigg\{-\tfrac{1}{2(\gfl + 2\eln)}\Big(\lppls_n\big(i,-2\eln  \big) - x_i \Big)^2 - \, \tfrac{1}{2(\gfr - \gfl)} \big( x_i - y_i \big)^2 \label{e.zint.new} \\ 
 & & \qquad \qquad \qquad \qquad \qquad \qquad  - \, \, \tfrac{1}{2(2\eln  - \gfr)} \Big( \lppls_n\big(i,2\eln  \big) - y_i \Big)^2 \bigg\} \,\cdot\,\BP_{1,k+1} \big( \bar{x},\bar{y},\bar{b} \big)
  \, \dd \bar{x} \dd \bar{y} \dd \bar{b} \, , \nonumber 
\end{eqnarray}
where the integral is over the outcome set and $\dd \bar{x} \dd \bar{y} \dd \bar{b}$ denotes integration with respect to the reference measure.

In other words, the witness of $\mc{F}$, in considering the eventuality that the $\mc{F}$-random data adopts the value $\big(\bar{x},\bar{y},\bar{b}\big)$, determines that this outcome occurs with density $h$ evaluated at this point. Two considerations, that may be labelled {\em kinetic} and {\em potential}, contribute to the formula for $h$: the first, represented by the omission of the factor $\BP_{1,k+1} \big( \bar{x},\bar{y},\bar{b} \big)$ in~(\ref{e.gdens.first.new}), expresses the Gaussian costs associated to the placement of the endpoint vectors $\bar{x}$ and $\bar{y}$ at $\gfl$ and $\gfr$.

The second, potential, factor is concerned with the need to check that the line ensemble that results from the placement of curves dictated by the triple~$\big(\bar{x},\bar{y},\bar{b}\big)$
observes the curve avoidance requirements of the ensemble.
 It is useful to categorize the avoidance conditions that are expressed by the equation $\BP_{1,k+1} \big( \bar{x},\bar{y},\bar{b} \big) = 1$. Before we do this, a word on our notation: we write $\BP_{1,k+1}$ to indicate that all of the top $k+1$ curves are implicated in the constraints; indeed, the concerned event, in~(\ref{e.intextavoid}), may be viewed as the intersection of an internal avoidance constraint involving the top $k$ curves (that are random for the witness), and an external constraint that stipulates avoidance of the $k\textsuperscript{th}$ curve with the non-random boundary condition $\mc{L}_n\big( k+1,\cdot \big):[-2\eln,2\eln] \to \R$.

For $A \subseteq [-2\eln ,2\eln ]$, write $\BP_{1,k+1}^A\big( \bar{x},\bar{y},\bar{b} \big)$ for the indicator function of the event 
\begin{eqnarray*}
& & \Big\{ \, \larg{x_i}{y_i}{b_i} \big(i,s\big) >  \larg{x_{i+1}}{y_{i+1}}{b_{i+1}} \big( i+1,s\big) \, \, \, \, \forall \, \, \big(s,i\big) \in A \times \intint{k- 1} \, \Big\} \\
&  \cap & \Big\{ \,  \larg{x_k}{y_k}{b_k} \big(k,s \big) >  \mc{L}_n \big( k+1,s \big)\, \, \, \, \forall \, \, s  \in A  \,  \Big\}  \, .
\end{eqnarray*}

The quantity $\BP_{1,k+1}\big( \bar{x},\bar{y}, \bar{b} \big)$ is thus alternatively denoted $\BP_{1,k+1}^{[-2\eln,2\eln]}\big( \bar{x},\bar{y}, \bar{b} \big)$.
It equals  the product 
$$
\BP_{1,k+1}^{[-2\eln,\gfl]} \cdot \BP_{1,k+1}^{[\gfl,\gfr]} \cdot \BP_{1,k+1}^{[\gfr,2\eln]} 
$$
evaluated at $\big( \bar{x},\bar{y}, \bar{b} \big)$.

We may think of the witness of $\mc{F}$ as testing the viability of the data $(\bar{x},\bar{y},\bar{b})$ by performing two checks: a {\em side intervals} test, which checks that both 
$\BP_{1,k+1}^{[-2\eln,\gfl]}$ and $\BP_{1,k+1}^{[\gfr,2\eln]}$ equal one; and a {\em middle interval} test, which checks that  $\BP_{1,k+1}^{[\gfl,\gfr]}$  also equals one. Each of these conditions is determined by $(\bar{x},\bar{y},\bar{b})$  and the $\mc{F}$-measurable data to which the witness is privy. In fact, the two criteria in the side intervals test are determined respectively by $\bar{x}$ and $\bar{y}$ alone (alongside the data in $\mc{F}$). We demonstrate this fact now, by providing an explicit characterization of when the two side interval subtests are met in terms of the values of $\bar{x}$ and $\bar{y}$.

\begin{lemma}\label{l.glv}
There exists $\gbarxmin \in \R_\geq^k$ such that $\bar{x} \in \R^k$ satisfies $\BP_{1,k+1}^{[-2\eln,\gfl]}(\bar{x},\bar{y},\bar{b}) =1$ if and only if $\bar{x} - \gbarxmin \in (0,\infty)^k_>$, i.e., this vector is a $k$-decreasing list of positive elements. In particular, the value of $(\bar{y},\bar{b}) \in \R^k_\geq \times \mc{C}_{0,0}\big([\gfl,\gfr],\R \big)^k$ plays no role in determining whether the left side interval test is met.
Similarly, there exists $\gbarymin \in \R_\geq^k$ such that $\bar{y} \in \R^k$ satisfies $\BP_{1,k+1}^{[\gfr,2\eln]}(\bar{x},\bar{y},\bar{b}) =1$ if and only if $\bar{y} - \gbarymin \in (0,\infty)^k_>$.
\end{lemma}
\noindent{\em Remark.} The lemma demonstrates that the admissible set of $\bar{x}$-locations is a set in $\R^k$ with a `lower' extreme point at $\gbarxmin$. We will encounter a similar circumstance also in regard to later $\sigma$-algebras that contain $\mc{F}$, and we will use this corner notation for them as well. The bar notation for $k$-vectors is being used for $\gbarxmin$ and $\gbarymin$ in a slightly altered form, so that it does not look ugly.

\medskip 
 
\begin{figure}[ht]
\begin{center}
\includegraphics[height=10cm]{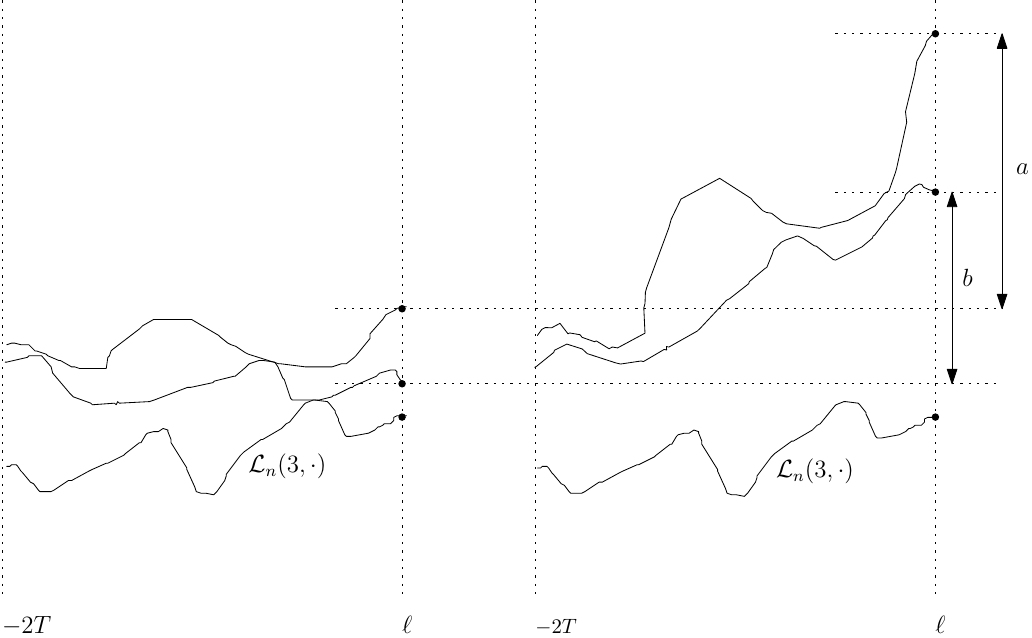}
\caption{The proof of Lemma~\ref{l.glv} is illustrated for $k=3$. The top two curves on the left are  $\larg{\gxmin_1}{y_1}{b_1} \big( 1,\cdot \big)$ and  $\larg{\gxmin_2}{y_2}{b_2} \big( 2,\cdot \big)$, where the $y$ and $b$ data is irrelevant. The top pair of beads over $\ell$ in the left sketch are located at $\gbarxmin$. 
On the right, this pair has been displaced to $\gbarxmin + (a,b)$, where $(a,b)$ represents a generic point in $(0,\infty)^2_>$; thus, the new pair verifies  $\BP_{1,3}^{[-2\eln,\gfl]}\big(\gbarxmin + (a,b)\big) =1$.} 
\label{f.leftcorner}
\end{center}
\end{figure}

\noindent{\bf Proof of Lemma~\ref{l.glv}.} The second statement is proved similarly to the first and we prove only the first, illustrating the proof with Figure~\ref{f.leftcorner}.
It is natural to think of specifying the components $\gxmin_i$, $i \in \intint{k}$, in decreasing order of the index~$i$. 
Thus take $i = k$ to begin. 
 Since the curve $\mc{L}_n\big( k+1,\cdot \big)$, and the $\mc{F}$-measurable bridge $\mc{L}_n^{[-2\eln,\gfl]}\big( k,\cdot \big):[-2\eln,\gfl] \to \R$ 
are almost surely continuous functions,  there is a unique value $\gxmin_k \in \R$ such that the curve
$$
\textrm{$\larg{\gxmin_k}{y_k}{b_k} \big( k,s \big)$ touches, but does not cross underneath,  $\mc{L}_n(k+1,s)$ on the interval $s \in [-2\eln,\gfl]$.} 
$$
(The values of $y_k$ and $b_k$ are irrelevant here.) Similarly, we may choose $\gxmin_{k-1} \in \R$ to be the unique value such that  
$$
\textrm{$\larg{\gxmin_{k-1}}{y_{k-1}}{b_{k-1}} \big( k-1,s \big)$ touches, but does not cross,  $\larg{\gxmin_k}{y_k}{b_k} \big( k,s \big)$  for $s \in [-2\eln,\gfl]$.} 
$$
Iteratively, we construct the vector $\gbarxmin$. Note that this vector specifies the values at $\gfl$ of a collection of curves $\Big\{ \larg{\gxmin_i}{y_i}{b_i} \big( i,s \big) : i \in\intint{k} \Big\}$ that do not cross each other for $s \in [-2\eln,\gfl]$; again, it makes no difference what the values of $y_i$ and $b_i$ are. In this way, we see that $\gbarxmin \in \R^k_\geq$.  

If $\bar{x} \in \R^k$ is to satisfy the non-touching condition that  $\BP_{1,k+1}^{[-2\eln,\gfl]}(\bar{x},\bar{y},\bar{b}) =1$, then $x_k > \gxmin_k$ is required, in order that the curves with indices $k$ and $k+1$ not touch. Moreover, if we begin by considering the curves' location dictated by the vector $\gbarxmin$ at $\gfl$, then  an upward push of $x_k - \gxmin_k$ made to curve $k$ must also be delivered to all curves of lower index, in order that no pair of consecutive curves among them begin to cross. After these equal pushes, the curves with indices $k-1$ and $k$ remain in contact, the higher requiring a further upward push, delivered via an increase in the value of $x_{k-1}$, with the lower indexed curves again being subjected to the same push via a similar increase in order to maintain non-crossing. Proceeding through the indices $i \in \intint{k}$ in decreasing order,
we see then that the difference vector $\bar{x} - \gbarxmin$ must be strictly increasing in its reverse-ordered components if and only if the condition $\BP_{1,k+1}^{[-2\eln ,\gfl]}(\bar{x},\bar{y},\bar{b}) =1$ is to be satisfied. \qed

\medskip

The middle interval test is more naturally analysed using the one-piece list presentation. In order to explain why, and to present a satisfying characterization of the overall condition $\BP_{1,k+1}^{[-2\eln,2\eln]} = 1$,
we begin by equipping the probability space carrying the law $\PP$ 
with some auxiliary random variables. 
The vectors 
$$
\ovbar{\mc{L}}_n(-2\eln) = \Big( \mc{L}_n\big(i,-2\eln  \big): i \in \intint{k} \Big) \, \, \, \, \textrm{and} \, \, \, \, 
\ovbar{\mc{L}}_n(2\eln) = \Big( \mc{L}_n\big(i, 2\eln  \big): i \in \intint{k} \Big)
$$
are $\mc{F}$-measurable random variables valued in $\R^k_>$.
Denote these two vectors temporarily by $\bar{u}$ and $\bar{v}$, the lower-case notation reflecting the deterministic status of these vectors in the eyes of the witness of~$\mc{F}$.
Recall that $\mc{B}_{k;\bar{u},\bar{v}}^{[-2\eln,2\eln]}$ denotes a collection $B(i,\cdot):[-2\eln,2\eln] \to \R$ of independent Brownian bridges indexed by $i \in \intint{k}$  with endpoints $B(i,-2\eln) = u_i$ and $B(i,2\eln) = v_i$. We augment the probability space that carries the law $\PP$ with copies of these bridge laws. 
(How this construction is carried out is irrelevant for our purpose, but for example one may take a single copy of the law $\mc{B}_{k;\bar{0},\bar{0}}^{[-2\eln,2\eln]}$ that is independent of~$\mc{F}$ and then affinely shift each of its~$k$ curves in an $\mc{F}$-determined manner in order to produce the desired distribution $\mc{B}_{k;\bar{u},\bar{v}}^{[-2\eln,2\eln]}$.) 
We will be concerned with these new processes 
only via their form on the subset $\intint{k} \times [\gfl,\gfr]$, and will denote them there by the symbol $\wien$. In this way, the witness of $\mc{F}$ constructs an ensemble $\wien:\intint{k} \times [\gfl,\gfr] \to \R$ with the marginal law on $[\gfl,\gfr]$ of the above bridge ensemble.
For reasons to be explained shortly, the ensemble~$\wien$ will be called the {\em Wiener candidate}. 
We may express the new ensemble in the two-piece list presentation, at the same time recalling a little notation:
\begin{itemize}
\item we write $\ovbar{\wien}(\gfl) = \big( \wien(i,\gfl):i \in \intint{k} \big)$ and $\ovbar{\wien}(\gfr) = \big( \wien(i,\gfr):i \in \intint{k} \big)$;
\item and $\wien^{[\gfl,\gfr]}:  \intint{k} \times [\gfl,\gfr] \to \R$ for the standard bridge ensemble on $[\gfl,\gfr]$ formed from $\wien$.
\end{itemize}

Note the relationship between this construction and the two-piece list Radon-Nikodym derivative discussion. Under the law $\PP_\mc{F}$,
the random triple of $k$-vectors $\big(\ovbar{\wien}(\gfl),\ovbar{\wien}(\gfr), \wien^{[\gfl,\gfr]}\big)$ has Radon-Nikodym derivative evaluated at $(\bar{x},\bar{y},\bar{b})$ with respect to the reference measure given by~(\ref{e.gdens.first.new}) where the factor of $\BP_{1,k + 1} \big( \bar{x},\bar{y}, \bar{b} \big)$ is omitted, (and where the normalization $Z$ is now specified by~(\ref{e.zint.new}) with this same factor also omitted). 
Kinetic and potential considerations dictate the law of the top $k$ curves of $\mc{L}_n$ for the witness of $\mc{F}$. The witness views the Wiener candidate $\wien:\intint{k} \times [\gfl,\gfr] \to \R$ -- or, alternatively represented, the triple $\big(\ovbar{\wien}(\gfl),\ovbar{\wien}(\gfr),\wien^{[\gfl,\gfr]}\big)$ -- 
as the random process 
that results when kinetic costs are considered but potential constraints are neglected.

For this reason, the conditional distribution under $\PP_\mc{F}$
of $\mc{L}_n:\intint{k} \times [\gfl,\gfr] \to \R$ equals the conditional law of the ensemble $\wien:\intint{k} \times [\gfl,\gfr] \to \R$ under $\PP_\mc{F}$
conditioned on the potential constraint that 
$$
\BP^{[-2\eln,2\eln]}_{1,k+1} \big( \ovbar{\wien}(\gfl),\ovbar{\wien}(\gfr), \wien^{[\gfl,\gfr]} \big) = 1 \, .
$$
When we wish to think of the question of whether this constraint is satisfied using the one-piece list presentation, we will instead write 
$$
\BP^{[-2\eln,2\eln]}_{1,k+1} \big( \wien \big) = 1 \, ,
$$
with the use of a single argument, rather than a triple, indicating the meaning. Similarly, of course, if we replace the superscript by a set $A \subset [-2\eln,2\eln]$. For example, the middle interval test
\begin{equation}\label{e.onepiecemiddle}
\BP_{1,k+1}^{[\gfl,\gfr]} \big( \wien \big) = 1   
\end{equation}
is more naturally expressed in the one-piece presentation, when it amounts to checking the ordering of the curves of $\wien$ and that $\wien(k,\cdot)$ exceeds $\mc{L}_n(k+1,\cdot)$ on $[\gfl,\gfr]$.

\begin{figure}[ht]
\begin{center}
\includegraphics[height=9cm]{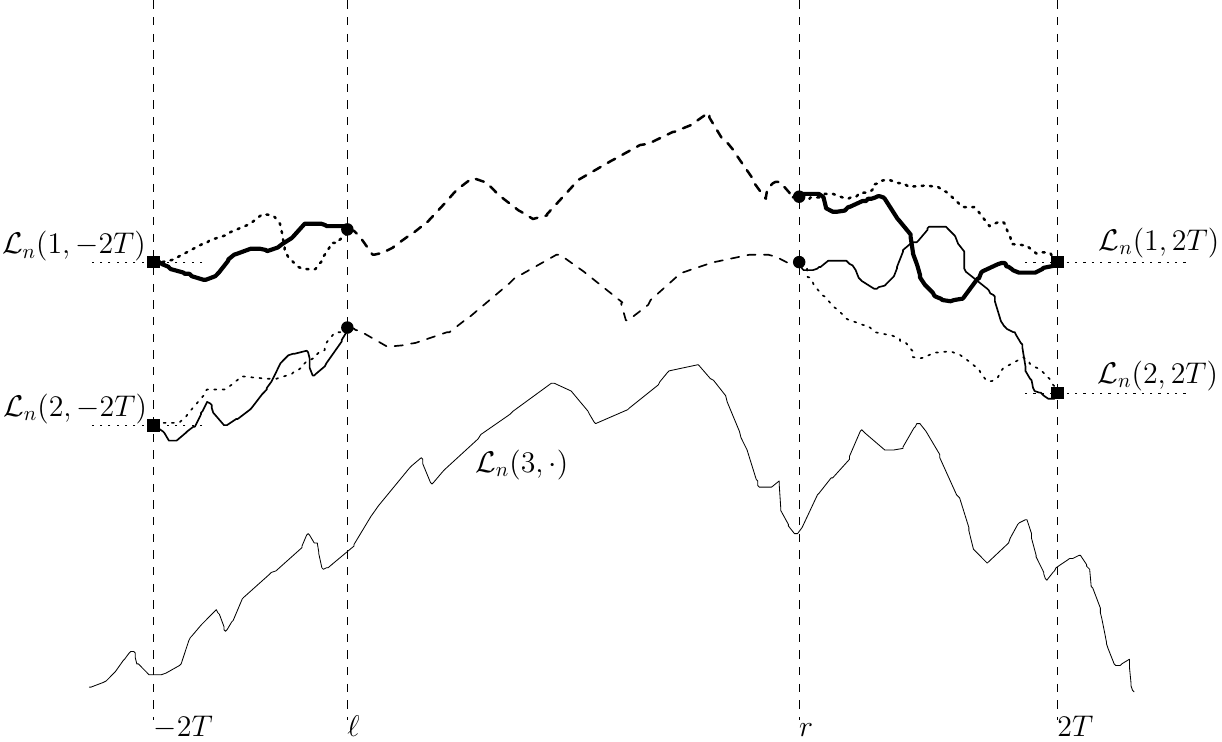}
\caption{The perspective of the witness of $\mathcal{F}$ and the construction of the Wiener candidate are depicted for $k=2$. The two curves on $[-2\eln,2\eln]$ that are dotted, dashed and then dotted again are samples of the law $\mc{B}_{2;\bar{u},\bar{v}}^{[-2\eln,2\eln]}$, with $(u_1,u_2) = \big( \mc{L}_n(1,-2\eln),  \mc{L}_n(2,-2\eln) \big)$ and $(v_1,v_2) = \big( \mc{L}_n(1,2\eln) ,  \mc{L}_n(2,2\eln) \big)$. The marginal of these curves on $[\ell,r]$, i.e., their dashed sections, forms the Wiener candidate $W$. 
The thicker solid black curve on $[-2\eln,\ell]$ is the affine translation of $\mc{L}_n^{[-2\eln,\ell]}$ with left endpoint $\mc{L}_n(1,-2\eln)$ and right endpoint $W(1,\ell)$. Similarly on the right, and for the thinner solid black second curve.
The black beads on the vertical line with $x$-coordinate $\ell$, located at heights $\big( W(1,\ell),W(2,\ell) \big)$, thus dictate the form of the
black curves to the left via affine translation subject to the fixed left endpoints; and similarly of course on the right.
Note that 
$\lefta \big(  \ovbar{\wien}(\gfl) \big) = 1$,
 $\mida ( \wien ) = 1$ and
  $\righta \big(  \ovbar{\wien}(\gfr) \big)  = 0$, and that the Wiener candidate is unsuccessful in this instance.}
\label{f.wienercandidate}
\end{center}
\end{figure}

That the conditional law of $\mc{L}_n: \intint{k} \times [\gfl,\gfr] \to \R$ under $\PP_\mc{F}$ may be obtained by conditioning $\wien$ helps to explain 
the reason for the Wiener candidate name. The witness of~$\mc{F}$ constructs $\wien$ and regards it as a candidate for the process $\mc{L}_n:\intint{k} \times [\gfl,\gfr] \to \R$:
the candidate is {\em successful} if the condition~$M_{1,k+1}^{[-2\eln,2\eln]}\big( \ovbar{\wien}(\gfl) , \ovbar{\wien}(\gfr) , \wien^{[\gfl,\gfr]} \big) = 1$ is met; indeed, under $\PP_\mc{F}$, the process $\wien$ conditioned on success  has the law of $\mc{L}_n$. See Figure~\ref{f.wienercandidate}.

We may summarise our understanding about how to check the outcome of the examination to which the candidate is subject by writing the examination indicator function in the form
$$
M_{1,k+1}^{[-2\eln,2\eln]}\big( \ovbar{\wien}(\gfl) , \ovbar{\wien}(\gfr) , \ovbar{\wien}^{[\gfl,\gfr]} \big)  = \lefta \big(  \ovbar{\wien}(\gfl) \big) \cdot \mida \big(  W \big) \cdot \righta \big(  \ovbar{\wien}(\gfr) \big)  \, ,   
$$
where the three right-hand factors are indicator functions whose meaning we now review.

The left and right side interval tests may be written
$$
\lefta \big(  \ovbar{\wien}(\gfl) \big)  = 
M_{1,k+1}^{[-2\eln,\gfl]}\big( \ovbar{\wien}(\gfl) , \ovbar{\wien}(\gfr) , \wien^{[\gfl,\gfr]} \big)
$$
and 
$$
\righta \big(  \ovbar{\wien}(\gfr) \big)  = 
M_{1,k+1}^{[\gfr,2\eln]}\big( \ovbar{\wien}(\gfl) , \ovbar{\wien}(\gfr) , \wien^{[\gfl,\gfr]} \big) \, .
$$
Note that Lemma~\ref{l.glv}
indeed implies that it is only the first and second argument in the respective triples that is used to determine the value of these right-hand sides. 

The middle interval test is naturally examined by using the one-piece list:
$$
  \mida ( \wien ) = 
\BP_{1,k+1}^{[\gfl,\gfr]} \big( \wien \big) \, .
$$

A final remark about the one-piece list perspective. For $i \in \intint{k}$, the witness of $\mc{F}$ may seek to reconstruct the curve $\mc{L}_n(i,\cdot)$ on the basis that the unknown data associated to this curve equals $f_i:[\gfl,\gfr] \to \R$. 
Denoting the reconstructed curve by $\mc{L}_n^{f_i}(i,\cdot):[-\xnmac,\infty) \to \R$, with the use of a single argument superscript indicating that the one-piece list presentation is being used, we have that 
\begin{equation}\label{e.reconcases.onepiece}
 \mc{L}_n^{f_i} \big( i, s \big) = \begin{cases}
  \, \mc{L}_n\big(i,s \big) \, & s \in [-\xnmac,-2\eln) \, , \\
 \, \mc{L}_n^{[-2\eln,\gfl]}\big(i,s\big) + 
 \tfrac{\gfl - s}{\gfl + 2\eln}  \lppls_n\big(i, - 2\eln \big)
 + \tfrac{s + 2\eln}{\gfl + 2\eln}  f_i(\gfl)   \, & s\in \big[-2\eln,\gfl \big] \, , \\
 \, f_i(s)  \, & s \in \big[\gfl ,\gfr \big] \, , \\
 \, \mc{L}_n^{[\gfr,2\eln]}\big(i,s\big) + 
 \tfrac{s - \gfr}{2\eln - \gfr}  \lppls_n\big(i,2\eln\big)
 + \tfrac{2\eln - s}{2\eln - \gfr}  f_i(\gfr)   \, & s\in \big[\gfr,2\eln \big]  \\
   \, \mc{L}_n\big(i,s \big) \, & s \in (2\eln,\infty) \, .
\end{cases}
\end{equation}

\section{Applications of the Wiener candidate approach}\label{s.wienercandidateapplications}

We are now ready to present the proofs of the curve closeness lower bound Theorem~\ref{t.airynt.lb} and the local maximal fluctuation Theorem~\ref{t.aestimate} (as well as Proposition~\ref{p.aestimateinference}).
These results are applications of missing closed middle reconstruction (via the Wiener candidate)
with the fixed parameter choice $\big(\eln,\gfl,\gfr \big) = \big(1,-1,1 \big)$.

There is a further similarity between the two proofs. 
Before we give the proofs, which appear in Sections~\ref{s.kcurvelowerbound} and~\ref{s.aestimate},
we describe this common element.
We introduce a {\em favourable} event~$F_t$, specified by a positive parameter $t$. This event may be detected by the witness of~$\mc{F}$, since it is
$\mc{F}$-measurable, and when the event does occur, the Wiener candidate behaves pleasantly under the law $\PP_{\mc{F}}$.

We now specify the favourable event; show in Lemma~\ref{l.favs} that it is typical; and, in Lemma~\ref{l.bcgood}, show that Wiener candidate success $M_{1,k+1}^{[-2,2]}(\wien) = 1$
is typical under $\PP_\mc{F}$ when the event occurs.

For $t > 0$, we define $F_t$ to be the event that
\begin{itemize}
\item $\mc{L}_n(i,-2) \in [-t,t]$ and $\mc{L}_n(i,2) \in [-t,t]$ for $i \in \intint{k}$;
\item $\gxmin_i \in [-t,t]$ and $\gymin_i \in [-t,t]$ for $i \in \intint{k}$;
\item and $\mc{L}_n(k+1,x) \leq t$ for $x \in [-1,1]$.
\end{itemize}
\begin{lemma}\label{l.favs}
For any $t > 0$, the event $F_t$ is $\mc{F}$-measurable. 
If $t \geq 2^{5/2}$,
then
\begin{equation}\label{e.favs}
 \PP \big( F_t^c \big) \leq 14 C_k \exp \big\{ - c_k t^{3/2}/8 \big\} 
\end{equation}
for all $n \geq k \vee (\rsc/3)^{-2(\phimac_1 \wedge \phimac_2)^{-1}} \vee 6^{2/\delta} \vee (t/2)^{1/\delta}$ (with $\delta = \phimac_1/2 \wedge \phimac_2/2 \wedge \phimac_3$).
Specifically, when  $t \geq \big( 8 c_k^{-1} \log (28 C_k) \big)^{2/3} \vee 2^{5/2}$, $\PP \big( F_t^c \big) \leq 1/2$.
\end{lemma}
\noindent{\bf Proof.} By ensemble ordering and Proposition~\ref{p.othercurves} applied with $x = \pm 2$ and $s = t - 2^{3/2}$, the lower tail conditions in the first bullet point fail with probability at most $4 C_k \exp \big\{ - c_k (t - 2^{3/2})^{3/2} \big\}$. The proposition applies because $t - 2^{3/2} \in [0,2n^\delta]$, as well as $2 \leq c/2 \cdot  n^\delta$. The upper tail conditions in this bullet point fail with probability at most $2\rsC \exp \big\{ -\rsc t^{3/2} \big\}$ due to ensemble ordering and $\rmreg(3)$
in light of $2 \leq \rsc n^{\phimac_2}$.

Regarding the next point, note that
$$
 \big[ \xmin_k , \xmin_1 \big] \subseteq \big[ \mc{L}_n(k+1,-1) , \mc{L}_n(1,-1) \big] \, .
$$ 
Bounds on the lower tail of $\mc{L}_n(k+1,-1)$ and on upper bound 
of $\mc{L}_n(1,-1)$ will yield an estimate on the probability of failure of the $\gfl$-condition in the second point. Similarly of course for the $\gfr$-condition. Arguing as in the preceding paragraph with $s = t - 2^{-1/2}$, we find that the probability that the event in this bullet point does not occur is at most
$$
4 C_k \exp \big\{ - c_k (t - 2^{3/2})^{3/2} \big\} + 2\rsC \exp \big\{ -\rsc t^{3/2} \big\} 
$$
since $t - 2^{-1/2} \in [1,2n^\delta]$.

The third point is treated by means of ensemble ordering and the `no big max' Lemma~\ref{l.nobigmax}. The failure probability is at most
$6 \rsC \exp \big\{ - \rsc t^{3/2}/8 \big\}$.

Using $C_k \geq \rsC$, $c_k \leq \rsc$ and $t \geq 2^{5/2}$, failure among any of the three points is seen to occur with probability at most $14 C_k \exp \big\{ - c_k t^{3/2}/8 \big\}$.
Since the final assertion of the lemma is an immediate consequence of~(\ref{e.favs}), this completes the lemma's proof. \qed

\medskip

\begin{lemma}\label{l.bcgood}
If $k \geq 1$ and  $t \geq 2 (\log 2)^{1/2}$, then
$$
\PP_{\mc{F}} \Big(  \,  M_{1,k+1}^{[-2,2]}(\wien) = 1  \, \Big) \geq 
 2^{-3k/2} \pi^{-k}   \exp \big\{ - 4k(k+2)^2 t^2 \big\}\cdot {\bf 1}_{F_t} 
  \, .
$$
\end{lemma}
\noindent{\bf Proof.} 
Set $\Delta \subset \R^k$ equal to the box $[2k,2k+1] \times \cdots \times [4,5] \times [2,3]$.
Let $G_t\big( \wien \big)$ denote the event that  $\ovbar{\wien}(-1)$ and $\ovbar{\wien}(1)$ belong to $t \cdot \Delta$. We will prove the lemma by deriving the stronger claim that 
\begin{equation}\label{e.strongerclaim}
\PP_{\mc{F}} \Big( G_t(\wien) \, , \,  M_{1,k+1}^{[-2,2]}(\wien) = 1 \Big) \geq 
 2^{-3k/2} \pi^{-k}   \exp \big\{ - 4k(k+2)^2 t^2 \big\}  \cdot {\bf 1}_{F_t} \, .
\end{equation}
To derive the claim, we begin by noting that 
\begin{equation}\label{e.gt}
\PP_{\mc{F}} \Big( G_t(\wien) \, , \,  M_{1,k+1}^{[-2,2]}(\wien) = 1 \Big) \cdot {\bf 1}_{F_t} = 
\PP_{\mc{F}} \big(   \, G_t(\wien) \, \big) \cdot
\PP_{\mc{F}} \Big(  \,  M_{1,k+1}^{[-2,2]}(\wien) = 1 \, \Big\vert \, G_t(\wien) \, \Big) \cdot {\bf 1}_{F_t}  \, .
\end{equation}
We are working with missing closed closed reconstruction with $\big(\eln,\gfl,\gfr\big)$ set equal to $\big( 1,-1,1 \big)$. Recall that, under $\PP_\mc{F}$, the Wiener candidate $W:\intint{k} \times [-1,1] \to \R$ has the marginal law on $[-1,1]$ of the Brownian bridge ensemble $\mc{B}_{k;\bar{\mc{L}}_n(-2),\bar{\mc{L}}_n(2)}^{[-2,2]}$. Applying Corollary~\ref{c.brbr} in regard to the curve $W(i,\cdot)$, $i \in \intint{k}$, with $l=2$, $a = -2$, $b=2$, $A_0 \subset [-t,t]$, $A_1 = A_2 = t\cdot[2k+2-2i,2k+3 - 2i]$
and $A_3 \subset [-t,t]$,
we learn that 
$$
\PP_{\mc{F}} \big(   \, G_t(\wien) \, \big)  \geq 2^{-k/2} \pi^{-k} \prod_{i=1}^k \exp \bigg\{ \frac{\big(t + t(2k+3 - 2i) \big)^2}{2} + \frac{t^2}{4}  +  \frac{\big(t + t(2k+3 - 2i) \big)^2}{2}  \bigg\}  \cdot {\bf 1}_{F_t} \, ,
$$
where we made use of the length of each interval $A_i$ being at least one (which is due to $t \geq 1$).
Thus,
\begin{equation}\label{e.bcgood1}
\PP_{\mc{F}} \big(   \, G_t(\wien) \, \big) \geq  2^{-k/2} \pi^{-k}   \exp \big\{ - 4k(k+2)^2 t^2 \big\}   \cdot {\bf 1}_{F_t}  \, .
\end{equation}
Consider now the conditional law 
$\PP_{\mc{F}} \big(  \, \cdot \, \big\vert \, G_t(\wien) \, \big)$. For the demand that  $M_{1,k+1}^{[-2,2]}(\wien) = 1$ to be met for this conditioned process, it is enough that each standard bridge $W^{[-1,1]}(i,\cdot)$ have absolute value whose supremum is less than $t/2$. Since these bridges independently have law $\mc{B}_{1;0,0}^{[-1,1]}$ under the measure in question, Lemma~\ref{l.maxfluc} implies that  
\begin{equation}\label{e.bcgood2}
\PP_{\mc{F}} \Big(  \,  M_{1,k+1}^{[-2,2]}(\wien) = 1 \, \Big\vert \, G_t(\wien) \, \Big) \geq \big( 1 - 2 e^{-t^2/4} \big)^k  \cdot {\bf 1}_{F_t} \geq 2^{-k}  \cdot {\bf 1}_{F_t}  
\end{equation}
since $t \geq 2 (\log 2)^{1/2}$.

Note that from~(\ref{e.bcgood1}) and~(\ref{e.bcgood2})
follow the claim~(\ref{e.strongerclaim}) and Lemma~\ref{l.bcgood}. \qed

\subsection{Proving Theorem~\ref{t.airynt.lb}, the lower bound on the $k$-curve closeness probability}\label{s.kcurvelowerbound}

 First, we reduce  
the theorem to the following assertion.
\begin{proposition}\label{p.airynt.lb}
For $\bar\phimac \in (0,\infty)^3$, $C,c > 0$ and $n \in \N$, let 
$$
\mc{L}_n:\intint{n} \times \big[-\xnmac,\infty\big) \to \R  
$$ 
be a    $\big(\bar\phimac,\rsc,\rsC\big)$-regular ensemble defined under the law~$\PP$. 
 Set $\delta = \phimac_1/2 \wedge \phimac_2/2 \wedge \phimac_3$, and, for $k \in \N$, $s_k = \big( 8 c_k^{-1} \log (28C_k) \big)^{2/3} \vee 2^{5/2}$. Then, for $\e \in \big( 0, k^{-2} s_k^{-1}/4 \big)$,
$$
  \PP \Big( \close\big( k ; \mc{L}_n , 0 , \e \big) \Big) 
  \geq    e^{- 52 s_k^2 k^3} \, \e^{k^2-1} 
$$ 
whenever  $n \geq k \vee (\rsc/3)^{-2 (\phimac_1 \wedge \phimac_2)^{-1}} \vee 6^{2/\delta} \vee (s_k/2)^{1/\delta}$. 
\end{proposition}
\noindent{\bf Proof of Theorem~\ref{t.airynt.lb}.}
We use parabolic invariance Lemma~\ref{l.parabolicinv}
and Proposition~\ref{p.airynt.lb} with the parameter choice $(\bar\phimac,c/2,C)$, while noting that the value of $c_k$ in~(\ref{e.littlec}) determined by $c/2$ is at least one-half of the value determined by $c$. \qed

\medskip

Preparing to prove Proposition~\ref{p.airynt.lb}, 
recall that the event $F_t$ is $\mc{F}$-measurable, and that we suggested at the beginning of Section~\ref{s.wienercandidateapplications}
that we might work by examining the law $\PP_\mc{F}$
for $\mc{F}$-data causing the occurrence of $F_t$.
Actually, this is not quite our approach (though this description is correct for the proof of Theorem~\ref{t.aestimate}). Rather, we now introduce a new event $H_{t}$ (which is not in fact~$\mc{F}$-measurable), establish in Lemma~\ref{l.bcgood.little} that the event is typical, and then finish the proof of  Proposition~\ref{p.airynt.lb} by showing that 
$k$-curve closeness at zero is typical when $H_t$ occurs.
We consider $F_t$ merely as a tool for proving that $H_t$ is typical.

As we did in the proof of Lemma~\ref{l.bcgood}, we let $\Delta \subset \R^k$ denote the box $[2k,2k+1] \times \cdots \times [4,5] \times [2,3]$.
For $s > 0$, let $H_s = H_{k,s}$ 
denote the event that
\begin{itemize}
\item the $k$-vectors $\ovbar{\mc{L}}_n(-1)$ and $\ovbar{\mc{L}}_n(1)$ belong to $s \cdot \Delta \,$;
\item and $\sup_{x \in [-1,1]} \mc{L}_n(k+1,x) \leq s$.
\end{itemize}
Developing the notation used in the proof of Lemma~\ref{l.bcgood}, and letting $G_t\big( \mc{L}_n \big)$ denote the event that  $\ovbar{\mc{L}}_n(-1)$ and $\ovbar{\mc{L}}_n(1)$ belong to~$t \cdot \Delta$, we  use the abbreviation $G_t = G_t\big( \mc{L}_n \big)$, and find that 
 $G_t \cap F_t \subseteq H_t$. Thus,
\begin{eqnarray}
\PP \big( H_t \big) \geq  \PP \big( G_t \cap F_t \big) = \EE \Big[ \PP_{\mc{F}} \big( G_t \big) \, {\bf 1}_{F_t} \Big]
 & = & \EE \bigg[ \PP_{\mc{F}} \Big( G_t(\wien) \, \Big\vert \,  M_{1,k+1}^{[-2,2]}(\wien) = 1 \Big)  \cdot {\bf 1}_{F_t} \, \bigg] \nonumber \\
  & \geq & \EE \bigg[ \PP_{\mc{F}} \Big( G_t(\wien) \, , \,  M_{1,k+1}^{[-2,2]}(\wien) = 1 \Big) \cdot {\bf 1}_{F_t} \, \bigg] \, . \label{e.htlb}
\end{eqnarray}

\begin{lemma}\label{l.bcgood.little}
If $t \geq \big( 8 c_k^{-1} \log (28 C_k) \big)^{2/3} \vee 2^{5/2}$, then
$$
\PP \big( H_{k,t} \big) \geq  \tfrac{1}{2} \cdot  2^{-3k/2} \pi^{-k} \, \exp \big\{ -  16   k^3 t^2 \big\}  
$$
provided that
 $n \geq k \vee (\rsc/3)^{-2(\phimac_1 \wedge \phimac_2)^{-1}} \vee 6^{2/\delta} \vee (t/2)^{1/\delta}$, where $\delta =  \phimac_1/2 \wedge \phimac_2/2 \wedge \phimac_3$.
\end{lemma}
\noindent{\bf Proof.} Note first that $\PP(F_t) \geq 1/2$, because we hypothesise the lower bound on $t > 0$
necessary to conclude this from  Lemma~\ref{l.favs}.
Using Lemma~\ref{l.bcgood} and $\PP(F_t) \geq 1/2$, we find that
$$
\PP \big( H_t \big) \geq \tfrac{1}{2} \cdot  2^{-k/2} \pi^{-k}   \exp \big\{ - 4k(k+2)^2 t^2 \big\}  \cdot  2^{-k} \, .
$$
Since $k \geq 2$, we obtain the lemma.  \qed

\medskip

\noindent{\bf Proof of Proposition~\ref{p.airynt.lb}.}
Consider a given instance of data $\ovbar{\mc{L}}_n(-1)$, $\ovbar{\mc{L}}_n(1)$ and $\mc{L}_n(k+1,\cdot):[-1,1] \to \R$.
Denote these pieces of data by $\bar{x}$, $\bar{y}$ and $f$.
Under the law $\PP$ conditionally on the observation of these data, the conditional distribution of $\mc{L}_n:\intint{k} \times [-1,1] \to \R$  is 
$$
\mc{B}_{k;\bar{x},\bar{y}}^{[-1,1]}\big( \cdot \, \big\vert \notouch^{[-1,1]}_f \big) \, .
$$
If the data $\big( \bar{x},\bar{y},f\big)$ is such that the event $H_{k,t}$ occurs for a value $t$ that is at least $2^{5/2}k$,
then this conditional distribution, after a downward shift of $\bar{x}$, $\bar{y}$ and $f$ by the amount $t + 2^{5/2}k$, satisfies the hypotheses of Proposition~\ref{p.expsimple}(2) with $K = 2k (t - 1)$. Applying this result, we find that
\begin{eqnarray*}
 & & \PP \Big( \, \close \big( \mc{L}_n ,0, \phi \big) \, \, \Big\vert \, \, H_{k,t}  \, \Big) \\ & \geq &
  \phi^{k^2-1} \cdot \big( 1 -  \pi^{-1/2}e^{-1} \big) \, \tfrac{1}{16}  (2k)^{-(k^2-1)}  \pi^{-(k-1)}  \exp \big\{ - 32 k^3   t^2 \big\}  \, , 
\end{eqnarray*}
where $\phi$ is assumed to be at most $\big( 4k^2(t-1) \big)^{-1}$.
We thus obtain Proposition~\ref{p.airynt.lb} by taking $\phi = \e$ and invoking Lemma~\ref{l.bcgood.little}, taking $t$ equal to the lowest value permitted in the proposition. \qed


In summary, the  curve closeness probability lower bound  Theorem~\ref{t.airynt.lb} has been proved
by reducing it to Proposition~\ref{p.airynt.lb}; and the latter result has been proved 
following the plan  
outlined in the second paragraph of Section~\ref{s.missingclosedmiddle}, the eventuality in question being $H_{k,t}$.

\subsection{Deriving Theorem~\ref{t.aestimate} on local maximal fluctuation}\label{s.aestimate}

 We continue to work with missing closed middle reconstruction with $\big(\eln,\gfl,\gfr\big) = \big( 1,-1,1 \big)$, and the Wiener candidate, in order to prove Theorem~\ref{t.aestimate} (and Proposition~\ref{p.aestimateinference}).  Our first application of the Wiener candidate approach, the proof of Theorem~\ref{t.airynt.lb}, made use of a specification of the parameter $t > 0$ that assured merely that the infimum over $n \in \N$ of the favourable event probability $\PP(F_t)$ was at least one-half.
As we now revisit this approach, we will make use of the same favourable event via Lemmas~\ref{l.favs} and~\ref{l.bcgood}, but instead specify $t > 0$ 
so that $\PP(F_t^c)$ 
decays to zero as the parameter $K > 0$ in  Theorem~\ref{t.aestimate} tends to infinity.

We begin the proof of Theorem~\ref{t.aestimate} by stating a result that will  emerge during the course of the derivation and which leads directly to Proposition~\ref{p.aestimateinference}.

\begin{proposition}\label{p.aestimate}
If $k \geq 1$, $x \in [-1,1]$, $\e \in (0,1/2)$ and $K \geq 2^{19/2} k^{1/2}(k+2)$, then,  with $t = 2^{-7} K (k+2)^{-1} k^{-1/2}$,
$$
 \PP \Big( \omega_k \big( \mc{L}_n , x,  \e \big)  \geq K\e^{1/2} \, , \, F_t \Big) \leq 
 2^{3k/2} \pi^{k}  k \cdot 60 K^{-1}  \exp \big\{ - 2^{-12}  K^2 \big\}   
$$
whenever  $n \geq k + 1$.
\end{proposition}
\noindent{\bf Proof of Proposition~\ref{p.aestimateinference}.}
By Lemma~\ref{l.parabolicinv}, it suffices to take $x = 0$. The result then follows from Proposition~\ref{p.aestimate}
and $\boundgood_{t + 2^{1/2}}(0) \subseteq F_t$ with $t = 2^{-7} K (k+2)^{-1} k^{-1/2}$. \qed

\medskip

\noindent{\bf Proof of Theorem~\ref{t.aestimate} and Proposition~\ref{p.aestimate}.}
In this argument, we will write  $\omega_k \big( E , x,  \delta \big)$ in place of  $\omega_{k,[x,x+\delta]} \big( E ,   \delta \big)$ since doing so simplifies the appearance of some expressions.

Note that 
\begin{eqnarray*}
 & & \PP \Big( \omega_k \big( \mc{L}_n , x,  \e \big) \geq K\e^{1/2} \, , \, F_t \Big) \\
 & \leq  & \EE \bigg[ \PP_{\mc{F}} \Big(  \omega_k \big( \wien , x,  \e \big) \geq K\e^{1/2} \, \Big\vert \,  M_{1,k+1}^{[-2\eln,2\eln]}(W) = 1 \Big) \cdot {\bf 1}_{F_t} \bigg] \\
 & \leq  & \EE \Bigg[ \frac{\PP_{\mc{F}} \Big(  \omega_k \big( \wien , x,  \e \big) \geq K\e^{1/2} \Big)}{\PP_{\mc{F}} \Big(  M_{1,k+1}^{[-2\eln,2\eln]}(\wien) = 1 \Big)} \cdot {\bf 1}_{F_t} \bigg]  \\
 & \leq &  2^{3k/2} \pi^{k}   \exp \big\{ 4k(k+2)^2 t^2 \big\} \EE \Big[ \PP_{\mc{F}} \Big(  \omega_k \big( \wien , x,  \e \big) \geq K\e^{1/2}\Big) \cdot {\bf 1}_{F_t} \Big]   \, , 
\end{eqnarray*}
where in the third inequality, we used~Lemma~\ref{l.bcgood}.

Recall that $\eln = 1$, $\gfl = -1$ and $\gfr = 1$. We leave the favourable event parameter $t > 0$ unspecified for now, but note that, when this event $F_t$ occurs, $\mc{L}_n \big( i,u \big) \in [-t,t]$ for $(i,u) \in \intint{k} \times \{ - 2  , 2 \}$. Thus, for any $t > 0$,
$$
  \PP_\mc{F} \Big(
 \omega_k \big( \wien , x,  \e \big) \geq K\e^{1/2}  
     \Big) \cdot   {\bf{1}}_{F_t}  \leq k \cdot \sup_{(u,v) \in [-t,t]^2} \mc{B}_{1;u,v}^{[-2,2]}\Big(  \omega_k \big( B , x,  \e \big) \geq K\e^{1/2} \Big) \, .
$$
For any $u,v \in \R$, Brownian bridges $B$ and $B'$ under the laws  $\mc{B}_{1;0,0}^{[-2,2]}$ and  $\mc{B}_{1;u,v}^{[-2,2]}$
may be coupled by affine shift; when they are, the bound
$\omega_1 \big( B' , x,  \e \big)
\leq \omega_1 \big( B , x,  \e \big) + \e  4^{-1} \vert v - u \vert$ results. This right-hand side is at most $\omega_1 \big( B , x,  \e \big) +  t \e/2$ when $u,v \in [-t,t]$.
Impose that $t\e \leq K \e^{1/2}$.
 Hence,  
\begin{eqnarray*}
 \PP_\mc{F} \Big(
 \omega_k \big( \wien , x,  \e \big) \geq K\e^{1/2}  
     \Big) \cdot   {\bf{1}}_{F_t} 
 &  \leq & k \cdot
 \mc{B}_{1;0,0}^{[-2,2]}\Big(  \omega_1 \big( B , x,  \e \big) \geq \tfrac{1}{2} \cdot K\e^{1/2} \Big) \\
 & = & k \cdot
 \mc{B}_{1;0,0}^{[0,1]}\Big(  \omega_1 \big( B ,  x/4 + 1/2  ,  \e/4 \big) \geq \tfrac{1}{4} \cdot K\e^{1/2} \Big) \, ,
\end{eqnarray*}
where since $x \in [-1,1]$, $x/4 + 1/2 \in [1/4,3/4]$.
Standard Brownian bridge $B:[0,1] \to \R$ may be represented in the form 
$$
B(t) \, = \, (1-t) \bmotion \big( \tfrac{t}{1-t} \big) \, \, ,  \, \, \, t \in [0,1] \, , 
$$
where $\bmotion:[0,\infty) \to \R$ is standard Brownian motion. Note that, when $r,s \in [0,7/8]$, $r < s$, 
$$
 \big\vert B(s) - B(r) \big\vert \, \leq \, \Big\vert  \bmotion \big( \tfrac{s}{1-s} \big) - \bmotion \big( \tfrac{r}{1-r} \big) \Big\vert + \big( s - r \big) \sup_{x \in [0,7/8]} \vert \bmotion \vert \big( \tfrac{x}{1-x} \big) \, .
$$
Since $[0,7/8]:x \to x/(1-x)$ has $64$ as a Lipschitz constant and $x/4 + 1/2 + \e/4 \leq 7/8$ since $\e \leq 1/2$, we find that
$$
 \omega_1 \big( B ,  x/4 + 1/2  ,  \e/4 \big)  \leq  \omega_1 \big(  X ,  \tfrac{1/2 + x/4}{1/2 - x/4}  ,  16 \e \big)  \, + \, \e/4 \cdot  \sup_{t \in [0,7]} \vert \bmotion  ( t ) \vert \, .
$$
Writing  $\PP$ for the probability measure associated to the Brownian motion $X$, we find that
\begin{eqnarray*}
 & & \mc{B}_{1;0,0}^{[0,1]}\Big(  \omega_1 \big( B ,  x/4 + 1/2  ,  \e/4 \big) \geq \tfrac{1}{4} \cdot K\e^{1/2} \Big) \\
 & \leq & \PP \Big(    \omega_1 \big(  X ,  \tfrac{1/2 + x/4}{1/2 - x/4}  ,  16 \e \big)  \geq  \tfrac{1}{8} \cdot K\e^{1/2} \Big)
 \, + \, \PP \Big(  \e/4 \cdot  \sup_{t \in [0,7]} \vert \bmotion  ( t ) \vert \geq  \tfrac{1}{8} \cdot K\e^{1/2} \Big) \\
  & \leq & 4 \, \nu_{0,16\e}\big(   \tfrac{1}{8} \cdot K\e^{1/2} , \infty \big) \, + \, 4 \, \nu_{0,7} \big( 2^{-1} K \e^{-1/2}   , \infty \big) \\
  & \leq & 4 (2\pi)^{-1/2} 32 K^{-1} \exp \big\{ - 2^{-1} (K/32)^2 \big\} + 4 (2 \pi)^{-1/2}  7^{1/2} \cdot 2 \e^{1/2} K^{-1} \exp \big\{ - 7^{-1} 8^{-1} K^2 \e^{-1} \big\}  \\
  & \leq & 60 K^{-1} \exp \big\{ - 2^{-11}  K^2 \big\} \, ,
\end{eqnarray*}
where the last inequality used $\e \leq 1$.

We find then that
$$
 \PP \Big( \omega_k \big( \mc{L}_n , x,  \e \big) \geq K\e^{1/2} \, , \, F_t \Big) \leq 
 2^{3k/2} \pi^{k}   \exp \big\{ 4k(k+2)^2 t^2 \big\} \cdot k \cdot 60 K^{-1} \exp \big\{ - 2^{-11}  K^2 \big\}
  \, . 
$$

We now set the favourable event parameter $t > 0$ so that  $4k(k+2)^2 t^2 =  2^{-12}  K^2$.
That is, $t = 2^{-7} K (k+2)^{-1} k^{-1/2}$.
Note that the condition that $t\e^{1/2} \leq K$, which we demanded earlier, is satisfied since $\e^{1/2} \leq 2^7 k^{1/2}(k+2)$ is implied by $k \geq 1$
and $\e < 1$. 
 We obtain 
$$
 \PP \Big( \omega_k \big( \mc{L}_n , x,  \e \big)  \geq K\e^{1/2} \, , \, F_t  \Big)  \leq  
 2^{3k/2} \pi^{k}  k \cdot 60 K^{-1}  \exp \big\{ - 2^{-12}  K^2 \big\}  \, ,
$$
which is Proposition~\ref{p.aestimate}.

Note further that
\begin{eqnarray*}
 \PP \Big( \omega_k \big( \mc{L}_n , x,  \e \big)  \geq
  K\e^{1/2} \Big) & \leq &  \PP \Big( \omega_k \big( \mc{L}_n , x,  \e \big)  \geq 
  K\e^{1/2} \, , \, F_t \Big) \, + \, \PP \big( F_t^c \big) \\
   & \leq & 
 2^{3k/2} \pi^{k}  k \cdot 60 K^{-1}  \exp \big\{ - 2^{-12}  K^2 \big\} \, + \,  14 C_k \exp \big\{ - c_k t^{3/2}/8 \big\} \\
  & \leq & 
 2^{3k/2} \pi^{k}  k \cdot 60 K^{-1}  \exp \big\{ - 2^{-12}  K^2 \big\} \\
 & & \qquad \qquad
  \, + \,  14 C_k \exp \big\{ - c_k 2^{-27/2}   K^{3/2}  k^{-3/4} (k+2)^{-3/2} \big\} \\
   & \leq & 
 \big( 2^{3k/2} \pi^{k}  k \cdot 60 K^{-1}     +   14 C_k \big) \exp \big\{ - c_k   k^{-9/4} 2^{-27/2}  3^{-3/2} K^{3/2}  \big\} \, ,
\end{eqnarray*}
where 
the second inequality uses Lemma~\ref{l.favs},
with the necessary bound  $t \geq 2^{5/2}$ taking the form  $K \geq 2^7 k^{1/2}(k+2) 2^{5/2}$ when  expressed in terms of $K$; and where, 
in the last inequality, we used $k \geq 1$ and $K \geq 2^{-3} c_k^2   k^{-9/2}  3^{-3}$ (which simply follows from $K \geq 1$ in view of $c_k \leq 1$) in the form $2^{-12}  K^2 \geq c_k   k^{-9/4} 2^{-27/2}  3^{-3/2} K^{3/2}$.

Also using $K \geq 1$, we obtain Theorem~\ref{t.aestimate}. \qed

\chapter{The jump ensemble method: foundations}\label{c.jumpensemblefoundations}

This chapter is devoted to introducing the general apparatus of this method; the next uses it to prove the $k$-curve closeness Theorem~\ref{t.airynt}(1) and the Brownian regularity Theorem~\ref{t.airytail.ln}.  
The present chapter has two sections. 
In Section~\ref{s.jumpensemblemethod},  we explain why missing closed middle reconstruction is inadequate for the purpose of proving such upper bounds as Theorem~\ref{t.airynt}(1). Faced with the need for a new technique, we refine this reconstruction approach  and so formulate the jump ensemble method, in the later part of  Section~\ref{s.jumpensemblemethod} and in Section~\ref{s.subpoly}.

\section{The jump ensemble method}\label{s.jumpensemblemethod}

In this section, we will present the apparatus of this general method for proving upper bounds on probabilities of events expressed in terms of the curves in Brownian Gibbs ensembles. 
\subsection{The need for a further method: how the Wiener candidate falls short}

We begin by advocating the need for a new approach 
by examining the $k$-curve closeness example. We will explain why the upper bound Theorem~\ref{t.airynt} may be expected to be harder to prove than the lower bound Theorem~\ref{t.airynt.lb}. This discussion will motivate the need for some of the apparatus of the jump ensemble method that we will specify later in Section~\ref{s.jumpensemblemethod}.

We may try to use missing closed middle reconstruction and the Wiener candidate~$\wien$, constructed under $\PP_\mc{F}$, in order to prove upper bounds on probabilities of events expressed in terms of the ensemble $\mc{L}_n$. Suppose that we set the reconstruction parameters $\eln$, $\gfl$ and $\gfr$ equal to $1$, $-1$ and $1$, as we did when we proved Theorem~\ref{t.airynt.lb} (and for that matter Theorem~\ref{t.aestimate}).
Recall that, under~$\PP_\mc{F}$, the conditional distribution of $\mc{L}_n:\intint{k} \times [-1,1] \to \R$ is given by the conditional law of the candidate ensemble $\wien:\intint{k} \times [-1,1] \to \R$ given that $M_{1,k+1}^{[-2,2]}(\wien) = 1$. 

The Wiener candidate is by its definition essentially a system of independent Brownian motions. In proving either the lower bound Theorem~\ref{t.airynt.lb} or the upper bound Theorem~\ref{t.airynt}, we seek to transmit information about the $k$-curve closeness probability for independent Brownian motions in Proposition~\ref{p.expsimple} to learn something comparable about the ensemble~$\mc{L}_n$. Conditioning as we do on the data specifying~$\mc{F}$, this transfer is a matter of specifying a reasonably probable {\em favourable} $\mc{F}$-measurable event~$F$ with the property that under the law $\PP_{\mc{F}}$
for a choice of $\mc{F}$-data causing the occurrence of~$F$, the conditional probability that  $M_{1,k+1}^{[-2,2]}(\wien) = 1$ is itself quite high. 

Of course, we saw this idea late in the preceding chapter, where we took the event $F$ to be $F_t$: this event is after all $\mc{F}$-measurable, and has a positive probability, uniformly in the parameter $n$, 
while Lemma~\ref{l.bcgood} records a lower bound on the conditional probability in question.

As we consider applying this approach to prove the upper bound, it is clear that it must be carried out in a more exacting way. From the form of Theorem~\ref{t.airynt}(1), it is apparent that we must specify a favourable event $F$ in terms of $\mc{F}$-data 
whose complement has a superpolynomial decay in $\e \searrow 0$ (or anyway a decay as fast as $\e^{k^2 - 1}$). This event specified, we must then argue that, under the associated $\PP_\mc{F}$ law, the condition $M_{1,k+1}^{[-2\eln,2\eln]}(\wien) = 1$ is met with at least a subpolynomially decaying $\e^{o(1)}$ probability as $\e \searrow 0$. If this scenario can be realized, then the probability of the putatively rare event $\close\big(\mc{L}_n,0,\e \big)$ may be gauged as follows:
\begin{equation}\label{e.wienergauge}
 \PP \Big( \close\big(\mc{L}_n,0,\e\big) \Big) \leq  \EE \bigg[ \PP_\mc{F} \Big(  \close(\wien,0,\e) \, \Big\vert \, M_{1,k+1}^{[-2\eln,2\eln]}(\wien) = 1 \Big) \cdot {\bf 1}_F \bigg] \, + \, \PP \big( F^c \big) \, .
\end{equation}
The latter term on the right-hand side is $O\big(\e^{k^2 - 1}\big)$ while the former is at most 
$$
 \EE \, \bigg[ \, \frac{\PP_\mc{F} \big( \close(\wien,0,\e) \, , \,  M_{1,k+1}^{[\gfl,\gfr]}(\wien) = 1 \big)}{\PP_{\mc{F}} \big(   M_{1,k+1}^{[-2\eln,2\eln]}(\wien) = 1 \big)} \cdot   {\bf 1}_F \, \bigg] \, .
$$
In the ratio here, the denominator behaves as $\e^{o(1)}$
as $\e \searrow 0$. The numerator concerns the $k$-curve closeness probability for a system of Brownian bridges conditioned on the avoidance constraints that comprise the middle interval test. As such,  Proposition~\ref{p.expsimple}(1) suggests that the numerator is at most $\e^{k^2 - 1 - o(1)}$. (One might object that this result requires that the entrance and exit data, in this case $\overline{W}(\gfl)$ and $\overline{W}(r)$, be well-spaced. However, it is plausible that this data is well-spaced fairly typically, and also that the bound cannot be significantly improved even when $\overline{W}(\gfl)$ and $\overline{W}(r)$ are constant vectors. Our expository purpose would not be advanced by examining this imprecision further.)

If this method can be implemented, then, we would learn that $\PP\big( \close(\mc{L}_n,0,\e) \big) \leq \e^{k^2 -1 - o(1)}$.

We now turn to considering how we may hope to carry out such an approach. The  parameters  $\eln > 0$, $\gfl \in [-\eln,0]$ and $\gfr \in [0,\eln]$ serve to specify the missing closed middle reconstruction $\sigma$-algebra~$\mc{F}$. 

Is the choice  $\eln = 1$, $\gfl = -1$
and $\gfr = 1$ that we have been considering
 realistic for the upper bound argument?
The favourable event $F$ must stipulate that the lower boundary condition $\mc{L}_n \big(k+1,\cdot \big): [\gfl,\gfr] \to \R$
not rise too high. After all, if $\mc{F}$-data that causes $F$ to occur is to result in at least a modest $\PP_{\mc{F}}$-probability of $M_{1,k+1}^{[-2\eln,2\eln]}(\wien) = 1$, then it must also entail such a probability for the weaker middle interval condition that  
 $M_{1,k+1}^{[\gfl,\gfr]}(\wien) = 1$. In the latter condition, we are demanding {\em inter alia} that Brownian bridges begun at unit-order locations at the presently chosen unit-order boundary times $-2\eln$ and $2\eln$
consistently rise above the curve $\mc{L}_n(k+1,\cdot)$ during $[\gfl,\gfr]$.

Since $\mc{L}_n$ is a regular ensemble,  it satisfies the one-point upper tail axiom $\rmreg(3)$ and is ordered (because it is Brownian Gibbs). For these two reasons, 
$$
 \PP \Big( \mc{L}_n \big(k+1,0 \big) \geq s \Big) \leq \rsC \exp \big\{ - \rsc s^{3/2} \big\}
$$
for $s \geq 1$.
In demanding that $\PP \big( F^c \big) = O \big(\e^{k^2 - 1}\big)$, we may thus incorporate into the definition of $F$ a demand on the $(k+1)\textsuperscript{st}$ curve in $\mc{L}_n$ no stronger than
$$
 \sup_{x \in [\gfl,\gfr]} \mc{L}_n \big( k+1,x \big) \leq  O \big( \log \e^{-1} \big)^{2/3} \, .
$$ 
The boundary data $\mc{L}_n \big( i,\pm 2\eln \big)$, $i \in \intint{k}$, are typically of unit-order (if $\eln = 1$ at least). Under $\PP_\mc{F}$, the~$k$ Brownian bridges on $[-2\eln,2\eln]$ used to form the candidate ensemble $\wien:\intint{k} \times [\gfl,\gfr] \to \R$ are being required to jump over a hill whose height is $O \big( \log \e^{-1} \big)^{2/3}$.
If $\eln$ is of unit order, the probability of this happening is
$\exp \big\{ - O( \log \e^{-1})^{4/3} \big\}$.
This four-thirds exponent is inadequate for our purpose: we want it to be less than one to carry out the proposed method.

We may improve matters by increasing the value of $\eln$. 
However, as we do so, we pay a price: the typical boundary data at $\pm 2\eln$ will no longer be at a unit-order location but instead at a parabolically determined location around $- O(\eln^2)$.

Suppose we write $\eln = \big( \log \e^{-1} \big)^\alpha$ for $\alpha > 0$. The interval $[-2\eln,2\eln]$ has length of order $\big( \log \e^{-1} \big)^\alpha$. The $k$ Wiener candidate curves  
defined on this interval begin and end at heights of order $- \big(\log \e^{-1}\big)^{2\alpha}$. They must clear a height of order $\big( \log \e^{-1} \big)^{2/3}$ on an interval of times including zero. The probability that they do so has order
$$
 \exp \bigg\{ - \tfrac{1}{\big( \log \e^{-1} \big)^\alpha} \Big( \big( \log \e^{-1} \big)^{2/3} + \big(\log \e^{-1}\big)^{2\alpha} \Big)^2 O(1) \bigg\} \, = \, 
 \exp \Big\{ - \big( \log \e^{-1} \big)^{(4/3 - \alpha) \vee 3\alpha} O(1) \Big\}   \, .
$$
This expression is maximized by setting $\alpha = 1/3$ when it takes the form $\exp \big\{ O ( \log \e^{-1}) \big\}$.

We were seeking a sub-polynomial decay in low~$\e$, for which such an expression would clearly qualify only if it took the form  $\exp \big\{ ( \log \e^{-1})^{1-\deltamac} \big\}$ for some $\deltamac > 0$; and this we did not find. Although the outcome that $\deltamac = 0$ is possible is ambiguous, it may seem unlikely, on the basis of these heuristics, that the method may work.

Indeed, a slightly closer review would show that it cannot work, and we will need to refine our approach in order to prove the upper bound Theorem~\ref{t.airynt}(1). We will allude to the problem identified here -- that $\deltamac$ cannot be made positive -- as the {\em high jump difficulty}.

Note the role of
the basic parameter $\e > 0$ in our discussion. For the $k$-curve closeness problem, its role is to specify up to a power the order of probability of $k$-curve $\e$-closeness for  a reference mutually avoiding Brownian bridge ensemble, this probability being $\e^{k^2 -1 + o(1)}$ by Proposition~\ref{p.expsimple}. Other upper bounds on Brownian Gibbs ensemble probabilities will be sought, such as in the proof of Theorem~\ref{t.airytail.ln}, a result which quantifies the assertion that the change of measure from mutually avoiding Brownian bridge law to Brownian Gibbs ensemble law transforms small probability events in the manner $\e \to \e^{1 - o(1)}$. Here, the parameter $\e > 0$ again plays the role of determining the order of magnitude of the probability of an event under study for the law in the transform's preimage.
The parameter $\e$ is embedded in the jump ensemble method that we are developing, and {\em the use of this symbol is henceforth reserved} for its use in the method; it always has a conceptual interpretation of the form we have just discussed.

\subsection{Setting the missing close middle reconstruction parameters}\label{s.elnfix}

Set back by the high jump difficulty we may seem to be, but we have learnt something useful: a sensible choice of the parameter $\eln > 0$. 
We now set, {\em for the rest of the paper}, the value of $\eln$ to be equal to $\const \big( \log \e^{-1} \big)^{1/3}$.
 The quantity $\const > 0$ is a constant, without dependence on $\e > 0$, though it depends on $k \geq 2$; in the general formulation of the jump ensemble method, we will take 
\begin{equation}\label{e.constvalue}
 \const  \geq \max \Big\{ k^{1/3} c_k^{-1/3} \big( 2^{-9/2} - 2^{-5} \big)^{-1/3} \, , \,  36(k^2 -1) \Big\} \, ,
\end{equation}
where the $c$-sequence is specified in~(\ref{e.littlec}). We
adopt equality in~(\ref{e.constvalue}) by default, but permit an increase in the value of~$\const$ as the need arises in  specific applications of the method.
On the parameter $\e > 0$, we require that 
\begin{equation}\label{e.epsilonupperbound}
\e < e^{-1} \wedge (17)^{-1/k} C_k^{-1/k} \const^{-1} \, ;
\end{equation}
in fact, the first condition $\e < e^{-1}$
is implied by the second for each $k \geq 2$, but we include it because we will sometimes invoke it directly.
The ensemble index $n \geq k$ will be supposed to satisfy
\begin{equation}\label{e.nlbone}
    n^{\phimac_1/4 \, \wedge \,  \phimac_2/4 \, \wedge \,  \phimac_3/2}   \geq    \big( \rsc/2 \wedge 2^{1/2} \big)^{-1} \const \big( \log \e^{-1} \big)^{1/3}  
\end{equation}
and    
\begin{equation}\label{e.nlbtwo}
n \geq k 
  \vee  (c/3)^{-2(\phimac_1 \wedge \phimac_2)^{-1}} \vee  6^{2/\delta} \, ,
\end{equation}
 where  $\delta = \phimac_1/2 \wedge  \phimac_2/2 \wedge \phimac_3$.

We emphasise that the preceding four hypotheses 
will be enforced throughout the use of the jump ensemble method in this article. The hypotheses are not explicitly recalled except for in the statements of theorems that are derived by the method. Specifically, we are supposing that~(\ref{e.nlbone}) holds. In applications, where $\bar\phimac = \big( 1/3,1/9,1/3 \big)$, this hypothesis asserts that $n \geq O(1) \big( \log \e^{-1} \big)^{12}$, where the $O(1)$ term permits $k$-dependence; in other words, $\e \geq \exp \big\{ - O(1) n^{1/12} \big\}$ -- so that extremely small probability events are not treated by the method uniformly in $n$. For example, this hypothesis transmits to Theorem~\ref{t.airytail.ln}, one of our principal results, which concerns Brownian regularity of regular ensembles; we will allude explicitly to this consideration again, when we prove this result at the end of  Chapter~\ref{c.jumpensembleapplications}.

We must also set the values of $\gfl \in [-\eln,0]$ and $\gfr \in [0,\eln]$. It would be natural enough to consider the choice $\gfl = -\eln$ and $\gfr = \eln$. In fact, we will make choices slightly closer to the origin than these, in a way that ensures that the lower boundary condition $\mc{L}_n \big( k+1,\cdot \big): [\gfl,\gfr] \to \R$ enjoys a little more regularity than would be the case with the choice $\gfl = -\eln$ and $\gfr = \eln$. 

To specify our choice of $(\gfl,\gfr)$, we introduce the {\em least concave majorant} $\cm_+:[-\eln,\eln] \to \R$ of the curve $\mc{L}_n(k+1,\cdot):[-\eln,\eln] \to \R$.

Define a random variable pair $(\mfl,\mfr)$ according to
\begin{eqnarray*}
  \mfl  & = & \inf \big\{ x \in [-\eln,\eln]: \cm_+'(x) \leq    4 \eln \big\} \\
\textrm{and}  \, \, \,   \mfr  & = & \sup \big\{ x \in [-\eln,\eln]: \cm_+'(x) \geq -  4 \eln \big\} \, ,  
\end{eqnarray*}
where the convention that $\inf \emptyset = \eln$ and 
$\sup \emptyset = - \eln$ is adopted.

We now specify, {\em for the rest of the paper}, our choice of $(\gfl,\gfr)$ to be equal to  $(\mfl,\mfr)$.
As such, we will be working with the missing closed middle  reconstruction whose side intervals are $[-2\eln,\mfl]$ and $[\mfr,2\eln]$ and whose middle interval is $[\mfl,\mfr]$. There would seem to be  formal problem that $\mfl > \mfr$ is possible
due to the cases where one of $\mfl = \eln$
or $\mfr = - \eln$ is forced by the above convention; but, in either of these cases, $\mfr = \mfl$. 
Thus, in fact, we always have $\mfl \leq \mfr$; and we will soon see that a stronger bound holds in our applications of the method. 

 Note that $\intint{k+1} \times (\mfl,\mfr)$ is not a stopping domain, because the form of $\mc{L}_n(k+1,\cdot)$ near the origin dictates the values of $\cm_+$ further away. 
 What is needed for our purpose is that $(\mfl,\mfr)$ is $\mc{F}$-measurable (which is true because the entire $(k+1)\textsuperscript{st}$ curve is known to the witness of~$\mc{F}$). 

(The use of the pair $(\mfl,\mfr)$ develops a technique in the construction~\cite{AiryLE} of the Airy line ensemble by which was proved the key technical Proposition~$3.5$ showing a uniform lower bound on acceptance probability for a prelimiting ensemble sequence. It would be of interest however to revisit this proof with merely the use of missing closed middle reconstruction.)

\subsection{Specifying the favourable event $\mpg$}\label{s.favevent}
Disregarding for a moment longer the identified high jump difficulty, 
we endeavour to specify in precise terms the form of the proposed method, now defining the favourable event $\mpg$: we set
$$
\mpg = \fev_1 \cap \fev_2 \cap \fev_3 \, ,
$$
where 
\begin{eqnarray*}
 \fev_1 & = & \Big\{ \mc{L}_n\big(i, x  \big) \in   \eln^2 \big[ - 2\sqrt{2} - 1 , - 2 \sqrt{2} + 1 \big] \, \, \textrm{for $(i,x) \in \intint{k} \times \big\{ -2\eln, 2\eln \big\}$} \,
 \Big\}  \, , \\
\fev_2 & = & \Big\{ - \eln^2 \leq  \mc{L}_n\big(k+1, x  \big)  \leq  \eln^2 \, \, \textrm{for $x \in [ -\eln, \eln ]$}   \Big\} \, , \\
\textrm{and} \, \, \, \fev_3 & = & \bigcap_{i \in \intint{k} } \Big\{  \xmin_i  \in \big[ - \eln^2 , \eln^2 \big]  \Big\} \cap
 \Big\{  \ymin_i  \in \big[ - \eln^2 , \eln^2 \big]  \Big\} \, .
\end{eqnarray*}

Note that, since $\mc{L}_n\big(k+1, \mfl \big) \geq \mc{L}_n\big(k+1, -\eln \big) +  4\eln \big( \mfl + \eln \big)$, the occurrence of $\mpg \subseteq \big\{ \mc{L}_n(k+1,-\eln) \geq - \eln^2 \big\} \cap \big\{ \mc{L}_n(k+1,\mfl) \leq \eln^2 \big\}$ entails that $\mfl \leq - \eln/2$; similarly, it ensures that $\mfr \geq  \eln/2$. Thus, $[-\eln/2,\eln/2]$ is always a subset of the middle interval $[\mfl,\mfr]$
whenever middle interval reconstruction is attempted.

The next lemma establishes that the favourable event $\mpg$
has the superpolynomial decay property that the upper bound method we have been proposing demands of it: $\PP \big( \mpg^c \big) = O(\e^{O(c_k \const^3)})$ as~$\e \searrow 0$, so that high choices of $\const$ give rise to any desired polynomial decay.

\begin{lemma}\label{l.glub.new}
We have that
$$
\PP \big( \mpg^c \big) \leq   \e^{2^{-5} c_k \const^3} \, .
$$
\end{lemma}
\noindent{\em Remark.} 
When applications are made, the value of $\const$ will be increased if necessary over the explicit value given in~(\ref{e.constvalue}) in order to ensure that the decay rate 
$\PP \big( \mpg^c \big)$ is suitable for the problem at hand.

\medskip

\noindent{\bf Proof of Lemma~\ref{l.glub.new}.} We must bound above $\PP(\fev_i^c)$ for each $i \in \intint{3}$.

\medskip

\noindent{\em Bounding  $\PP(\fev_1^c)$.}
By hypothesis, $\mc{L}_n$ is a regular ensemble. 
As such, 
Proposition~\ref{p.othercurves} 
 with $x = \pm 2\eln$ and $s = \eln^2$ 
provides an estimate needed to control the one-point lower tail, namely
$$
 \sum_{x \in \{ -2\eln,2\eln \}} \PP \Big( \, \mc{L}_n\big( k,x\big) + 2^{-1/2} x^2 \leq - \eln^2 \, \Big) \, \leq \,  2C_k \exp \big\{ - c_k  \eln^3 \big\} = 2C_k \e^{c_k \const^3} \, ;
$$
here are used 
$n \geq k 
  \vee  (c/3)^{-2(\phimac_1 \wedge \phimac_2)^{-1}} \vee  6^{2/\delta}$, $2\eln \leq 2^{-1}c n^\delta$ and  $\eln^2 \leq 2 n^\delta$, 
 where  $\delta = \phimac_1/2 \wedge  \phimac_2/2 \wedge \phimac_3$ -- hypotheses which hold due to~(\ref{e.nlbone}),~(\ref{e.nlbtwo}) and $c \leq 1$.

The one-point upper tail axiom $\rmreg(3)$ 
 with $z = \pm 2\eln$ and $s = \eln^2$  offers the control 
$$
  \sum_{z \in \{ -2\eln,2\eln \}} \PP \Big( \, \mc{L}_n\big( 1,z\big) + 2^{-1/2} z^2 \geq  s \, \Big) \, \leq \,   2 C \exp \big\{ - c  \eln^3 \big\} = 2C \e^{c \const^3}
$$
where we used $2\eln \leq c n^{\phimac_2}$ and $\eln^2 \geq 1$. 

From the two displayed bounds, as well as $C_k \geq C$ and $c_k \leq c$,  we learn that
$$
 \PP\big(\fev_1^c\big) \leq 4C_k \e^{c_k \const^3} \, .
$$

\noindent{\em Bounding $\PP(\fev_2^c)$.}
The probability that the lower-tail  condition in $\fev_2$ fails satisfies
\begin{eqnarray}
& &  \PP \,  \Big( \, \exists \, x \in [ -\eln, \eln ] : \mc{L}_n\big(k+1, x  \big)  < - \eln^2 \,  \Big) \label{e.fthreelb} \\
& \leq  & \PP \, \Big( \, \inf_{x \in [-\eln,\eln]} \big(  \mc{L}_n (k+1, x  ) + 2^{-1/2} x^2 \big) \leq - \big( 1 - 2^{-1/2} \big) \eln^2 \, \Big) \nonumber \\
& \leq & \eln^k \cdot \Cstrong_k \exp \big\{ - c_k (1 - 2^{-1/2})^{3/2} \eln^3 \big\} \leq \eln^k C_k \e^{(1 - 2^{-1/2})^{3/2} c_k \const^3} \, , \nonumber
\end{eqnarray}
where the second inequality is a consequence of Proposition~\ref{p.strongothercurves} applied at curve index $k+1$ with parameter choices~$y=0$, $t=\eln$ and $r=(1 - 2^{-1/2})\eln^2$. This proposition is a refinement of Proposition~\ref{p.othercurves} and appears in the appendix. The constant $\Cstrong_k$
is introduced in this proposition and, in view of~(\ref{e.formere}), satisfies $\Cstrong_k \leq C_k$
-- a bound used in the final inequality in the above display. The application of the proposition requires that certain upper and lower bounds on the parameters $t$ and $r$ be satisfied (as well as a trivial condition on $y$). It is readily checked that the hypotheses~(\ref{e.constvalue}) and~(\ref{e.nlbone}), alongside the definition~(\ref{e.littlec}) and $\e \leq e^{-1}$, ensure these; thus is the preceding display confirmed.

We will treat the upper-tail condition in $\fev_2$ by applying Proposition~\ref{p.nobigmax} in order to find that
\begin{equation}\label{e.elnsup}
\PP \Big( \sup_{x \in [-\eln,\eln]} \mc{L}_n \big( 1,x \big) \geq \eln^2 \Big) \leq  12 \eln \rsC \e^{2^{-9/2} \rsc \const^3} \, .
\end{equation}
For this application of Proposition~\ref{p.nobigmax}, recall that $\mc{L}_n$ is supposed to be a $\big(\bar\phimac,\rsc,\rsC\big)$-regular ensemble 
 for some $\bar\phimac \in (0,\infty)^3$. Set 
 $r = \eln$ and $t = \eln^2$. The proposition's hypotheses are  $\tfrac{1}{2} \rsc n^{\phimac_1 \wedge \phimac_2} \wedge 2^{1/2} n^{\phimac_3/2} \geq \eln \geq 2^{7/4}$ and $n \geq (2\rsc)^{-2(\phimac_1 \wedge \phimac_2)^{-1}}$; they hold by~(\ref{e.nlbone}) and~(\ref{e.nlbtwo}).
Applying the  proposition and using $\eln \geq 1$, we bound the probability in~(\ref{e.elnsup}) above by  $12\eln   \rsC \exp \big\{ - 2^{-9/2} \rsc  \eln^3 \big\}$, and thus derive~(\ref{e.elnsup}).

Note then that, if the upper tail condition in $\fev_2$ fails, the event in the displayed estimate occurs. Indeed, although the $\fev_2$-condition concerns the $(k+1)\textsuperscript{st}$ curve, rather than the first, the event becomes less probable if the $(k+1)\textsuperscript{st}$ curve is considered,
 because  the ensemble $\big\{ \mc{L}_n: n \in \N \big\}$ is ordered.

Combining (\ref{e.fthreelb}) with (\ref{e.elnsup}), we arrive at the bound
$$
\PP \big( \fev_2^c \big) \leq 
 12 \eln \rsC \e^{2^{-9/2} \rsc \const^3}   +  \eln^k C_k \e^{(1 - 2^{-1/2})^{3/2} c_k \const^3} \, .
$$

\noindent{\em Bounding $\PP(\fev_3^c)$.} Note that $\xmin_i \geq \mc{L}_n(k+1,\mfl)$
and $\ymin_i \geq \mc{L}_n(k+1,\mfr)$ for $i \in \intint{k}$. Since $\{\mfl,\mfr\} \subset [-\eln,\eln]$, we see that the lower tail conditions in the definition of $\fev_3$ are in fact implied by the conditions in $\fev_2$. Regarding the upper conditions, note that, for such $i$,
$\xmin_i \leq \mc{L}_n(1,\mfl)$ and 
$\ymin_i \leq \mc{L}_n(1,\mfr)$.
Since $\{ \mfl,\mfr \} \subset [-\eln,\eln]$, the occurrence of the event in~(\ref{e.elnsup}) entails that these conditions are satisfied. In this way, the right-hand side of the preceding display
 is seen to be an upper bound not merely on 
$\PP \big( \fev_2^c \big)$ but in fact also on 
$\PP \big( \fev_2^c \big) +
\PP \big( \fev_3^c \big)$.

\noindent{\em Gathering the estimates.}
In this way, we find that
$$
\PP \big( \mpg^c \big)  \leq  
4C_k \e^{c_k \const^3} + 12 \eln \rsC \e^{2^{-9/2} \rsc \const^3}   +  \eln^k C_k \e^{(1 - 2^{-1/2})^{3/2} c_k \const^3}  \, .
$$
Since $\eln \geq 1$, $C_k \geq \rsC$ and $c_k \leq \rsc$, we have that $\PP \big( \mpg^c \big) \leq 17 C_k \eln^k \e^{2^{-9/2} c_k \const^3}$.
From~(\ref{e.epsilonupperbound}), we have
$17 C_k \eln^k \leq \e^{-k}$;
note further that $\const \geq k^{1/3} c_k^{-1/3} \big( 2^{-9/2} - 2^{-5} \big)^{-1/3}$
yields $\e^{ 2^{-9/2} c_k \const^3} \leq \e^k \cdot \e^{c_k 2^{-5} \const^3}$. 
 Thus, we obtain Lemma~\ref{l.glub.new}. \qed

\subsection{More promising than the Wiener candidate: the jump ensemble}\label{s.morepromising}

Consider the law $\PP_\mc{F}$
for such $\mc{F}$-data that $\mpg$ occurs. We have understood that the candidate examination condition~$\BP^{[-2\eln,2\eln]}_{1,k+1}(\wien) = 1$
 may be satisfied only with $\PP_\mc{F}$-probability $\e^{O(1)}$ as $\e \searrow 0$, rather than with the desired $\e^{o(1)}$ probability. 

We must change our approach in order to solve the discussed high jump difficulty: we will consider a variant of the Wiener candidate ensemble $\wien:\intint{k} \times [\gfl,\gfr] \to \R$.
The new ensemble will be called the jump ensemble and denoted by
$J:\intint{k} \times [\gfl,\gfr] \to \R$. Under $\mpg$-satisfying choices of $\mc{F}$, it will verify~$\BP^{[-2\eln,2\eln]}_{1,k+1}(J) = 1$ with the sought subpolynomial $\e^{o(1)}$ $\PP_\mc{F}$-probability.  
 
The Wiener candidate has been found wanting because of the difficulty it encounters in passing the middle interval test~(\ref{e.onepiecemiddle}). 
We now specify a new test, the {\em jump test}, the passing of which will be a necessary condition for~(\ref{e.onepiecemiddle}) to hold. The new test will demand that the $k$ candidate curves jump over a large set of extreme points with $x$-coordinates in $[\mfl,\mfr]$ of the graph of the $(k+1)\textsuperscript{st}$ curve's concave majorant $\cm_+$,
so that the Wiener candidate who passes the jump test necessarily has an overall curve geometry that offers a viable prospect of success.
We will define the jump ensemble~$J$ by conditioning the Wiener candidate on passing the jump test, alongside the side intervals test; this ensemble will indeed resolve the high jump difficulty.  

Recall then that $\cm_+:[-\eln,\eln] \to \R$ is  the least concave majorant of the curve $\mc{L}_n(k+1,\cdot):[-\eln,\eln] \to \R$. 
Let $\xext \subset [\mfl,\mfr]$ denote the set of $x$-coordinates of extreme points of the closed set $\big\{ (x,y) : \mfl \leq x \leq \mfr \, , \, y \leq \cm^+(x) \big\}$.
 Note that $\xext$ consists  of the intersection with $[\mfl,\mfr]$ of the set of points of local non-constancy of $\cm_+'$; necessarily, $\{ \mfl,\mfr \} \in \xext$.   Introducing a parameter $\ipd \in [1,\mfr-\mfl)$, we let $\pole$ denote a subset of $\xext$  with the properties that 
\begin{itemize}
\item $\{ \mfl,\mfr \} \in \pole$;
\item any distinct elements $\pp_1,\pp_2 \in \pole$ satisfy $\vert \pp_1 - \pp_2 \vert > \ipd$;  
\item and, if $x \in \xext  \setminus \pole$, then some element $\pp \in \pole$ satisfies $\vert \pp - x \vert \leq \ipd$. 
\end{itemize}

In fact, several subsets of $\xext$ may satisfy these conditions. 
If this is the case, we may select~$\pole$ to be the subset among the choices of maximal cardinality that is maximal in the lexicographical ordering. The set $\pole$, rather than $\xext$, will be the focus of our attention; we will sometimes call it the {\em pole} set. 
The quantity $\ipd$ is the {\em inter-pole distance} parameter.
(The role of this parameter is determined by the application; if an understanding of the behaviour of for example $\mc{L}_n(1,\cdot)$ is sought on an interval of a certain length $l$, we would set $\ipd$ to equal $l$, or perhaps a bounded multiple thereof. In applications, we may think of $\ipd$ as being a constant, independent of $\e > 0$ and even of $k \geq 2$.)

We record for future use that
\begin{equation}\label{e.polecard}
  \vert P \vert \leq  2\eln   \, .
\end{equation}
Indeed, elements of $\pole$ are separated by a distance that exceeds $\ipd \geq 1$ and all lie in $[\mfl,\mfr] \subseteq [-\eln,\eln]$. 

We also define the {\em tent} map  $\tent:[\mfl,\mfr] \to \R$ to be the piecewise affine function which on any closed interval between consecutive pole set values $p_1$ and $p_2$ is equal to the affine function whose values at $p_1$ and $p_2$ are respectively $\mc{L}_n(k+1,p_1)$ and $\mc{L}_n(k+1,p_2)$. 
Note that $\tent$ is an $\mc{F}$-measurable function. 

We choose these names because we may think of a pole $\{ p \} \times \big(-\infty, \mc{L}_n(k+1,p)\big]$ being built over each element $p \in \pole$, so that the leftmost and rightmost poles are supported at $\mfl$ and $\mfr$. The graph of the tent map $\tent: [\mfl,\mfr] \to \R$
is thus propped up by the collection of poles. 

We say that the Wiener candidate passes the {\em jump} test if each of its curves clears all the pole tops, namely if 
$$
\wien(i,x) > \mc{L}_n\big(k+1,x \big) \ \, \, \textrm{for all $(i,x) \in \intint{k} \times \pole$} \, .
$$
The middle interval test is comprised of internal and external avoidance constraints for $\wien$ on $[\mfl,\mfr]$. These constraints collectively ensure that each $\wien$-curve exceeds $\mc{L}_n(k+1,\cdot)$ throughout $[\mfl,\mfr]$. The jump test checks that the latter condition holds only at points in $\pole$. Since $\pole \subseteq [\mfl,\mfr]$, the jump test is weaker than the middle interval test.

The jump test is fundamental to our method. As a practical matter, it is convenient to also introduce the {\em order-over-poles} test, which is passed when 
$$
\wien(i,x) > \wien \big(i+1,x \big) \ \, \, \textrm{for all $(i,x) \in \intint{k-1} \times \pole$} \, .
$$
When this test is passed, the candidate curves are in decreasing order above each pole.

In summary, then: under $\PP_\mc{F}$, the candidate~$\wien$'s overall success~$\BP_{1,k+1}^{[-2\eln,2\eln]}\big(\wien\big) = 1$ may be tested in four steps: 
\begin{itemize}
\item Test $1$ is the side intervals test: in one-piece list notation, $\BP_{1,k+1}^{[-2\eln,\mfl]\cup [\mfr,2\eln]} (\wien) = 1$;
\item Test $2$ is the jump test;
\item Test $3$ is the order-over-poles test; and
\item Test $4$ is  the middle interval test, i.e., 
$\BP_{1,k+1}^{[\mfl,\mfr]} (\wien) = 1$.  
\end{itemize}
Since Tests $2$ and~$3$ are weaker than Test $4$,~$\wien$ passes the test sequence if and only if it is successful.

The  indicator function of the event that is checked in Test $i$ will be denoted by $\test_i$, for $i \in \intint{4}$. We also write for example $T_{12} = T_1 \wedge T_2$; thus, $T_{1234}(W) = 1$ denotes the event that  $\BP_{1,k+1}^{[-2\eln,2\eln]}\big(\wien\big)  = 1$.

We now construct the jump ensemble~$J$ under the law~$\PP$. The ensemble $J:\intint{k} \times [\mfl,\mfr] \to \R$ is constructed so that, under $\PP_\mc{F}$, it has the conditional distribution of $\wien:\intint{k} \times [\mfl,\mfr] \to \R$ given that $\test_{12}(W) = 1$.

\begin{figure}[ht]
\begin{center}
\includegraphics[height=11cm]{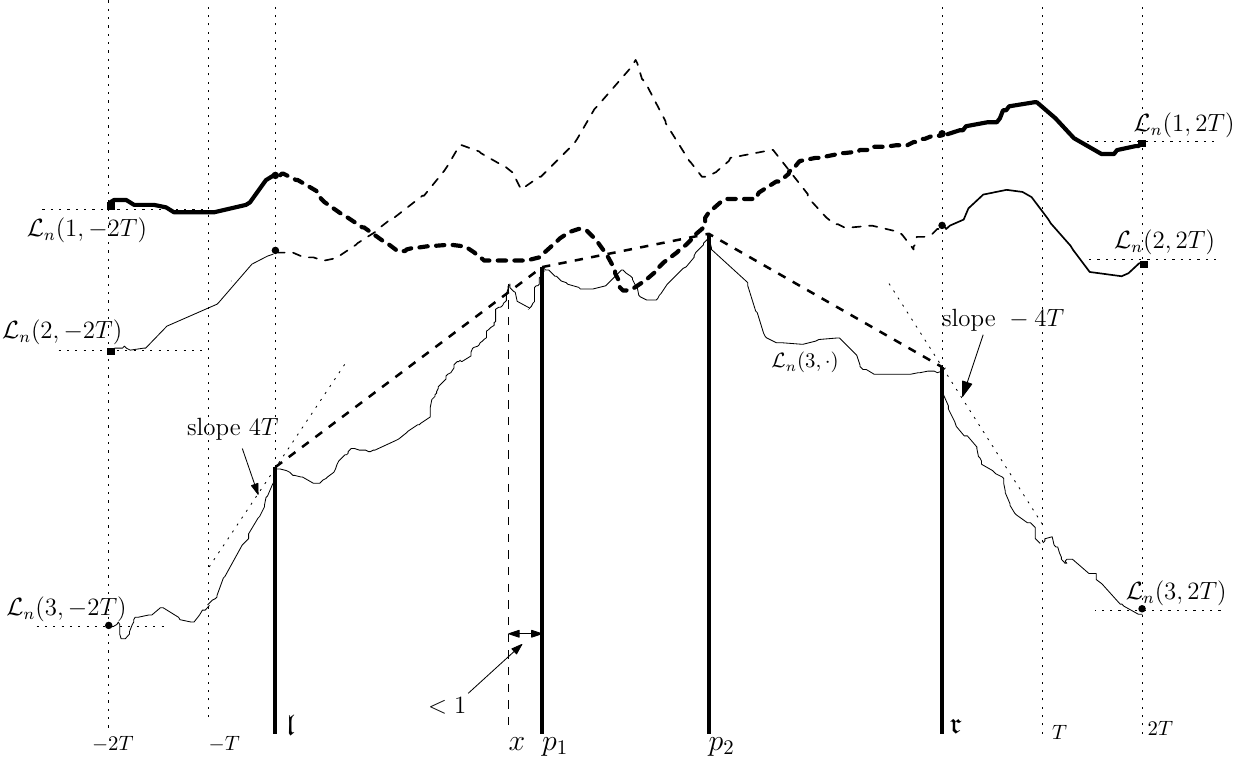}
\caption{The jump ensemble depicted with $k=2$ (and $\ipd = 1$).
The depiction is not to scale: the intervals $[-2\eln,-\eln]$ and  $[\eln,2\eln]$ are too short. 
The pole set~$P$  in this example equals $\{ \mfl, p_1 , p_2, \mfr \}$. The vertical poles are depicted in thick solid lines. The point $x$ is an element of $\xext$ but not of $P$, because $\vert x - p_1 \vert < 1$. (In fact, there are almost surely infinitely many elements of $\xext \setminus P$.) The piecewise affine dashed curve defined on $[\mfl,\mfr]$ is $\tent$. The rougher dashed curves are the jump ensemble $J:[1,2] \times [\mfl,\mfr] \to \R$. The jump ensemble fails the criterion $\test_4(J) = 1$ due to the meeting of the two $J$-curves and contact between $J(1,\cdot)$ and $\mc{L}_n(3,\cdot)$.} 
\label{f.jumpensemble}
\end{center}
\end{figure}

In this way, to the witness of~$\mc{F}$ the jump ensemble~$J$
is a candidate~$\wien$ that passes Tests~$1$ and~$2$. The construction is depicted in Figure~\ref{f.jumpensemble}.

\subsection{The jump ensemble as a halfway house}
The jump ensemble represents a halfway-house between the unadulterated Brownian randomness of the Wiener candidate $\wien$ and the desired conditional distribution under $\PP_\mc{F}$ of the actual line ensemble~$\mc{L}_n$. Under $\PP_\mc{F}$, the latter distribution is obtained as the law of the jump ensemble conditioned on $\test_4(J) = 1$. 

\subsubsection{The jump ensemble's promise realized: Proposition~\ref{p.vaultsuccess}} Our next result emphasises how the jump ensemble is a serious contender for passing the final test: 
provided that the highly typical $\mc{F}$-measurable event $\mpg$ occurs, this ensemble meets the demand that $\test_4(J) = 1$---equivalently, that $\test_{34}(J)=1$---with a probability that decays slowly, subpolynomially, as $\e \searrow 0$.
\begin{proposition}\label{p.vaultsuccess}
We have that
$$
 \PP_{\mc{F}} \Big( \test_4(J) = 1 \Big)  \, \geq \, 
 \exp \Big\{ -  3973 k^{7/2} \ipd^2 \const^2 \big( \log \e^{-1} \big)^{2/3} \Big\}  \cdot  {{\bf 1}}_{\mpg} \, .
$$
\end{proposition}

In a sense, we have solved the high jump difficulty: with a choice of $\deltamac = 1/3$, the $\PP_\mc{F}$-probability that
$J$ passes $\test_4(J) = 1$, and so renders to the witness of $\mc{F}$ a copy of the law $\mc{L}_n$, is at least $\exp \big\{ - O(1) (\log \e^{-1})^{1-\deltamac} \big\}$.

\subsection{The jump ensemble method: summary}
We may summarise the jump ensemble method for upper bounds by revisiting the
articulation in~(\ref{e.wienergauge}) of the proposed Wiener candidate method (that was found wanting for some such purposes), and formulating a counterpart inequality for the jump ensemble. Consider a general event $A$, expressed as a collection of $k$ curves on $[\mfl,\mfr]$. 
The analogue of~(\ref{e.wienergauge}) in
 the new formulation is
 $$
 \PP \big( \mc{L}_n \in A  \big) \leq  \EE \Big[ \PP_\mc{F} \big(  J \in A \, \big\vert \,  \test_4(J) = 1 \big) \cdot {\bf 1}_\mpg \Big] \, + \, \PP \big( \mpg^c \big) \, .
 $$
 On the left-hand side, we have abused notation to write $\mc{L}_n$ to indicate the restriction of this ensemble to $\intint{k} \times [\mfl,\mfr]$.
 Recall that the favourable event $\mpg$ has been specified in Section~\ref{s.favevent}.

Using Proposition~\ref{p.vaultsuccess} and Lemma~\ref{l.glub.new}, we see that
 $$
 \PP \big( \mc{L}_n \in A  \big) \leq \e^{-o(1)} \, \EE \Big[ \, \PP_\mc{F} \big(  J \in A \big) \cdot {\bf 1}_\mpg \, \Big] \, \, + \, \, \e^{\infty}  \, ,
 $$
 where in fact $\e^{-o(1)}$ is $\exp \big\{ O_k(1) (\log \e^{-1})^{2/3} \big\}$ (with $O_k(1)$ bounded for given $k$) and $\e^\infty$ indicates arbitrarily fast polynomial decay whose exponent is determined by choosing $\const > 0$ high enough. 
 
The scenario that we sought to realize in discussing~(\ref{e.wienergauge}) was not achievable, but now we have a new variant of that scenario. 
In spirit, the jump ensemble method works as follows. Consider an event $A$ of the type we are discussing,
whose Brownian bridge probability is of order $\e$. This means something to the effect that $\mc{B}_{k;\bar{u},\bar{v}}^{[-1,1]}(A) = \e$, where  $\bar{u}$ and $\bar{v}$ are say typical unit-order entrance and exit data. (We could also randomize over these vectors, choosing them to have say a Gaussian law conditional on being $k$-decreasing lists.)
Then, in order to  conclude that
 the ensemble probability $\PP \big( \mc{L}_n \in A \big)$
is at most $\e^{1-o(1)}$, it is enough to establish the jump ensemble estimate that the
$\PP_{\mc{F}}$-probability that the marginal process $J: \intint{k} \times [-1,1] \to \R$ lies in~$A$
is comparable (at most a large constant multiple would suffice) to the Brownian bridge probability $\mc{B}_{k;\bar{u},\bar{v}}^{[-1,1]}(A)$, for all instances of $\mc{F}$-measurable data for which $\mpg$ occurs.

\section{General tools for the jump ensemble method}\label{s.subpoly}

In this section, we complete our discussion of the general apparatus of the jump ensemble method by proving Proposition~\ref{p.vaultsuccess}.

\subsection{The strong jump ensemble}

As we begin a journey leading to the proof of Proposition~\ref{p.vaultsuccess}, we introduce a technically useful construct. The {\em strong jump} ensemble $\strongjump:\intint{k} \times [\mfl,\mfr] \to \R$ is specified under the law $\PP$ so that, under $\PP_\mc{F}$, it has the conditional distribution of $\wien:\intint{k} \times [\mfl,\mfr] \to \R$ given that $\test_{123}(\wien) = 1$. Thus, the jump ensemble $J$ further conditioned on passing the order-over-poles test is the strong version~$\strongjump$. The next two results, which concern the new ensemble, will be of value in the upcoming proof. 

\subsubsection{The Wiener candidate's prospects of promotion to the strong jump ensemble}
A lower bound on the $\PP_{\mc{F}}$-probability of this eventuality may be seen as part of the overall apparatus of the jump ensemble method. 
\begin{proposition}\label{p.wienerpromotion}
We have that 
$$
\PP_\mc{F} \big(  \test_{123}(\wien) = 1 \big) \geq \e^{36(k+2)^2 k \const^3}   \, 2^{-k} \, {\bf 1}_{\mpg}  \, . 
$$
\end{proposition}

\subsubsection{The jump ensemble's prospects of the same promotion}
The jump ensemble faces a less daunting challenge than does the Wiener candidate in seeking to attain strong jump status. Indeed, the lower bound offered by the next result decays at a much slower rate than does the preceding bound. The result's role in the proof of Proposition~\ref{p.vaultsuccess} will be clarified via the equality~(\ref{e.testguide}) with which we will later begin this proof.
\begin{lemma}\label{l.jumppromotion}
We have that 
$$
\PP_\mc{F} \big(  \test_3(J) = 1 \big) \geq  \exp \big\{ - \big( \log \e^{-1} \big)^{1/3} \cdot 2 D_k \,  k \log k \big\} \, . 
$$
\end{lemma}
We offer the short proof of this lemma; give the longer proof of Proposition~\ref{p.wienerpromotion}; and then turn to the proof of Proposition~\ref{p.vaultsuccess}.

{\bf Proof of Lemma~\ref{l.jumppromotion}.} Since the jump ensemble $J$ passes the side intervals' test, Lemma~\ref{l.glv} implies that its curves are in decreasing order above the leftmost and rightmost poles, $\mfl$ and $\mfr$. It suffices to prove Lemma~\ref{l.jumppromotion} in the strengthened form in which the left-hand event is further conditioned on any given value for the pair of decreasing $k$-vectors formed by $\overline{J}(\mfl)$
and  $\overline{J}(\mfr)$ and on the form of the sets $A_p$ of values adopted by $\overline{J}:\intint{k} \times \{ p\} \to \R$
at each pole $p \in P \setminus \{ \mfl, \mfr \}$. Each of the sets~$A_p$ may be supposed to have  cardinality $k$, since it is easy to see that coincidences among the values of the ensemble~$J$ at a given point occur with probability zero. Let $n = \vert P \vert - 2$ denote the number of sets $A_p$; note from~(\ref{e.polecard}) that $n \leq 2T = 2 D_k \big( \log \e^{-1}\big)^{1/3}$. In the random experiment formed by this conditioning, there are $(k!)^n$ possible outcomes for the passage of $J$ above the poles, each outcome being indexed by a product over the concerned poles~$p$ of a permutation on $\intint{k}$ that specifies the vertical order in which the curves in $J$ visit $A_p$. Among these is a unique outcome that realizes  $\test_3(J) = 1$---the outcome under which the curves in $J$ are in decreasing order over every pole. To prove Lemma~\ref{l.jumppromotion}, it will be enough, since $k! \leq k^k$, to argue that this outcome is at least as likely as any other. To any given outcome, we may associate the $k$-curve graph given by the union of the affine interpolations of the visits of the curves in $J$ to horizontal coordinates given by $p \in P$. The outcome that we seek to show is most likely is the one for which no curves intersect in the associated graph. 
To show the sought statement, it is thus enough to argue that the undoing of any given crossing in a graph---the replacement between consecutive poles of crossing affine segments with non-crossing affine segments that share the concerned endpoint pairs---can only increase the probability of the outcome in question. Consider a crossing configuration $(a,p) \to (d,p')$ and $(b,p) \to (c,p')$, and its non-crossing counterpart $(a,p) \to (c,p')$ and $(b,p) \to (d,p')$; here,   $a,b,c,d \in \R$ satisfy $a < b$ and $c < d$ (and $p,p' \in \R$ satisfy $p < p'$). 
Two pair-product Gaussian densities must be compared to see that it is the latter configuration that is the more probable; the bound  $bd + ac > ad + bc$ provides what is needed. Thus, we obtain  Lemma~\ref{l.jumppromotion}. \qed

\subsection{Trying for strong jump ensemble status: the Wiener candidate's flying leap}
Here we prove  Proposition~\ref{p.wienerpromotion}.
 Recall that, conditionally on $\mc{F}$,
the strong jump ensemble $\strongjump$ has the law of $\wien$ given $\test_{123}(\wien) = 1$.
 We begin by providing an explicit condition on $\wien$ that ensures that the three concerned tests are passed. Define the {\em flying leap} event 
\begin{eqnarray*}
\flyleap & = & \Big\{ \ovbar{\wien}(\mfl) \in \barxmin + 3 \eln^2 \cdot D \, \, , \, \, \ovbar{\wien}(\mfr)  \in \barymin + 3 \eln^2 \cdot D  \Big\} \\
   & & \cap \, \, \Big\{ \, \big\vert \wien^{[\mfl,\mfr]}\big( i, x \big) \big\vert < \eln^2 \, \, \, \forall \, (i,x) \in \intint{k} \times \big[\mfl,\mfr \big] \, \Big\} \, .
\end{eqnarray*}
Here, $D$ denotes the box $[2k-1,2k] \times \cdots \times [3,4] \times [1,2]$ and $3\eln^2 \cdot D$ the dilation by a factor of $3\eln^2$. 
We first check the inclusion
\begin{equation}\label{e.flyleapinc}
\mpg \cap \flyleap \, \subseteq \,  \Big\{ \test_{123}(\wien) =    1 \Big\} \, .
\end{equation}
To verify this, note first that Lemma~\ref{l.glv} implies that the side intervals test $\test_1(\wien) = 1$ is passed provided that $\ovbar{\wien}(\mfl) - \barxmin \in (0,\infty)^k_>$
and  $\ovbar{\wien}(\mfr) - \barymin \in (0,\infty)^k_>$. 
These latter conditions apply because $3\eln^2 \cdot D \subset (0,\infty)^k_>$. For the jump test $\test_2(\wien) = 1$ to be passed, it is sufficient that $\wien(i,x)$ exceed $\mc{L}_n(k+1,x)$ whenever $(i,x) \in \intint{k} \times [\mfl,\mfr]$. 
Taking $i = k$, we have that, for all such $x$, 
the occurrence of $\flyleap$ entails that
$\wien(i,x) > 3 \eln^2 + \xmin_k \wedge \ymin_k \,  - \eln^2$;
on $\mpg$, $\xmin_k \wedge \ymin_k \geq - \eln^2$.
Thus, on $\mpg \cap \flyleap$, 
$\wien(i,x) > \eln^2$. On the other hand,
 $\mpg$ entails that $\mc{L}_n(k+1,x) \leq \eln^2$ (since $[\mfl,\mfr] \subseteq [-\eln,\eln]$). The lower bound on~$\wien(i,x)$ is equally true for $i < k$: we have confirmed the condition $\test_2(\wien)=1$. Regarding Test~$3$, note that $\wien(i,x) - \wien(i+1,x) \geq 3T^2$ for $x \in \{ \mfl,\mfr\}$ and $i \in \intint{k-1}$ when $\flyleap$ occurs.
 We thus have $\wien(i,x) - \wien(i+1,x) \geq T^2$ for all $x \in [\mfl,\mfr]$ in view of the latter condition specifying $\flyleap$. So the $W$-curves are ordered, as the third test demands; we have demonstrated~(\ref{e.flyleapinc}).

In order to prove Proposition~\ref{p.wienerpromotion}, it is enough, 
in view of~(\ref{e.flyleapinc}), 
to verify that
\begin{equation}\label{e.flyingleap}
\PP_\mc{F} \big(  \flyleap \big) \geq \e^{36(k+2)^2 k \const^3}   \, 2^{-k} \, {\bf 1}_{\mpg}  \, . 
\end{equation}
We now do so.
Let $i \in \intint{k}$. When $\mpg$ occurs, $\xmin_i$ and $\ymin_i$ lie on the interval $[-\eln^2,\eln^2]$, while $\mc{L}_n(i,-2\eln)$ and $\mc{L}_n(i,2\eln)$ are at distance at most $\eln^2$ from $-2\sqrt{2} \eln^2$. Thus, the quantity  
\begin{eqnarray}
 & & \PP_\mc{F} \Big( \wien(i,\mfl)  \in \xmin_i + 3 \eln^2 \cdot \big[2k+1 - 2i,2(k+1 - i)\big]  \, , \label{e.flyleapbound} \\
  & & \qquad \qquad \qquad \qquad    \wien(i,\mfr)  \in \ymin_i + 2 \eln^2 \cdot \big[2k+1 - 2i,2(k+1 - i)\big] \Big) \nonumber
\end{eqnarray}
is seen to be at least
\begin{eqnarray*}
 & &  \inf \,
 \mc{B}_{1;x,y}^{[-2\eln,2\eln]} \Big( B(\mfl) \in a  + 3 \eln^2 \cdot \big[2k+1 - 2i,2(k+1 - i)\big] \, , \\
  & & \qquad \qquad \qquad \qquad \qquad \qquad   B(\mfr) \in b  + 3 \eln^2 \cdot \big[2k+1 - 2i,2(k+1 - i)\big] \Big) \, , 
\end{eqnarray*}
when $\mpg$ occurs.
The infimum is taken over choices of $x$ and $y$  in the interval $[-2\sqrt{2} - 1,-2\sqrt{2} + 1] \eln^2$ and of $a$ and $b$ in the interval $[-\eln^2,\eln^2]$. We may now use Corollary~\ref{c.brbr} to find a lower bound
on the expression whose infimum is taken. Using $3\eln^2 \geq 1$ in order to omit the product of the interval lengths in the corollary, we find that the expression is at least 
 $$
  Z^{-1} G_1(x,a) G_2(a,b) G_3(b,y)
  $$
where 
$$
G_1(x,a) = \inf_{u \in 3 \eln^2 \cdot \big[2k+1 - 2i,2(k+1 - i)\big]} g_{x,\mfl + 2\eln}(a + u) \, ,
$$  
$$
 G_2(a,b) =  \inf_{u,v \in 3 \eln^2 \cdot \big[2k+1 - 2i,2(k+1 - i)\big]} g_{a+u,\mfr - \mfl}(b+v)  \, ,
$$ 
 and 
$$  
 G_3(b,y) =   \inf_{v \in 3 \eln^2 \cdot \big[2k+1 - 2i,2(k+1 - i)\big]}  g_{b+v,2\eln - \mfr}(y)  \, .
$$   
The normalization $Z$ equals $g_{x,4\eln}(y)$; note that $Z^{-1} \geq 1$ since $4\eln (2\pi)^{1/2} \geq 1$. 

Write $G_1$ for the infimum of $G_1(x,a)$ over the stated choices of $(x,a)$; and similarly define~$G_2$ and~$G_3$. 
Note that
$$
 G_1 \geq 
 g_{(-2\sqrt{2} - 1)\eln^2,\mfl + 2\eln}\big( \eln^2 +  6 \eln^2 (k+1 - i) \big) \geq (2\pi)^{-1/2} \big( \tfrac{2}{3\eln} \big)^{1/2} \exp \Big\{ - \tfrac{1}{2\eln} \Big( \big( 6  (k+1 - i)  +  2\sqrt{2} + 2 \big) \eln^2 \Big)^2 \Big\} \, ,
$$
since $\eln \leq \mfl + 2\eln \leq 3\eln/2$. 
Thus,
$$
 G_1 \geq (2\pi)^{-1/2} \big( \tfrac{2}{3\eln} \big)^{1/2} \times \exp \big\{ - \tfrac{1}{2}  ( 6  k  +  5 )^2 \eln^3  \big\} \, \geq \, \exp \big\{ -  18 (   k  +  1 )^2 \eln^3  \big\} \, ;
$$
in the latter inequality, we bounded the pre-$\times$ term below by $\e^{5\const^3}$, using $\e^{30 \const^3 - 1} \leq \big( 3\pi\const \big)^{-3}$ to do so.

The quantity $G_3$ satisfies the same bound.  

Note that
$$
G_2 \geq g_{0,\mfr-\mfl}\big(5\eln^2\big) \geq (2\pi)^{-1/2} \big( \tfrac{1}{2\eln} \big)^{1/2} \cdot \exp \big\{ - \tfrac{1}{2\eln} 25 \eln^4 \big\}
$$
since $\eln \leq \mfr - \mfl \leq 2\eln$. Using
$\e^{3\const^3 - 1} \leq 2^{-6} \pi^{-3} \const^{-3}$,  we find that
$$
G_2 \geq \exp \big\{ -  13 \eln^3 \big\} \, .
$$

In this way, we find that the quantity in~(\ref{e.flyleapbound}) is at least
$$
\exp \big\{ - 36(k+1)^2 \eln^3 - 13 \eln^3 \big\} \geq \e^{36(k+2)^2 \const^3} \, .
$$
 We have found that the first event of the pair whose intersection constitutes $\flyleap$ has $\PP_{\mc{F}}$-probability at least  $\e^{36(k+2)^2 k \const^3}$ when $\mpg$ occurs, because the event is the intersection over $i \in \intint{k}$ of the event 
in~(\ref{e.flyleapbound}).   

Next, note that, for $i \in \intint{k}$,  
$$
 \PP_\mc{F} \Big( \sup_{x \in [\mfl,\mfr]}\big\vert \wien^{[\mfl,\mfr]}\big( i, x \big) \big\vert \leq \eln^2 \Big) = 1 - 2 \exp \big\{ - 2 \tfrac{\eln^4}{\mfr - \mfl} \big\} \, ; 
$$
since $\mpg \subseteq \{ \mfr - \mfl \geq \eln \}$ (and $\eln \geq 1$),  the right-hand side is at least $1 - 2 \exp \big\{ - 2 \eln^3 \big\} \geq 1 - 2 e^{-2} \geq 1/2$ when $\mpg$ occurs.

The bridge $\wien^{[\mfl,\mfr]}\big( i, \cdot \big)$ is independent of the values $\wien(i,\mfl)$ and $\wien(i,\mfr)$. Thus, we obtain~(\ref{e.flyingleap}) and Proposition~\ref{p.wienerpromotion}. \qed

\subsection{Subpolynomial success probability for the jump ensemble} 
Here we prove  Proposition~\ref{p.vaultsuccess}.  
This result offers a lower bound, valid for $\mc{F}$-data realizing $\mpg$, on the probability that a candidate who passes Test~$2$ will pass Test~$4$ (and thus also Test~$3$). Indeed, the equality
\begin{equation}\label{e.testguide}
 \PP_{\mc{F}} \big( \test_4(J) = 1 \big) \cdot {{\bf 1}}_\mpg  \,  = \,  \PP_{\mc{F}} \big( \test_4(\strongjump) = 1 \big) \PP_{\mc{F}} \big( \test_3(J) = 1 \big) \cdot {{\bf 1}}_\mpg
\end{equation}
guides the upcoming proof; given Lemma~\ref{l.jumppromotion}, we see that our task is to find a lower bound, valid for $\mc{F}$-data realizing~$\mpg$, on  the probability  $\PP_{\mc{F}} \big( \test_4(\strongjump) = 1 \big){{\bf 1}}_\mpg$
that the strong jump ensemble  passes the final Test~$4$.

To begin accomplishing this task,
we introduce an augmentation~$\mc{F}^1$ of the missing closed middle $\sigma$-algebra~$\mc{F}$,  and discuss its properties in a sequence of five lemmas.

The pole set $\pole$, which is determined by the curve $\mc{L}_n(k+1,\cdot):[-\eln,\eln] \to \R$, contains at most $2\eln$ elements, as we noted in~(\ref{e.polecard}). Let $\big\{ \pp_r: r \in \intint{\vert \pole \vert} \big\}$ be a list of $\pole$'s elements in increasing order.
In particular, $\pp_1 = \mfl$  and  $\pp_{\vert \pole \vert} = \mfr$.

The strong jump ensemble $\strongjump:\intint{k} \times [\mfl,\mfr] \to \R$ is specified by the data
\begin{enumerate}
\item the difference values $\strongjump\big(i,\pp_{r+1}\big) - \strongjump\big(i,\pp_r \big)$ for $(i,r) \in \intint{k} \times \intint{\vert \pole \vert - 1}$;
\item the standard bridges $\strongjump^{[\pp_r,\pp_{r+1}]}\big(i,\cdot\big):[\pp_r,\pp_{r+1}] \to \R$, $i \in \intint{k}$, $1 \leq r \leq \vert \pole \vert - 1$,
defined between each pair of consecutive pole set values;
\item and the values $\strongjump(i,\mfl)$ for $i \in \intint{k}$.
\end{enumerate}

Write $\mc{F}^1$ for the $\sigma$-algebra generated by $\mc{F}$ and the first item in the list, and let $\PP_{\mc{F}^1}$ denote conditional probability given $\mc{F}^1$. The $\sigma$-algebra generated by $\mc{F}^1$ and the third item in the list will be called $\mc{F}^{13}$, and the law $\PP_{\mc{F}^{13}}$ defined correspondingly.

We wish to note that the conditional distribution under~$\PP_{\mc{F}^{13}}$ of the item $(2)$ data coincides with this law under~$\PP_{\mc{F}^{1}}$: in each case, the product law of $\mc{B}_{k;\bar{0},\bar{0}}^{[p_r,p_{r+1}]}$ over $r \in \llbracket 0, \vert P \vert - 1 \rrbracket$.  
That this is so may be seen by noting that this product law, which is enjoyed under~$\PP_\mc{F}$ by the bridges associated to the Wiener candidate~$\wien$, is unperturbed by the conditioning on $\test_{123}(\wien) = 1$ that results in the ensemble~$\strongjump$.
Indeed, the conditions in these three tests are expressible in terms of the
$k$-vectors $\ovbar{W}(x)$ 
as $x$ varies over $P$, including its endpoints $\mfl$ and $\mfr$. The second item plays no role in determining the values of the $k$ curves in $\strongjump$ over any point in $P$ because these values are dictated by the first and third items.

Similarly, for $i \in \intint{k}$ given, the witness of $\mc{F}^1$ needs one piece of data to know the form of $\strongjump(i,\cdot):[-2\eln,\mfl] \cup P \cup [\mfr,2\eln] \to \R$: the value $\strongjump(i,\mfl)$.  She may reconstruct this curve supposing that this unknown value equals $z \in \R$. Writing $\strongjump^z(i,\cdot):[-2\eln,\mfl] \cup P \cup [\mfr,2\eln] \to \R$ for the reconstructed curve, we have $\strongjump^z(i,p) = z + \strongjump(i,p) - \strongjump(i,\mfl)$ for $p \in P$; note that the $\strongjump$-difference here is $\mc{F}^1$-measurable. 
The proposed value~$z$ also dictates a reconstructed form for $\mc{L}_n(i,s)$ at $s \in [-2\eln,\mfl]$, given by the second line of the formula in~(\ref{e.reconcases}) with $x_i$ taken equal to $z$. In the present context, we will record this formula by $\strongjump^z(i,s)$. (This is out of keeping with the domain of definition of the jump ensemble being $[\mfl,\mfr]$, which excludes such a choice of $s$, but other notational choices such as $\mc{L}_n^z$ clash with earlier usage.) 
Similarly, when $s \in [\mfr,2\eln]$, we set $\strongjump^z(i,s)$ equal to the fourth line in~(\ref{e.reconcases}) with $y_i$ set equal to $z + \strongjump(i,\mfr) - \strongjump(i,\mfl)$.

Let $h_{\mc{F}^1}: \R^k \to [0,\infty)$ denote the density  with respect to Lebesgue measure on $\R^k$ of the conditional law under $\PP_{\mc{F}^1}$ of the $k$-vector~$\overstrongjump(\mfl)$. (We prefer to write ${\overline{J}}_+$ in place of $\ovbar\strongjump$ in using this notation for $k$-vectors.) Using Lemma~\ref{l.brbr}, note that  $h_{\mc{F}^1}\big(\bar{x} \big)$ equals
\begin{equation}\label{e.gdensj.first.new}
   Z_{\mc{F}^1}^{-1}  \prod_{i=1}^k \, \exp\bigg\{-\tfrac{1}{2(\mfl + 2\eln)}\Big(\lppls_n\big(i,-2\eln \big) - x_i \Big)^2  - \tfrac{1}{2(2\eln - \mfr)} \Big( x_i + \strongjump(i,\mfr) - \strongjump(i,\mfl)  - \lppls_n\big(i,2\eln  \big)  \Big)^2 \bigg\} \,\cdot\,\mfone \big( \bar{x} \big) \, ,
\end{equation}
where 
$\mfone \big( \bar{x} \big)$ denotes the indicator function of the event
\begin{eqnarray*}
 & & \Big\{ \, \strongjump^{x_i} \big(i,s\big) >  \strongjump^{x_{i+1}} \big( i+1,s\big) \, \, \, \, \forall \, \, (s,i) \in \big( [-2\eln,\mfl] \cup [\mfr,2\eln] \big) \times \intint{k-1} \, \Big\} \\
 & \cap &  \Big\{ \,  \strongjump^{x_k} \big(k,s \big) >  \mc{L}_n \big( k+1,s \big)\, \, \, \, \forall \, \,  s  \in  [-2\eln,\mfl] \cup [\mfr,2\eln]    \,  \Big\} \\
  & \cap &  \Big\{ \,  \strongjump^{x_i} \big(i,s \big) >  \mc{L}_n \big( k+1,s \big)\, \, \, \, \forall \, \,  (s,i)  \in  \pole \times \intint{k}   \,  \Big\} 
  \, .
\end{eqnarray*}
The normalizing quantity $Z_{\mc{F}^1} \in (0,\infty)$ is $\mc{F}^1$-measurable.

\begin{figure}[ht]
\begin{center}
\includegraphics[height=11cm]{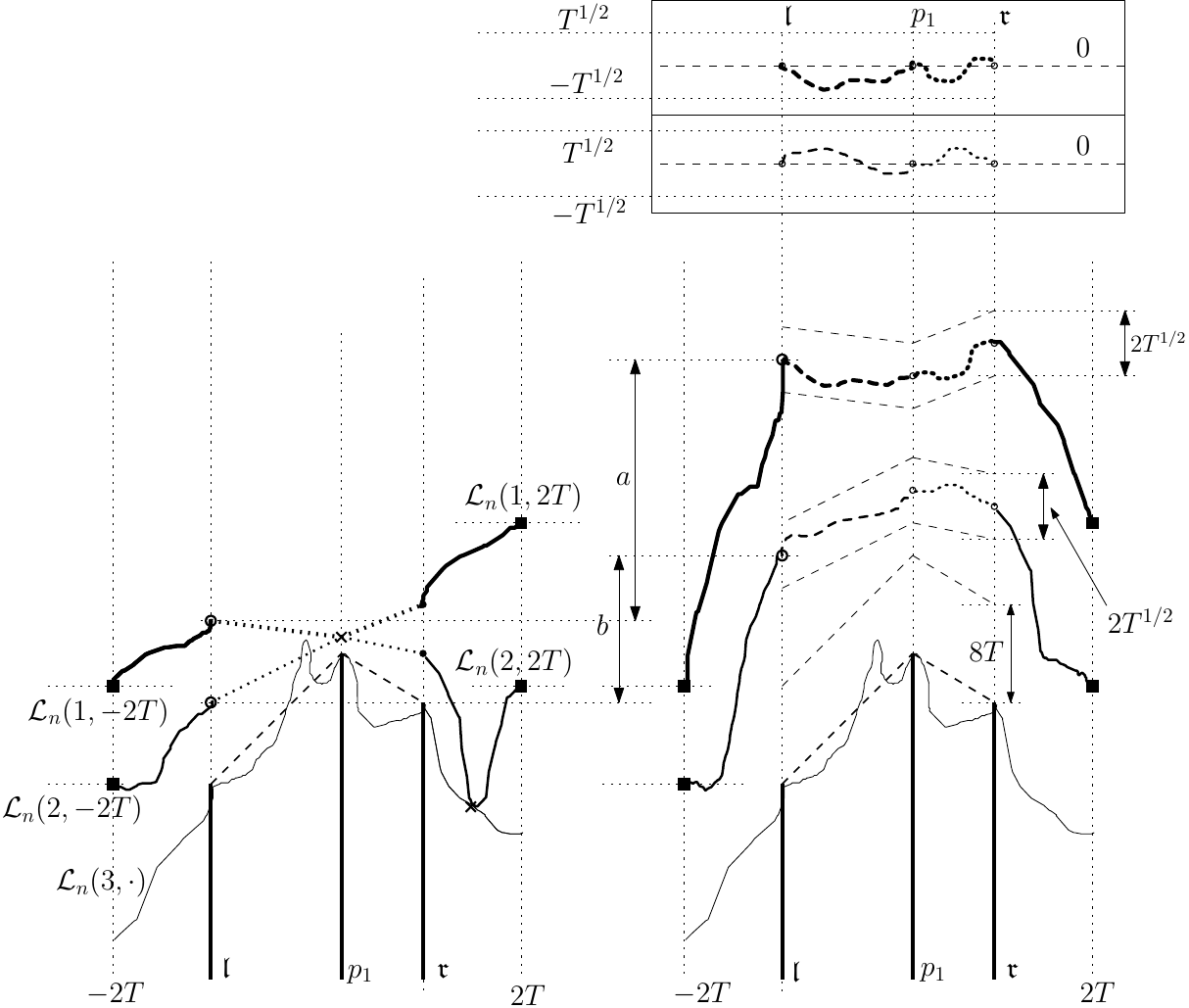}
\caption{The two sketches depict different aspects of the perspective of the witness of $\mc{F}^1$, in an example where $k=2$, $P = \{ \mfl,p_1,\mfr \}$ and $\ipd = 1$. 
The three vertical poles are shown in bold and the tent map they support is dashed in both sketches. In each sketch, the thickest curves correspond to the top curve index and the second thickest to the second curve. 
{\em Left sketch:}  
The vertical location of the pair of circular beads hanging over $\mfl$ 
represents a possible outcome for the value of $\big( \strongjump(1,\mfl),\strongjump(2,\mfl) \big)$ in the eyes of the witness; in this instance, this location is $\barqmin \in \R^2_>$. 
From each bead emanates a dotted piecewise affine curve on $[\mfl,\mfr]$. These two crooked rigid rods are $\mc{F}^1$-determined;  they may be translated vertically by the witness of $\mc{F}^1$ by her rolling the circular beads up or down the line at $\mfl$. In this way,  the pair dictates the values of $\strongjump(1,\cdot)$ and $\strongjump(2,\cdot)$ on~$P$. 
In the language of the proof of Lemma~\ref{l.hfl}, the $(1,2)$-contact point is at $p_1 \in P$, and the $(2,3)$-contact point is in $[\mfr,2\eln]$: see the two crosses. {\em Right sketch:} A sample of the witness's randomness is depicted. The pair of circles $\big( \strongjump(1,\mfl), \strongjump(2,\mfl) \big)$ have been translated from the left sketch location by $(a,b)$, where $b \geq 8\ipd\eln + \eln^{1/2}$ and $a \geq b + 2\eln^{1/2}$. In the upper-right boxes, we see the random bridges: for $i \in \{1,2\}$, the dashed $\strongjump^{[\mfl,p_1]}(i,\cdot)$ and the dotted $\strongjump^{[p_1,\mfr]}(i,\cdot)$, with the higher box containing the bolder $i=1$ curves.
This sample illustrates the proof of Proposition~\ref{p.vaultsuccess} and it verifies $\test_4(\strongjump)=1$.}
\label{f.fonetwosketches}
\end{center}
\end{figure}

We now present a counterpart for the witness of $\mc{F}^1$ to the `corner' Lemma~\ref{l.glv}. The  left sketch in Figure~\ref{f.fonetwosketches}, a few pages hence, may be consulted for an 
illustration of the new lemma's proof. 
\begin{lemma}\label{l.hfl}
There exists a $\mc{F}^1$-measurable random vector $\barqmin \in \R^k_>$ such that 
$$
\Big\{ \, \bar{x} \in \R^k : \mfone\big( \bar{x} \big) = 1 \, \Big\} \, = \, \barqmin \, + \, (0,\infty)^k_> \, .
$$
\end{lemma}
\noindent{\bf Proof.}  Recall the one-piece list reconstructed curve notation in~(\ref{e.reconcases.onepiece}). Given $\mc{F}^1$, and for any vector $\bar{x} \in \R^k$,
consider the $k+1$ curves on $[-2\eln,\mfl] \cup P \cup [\mfr,2\eln]$ the first $k$ of which are given by the
  reconstructed jump ensemble 
$\mc{L}_n^{x_i + \strongjump(i,\cdot) - \strongjump(i,\mfl)}\big(i,\cdot\big):[-2\eln,\mfl] \cup P \cup [\mfr,2\eln] \to \R$, $1 \leq i \leq k$; and the $(k+1)\textsuperscript{st}$ of which equals $\mc{L}_n(k+1,\cdot):[-2\eln,\mfl] \cup P \cup [\mfr,2\eln] \to \R$. Our task is to determine a condition on $\bar{x} \in \R^k$ that characterises this reconstructed ensemble being ordered.

We begin to show that this condition takes the form asserted by the lemma by identifying the value of the vector $\barqmin$.
It is specified so that for each consecutive pair of indices $(i,i+1)$ of this $(k+1)$-curve system, there exists a value of $x$ in the set $[-2\eln,\mfl) \cup \pole \cup (\mfr,2\eln]$ such that the two curves indexed by the pair are equal at~$x$; moreover, if the $i\textsuperscript{th}$ component of the vector is increased while the $(i+1)\textsuperscript{st}$ is held fixed, there is no such contact between the pair of curves. (This value of $x$ is in fact unique $\PP$-almost surely, but we will not use this fact. In this proof, we will call the value the $(i,i+1)$-contact point, in what we hope are the interests of exposition; however, our argument does not depend on the uniqueness that this phrase implicitly asserts.)

This description specifies the value of $\barqmin$. It is useful however to describe a more explicit means of determining this vector's value. Doing so offers an opportunity to consider the perspective of the witness of $\mc{F}^1$, a point of view that is valuable in understanding our use of the strong jump ensemble~$\strongjump$.

The  vector $\barqmin$ may be determined in decreasing order of its component index, similarly to the proof of Lemma~\ref{l.glv}. The highest indexed component  $\qmin_k$ equals the infimum of $q \in \R$ such that
\begin{itemize}
 \item $q \geq \xmin_k$; 
\item
 $q + \strongjump(k,\pp) - \strongjump(k,\mfl) \geq \mc{L}_n(k+1,\pp)$ for all $\pp \in \pole$;
 \item and $q + \strongjump(k,\mfr) - \strongjump(k,\mfl) \geq \ymin_k$.
\end{itemize} 
Indeed, Lemma~\ref{l.glv} shows that the first, and third, conditions correspond to the $(k,k+1)$-contact point lying in~$[-2\eln,\mfl]$, or $[\mfr,2\eln]$. The second clearly corresponds to this point lying in the pole set~$\pole$. Note also that the $\strongjump$-differences are $\mc{F}^1$-measurable while the data $\xmin_k$ and $\ymin_k$, as well as the form of $\mc{L}_n(k+1,\cdot)$, is measurable in the smaller $\sigma$-algebra~$\mc{F}$. Thus, the witness of~$\mc{F}^1$
is certainly equipped to determine~$\qmin_k$ according to this rule. 

At the generic step, $\qmin_i$ will be determined for $i \in \intint{k-1}$. The values $\qmin_j$ have been decided for $j \in \llbracket i+1,k \rrbracket$, and now 
 $\qmin_i$ is set equal to the infimum of $q \in \R$ such that
\begin{itemize}
 \item $q - \qmin_{i+1} \geq \xmin_i - \xmin_{i+1}$; 
\item
 $q + \strongjump(i,\pp) - \strongjump(i,\mfl) \geq \mc{L}_n(k+1,\pp)$ for all $\pp \in \pole$;
 \item and $\Big( q + \strongjump(i,\mfr) - \strongjump(i,\mfl) \Big)  - \Big( \qmin_{i+1}  + \strongjump(i+1,\mfr) - \strongjump(i+1,\mfl)  \Big) \geq  \ymin_i - \ymin_{i+1}$.
\end{itemize} 
The three cases correspond to the location of the $(i,i+1)$-contact point just as they did in the specification of $\qmin_k$. Take the third item, for example. If the witness of $\mc{F}^1$ considers the eventuality that the unknown $\strongjump(i,\mfl)$ equals $q$, the outcome that $\strongjump(i,\mfr)$ equals $q + \strongjump(i,\mfr) - \strongjump(i,\mfl)$ would be dictated; alongside the circumstance that $\strongjump(i+1,\mfl)$ equals $\qmin_{i+1}$, this would force $\strongjump(i,\mfr) - \strongjump(i+1,\mfr)$ to equal the third item left-hand side; thus, Lemma~\ref{l.glv} shows that equality in the third condition stipulates a point of contact between the $i\textsuperscript{th}$ and $(i+1)\textsuperscript{st}$ curves somewhere in $[\mfr,2\eln]$ and that there is no such contact when equality is replaced by strict inequality. 

That  $\bar{x} \in \R^k$ satisfies $\mfone(\bar{x})=1$ precisely when $\bar{x}- \barqmin \in (0,\infty)^k_>$ follows by the reasoning in the final paragraph of Lemma~\ref{l.glv}'s proof. \qed

\medskip

Let $\mcgone$ denote the $\mc{F}_1$-measurable event that
\begin{itemize}
\item $\strongjump(i,\pp)  - \strongjump(i,\mfl) \geq - 37 k^{3/2} \eln^2$ for $(i,\pp) \in \intint{k} \times \pole$;
\item $\strongjump(i,\mfr)  - \strongjump(i,\mfl) \leq  37 k^{3/2}  \eln^2$ for $i \in \llbracket 1,k \rrbracket$.
\end{itemize}

\begin{lemma}\label{l.estgone}
When $\mpg \cap \mcgone$ occurs, 
$$
 -  \eln^2 \leq \qmin_i \leq 74 k^{5/2} \eln^2 
$$
for each $i \in \intint{k}$.
\end{lemma}
\noindent{\bf Proof.} We begin by arguing that 
$$
 - \eln^2 \leq \qmin_k \leq \big( 37 k^{3/2} + 1 \big) \eln^2 
$$
on $\mpg \cap \mcgone$.
Recall that $\qmin_k$ is an infimum over $q \in \R$ satisfying the first set of three bullet-pointed lower bounds in the proof of Lemma~\ref{l.hfl}. The first of these conditions implies that $\qmin_k$ is at least $-\eln^2$ in light of $\mpg \subseteq \big\{ \xmin_k \geq -\eln^2 \big\}$.

To find an upper bound on $\qmin_k$, we must find for each bullet point an admissible value of~$q$. In the first case, that $\mpg$ entails $\xmin_k \leq \eln^2$ implies that $q$ is satisfactory when it equals $\eln^2$.  In the second and the third cases, $q =  \big( 37 k^{3/2} + 1 \big) \eln^2$ works.  In the second case, this is a consequence of the $\mcgone$ constraint that $\strongjump(k,\pp) - \strongjump(k,\mfl) \geq - 37 k^{3/2} \eln^2$ as well as $\pole \subset [\mfl,\mfr]$ and the $\mpg$ constraint that $\mc{L}_n\big(k+1,x\big) \leq \eln^2$ for $x \in [\mfl,\mfr]$. In the third, it follows from $\ymin_k \leq \eln^2$ on $\mpg$ and this same $\mcgone$ constraint.

We have established the base case $i = k$ of the assertion that
$$
 - \eln^2 \leq \qmin_i \leq \Big( 37 k^{3/2} + 1 + (k-i) \big(  74 k^{3/2} + 2  \big) \Big) \eln^2 \, ,
$$
which we now verify by an induction in decreasing order on the variable $i \in \intint{k}$.

Taking $i \in \intint{k-1}$ and assuming these bounds for index~$i+1$, we recall the three bullet point inequalities specifying $\qmin_i$ and note that the first of these already demonstrates that $\qmin_i \geq \qmin_{i+1}$ in view of $\xmin_i \geq \xmin_{i+1}$ (a fact which is due to $\barxmin \in \R^k_\geq$). Thus, we obtain the lower bound in the inductive hypothesis at index~$i$.
 
Regarding the upper bound, note that the three bullet point conditions are respectively satisfied by the following values of $q$:
\begin{itemize}
\item $q = \qmin_{i+1} + 2\eln^2$, since $\xmin_i \leq \eln^2$ and $\xmin_{i+1} \geq  - \eln^2$ on $\mpg$;
\item $q = \big(  37 k^{3/2} + 1  \big) \eln^2$, since $\strongjump(i,\pp) - \strongjump(i,\mfl) \geq - 37 k^{3/2} \eln^2$ on $\mcgone$ and 
$\mc{L}_n \big( k+1 , \pp \big) \leq \eln^2$ for $\pp \in \pole$ on $\mpg$;
\item and $q = \qmin_{i+1} + 2\eln^2 + 74 k^{3/2} \eln^2$,
since $\ymin_i - \ymin_{i+1} \leq 2\eln^2$ on $\mpg$, and 
$\strongjump(i,\mfr) - \strongjump(i,\mfl) \geq - 37 k^{3/2} \eln^2$ and $\strongjump(i+1,\mfr) - \strongjump(i+1,\mfl) \leq 37 k^{3/2} \eln^2$ on $\mcgone$.
\end{itemize}
Hence, $\qmin_i$ is seen to be at most the maximum of
$\big(  37 k^{3/2} + 1  \big) \eln^2$ and 
$\qmin_{i+1} + \big(  74 k^{3/2} + 2 \big) \eln^2$.

The inductive hypothesis upper bound at index $i$ thus follows from its counterpart at index~$i+1$.

The bounds stated in Lemma~\ref{l.estgone} follow directly from the inductively established bounds. \qed

\begin{lemma}\label{l.gprime}
$$
\PP_\mc{F} \big(  \mcgone^c  \big) \cdot {\bf 1}_\mpg   \leq \e^{15 k^3 \const^3}  \, . 
$$
\end{lemma}
\noindent{\bf Proof.}
Define the counterpart to the event $\mcgone$ for the Wiener candidate $W$, namely
$$
\mathsf{A} =  \bigcap_{(i,\pp) \in \intint{k} \times \pole} \Big\{ \wien(i,\pp)  - \wien(i,\mfl) \geq - 37k^{3/2}\eln^2 \Big\} \, \, \cap \, \,  \bigcap_{i \in \llbracket 1,k \rrbracket}  \Big\{ \wien(i,\mfr)  - \wien(i,\mfl) \leq  37k^{3/2}\eln^2  \Big\} \, .
$$
We claim that
\begin{equation}\label{e.acbound}
\PP_{\mc{F}} \big(  \mathsf{A}^c  \big) \, {\bf 1}_{\mpg} \leq \e^{160 k^3 \const^3}  \, . 
\end{equation}
To see this, note that, under $\PP_\mc{F}$, the Wiener candidate ensemble member $x \to \wien(i,x):[\mfl,\mfr] \to \R$, for given $i \in \intint{k}$, 
has the marginal law on $[\mfl,\mfr]$ of Brownian bridge $B$ distributed as $\mc{B}_{1;u_i,v_i}^{[-2\eln,2\eln]}$, where $u_i$ and $v_i$ each lie within $\eln^2$ of $-2\sqrt{2} \eln^2$ should the $\mc{F}$-measurable event $\mpg$ occur. Consequently, if $\mpg$ occurs, then if it is the case that $\big\vert \wien(i,x)  - \wien(i,\mfl) \big\vert >  37k^{3/2}\eln^2$ for any given (or indeed $\mc{F}$-measurable) choice of $x \in [\mfl,\mfr]$, it is easily seen that the standard bridge $[-2\eln,2\eln] \to \R:x \to B^{[-2\eln,2\eln]}(i,x)$ associated to $B$ must have maximum absolute value at least $\tfrac{1}{2} \big( 37k^{3/2}\eln^2 - 2\eln^2\big) \geq 18k^{3/2}\eln^2$ (since $k \geq 2$). This latter eventuality has probability at most $2\exp \big\{- \tfrac{18^2}{2} k^3 \eln^3 \big\}$ as we see by applying Lemma~\ref{l.maxfluc} with $r=18k^{3/2}\eln^2$ and $b-a = 4\eln$. There are $k \vert \pole \vert + k \leq k (2\eln + 2)$ such eventualities that may cause $\mathsf{A}$ to fail to occur. Recalling that $\eln = \const \big( \log \e^{-1} \big)^{1/3} \geq 1$,
$$
\PP_{\mc{F}} \big(  \mathsf{A}^c  \big) \, {\bf 1}_{\mpg} \leq 8 k \const  \big( \log \e^{-1} \big)^{1/3} \cdot \e^{162 k^3 \const^3} \, .
$$
The first factor on the right-hand side is at most $\e^{-\const^3}$ because 
$\e^{3\const^3 - 1} \leq \big( 8 k \const \big)^{-3}$.
We have verified~(\ref{e.acbound}).

Note that
$$
\PP_\mc{F} \big( \mcgone^c  \big) \, {\bf 1}_{\mpg}
= \frac{\PP_\mc{F}  \Big(  \mathsf{A}^c \, , \, \test_{123}(\wien)   = 1  \,  \Big)}{
\PP_\mc{F}  \big( \test_{123}(\wien)   = 1    \big)} \, {\bf 1}_{\mpg} \ \, . 
$$
By  Proposition~\ref{p.wienerpromotion} and~(\ref{e.acbound}), the right-hand side is when $\mpg$ occurs at most 
$$
\e^{-36 k (k+2)^2 \const^3} \, 2^{k} \cdot \e^{160 k^3 \const^3} \leq 2^k \e^{16 k^3 \const^3} \leq \e^{15 k^3 \const^3} \, ,
$$ 
since $k \geq 2$, $\e < 1/2$ and $\const \geq 1$. The proof of Lemma~\ref{l.gprime} is complete. \qed

\begin{lemma}\label{l.misigmai}
The conditional distribution under $\PP_{\mc{F}^1}$ of the third item vector  $\overstrongjump(\mfl)$
is the law of an independent sequence $\ovbar{N} = \big( N_i: i \in \intint{k} \big)$ of normal random variables conditionally on the event that $\ovbar{N} \in \barqmin + [0,\infty)^k_\geq$. 

The component random variables $N_i$ share a common, $\mc{F}$-measurable, variance $\sigma^2$ that satisfies 
$\sigma^2 \in \eln \cdot \big[1/2,3/4 \big]$  when $\mpg$ occurs. The mean of $N_i$, which we denote by $m(i)$, is $\mc{F}^1$-measurable.  When $\mpg \cap \mcgone$ occurs, it satisfies $\vert m(i) \vert \leq 38 k^{3/2} \eln^2$. 
\end{lemma}
\begin{lemma}\label{l.doublegauss}
Let $\ell_1 < a < b < \ell_2$ and $z,j \in \R$. Under the bridge law 
$\mc{B}_{1;0,z}^{[\ell_1,\ell_2]}$ conditioned on $B(b) - B(a) = j$, the conditional distribution of $B(a)$ is normal with mean $m$ and variance $\sigma^2$, where $m = \tfrac{z-j}{\ell_2 - b} \sigma^2$ and
$\sigma^{-2}  = (a - \ell_1)^{-1} + (\ell_2 - b)^{-1}$.
\end{lemma}
\noindent{\bf Proof.} 
Let $N_1$ and $N_2$ be  two independent normal random variables.
The first has  mean zero and variance $a - \ell_1$; the second, mean $z$ and variance $\ell_2 - b$. Under the conditioning in the lemma, the random variable $\big( B(a),B(b) \big)$ has the law of the pair $(N_1,N_2)$ conditioned on $N_2 = N_1 + j$.

Thus, the conditioned random variable
$B(a)$ adopts the value $x$ with density given up to normalization by 
$$
\exp \Big\{ - x^2 \tfrac{1}{2(a-\ell_1)} \Big\} \cdot \exp \Big\{ - (x + j - z)^2 \tfrac{1}{2(\ell_2 -b)} \Big\} \, .
$$
When it is normalized, this quantity equals $g_{m,\sigma^2}(x)$. \qed

\medskip

\noindent{\bf Proof of Lemma~\ref{l.misigmai}.} 
We apply Lemma~\ref{l.doublegauss}, setting $\ell_1 = -2\eln$, $a = \mfl$, $b = \mfr$, $\ell_2 = 2 \eln$, $z = \mc{L}_n(i,2\eln) - \mc{L}_n(i,-2\eln)$ and $j = \strongjump(i,\mfr) - \strongjump(i,\mfl)$. In doing so, we learn that 
$$
 \sigma^{-2} = \big( \mfl + 2\eln \big)^{-1} + \big( 2\eln - \mfr \big)^{-1} 
$$  
and
$$
 m(i) =  \frac{\mc{L}_n(i,2\eln) - \mc{L}_n(i,-2\eln) - j_i}{(\mfl + 2\eln)^{-1}(2\eln - \mfr) + 1}  \, ,
$$
where $j_i$ denotes $\strongjump(i,\mfr) - \strongjump(i,\mfl)$. 

The occurrence of $\mpg$ entails that $\mfl \in [-\eln,-\eln/2]$
and $\mfr \in [\eln/2,\eln]$, 
whence
$\eln/2 \leq \sigma^2 \leq 3\eln/4$.

The denominator in the expression for $m(i)$ is at least one. On $\mpg$, the quantities $\mc{L}_n(i,-2\eln)$ and $\mc{L}_n(i,2\eln)$ differ by at most $2\eln^2$;  and, on $\mcgone$, $\vert j_i \vert \leq 37k^{3/2} \eln^2$. That $\vert m(i) \vert \leq 38 k^{3/2} \eln^2$ follows from $k \geq 2$. \qed

\medskip

\noindent{\bf Proof of Proposition~\ref{p.vaultsuccess}.}
For $\bar{q} \in \R^k_>$ and $\alpha,\beta > 0$, let the {\em pair
separated} set $\pairsep_{\alpha,\beta,\bar{q}} \subseteq \R^k_>$ be given by 
$$
\pairsep_{\alpha,\beta,\bar{q}} = \Big\{ \bar{x} \in \R^k:  (x_i - q_i) - (x_{i+1}-q_{i+1}) \geq 2\alpha \, \, \, \textrm{for all $i \in \intint{k-1}$; and} \, \, \, x_k \geq q_k + \beta \Big\} \, .
$$ 

We will prove Proposition~\ref{p.vaultsuccess} by adopting the perspective of the witness of $\mc{F}^1$ and considering the eventuality that
$$
\overstrongjump(\mfl)  \in  \pairsep_{\alphapr,\betapr,\barqpr} \, \, \,  
\textrm{where} \, \, \, \big( \alphapr , \betapr , \barqpr \big)  = \big( \eln^{1/2} ,  8 \ipd \eln + \eln^{1/2} , \barqmin \big) \in (0,\infty) \times (0,\infty) \times \R^k_> \, . 
$$
Indeed, we set the value of $(\alpha',\beta',\bar{q}')$ in this way throughout the proof of this proposition (and thus for the remainder of Section~\ref{s.subpoly}).

There are two main steps for the proof: first, in Proposition~\ref{p.jpass}, we will argue that the above eventuality results in a subpolynomial decay in $\e$ for the conditional probability of $\test_4(\strongjump)=1$; then, in Proposition~\ref{p.jabovecorner}, we will find that the eventuality occurs with such a $\PP_{\mc{F}^1}$-probability. The right sketch in Figure~\ref{f.fonetwosketches} illustrates the ideas.

\begin{proposition}\label{p.jpass}
$$
 \PP_{\mc{F}^1} \Big( \test_4(\strongjump) = 1 \, \Big\vert \,  \overstrongjump(\mfl)  \in  \pspr \Big) \geq  \big( 1 - 2 e^{-1} \big)^{2\eln k}  
 \, . 
$$
\end{proposition}
The idea of the proof of this result: under the conditioning in question, $\strongjump(k,\mfl)$ has a margin of $8\ipd\eln + \eln^{1/2}$ over the minimum needed not to definitely violate an avoidance constraint in the eyes of the witness of~$\mc{F}^1$. The margin of $8\ipd\eln$ is needed because, as we shall see in Lemma~\ref{l.jumppole}, this is the amount by which $\mc{L}_n(k+1,\cdot):[\mfl,\mfr] \to \R$ may rise above the tent map; the extra margin of $\eln^{1/2}$, and an additional margin of $2\eln^{1/2}$ associated to the lower indices, permit the use of channels  of width $2\eln^{1/2}$ about the affine interpolation of consecutive jump ensemble curve values over the pole set, with the channels being disjoint from each other and from the graph of $\mc{L}_n(k+1,\cdot)$. This channel width is high enough that the channels may be shown via Lemma~\ref{l.mbfprob} to have a reasonable $\PP_{\mc{F}^1}$-probability of   
housing the ensemble's curves.

\medskip

\noindent{\bf Proof of Proposition~\ref{p.jpass}.} To begin the rigorous argument,
 recall the system of standard bridges $\strongjump^{[\pp_r,\pp_{r+1}]}\big(i,\cdot\big):[\pp_r,\pp_{r+1}] \to \R$, $i \in \intint{k}$, $0 \leq r \leq \vert \pole \vert - 1$, in the second item of the three-item list used in specifying~$\mc{F}^1$. Under the law $\PP_{\mc{F}^1}$ given that $\overstrongjump(\mfl)$  is an element of  $\pspr$, this system has the law of independent standard Brownian bridges on the respective intervals. It is straightforward to find a simple criterion for these bridges that is sufficient for the success of the final test condition $\test_4(\strongjump) = 1$ under this conditional law. Namely, the {\em modest bridge fluctuation} event $\mathsf{MBF}$ is given by
$$
\mathsf{MBF} = \Big\{ \big\vert \strongjump^{[\pp_r,\pp_{r+1}]}\big(i,x\big) \big\vert \leq \eln^{1/2}: i \in \intint{k} \, , \, x \in \big[\pp_r,\pp_{r+1} \big] \, , \, 0 \leq r \leq \vert \pole \vert - 1 \Big\} \, .
$$  
The second of the next two lemmas shows that $\mathsf{MBF}$ is such a criterion; the first that this event is not atypical: its probability decays to zero at most subpolynomially as $\e \searrow 0$. 
\begin{lemma}\label{l.mbfprob}
We have that
$$
\PP_{\mc{F}^{13}} \big( \mathsf{MBF} \big) \geq \big( 1 - 2e^{-1} \big)^{2\eln k} \, .
$$
\end{lemma}
\begin{lemma}\label{l.mbfinc}
The inclusion
$$
 \mathsf{MBF} \, \cap \, \Big\{ \, \overstrongjump(\mfl) \in \pspr \, \Big\} \,  \subseteq \, \big\{ \, \test_4(\strongjump) = 1 \, \big\}
$$
holds up to a $\PP$-null set.
\end{lemma}
Since Proposition~\ref{p.jpass} is a direct consequence of these two lemmas, its proof is complete subject to deriving them. \qed

\medskip

\noindent{\bf Proof of Lemma~\ref{l.mbfprob}.} We have noted that, under $\PP_{\mc{F}^{13}}$, the bridge $\strongjump^{[\pp_r,\pp_{r+1}]}\big(i,\cdot \big)$ has the standard bridge law $\mc{B}_{1;0,0}^{[\pp_r,\pp_{j+1}]}$, independently of the other bridges.  By Lemma~\ref{l.maxfluc}, this bridge has maximum exceeding $\eln^{1/2}$ with probability at most $\exp\big\{- 2\eln (\pp_{r+1} - \pp_r)^{-1} \big\}$; and likewise for the minimum being less than $-\eln^{1/2}$.

Since the pole set $\pole$ is a subset of $[\mfl,\mfr] \subseteq [-\eln,\eln]$,  $\pp_{r+1} - \pp_r \leq 2\eln$ for $0 \leq r \leq \vert \pole \vert  - 1$. Thus, $\big\vert J^{[\pp_r,\pp_{r+1}]} \big\vert \leq \eln^{1/2}$  with $\PP_{\mc{F}^{13}}$-probability at least $1 - 2e^{-1}$. By~(\ref{e.polecard}), the number of bridges is at most $k \vert \pole \vert \leq 2 \eln k$, whence we obtain Lemma~\ref{l.mbfprob}. \qed

\medskip

\noindent{\bf Proof of Lemma~\ref{l.mbfinc}.} Consider the law $\PP_{\mc{F}^1}$.
By Lemma~\ref{l.hfl}, the vector $\overstrongjump(\mfl)$ is a random element of the set $\barqmin + (0,\infty)^k_>$. If this random vector were to adopt the value $\barqmin$ (the minimal possible value, which is in fact only in the closure of the support of the vector's law),  then each of the vectors $\overstrongjump(x)$, for $x$ in the pole set $\pole$, would be elements of $\R^k_\geq$, with equality between consecutive components being possible, but not a strict violation of decreasing order.  When the event $\overstrongjump(\mfl)  \in  \pspr$ occurs, we may think of this vector of values being dynamically constructed as follows: from an initial vector equal to $\barqmin$, the value $\strongjump(k,\mfl)$ is determined by increasing from the initial value $\qmin_k$ by at least $8\ipd\eln + \eln^{1/2}$; at the same time, each lower indexed value $\strongjump(i,\mfl)$ receives the same upward push from an initial location of $\qmin_i$; $\strongjump(k-1,\mfl)$ is then further pushed up by at least $2\eln^{1/2}$, with lower indexed $\strongjump$-values likewise rising; and each $\strongjump(i,\mfl)$ in decreasing $i$ receives its own upward push of $2\eln^{1/2}$, similarly also forcing up lower indexed values, until the final outcome values of the vector are obtained.

The values of the $\strongjump$-vector at $\mfl$ dictate the corresponding values at  any $x \in \pole$. Indeed, we may view the $k$-vector of $\strongjump$-values above $x$ as being dynamically obtained from an initial state in correspondence with the choice that 
$\overstrongjump(\mfl) = \barqmin$ by the same process of upward pushes.

For $i \in \intint{k}$, consider the piecewise affine function $A(i,\cdot):[\mfl,\mfr] \to \R$ that adopts the values $\strongjump(i,x)$ for each $x \in \pole$. We may associate a corridor $\cdor_i \subset [\mfl,\mfr] \times \R$,
$$
\cdor_i  \, = \, \Big\{ \, (x,y) \in [\mfl,\mfr] \times \R : \big\vert  y - A(i,x) \big\vert \leq \eln^{1/2} \, \Big\} \, .
$$

Note then that the occurrence of the event $\overstrongjump(\mfl)  \in \pspr$ ensures that
\begin{itemize}
\item the corridors are ordered, in the sense that, for any $x \in [\mfl,\mfr]$ and $i \in \intint{k-1}$, the interval of $y$-coordinates such that $(x,y) \in \cdor_i$ lies strictly to the right of the corresponding interval for $\cdor_{i+1}$;
\item and the corridor $\cdor_k$ lies above the graph of the tent map by more than a distance of $8\ipd\eln$, in that the corridor's lower boundary function  $[\mfl,\mfr] \to \R: x \to A(i,x) - \eln^{1/2}$ strictly exceeds $[\mfl,\mfr] \to \R: x \to \tent(x) + 8\ipd\eln$.
\end{itemize}
The event $\mathsf{MBF}$ entails that each jump ensemble curve $\strongjump(i,\cdot):[\mfl,\mfr] \to \R$ has a graph that is a subset of the corridor~$\cdor_i$. When  $\overstrongjump(\mfl) \in \pspr$ also occurs, we see that $\strongjump(i,x) > \strongjump(i+1,x)$ for $(i,x) \in \intint{k-1} \times [\mfl,\mfr]$. Moreover, the two events' occurrence forces 
$\strongjump(k,x) > \tent(x) + 8\ipd\eln$ for $x \in [\mfl,\mfr]$. Since $\mc{L}_n(k+1,x) \leq \tent(x) + 8\ipd\eln$ for such $x$ by the next stated Lemma~\ref{l.jumppole}, we see that $\strongjump(k,x) > \mc{L}_n(k+1,x)$ for $x \in [\mfl,\mfr]$. We have verified that the occurrence of the two events ensures that $\test_4(\strongjump) = 1$; this completes the proof of Lemma~\ref{l.mbfinc}. \qed

\begin{lemma}\label{l.jumppole}
For each $x \in [\mfl,\mfr]$, $\mc{L}_n(k+1,x) \leq \tent(x) + 8\ipd \eln$.
\end{lemma}
\noindent{\bf Proof.}
Suppose the result fails at $x \in [\mfl,\mfr]$.
Let $x \in [\pp_1,\pp_2]$ where $\pp_1$ and $\pp_2$ are consecutive elements in~$\pole$.  Write $\ell_{\pp_1,\pp_2}:\R \to \R$ for the affine function whose graph interpolates the points $\big(\pp_1,\mc{L}_n(k+1,\pp_1) \big)$ and $\big(\pp_2,\mc{L}_n(k+1,\pp_2) \big)$.
Note that the graph of the function $[\mfl,\mfr]: \cdot \to \mc{L}_n(k+1,\cdot)$ intersects the region $\big\{ (u,v): v> \ell_{\pp_1,\pp_2}(u) + 8 \ipd \eln \big\}$. 
This forces the set of extreme points of the convex hull of this graph also to have such an intersection, at a location that we label $\big( y , \mc{L}_n(k+1,y) \big)$.
Note that $y$ is an element of the set $\xext$ of which $\pole$ is by definition a subset; furthermore, if we let $y \in [\pp_3,\pp_4]$ for consecutive pole set values $\pp_3$ and $\pp_4$, the definition of $\pole$ implies that $\vert \pp_3 - y \vert \wedge \vert \pp_4 - y \vert \leq \ipd$. Suppose that $\vert \pp_3 - y \vert \leq \ipd$; the other case is no different. 
Noting that $\tent(y)$ is bounded above by $\ell_{\pp_1,\pp_2}(y)$, 
and that the gradient of each of the planar line segments comprising the graph of $[\mfl,\mfr] \to \R: x \to \tent(x)$  is in absolute value at most $4\eln$,
we see that
$$
\mc{L}_n \big( k+1,y \big) > \ell_{\pp_1,\pp_2}(y) + 8 \ipd \eln \geq \tent(y) + 8\ipd \eln \geq \tent(\pp_3) - 4\eln \vert \pp_3 - y \vert + 8\ipd \eln \geq \tent(\pp_3) + 4\ipd \eln \, ,
$$
where the last inequality used $\vert \pp_3 - y \vert \leq \ipd$.
The points $\big(\pp_3, 
\mc{L}_n ( k+1, \pp_3 ) \big)$ and  $\big(y, 
\mc{L}_n ( k+1, y ) \big)$  lie on the graph of the convex hull of 
$[\mfl,\mfr]: \cdot \to \mc{L}_n(k+1,\cdot)$ (since they are extreme points of this graph). Thus, the absolute value of the gradient of the line that connects them is at most $4\eln$,
and~so 
$$
\mc{L}_n \big( k+1,y \big) \leq 
\mc{L}_n \big( k+1, \pp_3 \big) + 4\ipd \eln \, . 
$$
Since $\tent(\pp_3)  = 
\mc{L}_n \big( k+1, \pp_3 \big)$, we have arrived at a contradiction and thus complete the proof. \qed

\medskip

The elements needed in the proof of Proposition~\ref{p.jpass} having been provided, we reach in the next result the second stage  on the road to proving Proposition~\ref{p.vaultsuccess}.

\begin{proposition}\label{p.jabovecorner}
Recalling that  $\big( \alphapr , \betapr , \barqpr \big)  = \big( \eln^{1/2} ,  8 \ipd \eln + \eln^{1/2} , \barqmin \big)$,
$$
 \PP_{\mc{F}^1} \Big( \ovbar{J}(\mfl) \in \pspr \Big) \geq  
 \tfrac{1}{4} \exp \big\{ -  3970 k^{7/2} \ipd^2  \eln^2  \big\} \,  {\bf 1}_{\mpg \cap \mcgone} 
    \, . 
$$
\end{proposition}
\noindent{\em Remark.} It may be worth emphasising a basic feature of the upcoming proof, which is already apparent if we take $k=1$. (With this choice made, $\strongjump$ equals $J$ in law.) By Lemmas~\ref{l.hfl} and~\ref{l.misigmai}, the distribution of $J(1,\mfl)$ under~$\PP_{\mc{F}^1}$ is given by a Gaussian random variable~$N$, with mean of order~$\eln^2$ and variance of order~$\eln$, conditioned on $N \geq \qmin_1$, where $\qmin_1$ is itself of order~$\eln^2$. Note then that the conditional probability that 
the underlying Gaussian random variable~$N$ exceeds $\qmin_1$
behaves as $\exp \big\{ - O(1) (\eln^2)^2/\eln \big\} = \exp \big\{ - O(1) \eln^3 \big\} = \e^{O(1)}$, an unacceptably small term. However, when we discuss, as we need to, the conditional probability, given that $N$ is at least $\qmin_1$, that $N$ exceeds the sum of $\qmin_1$ and an extra margin $x$, we find that, if $x = o(\eln^2)$, this quantity behaves as
$$
 \frac{\exp \big\{ -  \big( \Theta(1) \eln^2 + x \big)^2 \eln^{-1} \big\}}{\exp \big\{ -  \big(\Theta(1)\eln^2 \big)^2 \eln^{-1} \big\}} \, = \,
  \exp \big\{ \Theta(1) \eln x \big\} \, ,
$$ 
which is the tolerable $\e^{o(1)}$. Note that there is a {\em cancellation of first order kinetic costs} in the displayed line, with the smaller cross term in the exponential left to dominate on the right-hand side. In the context of the upcoming proof, we will take $x = 8\ipd \eln + \eln^{1/2}$. Pointing out this cancellation is an almost trivial observation, but 
similar such cancellations will play important roles in several arguments.

\medskip

\noindent{\bf Proof of Proposition~\ref{p.jabovecorner}.}
Let $\bar{m} \in \R^k$ denote the vector $\big( m(i): i \in \intint{k} \big)$.
By Lemma~\ref{l.misigmai}, the probability in question is given by 
$\nu^k_{\bar{m},\sigma^2} \big( \pspr  \, \big\vert \, \barqmin + [0, \infty)^k_> \big)$, where the $\mc{F}^1$-measurable quantities $\bar{m}$ and $\sigma^2$ satisfy bounds stated in that lemma.

\begin{lemma}\label{l.bmalpha}
Let $\bar{q} \in \R^k_>$, and write $Q = \bar{q} + [0,\infty)^k_>$.
For $M,\alpha,\beta > 0$ arbitrary, write 
$$
B_{M,\alpha,\beta} = Q \cap \pairsep_{\alpha,\beta,\bar{q}}^c \cap \big\{ \bar{x} \in \R^k: x_1 \leq M \big\} \, .
$$

If $\bar{a} \in [-M,\infty)^k$, then, for any $\psi^2 > 0$,
$$
 \nu^k_{\bar{a},\psi^2} \big( B_{M,\alpha,\beta} \, \big\vert \, Q \big) \leq \big( 1 + \kappa \big)^{-1} \, ,
$$ 
where
$$
\kappa = \exp \Big\{ - k\psi^{-2} \Big( 2M\big( 2\alpha (k-1) + \beta \big) + \tfrac{1}{2} \big( 2\alpha(k-1) + \beta \big)^2 \Big) \Big\} \, .
$$
\end{lemma}
\noindent{\bf Proof.} 
Recall that $\vecint \in \R^k_>$ denotes the vector $(k-1,k-2,\cdots,0)$. Also write $\bar{1} \in \R^k$ 
for the constant vector whose components have value one.
For $\alpha,\beta > 0$, define the translate $\Phi = \Phi_{\alpha,\beta}: \R^k \to \R^k$ according to $\Phi(\bar{x}) = \bar{x} + 2\alpha \vecint + \beta \bar{1}$. Note that $\Phi(Q) \subseteq Q$ and $\Phi\big(Q \cap \pairsep_{\alpha,\beta,\bar{q}}^c\big) \subseteq \pairsep_{\alpha,\beta,\bar{q}}$. Thus, $B_{M,\alpha,\beta}$ and $\Phi\big(B_{M,\alpha,\beta}\big)$ are disjoint subsets of $Q$, so that
$$
 \nu^k_{\bar{a},\psi^2} \big( B_{M,\alpha,\beta} \, \big\vert \, Q \big)  \leq 
 \frac{\nu^k_{\bar{a},\psi^2} \big( B_{M,\alpha,\beta}  \big)}{\nu^k_{\bar{a},\psi^2} \big( B_{M,\alpha,\beta}  \big) \, + \, \nu^k_{\bar{a},\psi^2} \big( \Phi ( B_{M,\alpha,\beta} ) \big)} \,. 
 $$
It is enough then to establish that 
\begin{equation}\label{e.jacobianphi}
 \nu^k_{\bar{a},\psi^2} \big( \Phi ( B_{M,\alpha,\beta} ) \big)
 \geq  \kappa \cdot \nu^k_{\bar{a},\psi^2} \big(  B_{M,\alpha,\beta}  \big) \, ,
\end{equation}
and this we now do. 
Since the Jacobian of the translation $\Phi$ equals one, we have that
\begin{eqnarray*}
 & & \nu^k_{\bar{a},\psi^2} \big( \Phi ( B_{M,\alpha,\beta} ) \big) 
  =   \int_{\Phi ( B_{M,\alpha,\beta} )} g^k_{\bar{a},\psi^2}(\bar{x}) \dd \bar{x}  
  =  \int_{B_{M,\alpha,\beta}} g^k_{\bar{a},\psi^2}\big(\Phi(\bar{x})\big) \dd \bar{x} \\
  \\
 & \geq & \kappa_0 \int_{B_{M,\alpha,\beta}} g^k_{\bar{a},\psi^2}\big(\bar{x}\big) \dd \bar{x} = 
\kappa_0 \cdot \nu^k_{\bar{a},\psi^2} \big( B_{M,\alpha,\beta}  \big) \, ,
\end{eqnarray*}
where $\kappa_0$ denotes the infimum over $\bar{x} \in B_{M,\alpha,\beta}$ of 
 $\frac{g^k_{\bar{a},\psi^2}\big(\Phi(\bar{x})\big)}{g^k_{\bar{a},\psi^2}\big(\bar{x}\big)}$. For such $\bar{x}$, this last ratio equals
\begin{eqnarray*}
 & & \prod_{j=1}^k \exp \Big\{ - \tfrac{1}{2} \psi^{-2} \big( x_j - a_j  + 2\alpha (k-j) + \beta \big)^2 \Big\} \cdot  \exp \Big\{  \tfrac{1}{2} \psi^{-2} \big( x_j - a_j \big)^2 \Big\} \\
 & = &
 \prod_{j=1}^k \, \exp \Big\{ - \psi^{-2} \Big( ( x_j - a_j ) \cdot \big( 2\alpha(k-j) + \beta  \big) + \tfrac{1}{2} \big( 2 \alpha(k-j) + \beta \big)^2 \Big) \Big\} \, .
\end{eqnarray*}
The set  $B_{M,\alpha,\beta}$ is comprised of decreasing $k$-vectors whose first component is at most $M$. Thus, $x_j \leq M$ for $j \in \intint{k}$.
Also using $a_j \geq -M$, we find that $\kappa_0$ is at least $\kappa$, as we sought to show.
We have verified~(\ref{e.jacobianphi}) and thus completed the proof of Lemma~\ref{l.bmalpha}. \qed

\medskip

To prove Proposition~\ref{p.jabovecorner}, recall that $\barqpr = \barqmin$ and note that
\begin{eqnarray*}
& & 
 \PP_{\mc{F}^1} \Big( \overstrongjump(\mfl) \in \pspr \Big) \cdot {\bf 1}_{\mpg \cap \mcgone} 
  =  \nu^k_{\bar{m},\sigma^2} \Big( \pspr \, \Big\vert \, \barqpr + [0, \infty)^k_> \Big)  \cdot {\bf 1}_{\mpg \cap \mcgone}  \\
 & \geq &  \nu^k_{\bar{m},\sigma^2} \Big( (-\infty,M]^k \cap  \pspr \, \Big\vert \, \barqpr + [0, \infty)^k_> \Big)  \cdot {\bf 1}_{\mpg \cap \mcgone} \\
 & \geq &  \bigg( 1 - \nu^k_{\bar{m},\sigma^2} \Big( (-\infty,M]^k \cap \pspr^c \, \Big\vert \, \barqpr + [0, \infty)^k_> \Big) \\
  & & \qquad \qquad \qquad \qquad \qquad \qquad - \, \, \, \,  \nu^k_{\bar{m},\sigma^2} \Big( (M,\infty) \times \R^{k-1} \, \Big\vert \, \barqpr + [0, \infty)^k_> \Big) \bigg) \,  \cdot \, {\bf 1}_{\mpg \cap \mcgone} \, ;
\end{eqnarray*}
the latter inequality depends on $\barqmin + [0, \infty)^k_>  \subset \R^k_>$.

We now apply Lemma~\ref{l.bmalpha}
with the choices $\psi = \sigma$, $M = \big( 38 k^{3/2} + 144k \big) \eln^2$, $\alpha = \eln^{1/2}$, $\beta = 8\ipd \eln + \eln^{1/2}$ and $\bar{q} = \barqmin$, recalling that $\sigma^{-2} \leq 2 \eln^{-1}$ (and that $\ipd \geq 1$, $k \geq 2$ and $\eln \geq 1$); and also using 
the bound $M \leq 140 k^{3/2} \eln^2$ (which is valid since $k \geq 2$).
Doing so while also making use of the next 
Lemma~\ref{l.eccubed}, and $k \geq 2$, we find that the bracketed expression in the final double-line of the last display is at least 
$$
  \tfrac{1}{2} \exp \big\{ -  3970 \ipd^2 k^{7/2}  \eln^2  \big\} \, - \, \e^{1054 k^2 \const^3} \, .
$$
Of the two terms in this difference, the second has value at most one-half of the first, due to $\e < e^{-1}$, $\const \geq 1$, $k \geq 2$ and $\const \big( \log \e^{-1} \big)^{1/3} \geq \tfrac{3971}{1054} \ipd^2 k^{3/2}$. The proof of Proposition~\ref{p.jabovecorner} is complete. \qed

\begin{lemma}\label{l.eccubed}
When the event $\mpg \cap \mcgone$ occurs, 
$$
 \nu^k_{\bar{m},\sigma^2} \Big( (M,\infty) \times \R^{k-1} \, \Big\vert \, \barqmin + [0,\infty)^k_> \Big) \leq \e^{1054 k^2 \const^3} \, ,
$$
provided that $M \in \R$ verifies $M - m(1) \geq (12)^2 k \eln^2$.
\end{lemma}
\noindent{\bf Proof.} We will use the next result.
\begin{lemma}\label{l.elnkcubed}
On the event $\mpg \cap \mcgone$, 
$$
 \nu^k_{\bar{m},\sigma^2} \Big( \barqmin + [0,\infty)^k_> \Big) \geq  \e^{12770 k^2 \const^3} \, .
$$
\end{lemma}
\noindent{\bf Proof.}
Set $B \subset \R^k$ equal to the box $[k-1,k] \times [k-2,k-1] \times \cdots \times [0,1]$. Note that
\begin{eqnarray*}
& &   \nu^k_{\bar{m},\sigma^2} \Big( \barqmin + [0,\infty)^k_> \Big) \geq 
 \nu^k_{\bar{m},\sigma^2} \Big( \barqmin + B \Big) \\
 & = & \prod_{j=1}^k \nu_{m(j),\sigma^2} \big( \qmin_j + [j-1,j] \big) =  \prod_{j=1}^k \int_{\qmin_j + j-1}^{\qmin_j + j} g_{m(j),\sigma^2}(x) \dd x  \\
 & \geq & \big( 2 \pi \sigma^2 \big)^{-k/2} \prod_{j=1}^k \exp \Big\{  - \tfrac{1}{2} \sigma^{-2} \big( \vert \qmin_j - m(j) \vert + j \big)^2 \Big\} \\
 & \geq & \big( 3\pi\eln/2 \big)^{-k/2} \cdot \exp \big\{ - (113)^2 k^5 \eln^3 \big\}  \geq  \e^k \cdot \e^{(113)^2  k^5 \const^3}  \geq \e^{12770 k^2 \const^3} \, ,
\end{eqnarray*}
since $\eln/2 \leq \sigma^2 \leq 3\eln/4$ and $\vert m(j) \vert \leq 38 k^{3/2} \eln^2$ by Lemma~\ref{l.misigmai} and $\vert \qmin_j \vert \leq 74 k^{5/2} \eln^2$ by Lemma~\ref{l.estgone}. The penultimate inequality uses $\e < \big( \tfrac{2}{3 \pi \const} \big)^{1/3}$. \qed

\medskip

Note also that
\begin{eqnarray}
 & & \nu_{m(1),\sigma^2}(M,\infty) \leq (2\pi)^{-1/2} \cdot \tfrac{\sigma}{M - m(1)} \cdot \exp \Big\{ - \tfrac{1}{2} \sigma^{-2} \big(M - m(1) \big)^2 \Big\} \nonumber \\
 & \leq & \tfrac{1}{2} 3^{1/2}(2\pi)^{-1/2} (12)^{-2} k^{-1} \eln^{-3/2} \exp \big\{ - \tfrac{2}{3} \, (12)^4 k^2 \eln^3 \big\} \leq \e^{13824 k^2 \const^3} \, , \label{e.eccubedub}
 \end{eqnarray}
where we used $\eln/2 \leq \sigma^2 \leq 3\eln/4$ and $M - m(1) \geq (12)^2 k \eln^2$. 

Since 
$$
\nu^k_{\bar{m},\sigma^2} \Big( (M,\infty) \times \R^{k-1} \, \Big\vert \,  \barqmin + [0,\infty)^k_>   \Big) \leq \frac{
\nu^1_{m(1),\sigma^2} \big( M,\infty \big)}{
\nu^k_{\bar{m},\sigma^2} \big(  \barqmin + [0,\infty)^k_>   \big)} \, ,
$$
we find from Lemma~\ref{l.elnkcubed} and~(\ref{e.eccubedub}) that
$$
\nu^k_{\bar{m},\sigma^2} \Big( (M,\infty) \times \R^{k-1} \, \Big\vert \, \barqmin + [0,\infty)^k_>   \Big)  \leq  \e^{-12770 k^2 \const^3}  \cdot   \e^{13824 k^2 \const^3} \, .
$$
 This completes the proof of Lemma~\ref{l.eccubed}. \qed

\medskip

We may now complete the proof of Proposition~\ref{p.vaultsuccess} by applying the guiding equality~(\ref{e.testguide}). Setting $q_k = \exp \big\{ - 2D_k \big( \log \e^{-1} \big)^{1/3} k \log k \big\}$
equal to the lower bound on  $\PP_{\mc{F}} \big( \test_3(J) = 1 \big)$ provided by Lemma~\ref{l.jumppromotion},
note that
\begin{eqnarray*}
 & & \PP_{\mc{F}} \big( \test_4(J) = 1 \big) \cdot {{\bf 1}}_\mpg  \,  = \,  \PP_{\mc{F}} \big( \test_4(\strongjump) = 1 \big) \PP_{\mc{F}} \big( \test_3(J) = 1 \big) \cdot {{\bf 1}}_\mpg  \\ 
& \geq &  \EE_{\mc{F}} \Big[ \PP_{\mc{F}^1} \big( \test_4(\strongjump) = 1  \big) {\bf 1}_{\mcgone} \Big]  q_k \cdot {{\bf 1}}_{\mpg} \\ 
& \geq & \EE_{\mc{F}} \bigg[ \PP_{\mc{F}^1} \Big( \test_4(\strongjump) = 1 \, , \,  \overstrongjump(\mfl) \in \pairsep_{\eln^{1/2}, 8 \eln + \eln^{1/2},\barqmin}  \Big)  {\bf 1}_{\mcgone} \bigg]  q_k \cdot {{\bf 1}}_{\mpg} \\
 & \geq &  \tfrac{1}{4} \exp \Big\{ -  3970  k^{7/2}\ipd^2 \const^2 \big( \log \e^{-1} \big)^{2/3} 
 \Big\}  \cdot \big( 1 - 2e^{-1} \big)^{2k \const (\log \e^{-1})^{1/3}}  q_k \cdot \PP_{\mc{F}}\big( \mcgone \big) \cdot {{\bf 1}}_{\mpg} \\
   & \geq &   \exp \Big\{ -  3973  k^{7/2}\ipd^2 \const^2 \big( \log \e^{-1} \big)^{2/3} \Big\}   \cdot {{\bf 1}}_{\mpg} \, ,
\end{eqnarray*}
where we used Propositions~\ref{p.jpass} and \ref{p.jabovecorner} in the third inequality.
In the fourth inequality, we made use of  $\PP_{\mc{F}} \big(  \mcgone \big) \geq 2^{-1} \cdot {{\bf 1}}_\mpg$, which follows from Lemma~\ref{l.gprime} and $\e \leq 1/2$; as well as  $k \geq 2$, $\ipd \geq 1$, $\const \geq 1$ and   $\e < e^{-1}$.
Thus is the proof of Proposition~\ref{p.vaultsuccess} completed. 
\qed



\chapter{The jump ensemble method: applications}\label{c.jumpensembleapplications}

In this chapter, our new method, and related resampling ideas, are applied to prove the more subtle of the technical results in the article. There are three sections.

Section~\ref{s.closenessupper} is devoted to the proof of Theorem~\ref{t.airynt}(1), concerning the probability of $k$-curve closeness at a given point, via the jump ensemble method and other reconstruction techniques (the {\em snap up} and {\em swing through} arguments). 

Theorem~\ref{t.airynt}(2), concerning $k$-curve closeness at a general location, is a consequence of the one-point version, Theorem~\ref{t.airynt}(1), and the locally Brownian nature of ensemble curves. This Brownian nature, in essence Theorem~\ref{t.weakbound}, is established in Section~\ref{s.closenessgeneral}, where the proof of Theorem~\ref{t.airynt}(2)
is given.
 
The jump ensemble method is employed again in Section~\ref{s.brownianregularity} to obtain the principal Brownian regularity results: the Radon-Nikodym moment bound Theorem~\ref{t.rnbound} for the Airy line ensemble; Theorem~\ref{t.airytail}, whose first part reformulates Theorem~\ref{t.rnbound} in terms of the deformation of the probability of unlikely events; and Theorem~\ref{t.airytail.ln}, which is the finite-$n$, regular ensemble, counterpart to Theorem~\ref{t.airytail}(1).

\section{Upper bound on the probability of curve closeness over a given point}\label{s.closenessupper}

This section is devoted to the proof of Theorem~\ref{t.airynt}(1). 
This result concerns the $k$-curve closeness probability over a given point in the case that this point is permitted to lie in a rather long interval. We begin by reducing via Lemma~\ref{l.parabolicinv} to a counterpart result in which this interval is much shorter.

For this section  and the next, we let the value of $\const$ be specified to be 
\begin{equation}\label{e.constvalue.new}
 \const  = \max \Big\{ k^{1/3} c_k^{-1/3} \big( 2^{-9/2} - 2^{-5} \big)^{-1/3} \, , \,  36(k^2 -1) \, , \, 2^{5/3} c^{-1/3}_k (k^2 -1)^{1/3} \Big\} \, ,
\end{equation}
so that  the value of $\const$ is increased from the expression in~(\ref{e.constvalue}) if necessary so that it is at least  $2^{5/3} c^{-1/3}_k (k^2 -1)^{1/3}$. The new condition ensures that Lemma~\ref{l.glub.new} implies that $\PP \big( \mpg^c \big) \leq \e^{k^2 - 1}$.

\begin{theorem}\label{t.airynt.one}
For $\bar\phimac \in (0,\infty)^3$, $C,c > 0$ and $n \in \N$, let 
$$
\mc{L}_n:\intint{n} \times \big[-\xnmac,\infty\big) \to \R  
$$ 
be a    $\big(\bar\phimac,\rsc,\rsC\big)$-regular ensemble defined under the law~$\PP$. 
Let $k \in \N$ satisfy $k \geq 2$ and let $\e > 0$ satisfy~(\ref{e.epsilonupperbound}).
For $n \in \N$ satisfying $n \geq k 
  \vee  (c/3)^{-2(\phimac_1 \wedge \phimac_2)^{-1}} \vee  6^{2/\delta}$ and~(\ref{e.nlbone}), 
 the bound
$$
  \PP \Big( \close\big( k ; \mc{L}_n , \fa , \e \big)  \Big) 
  \leq     10^6 \exp \Big\{  4962 k^{7/2}  \const^{5/2} \big( \log \e^{-1} \big)^{5/6} \Big\}  \, \e^{k^2-1} 
$$ 
holds 
for any given  $\fa \in \R$ for which $\vert \fa \vert \leq \const \big( \log \e^{-1} \big)^{1/3}/2$.
\end{theorem}
\noindent{\bf Proof of Theorem~\ref{t.airynt}(1).}
We simply consider the choice $\fa = 0$ in Theorem~\ref{t.airynt.one} and apply the parabolic invariance Lemma~\ref{l.parabolicinv}. The value of the parameter $c$ drops by a factor of two, which in view of (\ref{e.constvalue}) entails the replacement of $\const^{5/2}$ by $2^{5/6} \const^{5/2}$.    \qed

\subsection{Upper bound on curve closeness high above the tent map}\label{s.abovepar}
It is our aim to use the jump ensemble method to find the upper bound on $k$-curve closeness in $\mc{L}_n$ stated in Theorem~\ref{t.airynt.one}. This involves finding an upper bound on such an occurrence in the jump ensemble $J$, something that is easier to do if   the vertical coordinate of the location at which the $k$ curves gather is comfortably above the lower boundary condition (by which we mean at a significant  distance above the tent map). In this subsection, we carry out the method in this circumstance, with the order of this distance being $(\log \e^{-1})^{1/2}$. In the next, we use different arguments to treat the probability of curve closeness at lower heights, and thus obtain Theorem~\ref{t.airynt.one}.
 
 The real-valued parameter $x_0$ is given by Theorem~\ref{t.airynt.one} throughout the upcoming proof. Recalling from Section~\ref{s.elnfix} that $\eln$ has been fixed to be  $\const \big( \log \e^{-1} \big)^{1/3}$, $\vert x_0 \vert$ is thus supposed to be at most $\eln/2$. 
 
 \begin{definition}\label{d.abovepar}
Let  $X$ denote one of the ensembles $J$ or $\mc{L}_n$. Define
$$
\abovepar\big(X,\fa\big) =  \abovepar\Big(X,\fa,15 k \const \big(\log \e^{-1}\big)^{1/2}\Big)
$$
to be the event that
$$ 
 X(i,\fa) - \tent(\fa) \geq 15 k \const (\log \e^{-1})^{1/2} \, \, \, \, \textrm{for $i \in \intint{k}$} \, .
$$
\end{definition}
  
\begin{proposition}\label{p.rare}
We have that
$$
\PP_\mc{F} \Big( \, 
 \test_4 ( J ) = 1   \, , \, \abovepar\big(J,\fa\big) \, , \, 
  \close(J,\fa,\e)  \,  \Big) \cdot {\bf{1}}_{\mpg}  
  \leq   \big( 10^4 \const^{k^2} + 2k \big) \e^{k^2-1} \big(\log \e^{-1}\big)^{k^2/2} \, .
$$
\end{proposition}
(Here, and similarly later, by $\close(J,\fa,\e)$ is understood  $\close(k;J,\fa,\e)$: recall the comment that ends Definition~\ref{d.closeenc}.)

Alongside Proposition~\ref{p.vaultsuccess}, 
Proposition~\ref{p.rare} provides the ingredient demanded when we seek to prove an upper bound on $k$-curve closeness high above the lower boundary condition using the jump ensemble as a candidate for $\mc{L}_n$. Indeed, we may now promptly infer the following on the basis of these two inputs.
\begin{proposition}\label{p.labovepar}
\begin{eqnarray*}
  & & 
\PP \Big( \, 
  \close\big(\mc{L}_n,\fa,\e\big) \cap    \abovepar\big(\mc{L}_n,\fa \big)   \, \Big) \\
   & \leq & 
    2  \big( 10^4 \const^{k^2} + 2k \big) \exp \Big\{  3973 k^{7/2} \ipd^2 \const^2 \big( \log \e^{-1} \big)^{2/3} \Big\}  \big(\log \e^{-1}\big)^{k^2/2}   \e^{k^2-1}  \, .
\end{eqnarray*}
\end{proposition}
\noindent{\bf Proof.} Under $\PP_\mc{F}$, the conditional law of $\mc{L}_n:\intint{k} \times [\mfl,\mfr] \to \R$ coincides with the conditional law of $J:\intint{k} \times [\mfl,\mfr] \to \R$ given that $\test_4(J) = 1$. For this reason, it follows from Propositions~\ref{p.vaultsuccess} and~\ref{p.rare} that
\begin{eqnarray*}
 & & \PP_\mc{F} \Big( \, 
   \abovepar\big(\mc{L}_n,\fa \big) \cap   
  \close\big(\mc{L}_n,\fa,\e\big) \, \Big) \cdot {\bf{1}}_{\mpg}  \\
 & \leq &    \big( 10^4 \const^{k^2} + 2k \big) \exp \Big\{  3973 k^{7/2} \ipd^2 \const^2 \big( \log \e^{-1} \big)^{2/3} \Big\}  \big(\log \e^{-1}\big)^{k^2/2}   \e^{k^2-1}
   \, .
\end{eqnarray*}

We may now write 
$$
    \PP \Big( \, 
  \close\big(\mc{L}_n,\fa,\e\big) \cap    \abovepar\big(\mc{L}_n,\fa\big))  \, \Big) \,
    \leq  \, \EE \, \bigg[ \, \PP_{\mc{F}}  \Big( \, 
  \close\big(\mc{L}_n,\fa,\e\big) \cap    \abovepar\big(\mc{L}_n,\fa\big))  \, \Big) {\bf 1}_{\mpg} \, \bigg] \, + \, \PP\big(\mpg^c\big) 
$$
and apply the last inequality and Lemma~\ref{l.glub.new} to deduce Proposition~\ref{p.labovepar}. \qed

\medskip

The remainder of Section~\ref{s.abovepar} is devoted to the next proof. 

\medskip

\noindent{\bf Proof of Proposition~\ref{p.rare}.}  
Set
$$
\smalljfluc = \Big\{ 
J(i,x) \geq \tent(\fa) + (\log \e^{-1})^{1/2} \, \, \, \textrm{for $(i,x) \in \intint{k} \times \{ \fa - 1, \fa + 1 \}$}  \Big\} \, .
$$
Why is this event so named? We will use it alongside the event $\abovepar\big(J,\fa \big)$ that $J$-curves at $\fa$ exceed the tent map by $15 k \const (\log \e^{-1})^{1/2}$. When the latter event occurs, the new event entails that the curves in $J$ not fluctuate too much between $\fa$ and the neighbouring times~$\fa - 1$ and~$\fa+1$.

 Note that 
\begin{eqnarray*}
& & \PP_\mc{F} \Big( \, 
 \test_4 ( J ) = 1   \, , \, \abovepar\big(J,\fa\big) \, , \, 
  \close(J,\fa,\e) \, , \, 
\smalljfluc  \Big) \cdot {\bf{1}}_{\mpg}  \\
& \leq & \sup \, \PP_\mc{F} \Big( \, 
 \notouch^{[\fa - 1,\fa+1]} ( J )   \, , \,
  \close(J,\fa,\e) \, \Big\vert \, 
 \ovbar{J}(\fa - 1) = \bar{x} \, , \, \ovbar{J}(\fa +1 ) = \bar{y} \Big) \, ,
\end{eqnarray*}
where the supremum is taken over all choices of $\bar{x}$ and $\bar{y}$ in the set $\big(  \tent(\fa) + (\log \e^{-1})^{1/2} , \infty \big)^k$.
(Conditioning the ensemble $J$ to assume given values is a singular conditioning which nonetheless has an unambiguous meaning because, under $\PP_\mc{F}$, this ensemble is given as the marginal on $[\mfl,\mfr]$ of a Brownian bridge ensemble conditioned on an event of positive probability.)
 For any given such choice of $\bar{x}$ and $\bar{y}$, we have that  
\begin{eqnarray}
& &  \PP_\mc{F} \Big( \, 
 \notouch^{[\fa - 1,\fa+1]} ( J )   \, , \,
  \close(J,\fa,\e) \, \Big\vert \, 
 \ovbar{J}(\fa - 1) = \bar{x} \, , \, \ovbar{J}(\fa +1 ) = \bar{y} \Big) \cdot {{\bf 1}}_{\mpg} \nonumber \\
  & \leq & \frac{\mc{B}_{k;\bar{x},\bar{y}}^{[\fa - 1,\fa+1]}  \Big( \, 
 \notouch^{[\fa - 1,\fa+1]} ( B )   \, , \,
  \close(B,\fa,\e) \, \Big)}{\mc{B}_{k;\bar{x},\bar{y}}^{[\fa - 1,\fa+1]}  \Big(  B(i,\pp) \geq \mc{L}_n(k+1,\pp) \, \, \, \textrm{for $i \in \intint{k}$ and $\pp \in \pole \cap [\fa - 1, \fa + 1]$} \, \Big)}  \cdot {{\bf 1}}_{\mpg} \nonumber \\
  & \leq &  2 \, \mc{B}_{k;\bar{x},\bar{y}}^{[\fa - 1,\fa+1]}  \Big( \, 
 \notouch^{[\fa - 1,\fa+1]} ( B )   \, , \,
  \close(B,\fa,\e) \, \Big) \, . \label{e.denomlb}
\end{eqnarray}
The inequality $\geq$ in the denominator in the middle line is permitted to be non-strict because Brownian bridge has zero probability of hitting any given point.
That this denominator is least one-half whenever $\mpg$ occurs should be justified.
By Lemma~\ref{l.monotonetwo}, this denominator is minimized by taking both of the vectors $\bar{x}$ and $\bar{y}$ equal to the constant vector (that we denote by $\bar{z}_0$) whose components equal  $\tent(\fa) + (\log \e^{-1})^{1/2}$; since  $\mc{L}_n(k+1,\pp) = \tent(\pp) \leq \tent(\fa) + \numcthree \eln$ for $\pp \in \pole \cap [\fa - 1, \fa + 1]$ on the event $\mpg$ (whose occurrence merely ensures that $\fa$, being at most $\eln/2$ in absolute value, is an element of the tent map domain $[\mfl,\mfr]$), we see by Lemma~\ref{l.maxfluc} that the denominator is at least
\begin{eqnarray*}
& &
 \mc{B}_{k;\bar{z}_0,\bar{z}_0}^{[\fa - 1,\fa + 1]} \Big( \inf_{(i,x) \in \intint{k} \times [\fa - 1, \fa +1]} B(i,x) \geq \tent(\fa)  + \numcthree \eln \Big) \\
  & \geq & 1 \, - \, k \exp \bigg\{ - \Big( \big(\log \e^{-1}\big)^{1/2} - \numcthree \const  \big(\log \e^{-1}\big)^{1/3}  \Big)^{2} \bigg\}
   \geq 1/2 \,  
\end{eqnarray*}
 where we used that $\log \e^{-1} \geq (8 \const)^6 \vee 4 \log (2k)$. The inequality in line~(\ref{e.denomlb}) has been justified.

Noting that $\const \geq 36(k^2 -1)$ and $\e < k^{-3}/6$, we may apply Lemma~\ref{l.browniannotouch} with $\phi = \e$ to the probability in~(\ref{e.denomlb}), finding that
\begin{eqnarray}
& & 
 \PP_\mc{F} \Big( \, 
 \test_4 ( J ) = 1   \, , \, \abovepar\big(J,\fa\big) \, , \, 
  \close(J,\fa,\e) \, , \, 
\smalljfluc  \Big) \cdot  {\bf{1}}_{\mpg} \nonumber \\ 
 & \leq &  10^4 \, \const^{k^2} \e^{k^2-1} \big(\log \e^{-1}\big)^{k^2/2} \, . \label{e.smalljfluc}
\end{eqnarray}

Set $R = 5k \const$.
We will now argue that
\begin{equation}\label{e.nosmalljfluc}
 \PP_\mc{F} \Big( \, \neg \, \smalljfluc \, , \, \abovepar\big(J,\fa \big) \Big) \cdot {\bf{1}}_{\mpg} \leq 2k \e^{R^2/2} \, .
\end{equation}

To this end, let $i \in \intint{k}$. Recall that $\mpg$ entails that
$\mc{L}_n \big( i,-2\eln \big) \geq - \big( 2\sqrt{2} + 1 \big) \eln^2$. If this event occurs, we thus see in light of Lemmas~\ref{l.monotoneone} and~\ref{l.monotonetwo} that, under $\PP_\mc{F}$ given that $J(i,\fa)$ equals $h \in \R$, the conditional distribution of $J(i,\cdot)$ on $[\mfl,\fa]$ stochastically dominates the marginal of $\mc{B}_{1,-(2\sqrt{2} +1)\eln^2,h}$ on $[\mfl,\fa]$. 
     Recalling Definition~\ref{d.abovepar} and that $R = 5k \const$, we thus find that
\begin{eqnarray*}
& & 
 \PP_\mc{F} \Big( \, J(i,\fa - 1) <  \tent(\fa) + (\log \e^{-1})^{1/2} \, , \, \abovepar\big(J,\fa\big) \Big) \cdot {\bf{1}}_{\mpg} \\ 
 & \leq & \sup_{r \geq 0} \mc{B}_{1;-(2\sqrt{2} + 1)\eln^2,\tent(\fa) + 3R  (\log \e^{-1})^{1/2} + r}^{[-2\eln,\fa]}
 \Big( B(1,\fa - 1) <  \tent(\fa) +   (\log \e^{-1})^{1/2} \Big) \, .
\end{eqnarray*}
The right-hand side equals
$$ 
  \mc{B}_{1;-(2\sqrt{2} + 1)\eln^2,\tent(\fa) + 3R  (\log \e^{-1})^{1/2}}^{[-2\eln,\fa]}
 \Big( B(1,\fa - 1) <  \tent(\fa) +   (\log \e^{-1})^{1/2} \Big)
$$
by another use of Lemma~\ref{l.monotonetwo}. 
This last expression is increasing when regarded as a function of the variable $\tent(\fa)$. Since the occurrence of $\mpg$
 entails that $\mc{L}_n\big(k+1,x\big)$, and thus also $\tent(x)$, is at most $\eln^2$ whenever $x \in [\mfl,\mfr]$, we see that
\begin{eqnarray}
& & 
 \PP_\mc{F} \Big( \, J(i,\fa - 1) <  \tent(\fa) + (\log \e^{-1})^{1/2} \, , \, \abovepar\big(J,\fa\big) \Big) \cdot {\bf{1}}_{\mpg} \nonumber \\ 
 & \leq & 
 \mc{B}_{1;-(2\sqrt{2} + 2)\eln^2,3R  (\log \e^{-1})^{1/2}}^{[-2\eln,\fa]}
 \Big( B(1,\fa - 1) <  (\log \e^{-1})^{1/2} \Big) \, . \label{e.bonefa}
 \end{eqnarray}
In the last line, the random variable $B(1,\fa-1)$  is normally distributed. Its mean is less than $3R  (\log \e^{-1})^{1/2}$ by 
$$
 \frac{1}{\fa + 2\eln} \Big( 3 R  (\log \e^{-1})^{1/2}  + 2(\sqrt{2} + 1)\const^2 (\log \e^{-1})^{2/3} \Big)
$$ 
which quantity is, in view of the bounds $\fa \geq - \eln$, $\e < e^{-1}$, $\const \geq 6$ and $R \geq 4(\sqrt{2} + 1)\const$, at most $R  (\log \e^{-1})^{1/2}$. The mean is thus at least $2R  (\log \e^{-1})^{1/2}$. The event in~(\ref{e.bonefa}) entails that $B(1,\fa)$ drops below its mean by at least $(2R - 1)(\log \e^{-1})^{1/2}$. Since $R \geq 1$, we find that the expression~(\ref{e.bonefa}) is at most
\begin{eqnarray*}
& &  \mc{B}_{1;0,0}^{[-2\eln,\fa]}
 \Big( B(1,\fa - 1) <  -R(\log \e^{-1})^{1/2} \Big) \\
  & \leq & \nu_{0,1} \big( R(\log \e^{-1})^{1/2} , \infty \big) \leq (2\pi)^{-1/2} R^{-1} \big( \log \e^{-1} \big)^{-1/2} \exp \big\{ - \tfrac{1}{2} R^2 \log \e^{-1} \big\} \leq \e^{R^2/2} \, ,
\end{eqnarray*}
where lastly we used $R (\log \e^{-1})^{1/2} \geq 1$.
We now sum over $i \in \intint{k}$ to learn that, when $\mpg$ occurs, the probability
$$
 \PP_\mc{F} \bigg( \, \Big\{ \, \exists \, i \in \intint{k}: J(i,\fa- 1) < \tent(\fa) -  (\log \e^{-1})^{1/2} \Big\} \, \cap \,  \abovepar\big(J,\fa\big)  \, \bigg)
$$
is at most $k \e^{R^2/2}$. 
A similar bound holds with $\fa + 1$ in place of $\fa - 1$ and thus we obtain~(\ref{e.nosmalljfluc}).

Combining~(\ref{e.smalljfluc}) and~(\ref{e.nosmalljfluc}),
\begin{eqnarray*}
 &  & \PP_\mc{F} \Big( \, 
 \test_4 ( J ) = 1   \, , \, \abovepar\big(J,\fa\big) \, , \, 
  \close(J,\fa,\e)  \,  \Big) \cdot {\bf{1}}_{\mpg} \\
  & \leq &  10^4 \const^{k^2} \e^{k^2-1} \big(\log \e^{-1}\big)^{k^2/2} + 2k \e^{R^2/2} \, .
\end{eqnarray*}
Noting that $R^2/2 \geq k^2-1$, we complete the proof of Proposition~\ref{p.rare}. \qed

\subsection{Curve closeness at low height: the  snap up and the swing through}\label{s.lowtohigh}

\begin{figure}[ht]
\begin{center}
\includegraphics[height=12cm]{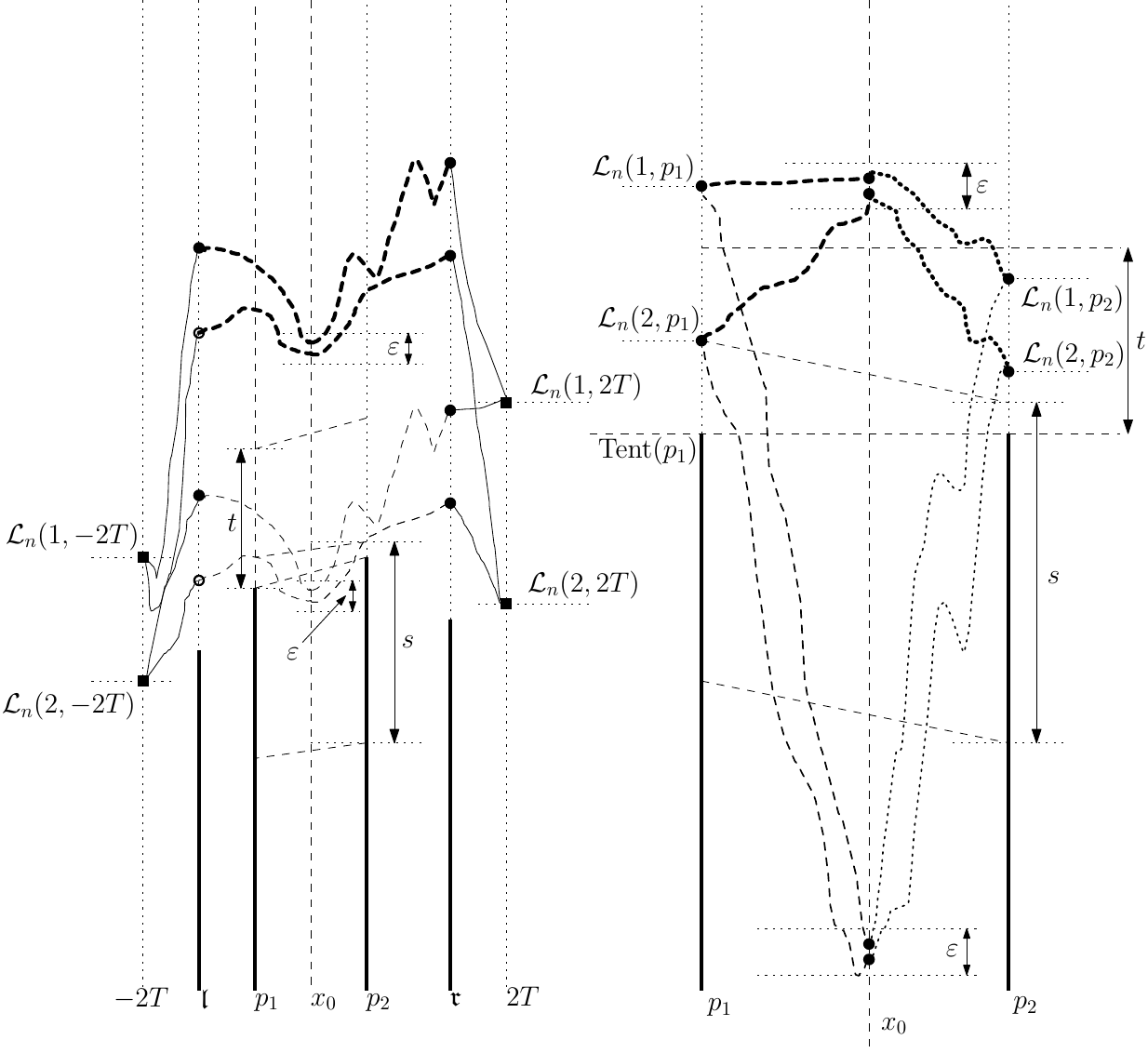}
\caption{The snap up on the left and the swing through on the right, with $k$ equal to two. Vertical poles are depicted in thick solid lines. In each sketch, $s = \big( 2(k-1) + 4\eln \big)(2\eln)^{1/2}$ is the big drop distance and $t = 15 k \const \big( \log \e^{-1} \big)^{1/2}$ is the smaller above-tent margin specifying the event~$\abovepar(\mc{L}_n,\fa)$. {\em Left:} The snap up is used to prove Lemma~\ref{l.apgt}. The clear circular bead over $\mfl$ represents all that is random to the witness of $\mc{G}$. As it rolls up or down, it dictates a joint displacement of the top two curves on $[\mfl,\mfr]$. 
An original configuration of these curves (dashed, thin), in which the governing bead is lower, is snapped up by pushing the bead higher to obtain a new configuration (dashed, thick) high above the tent map.
{\em Right:} The upcoming proof of Lemma~\ref{l.lowtohigh} is illustrated. The point $\fa$ lies in $[p_1,p_2]$ with $p_1,p_2 \in P$ consecutive pole set elements.  The original configuration of the top two curves on $[p_1,p_2]$ (dashed-dotted, thin) is $\e$-close at $\fa$ and has a big drop. The kinetic energy bound up in these curves' displacements during $[p_1,\fa]$ and $[\fa,p_2]$ is unleased in the swung through configuration (dashed-dotted, thick) that as a result reaches high above the tent map at $\fa$ while maintaining $\e$-closeness at $\fa$.}
\label{f.snapupswingthrough}
\end{center}
\end{figure}

   To obtain Theorem~\ref{t.airynt.one}, to which we have reduced Theorem~\ref{t.airynt}(1), it  remains to address the probability of $k$-curveness closeness at lower locations than those with which the event $\abovepar\big(\mc{L}_n,\fa\big)$ is concerned. This will be the subject of the present subsection.
   We want a version of Proposition~\ref{p.labovepar}
without the presence of the $\abovepar$ event. We will use a {\em snap up} to prove this new version: if the top $k$ curves in $\mc{L}_n$ at $\fa$ come close, they do so at a random height that in this paragraph of overview we will denote simply by  $h$; (the value of $h$ is essentially the same for each of the $k$ curves, because they are close). If $h$ is low relative to $\tent(\fa)$,  we will argue that this system of $k$ curves in a locale of $\fa$ can collectively be snapped up; that is, all these curve pieces can be raised by a quite large common amount. By `quite large', we mean at least the sum of two quantities: the distance below the tent map $\tent(\fa) - h$, and the above tent margin $15k \const (\log \e^{-1})^{1/2}$ that appears in the $\abovepar$ event Definition~\ref{d.abovepar}.  This operation thus forces the occurrence of $\abovepar\big(\mc{L}_n,\fa\big)$, while preserving the occurrence of the $k$-curve closeness event $\close\big(\mc{L}_n,\fa,\e \big)$. We will seek to argue that the snapped up configuration is roughly as probable as the original: in this way, we will be able to invoke Proposition~\ref{p.labovepar} to strengthen this proposition so that the presence of the $\abovepar$ event is dropped from its statement (at the expense of a manageable increase in the right-hand side). 
We will succeed in arguing that the new configuration has a comparable probability to the old one only in the case that 
the below tent distance $\tent(\fa) - h$ is not too large,
which is to say, at  most $O\big(\eln^{3/2}\big) = O\big( (\log \e^{-1})^{1/2} \big)$. We will specify an event $\bd$ that the below tent distance exceeds this order. When it fails to occur, we will use the snap up; when it does, we will apply a different technique, the {\em swing through}, to argue that the occurrence of $\abovepar$ may be forced with comparable probability. The main concepts in the two techniques are illustrated in Figure~\ref{f.snapupswingthrough}.

The snap and swing arguments make no use of the jump ensemble, though they do to a degree build on the more basic apparatus of the missing closed middle $\sigma$-algebra $\mc{F}$.  
The natural perspective of the snap up argument is the viewpoint of the witness of a new $\sigma$-algebra~$\mc{G}$ which we begin by specifying. It will contain $\mc{F}$. For the new witness, the top $k$ curves in $\mc{L}_n$ near $\fa$ will be in a determined configuration modulo a height shift. The only randomness for this witness will be a one-dimensional random variable that determines this height. 

Recall that $\mc{F}$ is generated by the ensemble $\mc{L}_n$
outside of the index set $\intint{k} \times \big(-2\eln,2\eln \big)$ alongside the side interval standard bridges $I \to \R: x \to \mc{L}^I_n(i,x)$, $i \in \intint{k}$, with $I$ equal to $[-2\eln,\mfl]$ and $[\mfr,2\eln]$. The information absent from $\mc{F}$ may be gathered in the following form:


\begin{itemize}
\item the shifted curves $[\mfl,\mfr] \to \R: x \to \mc{L}_n(i,x) - \mc{L}_n(i,\mfl)$, indexed by $i \in \intint{k}$;
\item the differences $\mc{L}_n(i,\mfl) - \mc{L}_n(i+1,\mfl)$ for $i \in \intint{k-1}$;
\item and the value $\mc{L}_n(k,\mfl)$.
\end{itemize}


In this section, we denote by $\mc{G}$ the $\sigma$-algebra generated by $\mc{F}$ and the first two items in this list.
To the witness of~$\mc{G}$, whose perspective is represented by the conditional law $\PP_\mc{G}$, only the $\R$-valued quantity $\mc{L}_n\big(k,\mfl\big)$ is random. 
If this random quantity adopts the value $u \in \R$, the ensemble may be reconstructed in its entirety by the witness.
Denoting the reconstructed ensemble by $\mc{L}_n^u$, we have that for $(i,x) \in \intint{k} \times [\mfl,\mfr]$,
$$
\mc{L}_n^u(i,x) \, = \, u \, + \, \Big( \mc{L}_n\big(i,\mfl\big) - \mc{L}_n\big(k,\mfl\big) \Big) \, + \Big( \mc{L}_n\big(i,x\big) - \mc{L}_n\big(i,\mfl\big)  \Big) \, .
$$
The two bracketed terms are determined by the data in the second and first items in our new list; thus, the witness of~$\mc{G}$ may compute $\mc{L}_n^u(i,x)$ from these values. 

For $(i,x) \in \intint{k} \times \big( [-2\eln,\mfl] \cup [\mfr,2\eln] \big)$, we instead have
\begin{equation*}
 \mc{L}_n^u \big( i, x \big) = \begin{cases}
 \, \mc{L}_n^{[-2\eln,\mfl]}\big(i,x\big) + 
 \tfrac{\mfl - x}{\mfl + 2\eln}  \lppls_n\big(i, - 2\eln \big)
 + \tfrac{x + 2\eln}{\mfl + 2\eln}   \mc{L}_n^u\big(i,\mfl\big)     \, & x \in \big[-2\eln,\mfl \big] \, , \\
 \, \mc{L}_n^{[\mfr,2\eln]}\big(i,x\big) + 
 \tfrac{x - \mfr}{2\eln - \mfr}  \lppls_n\big(i,2\eln\big)
 + \tfrac{2\eln - x}{2\eln - \mfr}  \mc{L}_n^u\big(i,\mfr\big)   \, & x\in \big[\mfr,2\eln \big]  \, ,
\end{cases}
\end{equation*}
similarly as we did in the $\mc{F}$-case in~(\ref{e.reconcases}). Note that the values  $\mc{L}_n^u\big(i,\mfl\big)$   and  $\mc{L}_n^u\big(i,\mfr\big)$  used in this formula have already been defined.

Our notation $\mc{L}_n^u$ clashes with that for the reconstruction made in Section~\ref{s.missingclosedmiddle} by the witness of $\mc{F}$. Embellishing the notation with the symbol~$\mc{G}$ would relieve this difficulty. Since the new notation is employed only in this subsection, we hope that by omitting the $\mc{G}$ symbol, there is little danger of confusion.

Let $h_{\mc{G}}: \R \to [0,\infty)$ denote the density  with respect to Lebesgue measure on $\R$ of the conditional law under $\PP_{\mc{G}}$ of the random variable 
$\mc{L}_n \big(k,\mfl \big)$. By Lemma~\ref{l.brbr}, $h_{\mc{G}}(u)$ equals
\begin{equation}\label{e.gdensj.three}
 Z_{\mc{G}}^{-1}  \prod_{i=1}^k \, \exp\bigg\{-\tfrac{1}{2(\mfl + 2\eln)}\Big(\lppls_n \big(i,-2\eln \big) - \mc{L}_n^u\big(i,\mfl\big)
  \Big)^2 \, - \, \tfrac{1}{2(2\eln - \mfr)} \Big(  \mc{L}_n^u\big(i,\mfr\big)   - \lppls_n\big(i,2\eln  \big)  \Big)^2 \bigg\} \,\cdot\,\mfthree ( u ) \, ,
\end{equation}
where 
$\mfthree ( u )$ denotes the indicator function of the event
$$
  \Big\{ \, \mc{L}_n^{u} \big(k,x\big) >  \mc{L}_n \big( k+1,x\big) \, \, \, \, \forall \, \, x \in  [-2\eln,2\eln] \, \Big\} \, .
$$
The $\mc{G}$-measurable $Z_\mc{G} \in (0,\infty)$ 
is of course a normalization.

There are only two kinetic terms in the product in~(\ref{e.gdensj.three}): the third, associated to the middle interval~$[\mfl,\mfr]$, is absent because it is independent of the value of~$u \in \R$.

\begin{lemma}\label{l.hflthree}
There exists a $\mc{G}$-measurable random variable $\cornthree \in \R$ such that 
$$
\big\{ \, x \in \R : \mfthree ( x ) = 1 \, \big\} \, = \, \big( \cornthree \, , \, \infty \big) \, .
$$
\end{lemma}
\noindent{\bf Proof.} Unsurprisingly, this result has a similar proof to Lemma~\ref{l.glv}'s.  As $u$ falls, so does the curve $\mc{L}_n^u\big(k,\cdot\big):[-2\eln,2\eln] \to \R$, until a value of $u$ is reached at which the curve makes contact with $[-2\eln,2\eln] \to \R: u \to \mc{L}_n(k+1,u)$. The value of $u$ at which contact is made is $\cornthree$. The location of contact may lie in $[-2\eln,\mfl]$, $[\mfl,\mfr]$ or $[\mfr,2\eln]$. Indeed, we may specify a formula for $\cornthree$ by considering each of these three eventualities: $\cornthree$ equals the infimum of $q \in \R$ such that
\begin{itemize}
 \item $q \geq \xmin_k$; 
\item
 $q + \mc{L}_n(k,x) - \mc{L}_n(k,\mfl) \geq \mc{L}_n(k+1,x)$ for all $x \in [\mfl,\mfr]$;
 \item and $q + \mc{L}_n(k,\mfr) - \mc{L}_n(k,\mfl) \geq \ymin_k$.
\end{itemize} 

That  $u \in \R$ satisfies $\mfthree(u)=1$ precisely when $u >  \cornthree$ follows by the reasoning in the final paragraph of Lemma~\ref{l.glv}'s proof, but without the need of iteratively making further upward displacements of lower indexed curves.   \qed

\medskip

Define the $\mc{G}$-measurable event
$$
 \mpgg  =  \Big\{ \big\vert  \mc{L}_n(i,\mfr) - \mc{L}_n(i,\mfl) \big\vert  \leq (2 + 2^{-1/2}) \eln^2 \textrm{ for $i \in \intint{k}$} \Big\} \cap \big\{ \cornthree \leq \eln^2 \big\} 
 \, .
$$

\begin{lemma}\label{l.msigma}
The conditional distribution of $\mc{L}_n(k,\mfl)$ under~$\PP_\mc{G}$ is that of a normal random variable~$N$, of mean $m$ and variance $\sigma^2$, conditionally on the event that $N \geq \cornthree$. We have that 
$m = k^{-1} \sum_{i=1}^k m(i)$, where $m(i)$ equals
$$
 \big( \mc{L}_n(k,\mfl) -  \mc{L}_n(i,\mfl) \big) \, + \, \frac{\mfl + 2\eln}{4\eln + \mfl - \mfr} \Big(  \mc{L}_n(i,2\eln) - \big( \mc{L}_n(i,\mfr) - \mc{L}_n(i,\mfl) \big)  \Big) \, +  \, \frac{2 \eln - \mfr}{4\eln + \mfl - \mfr} \mc{L}_n(i,-2\eln) 
$$
 and
$$
\sigma^{-2} = k \Big( \big( \mfl + 2\eln \big)^{-1} + \big( 2\eln - \mfr \big)^{-1} \Big) \, .
$$
The occurrence of the event $\mpg \cap \mpgg$ entails that $m \geq - 12 \eln^2$ and $\sigma^2 \in \big[ \tfrac{1}{2} , \tfrac{3}{4} \big] \cdot \eln k^{-1}$.
\end{lemma}
Note that the first, bracketed, term in the expression for $m(i)$ equals $-\sum_{j=i}^{k-1} \big( \mc{L}_n(j,\mfl) - \mc{L}_n(j+1,\mfl) \big)$ and is thus $\mc{G}$-measurable; and so is  
$\mc{L}_n(i,\mfr) - \mc{L}_n(i,\mfl)$.  The remaining random parameters, $\mc{L}_n(i,2\eln)$, $\mc{L}_n(i,-2\eln)$, $\mfl$ and $\mfr$, are measurable in the smaller $\sigma$-algebra~$\mc{F}$.
Thus, $m(i)$  is $\mc{G}$-measurable and $\sigma^2$ is $\mc{F}$-measurable. 

\begin{lemma}\label{l.triplegauss}
Let $\ell_1 < a < b < \ell_2$, and let $\bar{y}$, $\bar{z}$,  $\bar{\jmath}$ and $\bar{r}$ be elements of $\R^k$. Under the law 
$\mc{B}_{k;\bar{y},\bar{z}}^{[\ell_1,\ell_2]}$ conditioned on $\ovbar{B}(b) - \ovbar{B}(a) = \bar{\jmath}$ and $B(i,a) - B(i+1,a) = r_i - r_{i+1}$ for $i \in \intint{k-1}$, the conditional distribution of $B(k,a)$ is normal, of mean and variance that we denote by  $m_0$ and  $\zeta^2$. 
If we write $\hat{x} = k^{-1} \sum_{i=1}^k x_i$ for any $\bar{x} \in \R^k$, then 
$$
m_0 = r_k - \hat{r} \, + \, \frac{a - \ell_1}{(a - \ell_1) + (\ell_2 - b)} \big( \hat{z} - \hat{\jmath}  \big) \, +  \, \frac{\ell_2 - b}{(a - \ell_1) + (\ell_2 - b)} \hat{y} 
$$
and
$$
\zeta^{-2}  =  k \Big( (a - \ell_1)^{-1} + (\ell_2 - b)^{-1} \Big) \, .
$$
\end{lemma}
\noindent{\bf Proof.} 
The conditioned random variable
$B(k,a)$ adopts the value $x$ with density given up to normalization by 
$$
\prod_{i=1}^k \exp \Big\{ -  \tfrac{1}{2(a-\ell_1)} \big( x + r_i - r_k - y_i \big)^2 \Big\} \cdot \exp \Big\{ - \tfrac{1}{2(\ell_2 -b)} \big(x + r_i - r_k  + j_i - z_i \big)^2  \Big\} \, ;
$$
again up to normalization, this quantity equals
$$
 \exp \Big\{ - \tfrac{k}{2} \big( (a - \ell_1)^{-1} + (\ell_2 - b)^{-1} \big) \big( x - m \big)^2  \Big\} \, ,
$$
where $m_0$ satisfies
$$
 m_0 k \big( (a - \ell_1) + (\ell_2 - b) \big)  \, = \,
 - \sum_{i=1}^k \Big( (\ell_2 - b)(r_i - r_k - y_i) \, + \, (a-\ell_1)(r_i - r_k + j_i - z_i) \Big) \, .
$$
Rearranging yields the stated formulas for $m_0$ and the reciprocal of $\zeta^2$.
 \qed

\medskip

\noindent{\bf Proof of Lemma~\ref{l.msigma}.} 
We apply Lemma~\ref{l.triplegauss}, setting $\ell_1 = -2\eln$, $a = \mfl$, $b = \mfr$, $\ell_2 = 2 \eln$, 
$\bar{y} = \ovbar{\mc{L}}_n(-2\eln)$, 
$\bar{z} = \ovbar{\mc{L}}_n(2\eln)$ and
$\bar{\jmath} =  \ovbar{\mc{L}}_n(\mfr) -  \ovbar{\mc{L}}_n(\mfl)$ 
and $\bar{r} = \ovbar{\mc{L}}_n(\mfl)$. 

The occurrence of $\mpg$ entails that $\mfl \in [-\eln,-\eln/2]$
and $\mfr \in [\eln/2,\eln]$, 
whence
$\eln k^{-1}/2 \leq \sigma^2 \leq 3\eln k^{-1}/4$.

The quantity $m(i)$ is a sum of three terms. When $\mpg$ occurs, these terms satisfy certain bounds:
\begin{itemize}
\item the first term satisfies 
$$
\mc{L}_n(k,\mfl) - \mc{L}_n(i,\mfl) \geq - (2^{-1/2} + 1) \eln^2  - \eln^2 = - \big( 2 + 2^{-1/2} \big) \eln^2 \, , 
$$
since $\mfl \in [-\eln,0]$;
\item the second term is at least $-\tfrac{3}{4 \sqrt{2}} \big( 5 + 3\sqrt{2} \big) \eln^2$,
since $\big\vert \mc{L}_n(i,\mfr) - \mc{L}_n(i,\mfl) \big\vert \leq   \big( 2 + 2^{-1/2} \big) \eln^2$ on $\mpgg$
and $\mc{L}(i,2\eln) \in -2\sqrt{2}\eln^2 + \eln^2 \cdot [-1,1]$ on $\mpg$,
 while 
$$
 \frac{\mfl + 2\eln}{4\eln + \mfl - \mfr} \leq \frac{3\eln/2}{2\eln} = \frac{3}{4} \, \, ;
$$
\item and the third term satisfies $\mc{L}_n(i,-2\eln) \geq \big(-2\sqrt{2} - 1 \big) \eln^2$.
\end{itemize}
Note that $\sigma^{-2} \in k \eln^{-1} \cdot [4/3,2]$ since $2\eln - \mfr$ and $\mfl + 2\eln$ both lie in the interval $[\eln,3\eln/2]$. \qed

\begin{lemma}\label{l.gthree}
$$
 \PP \big( \mpg^c_\mc{G} \big) \leq \e^{k^2 - 1}  \, .
$$
\end{lemma}
\noindent{\bf Proof.} 
We begin by noting that the occurrence of $\mpg_\mc{G}$ is forced by the conditions
$$
    \inf_{x \in [-\eln,\eln]} \mc{L}_n(k+1,x) \geq - \eln^2 \, \, \, \textrm{and} \, \, \, 
  \sup_{x \in [-\eln,\eln]} \mc{L}_n(1,x) \leq  \eln^2 \, .
  $$
Indeed, $\mfl,\mfr \in [-\eln,\eln]$ and $\mc{L}_n$ being ordered imply that the conditions entail the stronger-than-required $\big\vert \mc{L}_n(i,\mfr) - \mc{L}_n(i,\mfl) \big\vert \leq 2 \eln^2$
for $i \in \intint{k}$; and, in view of  $\cornthree \leq \mc{L}_n(k,\mfl) \leq \mc{L}_n(1,\mfl)$,  the latter condition ensures that $\cornthree \leq \eln^2$. 

Note however that the failure of the displayed conditions entails the occurrence of one or other of the events whose probabilities are bounded above in the analysis of $\PP \big( \fev_2^c \big)$ in the proof of Lemma~\ref{l.glub.new}. Given the form of this lemma, we find then that  $\PP \big( \mpg^c_\mc{G} \big) \leq \e^{2^{-5}\rsc_k \const^3}$.
The imposition made on $\const$ in (\ref{e.constvalue.new}) that $2^{-5}\rsc_k \const^3 \geq k^2 - 1$ then yields Lemma~\ref{l.gthree}. \qed

 
\medskip

We are about to present our snap up assertion, Lemma~\ref{l.apgt}, a result to the effect that the witness of $\mc{G}$ typically has a reasonable, $\e^{o(1)}$-order, probability of observing a value for the random height $\mc{L}_n(k,\mfl)$ which is high enough compared to its minimum $\cornthree$ that the $\abovepar\big(\mc{L}_n,\fa\big)$ event occurs. As we mentioned in overview,  the snap up order of magnitude is great enough to deliver the occurrence of this  $\abovepar \big( \mc{L}_n,\fa, 15 k \const (\log \e^{-1})^{1/2} \big)$ event with the requisite probability only if we impose a further `no big drop' condition on $\mc{G}$-measurable data. 
 
To describe the necessary condition, let $r \in \intint{\vert \pole \vert - 1}$ satisfy $\fa \in [\pp_r,\pp_{r+1}]$
and define the 
event 
$$
\nbd\big(\mc{L}_n,\fa\big) =  \Big\{ \,  \mc{L}_n^{[\pp_r,\pp_{r+1}]}(k,x_0)  \geq - \big( 2 (k-1)  + 4 \eln \big) (2\eln)^{1/2} \, \Big\} \, .
$$
This is the event that $\mc{L}_n(k,\fa)$ does not drop by more than the stated quantity below the affine interpolation of the values of $\mc{L}_n(k,\cdot)$ at the pole set values $p_r$ and $p_{r+1}$ that neighbour~$\fa$. We also write $\bd \big(\mc{L}_n,\fa\big)$ for the complementary event $\neg \, \nbd\big(\mc{L}_n,\fa\big)$. 

When the witness of~$\mc{G}$  observes the $\mc{G}$-measurable event $\nbd\big(\mc{L}_n,\fa\big)$ alongside certain other, typical, data, she deploys the snap up, as we now explain in detail.

\begin{lemma}\label{l.apgt}
$$
\PP_{\mc{G}} \Big(  
   \abovepar \big( \mc{L}_n,\fa \big) \Big)  \geq \tfrac{1}{4}  \, \exp \big\{ - 987 k^3 \eln^{5/2} \big\} \cdot {\bf 1}_{\mpg \cap \mpgg \cap \nbd(\mc{L}_n,\fa)} \, .
$$
\end{lemma}
\noindent{\bf Proof.} We claim that
$$
   \Big\{ \mc{L}_n(k,\mfl) \geq \cornthree  + \big( 2 (k-1)  + 4 \eln \big) (2\eln)^{1/2} + 15 k \const (\log \e^{-1})^{1/2} \Big\} \, \cap \,  \nbd\big(\mc{L}_n,\fa\big) 
$$
is a subset of $\abovepar\big( \mc{L}_n,\fa \big)$.
To verify this, 
consider an instance of the data specifying $\mc{G}$ for which $\nbd\big(\mc{L}_n,\fa\big)$ occurs. Under the law $\PP_{\mc{G}}$, the lowest value that the random variable $\mc{L}_n(k,\mfl)$ may adopt is $\cornthree$, in which case, the maximal drop 
$\tent(\fa) - \mc{L}_n(k,\fa)$ 
of the $k\textsuperscript{th}$ curve in $\mc{L}_n$ below the tent map at $\fa$ would be at most the quantity
$$
 \big( 2 (k-1)  +  4 \eln  \big) (2\eln)^{1/2} \, .
$$
 By raising $\mc{L}(k,\mfl)$ above this minimum value by this quantity, and then by a further distance of $15 k \const (\log \e^{-1})^{1/2}$, the occurrence of the event $\abovepar\big(\mc{L}_n,\fa, 15 k \const (\log \e^{-1})^{1/2} \big)$ is assured. 

It is enough then  to confirm the lemma's statement with the event  $\abovepar\big(\mc{L}_n,\fa \big)$ replaced by 
$$
\Big\{ \mc{L}_n(k,\mfl) \geq \cornthree + \big( 2 (k-1)  +  4 \eln  \big) (2\eln)^{1/2} + 15 k \const (\log \e^{-1})^{1/2} \Big\} \, .
$$
We remarked on the cancellation of first order kinetic costs after Proposition~\ref{p.jabovecorner}, and we now exploit a similar cancellation. Indeed, note that, for $y > 0$, on $\mpg \cap \mpgg$,
\begin{eqnarray}
 & & \PP_{\mc{G}} \Big(  
  \mc{L}_n(k,\mfl) \geq \cornthree + y   \Big) 
  =  \frac{\nu_{m,\sigma^2}\big( \cornthree + y \, , \, \infty \big)}{\nu_{m,\sigma^2}\big( \cornthree \, , \, \infty \big)} \nonumber \\ 
 & \geq & \frac{\nu_{-12\eln^2,\sigma^2}\big( \eln^2 + y \, , \, \infty \big)}{\nu_{-12\eln^2,\sigma^2}\big( \eln^2  \, , \, \infty \big)}  
  =   \frac{\nu_{0,1}\big( 13 \eln^2 \sigma^{-1}  +  y \sigma^{-1} \, , \, \infty \big)}{\nu_{0,1} \big(  13 \eln^2 \sigma^{-1}  \, , \, \infty \big)}  \nonumber \\
  & \geq & \frac{ 13  \eln^2 \sigma^{-1}}{2 \big(  13 \eln^2 \sigma^{-1} +  y \sigma^{-1} \big)} \frac{g_{0,1}\big( 13  \eln^2 \sigma^{-1} +  y \sigma^{-1} \big)}{g_{0,1} \big(  13 \eln^2 \sigma^{-1}   \big)}  \nonumber \\
    & \geq & \tfrac{1}{4} \cdot \exp \big\{ - 13 \eln^2 \sigma^{-2} y  -  y^2 \sigma^{-2}/2 \big\} \geq  \tfrac{1}{4} \cdot \exp \big\{ - 26 k \eln  y  -  k   y^2 \eln^{-1} \big\} \, , \label{e.regfin}
\end{eqnarray}
the first inequality due to Lemmas~\ref{l.monotonenormal}  and~\ref{l.msigma};
the second to Lemma~\ref{l.gaussiantail} via
$$ 
 13 \eln^2 \sigma^{-1} +  y \sigma^{-1}  \geq 1 \, ,
$$
which certainly follows from the bound $\sigma^2 \leq \tfrac{3}{4} \eln k^{-1}$;
the third to $y \leq 13\eln^2$, and the fourth to $\sigma^2 \geq \eln k^{-1}/2$. We now take 
$$
y =  \Big( 2 (k-1)  +  4 \eln  \Big) (2\eln)^{1/2} + 15 k \const (\log \e^{-1})^{1/2} \, .
$$ 
Using $\eln \geq 1$, $k \geq 2$ and $\const \geq 1$, we see that $y \leq 21k\eln^{3/2}$; thus, $y$ is indeed at most $13\eln^2$, since
$\e < e^{-18 \const^{-3} k^6}$. Applying $y \leq 21k\eln^{3/2}$ and $\eln \geq 1$ to find a lower bound on the expression in line~(\ref{e.regfin}) completes the proof of Lemma~\ref{l.apgt}. \qed

\medskip

Of course, the witness of $\mc{G}$ may hold $\mc{G}$-data that dictates the occurrence of $\bd$.  The swing through will replace  the snap up in this case. We now state the resulting conclusion, Lemma~\ref{l.lowtohigh}. We then close out the proof of Theorem~\ref{t.airynt.one} before providing the proof of the lemma.

\begin{lemma}\label{l.lowtohigh}
$$
  \PP \Big(    
  \close\big(\mc{L}_n,\fa,\e\big) \cap  \bd\big(\mc{L}_n,\fa\big) \cap \mpg  \Big) 
   \leq 
  \PP \Big(    
  \close\big(\mc{L}_n,\fa,\e\big) \cap    \abovepar\big(\mc{L}_n,\fa\big)   \Big) \, .
$$
\end{lemma} 

\noindent{\bf Proof of Theorem~\ref{t.airynt.one}.}
Applying Lemma~\ref{l.apgt}, we find that
\begin{eqnarray*}
 & &  \PP \Big( \, 
  \close\big(\mc{L}_n,\fa,\e\big) \cap    \abovepar\big(\mc{L}_n,\fa\big)   \, \Big) 
  =  \EE \, \Big[ \, \PP_{\mc{G}} \big(  
   \abovepar ( \mc{L}_n,\fa )   \big)  \cdot {\bf{1}}_{\close\big(\mc{L}_n,\fa,\e\big)} \, \Big] \\
   & \geq & \tfrac{1}{4}  \, \exp \big\{ - 987 k^3 \eln^{5/2} \big\} \,   \PP \Big(    
  \close\big(\mc{L}_n,\fa,\e\big) \cap \mpg \cap \mpgg  \cap \nbd\big(\mc{L}_n,\fa\big) \Big)  \, .
\end{eqnarray*}

We set the inter-pole distance $\ipd$ used to specify the jump ensemble equal to one. Using Proposition~\ref{p.labovepar} and Lemma~\ref{l.lowtohigh}, we find that $\PP \big(    
  \close(\mc{L}_n,\fa,\e)   \big)$, being equal to the sum of the probabilities of two disjoint events,
$$
        \PP \Big(    
  \close\big(\mc{L}_n,\fa,\e\big) \cap \nbd\big(\mc{L}_n,\fa\big)  \Big) +  \PP \Big(    
  \close\big(\mc{L}_n,\fa,\e\big) \cap  \bd\big(\mc{L}_n,\fa\big)   \Big) \, ,
  $$
  is at most
  $$
   \Big( 4  \, \exp \big\{ 987 k^3 \eln^{5/2} \big\}  + 1 \Big)  \cdot 2  \big( 10^4 \const^{k^2} + 2k \big) \exp \Big\{  3973 k^{7/2} \const^2 \big( \log \e^{-1} \big)^{2/3} \Big\}  \big(\log \e^{-1}\big)^{k^2/2}   \e^{k^2-1} 
  \, \, + \, 2\e^{k^2-1} \, ,
$$
where the $ 2\e^{k^2-1}$ term is composed of two parts, contributed by Lemmas~\ref{l.glub.new} and~\ref{l.gthree}.

Using $\e < e^{-1}$ and $\const \geq 1$, the last display is bounded above by 
$$
2  \cdot    7  \big( 10^4 \const^{k^2} + 2k \big) \exp \Big\{  4960 k^{7/2}  \const^{5/2} \big( \log \e^{-1} \big)^{5/6} \Big\}  \big(\log \e^{-1}\big)^{k^2/2}   \e^{k^2-1}  
$$
and thus by 
$$
14  \big( 10^4 + 1 \big) \cdot  \exp \Big\{  4961 k^{7/2} \ipd^2 \const^{5/2} \big( \log \e^{-1} \big)^{5/6} \Big\}  \big(\log \e^{-1}\big)^{k^2/2}   \e^{k^2-1}  \, .
$$
Since $\e < e^{-2^{3/2}}$ ensures that $x \leq e^{x^{5/6}}$ where $x = \log \e^{-1}$, we use $\const \geq 1$ to obtain the upper bound stated in Theorem~\ref{t.airynt.one}. \qed

\medskip

\noindent{\bf Proof of Lemma~\ref{l.lowtohigh}.}
Recall that $r \in \intint{\vert \pole \vert - 1}$ satisfies $\fa \in [\pp_r,\pp_{r+1}]$. It is our task to explain why, conditionally on the top $k$ curves $\mc{L}_n$ collecting closely together at $\fa \in [p_r,p_{r+1}]$, it is not likely that the locale of this meeting is so low that the $k\textsuperscript{th}$ of these curves (and by ensemble ordering, the higher $k-1$ as well), must drop precipitously from time~$p_r$ to reach the locale at time $\fa$ before rising again rapidly as time $p_{r+1}$ is reached. (At least, this description is correct in coordinates chosen so that the $k\textsuperscript{th}$ curve has equal height at $p_r$ and $p_{r+1}$.)  

Such downward swooping on the part of the $k$ curves during $[p_r,p_{r+1}]$ is presumably kinetically costly. We will release this kinetic energy via the swing through in order to find an alternative configuration of comparable probability in which $k$-curve closeness at $\fa$ is maintained but where $\abovepar(\mc{L}_n,\fa)$ also occurs. 
To set up the swing through technique,
we introduce another $\sigma$-algebra that we label $\mc{H}[x_0]$. 

Note that the 
 ensemble $\mc{L}_n$ is entirely specified by the data:
\begin{itemize}
\item the curves $\mc{L}_n$ on $\intint{k} \times (\pp_r,\pp_{r+1})^c$ and $\llbracket k + 1 , n \rrbracket \times (-\xnmac,\infty)$;
\item the standard bridges $\mc{L}_n^{[\pp_r,\fa]}: \intint{k} \times [\pp_r,\fa] \to \R$ and $\mc{L}_n^{[\fa,\pp_{r+1}]}: \intint{k} \times [\pp_r,\fa] \to \R$;
\item the differences $\mc{L}_n(i,\fa) - \mc{L}_n(i+1,\fa)$, $i \in \intint{k-1}$;
\item and the value $\mc{L}_n(k,\fa)$. 
\end{itemize}
Let $\mcfa$ denote the $\sigma$-algebra generated by the first three items. Note first that $\mc{F} \subset \mcfa$. Only the final one-dimensional piece of data remains random to the witness of $\mcfa$.
This witness may be depicted as holding a bead at an unknown height $\mc{L}_n(k,\fa)$ on a vertical rod at coordinate $\fa$; as the bead is pushed up or down, it forces the values of the lower indexed curves at $\fa$ in lockstep; these placements in turn force the form of the curves on the intervals $[p_r,\fa]$ and $[\fa,p_{r+1}]$ in accord with the fixed values at $p_r$ and $p_{r+1}$ by means of affine displacement. 
  
Note that since the differences are fixed in the third bullet point, the event  $\close(\mc{L}_n,\fa,\e)$ is $\mcfa$-measurable. 

Under $\PP_{\mcfa}$, the occurrence of $\bd\big(\mc{L}_n,\fa\big)$
is characterized by the condition that the random quantity $\mc{L}_n(k,\fa)$ is at most an $\mcfa$-measurable quantity $R : = \ell(\fa) - \big( 2(k-1) + 4 \eln \big)(2\eln)^{1/2}$.
Here, $\ell:[p_r,p_{r+1}] \to \R$ denotes the affine interpolation of $p_r \to \mc{L}_n\big(k,p_r\big)$ and  $p_{r+1} \to \mc{L}_n\big(k,p_{r+1}\big)$. 
 Under this same measure, on the other hand, the event $\abovepar \big( \mc{L}_n,\fa , s \big)$ occurs precisely when   $\mc{L}_n(k,\fa)$ is at least $\tent(\fa) + s$ (for any $s > 0$).

Consider now the law $\PP_{\mcfa}$. By the argument that leads to Lemma~\ref{l.hflthree}, the conditional distribution of the random variable $\mc{L}_n(k,\fa)$ 
takes the form of a Gaussian random variable conditioned to exceed a certain $\mcfa$-quantity that we label $\cornfour$. (We use the shorthand $\mc{H} = \mcfa$ in so doing.) This Gaussian has an $\mcfa$-measurable mean and variance that we denote by $m$ and $\sigma^2$. 
With a view to making some inferences about the values of these parameters, we note that the kinetic energy associated to the Gaussian equals
\begin{eqnarray*}
 & &
 \tfrac{1}{2} \, \sum_{i=1}^k \bigg( \, \tfrac{1}{\fa - \pp_r} \, \Big( x \, + \ \mc{L}_n(i,\fa) - \mc{L}_n(k,\fa) -  \mc{L}_n\big(i,\pp_r\big) \Big)^2 \\
  & &  \qquad \qquad \qquad \qquad \qquad \qquad  + \, \,  \tfrac{1}{\pp_{r+1} - \fa} \, \Big( x \, +  \, \mc{L}_n(i,\fa) - \mc{L}_n(k,\fa) -  \mc{L}_n\big(i,\pp_{r+1}\big) \Big)^2  \, \bigg) \, ,
\end{eqnarray*}
which up to normalization is given by 
\begin{eqnarray*}
 & & \tfrac{k}{2} \cdot \frac{\pp_{r+1} - \pp_r}{(\fa - \pp_r)(\pp_{r+1} - \fa)} \, \Bigg(  x \, + \, \tfrac{1}{k} \sum_{i=1}^k \bigg( \, \Big( \mc{L}_n(i,\fa) - \mc{L}_n(k,\fa) -  \mc{L}_n\big(i,\pp_r\big) \Big) \tfrac{\pp_{r+1} - \fa}{\pp_{r+1} - \pp_r}  \\
  & & \qquad \qquad \qquad \qquad \qquad \qquad \qquad \qquad \qquad  \,  + \, 
  \Big( \mc{L}_n(i,\fa) - \mc{L}_n(k,\fa) -  \mc{L}_n\big(i,\pp_{r+1}\big) \Big) \tfrac{\fa - \pp_r}{\pp_{r+1} - \pp_r} \bigg) \Bigg)^2 \, .
\end{eqnarray*}

We find then that 
$$
 m = -  \tfrac{1}{k} \sum_{i=1}^k  \big( \mc{L}_n(i,\fa) - \mc{L}_n(k,\fa) \big) \, + \,  \tfrac{1}{k} \sum_{i=1}^k  \Big(  \tfrac{\pp_{r+1} - \fa}{\pp_{r+1} - \pp_r}    \mc{L}_n(i,\pp_r)  \, + \,  \tfrac{\fa - \pp_r}{\pp_{r+1} - \pp_r}  \mc{L}_n\big(i,\pp_{r+1}\big)   \Big)   \, ,
$$
and that $\sigma^2 = k^{-1} \tfrac{(\fa - \pp_r)(\pp_{r+1} - \fa)}{\pp_{r+1} - \pp_r}$.

In this formula for $m$, the summand in the latter sum is at least the affine interpolation $\ell(\fa)$ defined a few moments ago. On the event $\close(\mc{L}_n,\fa,\e)$, the earlier summand, $\mc{L}_n(i,\fa) - \mc{L}_n(k,\fa)$,  lies in $(0,\e]$ for each $i \in \intint{k-1}$; thus $m \geq \ell(\fa) - \e$ 
on this event. Note also that $\ell(\fa) \geq \tent(\fa)$ because $\ell$ and $\tent$ are both affine on $[p_r,p_{r+1}]$, with $\ell$ having the higher boundary data since $\tent(p)$ equals $\mc{L}_n(k+1,p)$ for $p \in \{ p_r,p_{r+1}\} \subset P$.

With these things observed, we may note that, for any specification of $\mcfa$ that realizes~$\close(\mc{L}_n,\fa,\e)$, 
\begin{eqnarray*}
\PP_{\mcfa} \big( \abovepar\big(\mc{L}_n,\fa\big)) \big)
 & = & \nu_{m,\sigma^2} \Big( \tent(\fa)+ 15 k \const (\log \e^{-1})^{1/2}  , \infty \, \Big\vert \, \cornfour , \infty  \Big) \\
 & \geq & \nu_{m,\sigma^2} \big( \tent(\fa) + 15 k \const (\log \e^{-1})^{1/2}  , \infty \big) \\
 & \geq & \nu_{ \tent(\fa)-\e,\sigma^2} \big(  \tent(\fa)   + 15 k \const (\log \e^{-1})^{1/2}  , \infty \big) \\
 & = & \nu_{0,\sigma^2} \big(  15 k \const (\log \e^{-1})^{1/2}  + \e, \infty  \big) \, ,
\end{eqnarray*}
where recall that $\abovepar\big(\mc{L}_n,\fa\big)) =  \abovepar \big( \mc{L}_n,\fa, 15 k \const (\log \e^{-1})^{1/2} \big)$. (In the second inequality, we used $\ell(\fa) \geq \tent(\fa)$.)
We also have that, for such a specification,
\begin{eqnarray*}
\PP_{\mcfa} \Big(  \bd(\mc{L}_n,\fa) \Big)
  & = & \nu_{m,\sigma^2} \big( -\infty, R \, \big\vert \, \cornfour , \infty  \big)
   \leq  \nu_{m,\sigma^2} \big( -\infty, R   \big) \\
   & \leq & \nu_{ \ell(\fa)-\e,\sigma^2} \Big( -\infty \, , \, \ell(\fa) -  \big( 2(k-1) + 4 \eln \big)(2\eln)^{1/2}   \Big) \\
   & = & \nu_{0,\sigma^2} \Big( \big( 2(k-1) + 4 \eln \big)(2\eln)^{1/2}  - \e, \infty  \Big) \, . 
\end{eqnarray*}
That
$$
   15 k \const (\log \e^{-1})^{1/2} + \e \leq \big(  2(k-1) + 4 \eln \big)(2\eln)^{1/2}  - \e  
$$
is due to $\const \geq 1$, $\e < e^{-1}$ and $\const \geq \tfrac{1}{32} (15 k +1)^2$;
we learn that
$$
\PP_{\mcfa} \big( \abovepar\big(\mc{L}_n,\fa\big)  \big) \geq 
\PP_{\mcfa} \Big( \bd\big(\mc{L}_n,\fa\big) \Big) \, .
$$
That $\mc{F} \subset \mcfa$ implies that $\mpg$ is $\mcfa$-measurable. Thus,
\begin{eqnarray*}
& &  \PP \Big(    
  \close\big(\mc{L}_n,\fa,\e\big) \cap  \bd\big(\mc{L}_n,\fa\big)   \cap \mpg \Big)  \\
   & \leq & \EE   \, \bigg[ \,
\PP_{\mcfa} \Big( \bd\big(\mc{L}_n,\fa\big)  \Big) \cdot
{\bf 1}_{\close(\mc{L}_n,\fa,\e) \cap \mpg}  \, \bigg]   \\
   & \leq &     \EE   \, \bigg[ \,
\PP_{\mcfa} \big(  \abovepar(\mc{L}_n,\fa) \big) \cdot 
{\bf 1}_{\close(\mc{L}_n,\fa,\e) \cap \mpg}  \, \bigg]   \leq   \PP \Big(    
  \close\big(\mc{L}_n,\fa,\e\big) \cap  \abovepar\big(\mc{L}_n,\fa\big))  \Big) \, . 
\end{eqnarray*}
This completes the proof of Lemma~\ref{l.lowtohigh}. \qed




\section{Closeness of curves at a general location}\label{s.closenessgeneral}

In this section, we prove Theorem~\ref{t.airynt}(2) and Theorem~\ref{t.weakbound}. In view of  Theorem~\ref{t.airynt}'s first part, the rough form of its second part is plausible: since the curves of $\mc{L}_n$ are expected to be locally Brownian, the failure of $k$-curve $\e$-closeness at any given point may be expected to dictate this eventuality in a neighbourhood of width $\e^2$ about the point. Since $\e^{-2}$ such neighbourhoods are needed to cover a unit-order interval, the exponent value changes from $k^2 -1$ and $k^2 -3$ between the theorem's first and second part. We will follow this proof approach, and for this reason, we begin by showing that the local fluctuation of ensemble curves is Gaussian in magnitude.

Recall that jump ensemble method parameter conditions
~(\ref{e.constvalue}),
~(\ref{e.epsilonupperbound}),~(\ref{e.nlbone}) and~(\ref{e.nlbtwo})
  remain in force. Recall also the $k$-line modulus of continuity Definition~\ref{d.klinemod}.

\begin{proposition}\label{p.raremod}
Take the inter-pole distance parameter $\ipd$ equal to one. 
Let $\phi \in \big( 0 , \tfrac{1}{4} \eln^{-1} \wedge C_1^2 \eln^{-2} \big)$.
Then, for any $\Cone \geq 96$,
$$
 \PP_\mc{F} \Big( \omega_{k,[\mfl,\mfr]}\big(J, \phi \big) >  \Cone \phi^{1/2} \big( \log \phi^{-1} \big)^{1/2}   \Big) \cdot {\bf{1}}_{\mpg} \, \leq \, 
 \e^{-36(k+2)^2 k \const^3} 2^k \cdot 3k \cdot \phi^{\tfrac{\Cone^2}{4608}}    \, .
$$
\end{proposition}

\begin{corollary}\label{c.lraremod}
Let $\phi \in \big( \e^8 , \tfrac{1}{4} \eln^{-1}  \wedge C_1^2 \eln^{-2} \big)$.
Then, for any $\Cone \geq 96$,
\begin{eqnarray*}
 & & 
 \PP \Big( \omega_{k, I}\big(\mc{L}_n, \phi \big) > 
 \Cone \phi^{1/2} \big( \log \phi^{-1} \big)^{1/2}  
   \Big) \\
 & \leq &  \phi^{\tfrac{\Cone^2}{4608}}  \cdot 
 \e^{-36(k+2)^2 k \const^3} 2^k \cdot 3k \cdot
8 \exp \Big\{  3973 k^{7/2} \const^2 \big( \log \e^{-1} \big)^{2/3} \Big\} 
  \, + \,  \PP \big( \mpg^c \big) \, , 
\end{eqnarray*} 
where the interval $I$ is given by $I = \tfrac{1}{4} \const \big(\log \phi^{-1}\big)^{1/3} \cdot [-1,1]$.
\end{corollary}
\noindent{\bf Proof.}
Note that, under $\PP_\mc{F}$, the marginal law of $\mc{L}_n$ on $\intint{k} \times [\mfl,\mfr]$ coincides with that of $J:\intint{k} \times [\mfl,\mfr] \to \R$ given $\test_4(J) =1$.
Since $\ipd = 1$, Propositions~\ref{p.raremod} and~\ref{p.vaultsuccess} thus imply that
\begin{eqnarray*}
 & &  \PP_\mc{F} \Big( \omega_{k,[\mfl,\mfr]}\big(\mc{L}_n, \phi \big) >  \Cone \phi^{1/2} \big( \log \phi^{-1} \big)^{1/2} 
  \Big) {\bf{1}}_{\mpg} \\
 & \leq &  \phi^{\tfrac{\Cone^2}{4608}}  \cdot 
 \e^{-36(k+2)^2 k \const^3} 2^k \cdot 3k \cdot
 \exp \Big\{  3973 k^{7/2} \const^2 \big( \log \e^{-1} \big)^{2/3} \Big\}  \, . 
\end{eqnarray*}
Recalling that $\mpg \subseteq \big\{ [-\eln/2,\eln/2] \subseteq [\mfl,\mfr] \big\}$ and $\phi \geq \e^8$ implies the result. \qed

\medskip

\noindent{\bf Proof of Theorem~\ref{t.weakbound}.}
We begin by recalling the parameter $\const$
that was introduced as part of the jump ensemble apparatus in Section~\ref{s.elnfix}.
We mentioned there that we would permit an increase in the value of this parameter according to the requirements of the given application.
 
Since the failure probability of $\mpg$
is known by Lemma~\ref{l.glub.new}
to be at most  $\e^{2^{-5} c_k \const^3}$, while the present application demands a conclusion of the form $\e^{O(1)K^2}$, we see that we must increase $\const$ if necessary so that these last two expressions are roughly equal.

The hypotheses of the present theorem already imply that
$\Cwb \geq \sqrt{3973} (k+2)k^{3/4}  \sqrt{3 \cdot 4608} \const^{3/2}$. This theorem presents a given choice of $\const$, and it hypotheses a lower bound on $K$ in terms of $\const$. We cope with the need to increase $\const$ as described for the purpose of proving the theorem by permitting a notational abuse, for the purposes only of the present proof, in which the value of $\const$ is increased so that in fact $\Cwb$ is actually equal to $\sqrt{3973} (k+2)k^{3/4}  \sqrt{3 \cdot 4608} \const^{3/2}$. This abuse would be problematic for hypotheses of Theorem~\ref{t.weakbound}
which involve bounds expressed in terms of $\const$. In this regard, note that the condition~(\ref{e.epsilonandn})
is~(\ref{e.nlbone}) reexpressed in terms of $K$, rather than $\const$, in view of 
$\const = \big( 3973 \cdot 3 \cdot 4608 \big)^{-1/3} (k+2)^{-2/3}k^{-1/2} K^{2/3}$; and similarly is the upper bound on $\e$ expressed.

We have then that  $\tfrac{1}{3 \cdot 4608} \Cwb^2 \geq 36(k+2)^2 k \const^3$ and, alongside $\e < e^{-1}$ and $\const \geq 1$, that the quantity $\tfrac{1}{3 \cdot 4608} \Cwb^2  \big( \log \e^{-1} \big)^{1/3}$ is at least $3973 k^{7/2} \const^2$. 
We now apply Corollary~\ref{c.lraremod} with $\phi = \e$ and with $\Cone = \Cwb$, using Lemma~\ref{l.glub.new} to bound~$\PP \big( \mpg^c \big)$. We learn that 
$$
 \PP \Big( \omega_{k, I}\big(\mc{L}_n, \e \big) > 
 \Cwb \e^{1/2} \big( \log \e^{-1} \big)^{1/2}  
   \Big)  \leq   24k \cdot 2^k \e^{\tfrac{\Cwb^2}{3 \cdot 4608}} + \e^{2^{-5} c_k \const^3} \, .
$$
We bound $\e^{2^{-5} c_k \const^3}$ above 
in terms of $\Cwb$ by the expression
$\e^{  10^{-10} c_k (k+2)^{-2}k^{-3/2} K^{2} }$.

The last displayed expression is thus seen to be at most
$$
 \big( 24 k 2^k + 1 \big) 
\e^{  (6 \cdot 10^9)^{-1} c_k (k+2)^{-2}k^{-3/2} K^{2} } 
  \leq 25 k 2^k \, 
\e^{ 10^{-10} c_k (k+2)^{-2}k^{-3/2} K^{2} }   \, .
$$
This completes the proof.  \qed

\medskip

For the proof of Proposition~\ref{p.raremod}, a lemma is needed.

\begin{lemma}\label{l.bbmodcon}
Let $\delta \in (0,1/6)$ and $R \geq 24 \sqrt{2} \, \delta^{1/2} \big( \log \delta^{-1} \big)^{1/2}$. Then
$$
 \mc{B}_{1;0,0}^{[0,1]} \Big( \omega(B,\delta) > R \Big) \leq
  3 \exp \big\{ - \tfrac{1}{1152} R^2 \delta^{-1} \big\} \, ,
$$
where $\omega = \omega_{1,[0,1]}$ denotes the modulus of continuity of a single curve.
\end{lemma}
\noindent{\bf Proof.} Our argument has similarities to the proof of Theorem~\ref{t.aestimate}. We mentioned there that standard Brownian bridge $B:[0,1] \to \R$ may be represented in the form 
$$
B(t) \, = \, (1-t) \bmotion \big( \tfrac{t}{1-t} \big) \, \, ,  \, \, \, t \in [0,1] \, , 
$$
where $\bmotion:[0,\infty) \to \R$ is standard Brownian motion. Note that, when $r,s \in [0,2/3]$, $r < s$, 
$$
 \big\vert B(s) - B(r) \big\vert \, \leq \, \Big\vert  \bmotion \big( \tfrac{s}{1-s} \big) - \bmotion \big( \tfrac{r}{1-r} \big) \Big\vert + \big( s - r \big) \sup_{x \in [0,2/3]} \vert \bmotion \vert \big( \tfrac{x}{1-x} \big) \, .
$$
Since $[0,2/3]:x \to x/(1-x)$ has nine as a Lipschitz constant, we find that
$$
 \omega_{1,[0,2/3]} (B,\delta) \, \leq \, \omega_{1,[0,2]} (\bmotion,9\delta) \, + \, \delta  \sup_{t \in [0,2]} \vert \bmotion  ( t ) \vert \, .
$$
Note that, for any parameter $\phi > 0$, interval $I \subseteq \R$ of length at least $\phi$, and function $f:I \to \R$, 
$$
\omega_{1,I} (f,\phi) \, \leq \, 2 \sup \vert f(x+y) - f(x) \vert \, ,
$$ 
where the supremum is taken over pairs $(x,y) \in \phi \Z \times [0,2\phi]$ such that $x$ and $x+y$ are elements of~$I$. From this, we find that, if we write
$$
X = 2 \max_{i \in \llbracket 0, 2/(9\delta) \rrbracket} \sup_{t \in [0,18\delta]} \big\vert \bmotion(9\delta i + t) -  \bmotion(9\delta i )  \big\vert  \, + \, \delta  \sup_{t \in [0,2]} \vert \bmotion ( t ) \vert \, ,
$$
then $X$ is an upper bound on $\omega_{1,[0,2/3]} (B,\delta)$. A random variable sharing the law of $X$ offers an upper bound on $\omega_{1,[1/3,1]} (B,\delta)$, from which it follows in light of $\delta \in (0,1/6)$ and stationarity of Brownian motion increments that
$$
 \mc{B}_{1;0,0}^{[0,1]} \Big( \omega(B,\delta) > R \Big) \leq 
 2 \, \big( \tfrac{2}{9\delta} + 1 \big)  \, \PP \Big(  2 \sup_{t \in [0,18\delta]}  \big\vert \bmotion ( t ) \big\vert  > R/2  \Big) \, + \,   \, 2 \, \PP \Big( \delta  \sup_{t \in [0,2]} \vert \bmotion ( t ) \vert  > R/2 \Big) \, ;
$$
here, we write $\PP$ for the probability measure associated to the Brownian motion $\bmotion$. Using $\delta \leq 1$, the $\bmotion \to -\bmotion$ symmetry of Brownian motion, and the reflection principle, followed by the standard upper bound on the tail of the Gaussian distribution,  we find that
\begin{eqnarray*}
 \mc{B}_{1;0,0}^{[0,1]} \Big( \omega(B,\delta) > R \Big) & \leq & 
 \tfrac{22}{9} \delta^{-1} \cdot 4  \, \PP \big(  \bmotion ( 18\delta )  > R/4  \big) \, +  \,  8 \, \PP \big( \bmotion ( 2 )  > \tfrac{1}{2} R \delta^{-1} \big) \\
 & \leq & \tfrac{88\cdot 4}{3} \pi^{-1/2}  \delta^{-1/2} R^{-1} \cdot \exp \big\{ - \tfrac{R^2}{32 \cdot 18 \delta} \big\} \, + \, 8 (2\pi)^{-1/2} \cdot 2\sqrt{2} R^{-1} \delta \cdot \exp \big\{ - \tfrac{R^2}{16 \delta^2} \big\} \, .
\end{eqnarray*}
We make use of our assumption that $R^2 \geq 2 \cdot 32 \cdot 18 \delta \log \delta^{-1}$ (and also use $\delta \leq e^{-1}$) in order to find that 
the latter expression is at most 
$$ 
 \tfrac{88\cdot 4}{3} \pi^{-1/2} \big( 2 \cdot 32 \cdot 18 \big)^{-1/2} \exp \big\{ - \tfrac{R^2}{2 \cdot 32 \cdot 18 \delta} \big\} \, + \, 4 \sqrt{2} \, \pi^{-1/2} \cdot \tfrac{2}{(32 \cdot 18)^{1/2}} \delta^{1/2} \, \exp \big\{ - \tfrac{R^2}{16 \delta^2} \big\} \leq 3  \, \exp \big\{ - \tfrac{1}{1152} R^2 \delta^{-1} \big\} \, .
$$ 

This completes the proof of Lemma~\ref{l.bbmodcon}. \qed

\medskip

\noindent{\bf Proof of Proposition~\ref{p.raremod}.}
Recall that, if $\mpg$ occurs, then, under $\PP_\mc{F}$, the conditional distribution of the jump ensemble $J:\intint{k} \times [\mfl,\mfr] \to \R$
is obtained by conditioning the Wiener candidate  $\wien:\intint{k} \times [\mfl,\mfr] \to \R$ on success in the test conditions $\test_{12}(\wien) = 1$. Recall also that, under $\PP_\mc{F}$, $\wien$ is the marginal on $\intint{k} \times [\mfl,\mfr]$ of a bridge ensemble with law $\mc{B}_{k;\ovbar{\mc{L}}_n(-2\eln),\ovbar{\mc{L}}_n(2\eln)}^{[-2\eln,2\eln]}$.
For our present purpose, it will be a convenient notational abuse to allow $\wien$ to denote this bridge ensemble of which it is in reality a marginal. In this way, $\wien$ is in this argument defined on $\intint{k} \times [-2\eln,2\eln]$.  We have then that
\begin{eqnarray}
 & &  \PP_\mc{F} \Big( \omega_{k,[\mfl,\mfr]}\big(J, \phi \big) >  \Cone \phi^{1/2} \big( \log \phi^{-1} \big)^{1/2}   \Big) \cdot {\bf{1}}_{\mpg} \label{e.pfmod} \\
 & \leq & 
 \PP_\mc{F} \Big( \omega_{k,[-2\eln,2\eln]}\big(\wien, \phi \big) >  \Cone \phi^{1/2} \big( \log \phi^{-1} \big)^{1/2}   \Big) \cdot  \PP_\mc{F} \Big( \test_{12}(\wien) = 1 \Big)^{-1} \cdot {\bf{1}}_{\mpg} \nonumber \\
 & \leq & \e^{-36(k+2)^2 k \const^3} 2^k \,
 \PP_\mc{F} \Big( \omega_{k,[-2\eln,2\eln]}\big(\wien, \phi \big) >  \Cone \phi^{1/2} \big( \log \phi^{-1} \big)^{1/2}   \Big) \cdot   {\bf{1}}_{\mpg} \nonumber
\end{eqnarray}
where the second inequality 
is due to Proposition~\ref{p.wienerpromotion}.

Recall that, when $\mpg$ occurs, $\mc{L}_n \big( i,x \big) \in - 2\sqrt{2} \eln^2 + [-\eln^2,\eln^2]$ for $(i,x) \in \intint{k} \times \{ - 2\eln , 2\eln \}$. Thus, for any $r > 0$,
$$
  \PP_\mc{F} \Big( \omega_{k,[-2\eln,2\eln]}\big(\wien, \phi \big) >   r   \Big) \cdot   {\bf{1}}_{\mpg}  \leq k \cdot \sup_{(x,y) \in [0,2\eln^2]^2} \mc{B}_{1;x,y}^{[-2\eln,2\eln]}\Big( \omega_{1,[-2\eln,2\eln]} \big(\wien, \phi \big) >  r  \Big) \, .
$$
For any $x,y \in \R$, Brownian bridges $B$ and $B'$ under the laws  $\mc{B}_{1;0,0}^{[-2\eln,2\eln]}$ and  $\mc{B}_{1;x,y}^{[-2\eln,2\eln]}$
may be coupled by affine shift; when they are, the processes' moduli of continuity satisfy $\omega(B',\delta) \leq \omega(B,\delta) + \delta  (4\eln)^{-1} \vert y - x \vert$. This right-hand side is at most $\omega(B,\delta) + \delta  \eln/2$ when $x,y \in [0,2\eln^2]$. Choose $\delta = \phi$ and note that 
$$
 \frac{\delta \eln}{2} \, \leq \, \frac{\Cone}{2} \cdot \phi^{1/2} \big( \log \phi^{-1} \big)^{1/2}
$$
holds since $\phi \leq e^{-1} \wedge \Cone^2 \const^{-2} \big( \log \e^{-1} \big)^{-2/3}$. Hence,
\begin{eqnarray*}
 & &  
 \PP_\mc{F} \Big( \omega_{k,[-2\eln,2\eln]}\big(\wien, \phi \big) >  \Cone \phi^{1/2} \big( \log \phi^{-1} \big)^{1/2}   \Big) \cdot   {\bf{1}}_{\mpg} \\
 & \leq & k \cdot \mc{B}_{1;0,0}^{[-2\eln,2\eln]}\Big( \omega_{1,[-2\eln,2\eln]} \big(B, \phi \big) >   \Cone/2 \cdot  \phi^{1/2} \big( \log \phi^{-1} \big)^{1/2}  \Big) \\ 
 & \leq & k \cdot \mc{B}_{1;0,0}^{[0,1]}\Big( \omega_{1,[0,1]} \big(B, \tfrac{\phi}{4\eln} \big) > (4\eln)^{-1/2} \cdot  \Cone/2 \cdot  \phi^{1/2} \big( \log \phi^{-1} \big)^{1/2}  \Big) \\ 
 & \leq & 3k \cdot \phi^{\tfrac{\Cone^2}{4608}} \, ,
\end{eqnarray*}
where we take $\delta = \tfrac{\phi}{4\eln}$ and $R = \tfrac{1}{4} \cdot \eln^{-1/2} \Cone \cdot  \phi^{1/2} \big( \log \phi^{-1} \big)^{1/2}$ in Lemma~\ref{l.bbmodcon} to obtain the final inequality. 
The hypothesis of the lemma that $R^2 \geq 1152 \delta \log \delta^{-1}$
is satisfied because $\Cone^2 \log \phi^{-1} \geq 4608 \big( \log (4\eln) + \log \phi^{-1} \big)$ holds due to 
$\phi^{-1} \geq 4\eln$ and $\Cone^2 \geq 2\cdot 4608$.

Applying this bound to~(\ref{e.pfmod}) alongside $\eln = \const \big( \log \e^{-1} \big)^{1/3}$ proves Proposition~\ref{p.raremod}. \qed


Turning next to the task of deriving Theorem~\ref{t.airynt}(2),
we begin by reducing to the assertion made for an interval centred at $y = 0$.

\begin{proposition}\label{p.airynt.two}
For $\bar\phimac \in (0,\infty)^3$, $C,c > 0$ and $n \in \N$, let 
$$
\mc{L}_n:\intint{n} \times \big[-\xnmac,\infty\big) \to \R  
$$ 
be a    $\big(\bar\phimac,\rsc,\rsC\big)$-regular ensemble defined under the law~$\PP$. For $k \geq 2$,
 let $\const$ be given by~(\ref{e.constvalue.new}). Let $\e > 0$ satisfy the bound~(\ref{e.epsilonupperbound}). 
For $n,\in \N$ satisfying $n \geq k 
  \vee  (c/3)^{-2\phimac_2^{-1}} \vee  6^{2/\delta}$ and 
    $n^{\phimac_1/4 \, \wedge  \, \phimac_2/4 \, \wedge \, \phimac_3/2}   \geq    \big( \rsc/2 \wedge 2^{1/2} \big)^{-1} \const \big( \log \e^{-1} \big)^{1/3}$, 
 the following bound holds:
\begin{eqnarray*}
 & & \PP \Big(  \exists \, x \in \R \, , \, \vert x \vert \leq \tfrac{1}{4} \const \big( \log \e^{-1} \big)^{1/3}:  \close\big( k ; \mc{L}_n , x , \e \big) \Big) \\
 & \leq &
\e^{k^2 - 3} \cdot 10^{44} \,  2^{6k^2} \const^{18} \exp \Big\{  4963 k^{7/2} \const^{5/2} \big( \log \e^{-1} \big)^{5/6} \Big\}  \, .
\end{eqnarray*}
\end{proposition}
\noindent{\bf Proof of Theorem~\ref{t.airynt}(2).} Proposition~\ref{p.airynt.two} 
implies the result via the parabolic invariance Lemma~\ref{l.parabolicinv}, with the value of $c$ decreasing by a factor of two, resulting in the replacement of factors of $\const$ by $2^{1/3} \const$. \qed

\medskip

\noindent{\bf Proof of Proposition~\ref{p.airynt.two}.}
It suffices to prove
\begin{eqnarray}
 & & \PP \Big(  \exists \, x \in \R \, , \, \vert x \vert \leq \tfrac{1}{4} \const \big( \log \e^{-1} \big)^{1/3}:  \lppls_n\big(1,x \big) \leq \lppls_n\big(k,x \big) + \e/2 \Big) \label{t.airynt.tworeduce} \\
 & \leq &
\e^{k^2 - 3} \cdot 10^{44}  k^{18} \const^{18} \exp \Big\{  4963 k^{7/2}  \const^{5/2} \big( \log \e^{-1} \big)^{5/6} \Big\} \, ; \nonumber
\end{eqnarray}
indeed, the result then follows by replacing $\e/2$ by $\e$, this replacement leading to an additional factor of  $2^{k^2 - 3}$ on the right-hand side, and applying $k^{18} \leq 2^{5k^2}$ when $k \geq 2$.

Note that
\begin{eqnarray*}
 & & \Big\{ \exists \, x \in \big[ -\eln/4,\eln/4 \big]: \close\big(\mc{L}_n,x,\e/2\big) \Big\} \, \cap \, \Big\{ \omega_{k,[-\eln/4,\eln/4]}\big( \mc{L}_n , \phi \big) \leq \e/2 \Big\} \\
 & \subseteq & 
\Big\{ \exists \, x \in \phi \Z \cap \big[ -\eln/4,\eln/4 \big]: \close\big(\mc{L}_n,x,\e\big) \Big\} \, .
\end{eqnarray*}
Set
\begin{equation}\label{e.ephi}
 \e/2 =  \Cone \phi^{1/2} \big( \log \phi^{-1} \big)^{1/2} \, .
\end{equation}
Note that $\phi < e^{-1}$ implies that $\phi \leq \tfrac{1}{4} \e^2 \Cone^{-2}$.
We also claim that
\begin{equation}\label{e.phicone}
\phi \geq \tfrac{1}{4} \e^2 \Cone^{-2} \Big( \log \big( 16 \Cone^4 \e^{-4} \big) \Big)^{-1} \, .
\end{equation}
To verify this, note first that $\log \phi^{-1} \leq \phi^{-1/2}$ when $\phi \leq e^{-1}$. Squaring~(\ref{e.ephi}), applying this last bound and squaring again, we find that
$\phi^{-1} \leq 16 \Cone^4 \e^{-4}$. Returning to the square of~(\ref{e.ephi}) with this bound, we obtain~(\ref{e.phicone}).

Thus, by Theorem~\ref{t.airynt.one} and Corollary~\ref{c.lraremod},
\begin{eqnarray}
 & & 
 \PP \Big( \exists \, x \in \big[ -\eln/4,\eln/4 \big]: \close\big(\mc{L}_n,x,\e/2\big)   \Big) \nonumber \\
 &  \leq &   
 \phi^{\tfrac{\Cone^2}{4608}}  \cdot 
 \e^{-36(k+2)^2 k \const^3} 2^k \cdot 3k \cdot
8 \exp \Big\{  3973 k^{7/2} \const^2 \big( \log \e^{-1} \big)^{2/3} \Big\} 
  \, + \, \e^{k^2-1} \label{e.firstterm} \\
 & & \, + \, \big(\eln/2 + 1 \big) \phi^{-1} \cdot  10^6 \exp \Big\{  4962 k^{7/2} \const^{5/2} \big( \log \e^{-1} \big)^{5/6} \Big\}  \, \e^{k^2-1} \, , \label{e.thirdterm}
 \end{eqnarray}
 where the first instance of $\e^{k^2-1}$ appears via Lemma~\ref{l.glub.new}
 in light of the comment following~(\ref{e.constvalue.new}). 

Choose $\Cone > 0$ so that 
 $\tfrac{\Cone^2}{2304} -36(k+2)^2 k \const^3  = k^2 - 3$. Using  $\phi \leq \tfrac{1}{4} \e^2 \Cone^{-2}$, as well as $k \geq 2$, $\e < e^{-1}$ and $\const \geq 1$, we find that
the term in line~(\ref{e.firstterm})  is then bounded above by
\begin{eqnarray*}
 & & \e^{\tfrac{\Cone^2}{2304} -36(k+2)^2 k \const^3}  \big( 4 \Cone^2 \big)^{-\tfrac{1}{4608}\Cone^2} \cdot  2^k \cdot 3k \cdot
8 \exp \Big\{  3973 k^{7/2} \const^2 \big( \log \e^{-1} \big)^{2/3} \Big\} \\
 &  \leq &  \e^{k^2 - 3} \cdot 
 2^k \cdot 3k \cdot
8 \exp \Big\{  3973 k^{7/2} \const^2 \big( \log \e^{-1} \big)^{2/3} \Big\} 
   \leq   \e^{k^2 - 3} \cdot 
  \exp \Big\{  3974 k^{7/2} \const^2 \big( \log \e^{-1} \big)^{2/3} \Big\}  \, .
\end{eqnarray*}
 
Using~(\ref{e.phicone}), the term is line~(\ref{e.thirdterm}) is  found to be at most
\begin{eqnarray*}
 & &  (\eln + 1) \cdot 4 \e^{-2} \Cone^2 \log \big( 16 \Cone^4 \e^{-4} \big) \cdot  10^6 \exp \Big\{  4962 k^{7/2}  \const^{5/2} \big( \log \e^{-1} \big)^{5/6} \Big\}  \, \e^{k^2-1}   \\
 &  \leq &  \e^{k^2 - 3}  \cdot 2 \const \big( \log \e^{-1} \big)^{1/3} \cdot 4  \Cone^2 \cdot  16 \Cone^4 \cdot 4 \log \big( \e^{-1} \big) \cdot  10^6 \exp \Big\{  4962 k^{7/2} \const^{5/2} \big( \log \e^{-1} \big)^{5/6} \Big\}  \bigg) \\
 &  \leq &  \e^{k^2 - 3} \cdot 512 \cdot 10^6 \cdot  (663552)^6 k^{18} \const^{18} \exp \Big\{  4963 k^{7/2}  \const^{5/2} \big( \log \e^{-1} \big)^{5/6} \Big\}  \, ,
\end{eqnarray*}
where in the first displayed inequality we used $\eln \geq 1$; and also
 the bound $\log (ar) \leq a\log r$, which is valid when $a \geq 1$ and $r \geq 2$, was applied in the case that $a = 16 \Cone^4$ and $r = \e^{-4}$. In the second inequality, we used that our choice of $\Cone$ satisfies 
$\Cone^2 \leq 2304 \cdot 2 \cdot 36 \cdot 4 k^3 \const^3$.
 We also used   $\big( \log \e^{-1} \big)^{4/3} \leq \exp \big\{ (\log \e^{-1})^{5/6} \big\}$ if $\e < e^{-8}$, and $\const \geq 1$.

Using these bounds alongside the inequality~(\ref{e.firstterm}), we find that, since $\e < 1$,  
 $$
 \PP \Big( \exists \, x \in \big[ -\eln/4,\eln/4 \big]: \close\big(\mc{L}_n,x,\e/2\big)   \Big) \leq  \e^{k^2 - 3} \cdot 10^{44}  k^{18} \const^{18} \exp \Big\{  4963 k^{7/2}  \const^{5/2} \big( \log \e^{-1} \big)^{5/6} \Big\}  \, .
$$

This completes the proof of~(\ref{t.airynt.tworeduce}) and thus of Proposition~\ref{p.airynt.two}. \qed

\section{Brownian bridge regularity of regular ensembles}\label{s.brownianregularity}

In this section, we will prove our principal results -- notably Theorems~\ref{t.rnbound} and~\ref{t.airytail.ln} -- concerning Brownian regularity of standard bridges derived from curves in regular ensembles.
Throughout the section, 
the value of $\const$ be specified to satisfy the general condition~(\ref{e.constvalue}) and as well a new constraint, equalling
\begin{equation}\label{e.constvalue.br}
 \const  = \max \Big\{ k^{1/3} c_k^{-1/3} \big( 2^{-9/2} - 2^{-5} \big)^{-1/3} \, , \,  36(k^2 -1) \, , \, 16 c_k^{-1} \Big\} \, .
\end{equation}



As we prepare for the proofs of these theorems, we first mention that the parabolic invariance Lemma~\ref{l.parabolicinv} in essence
allows us to reduce immediately to the case where $K=0$ in the theorems' statements (see the actual proofs later in the section for a precise explanation in this regard). As we develop the tools that we will use, we will make this choice.

\begin{proposition}\label{p.vaultbrownian}
Let $\ipdval \in [1,\eln/2)$. For any $i \in \intint{k}$ and $A \subseteq \mc{C}_{0,0}\big([0,\ipdval],\R\big)$ any measurable collection of standard bridges on $[0,\ipdval]$,
$$
 \PP_{\mc{F}} \Big( J^{[0,\ipdval]}(i,\cdot) \in A \Big)   {{\bf 1}}_{\mpg} \leq 2101 \ipdval^{1/2} \const^2 \log \e^{-1} \cdot    \exp \big\{ 54 \ipdval \const^2 \big( \log \e^{-1} \big)^{5/6} \big\} \cdot \mc{B}_{1;0,0}^{[0,\ipdval]} \Big(  B \in A \Big) \, + \, \e^{\const^2/2} \, . 
$$
\end{proposition}
The next result follows since $\e < e^{-1}$ and $\const \geq 1$.
\begin{corollary}\label{c.vaultbrownian}
Let $\ipdval \in [1,\eln/2)$. 
If the set~$A$ in the proposition satisfies  $\mc{B}_{1;0,0}^{[0,\ipdval]} \big(  B \in A \big) \geq \e^{\const^2/2}$, then, for $i \in \intint{k}$, 
$$
 \PP_{\mc{F}} \Big( J^{[0,\ipdval]}(i,\cdot) \in A \Big)   {{\bf 1}}_{\mpg} \leq 
 2102 \ipdval^{1/2} \const^2 \log \e^{-1} \cdot    \exp \big\{ 54 \ipdval \const^2 \big( \log \e^{-1} \big)^{5/6} \big\} \cdot \mc{B}_{1;0,0}^{[0,\ipdval]} \Big(  B \in A \Big) \, . 
$$
\end{corollary}
\noindent{\em Remarks: (1).} Throughout this section, we set the inter-pole distance parameter $\ipd$ equal to $\ipdval$. This assignation reflects our aim in proving Theorem~\ref{t.airytail.ln} of understanding the behaviour of the process $\mc{L}_n^{[K,K+\ipdval]}(k,\cdot)$ on an interval of length~$\ipdval$.\\
\noindent{\em (2).}
We remind the reader that the jump ensemble method's general hypotheses on the parameters~$\e$ and $n$ -- namely~(\ref{e.epsilonupperbound}),~(\ref{e.nlbone}) and~(\ref{e.nlbtwo}) -- remain in force; recall further that an $n$-dependence is suppressed in the notation $J$ for the jump ensemble.

\medskip 
 
\noindent{\bf Proof of Proposition~\ref{p.vaultbrownian}.}
Elements of the pole set $\pole$ are separated from each other by gaps whose distance exceeds $\ipd = \ipdval$; thus, there is at most one element of $\pole \cap [0,\ipdval]$. If this element exists, we label it~$\pp$. When $\mpg$, and thus $\big\{ [-\eln/2,\eln/2] \subseteq [\mfl,\mfr]  \big\}$, occurs, and $\pp$ exists, this element is neither the greatest nor the least member of~$\pole$; in this case, we write $\pp^-$ and $\pp^+$ for the elements of $\pole$ that precede and follow $\pp$ in the increasing order.


On the event that $\pole \cap [0,\ipdval] \not= \emptyset$, set
\begin{equation}\label{e.sigma12}
\sigma_1^2 = \frac{- \pp^- \!\cdot \pp }{\pp - \pp^-} \, \, \, \, \, \textrm{and} \, \, \, \, \, \sigma_2^2 = \frac{(\ipdval - \pp)(\pp^+ - \ipdval)}{\pp^+ - \pp} \, .
\end{equation}

For later use, we further define on the same event
\begin{equation}\label{e.sigma34}
   \sigma_3^2 = \frac{\pp (\ipdval - \pp )}{\ipdval} \, \, \, \, \, \textrm{and} \, \, \, \, \,  \sigma_4^2 = \big(1-\pp \ipdval^{-1}\big)^2 \sigma_1^2 + \pp^2 \ipdval^{-2} \sigma_2^2 \, .
\end{equation}

Recall from Section~\ref{s.morepromising} that $\tent:[\mfl,\mfr] \to \R$ is the $\mc{F}$-measurable {\em tent} map, which piecewise linearly interpolates the pole tops indexed by poles in $P \subset [\mfl,\mfr]$. 

The next lemma provides a Gaussian upper bound on the density of the location of $J(i,x)$ relative to the tent map in the event that this curve drops below the tent map at a given location $x$.
\begin{lemma}\label{l.belowtent}
Let $i \in \intint{k}$. 
\begin{enumerate}
\item
For $x \in \{ 0, \ipdval \}$, write  $h^{\mc{F}}_{1,x}:\R \to [0,\infty)$ for the $\mc{F}$-measurable random variable that, when $\mpg$ occurs, equals
the density under the law $\PP_\mc{F}$ 
of 
$$
J(i,x) - \tent(x) 
\, ; 
$$
(and equals zero when $\neg \, \mpg$ occurs).

Then, for any $s \in (-\infty,0)$, 
$$
 h^{\mc{F}}_{1;x}(s) \cdot {\bf 1}_{\mpg \cap \{ \pole \cap [0, \ipdval ] \not= \emptyset \}} \, \leq \, g_{0,\sigma^2}(s)  \, ,
$$
where $\sigma$ equals $\sigma_1$ or $\sigma_2$ according to whether $x$ equals $0$ or $\ipdval$.
\item
Write $h^{\mc{F}}_{2;0,\ipdval}:\R^2 \to [0,\infty)$ for the $\mc{F}$-measurable random variable that, when $\mpg$ occurs, equals
the density under the law $\PP_\mc{F}$ 
of the pair 
$$
\Big( J(i,0) - \tent(0) \, , \, J(i,\ipdval) - \tent(\ipdval) \Big)
$$
(and equals $(0,0)$ when $\neg \,  \mpg$ occurs). 

Then, for any $(s,t) \in (-\infty,0)^2$, 
$$
 h^{\mc{F}}_{2;0,\ipdval}(s,t) \cdot {\bf 1}_{\mpg \cap \{ \pole \cap [0,\ipdval] \not= \emptyset \}} \, \leq \, g_{0,\sigma_1^2}(s) g_{0,\sigma_2^2}(t) \, .
$$
\end{enumerate}
\end{lemma}
\noindent{\bf Proof.} The first statement follows by the reasoning that proves the second and its proof is omitted.
Recall that, when $\pole$ intersects $[0,\ipdval]$, 
$\pp$ is the unique element of intersection, and $\pp^-$ and $\pp^+$ the adjacent elements of $\pole$. Let $\mc{F}[i;\pp^-,\pp,\pp^+]$
denote the $\sigma$-algebra generated by $\mc{F}$ and the random variables $J(i,x)$ for $x \in \{ \pp^-,\pp,\pp^+ \}$. (These random variables provide extra information only when $\pole$ intersects $[0,\ipdval]$.) The density $h^{\mc{F}}_{2;0,\ipdval}(s,t)$ has a counterpart under the augmented $\sigma$-algebra, and indeed it is sufficient to argue that 
$$
 h^{\mc{F}[i;\pp^-,\pp,\pp^+]}_{2;0,\ipdval}(s,\ipdval) \cdot {\bf 1}_{\mpg \cap \{ \pole \cap [0,\ipdval] \not= \emptyset \}} \, \leq \, g_{0,\sigma_1^2}(s) g_{0,\sigma_2^2}(t) \, ,
$$
since then Lemma~\ref{l.belowtent}(2) will arise by averaging.

Under the law $\PP_{\mc{F}[i;\pp^-,\pp,\pp^+]}$, the processes $J(i,\cdot)$ on $[\pp^-,\pp]$ and $[\pp,\pp^+]$ are conditionally independent. Supposing that the data in $\mc{F}[i;\pp^-,\pp,\pp^+]$ causes $\mpg \cap \{ \pole \cap [0,\ipdval] \not= \emptyset \}$ to occur, it is thus enough to argue that 
\begin{itemize}
\item the conditional density of $J(i,0)- \tent(0)$ at $s \leq 0$ is at most $g_{0,\sigma_1^2}(s)$;
\item and 
the conditional density of $J(i,\ipdval)- \tent(\ipdval)$ at $t \leq 0$ is at most $g_{0,\sigma_2^2}(t)$.
\end{itemize}
These statements are quite straightforward to verify.
Indeed, the conditional law under $\mc{F}[i;\pp^-,\pp,\pp^+]$ of $J(i,0)$ is normal with mean
$\tfrac{\pp}{\pp - \pp^-} J(i,\pp^-) + \tfrac{- \pp^-}{\pp - \pp^-} J(i,\pp)$ and variance $\sigma_1^2$. Note that $J(i,\pp^-) \geq \tent(\pp^-)$ and $J(i,\pp) \geq \tent(\pp)$ since $\pp^-,\pp \in \pole$, and that $\tent$ is affine on the interval between consecutive pole set elements $\pp^-$ and $\pp$; thus, we see that this mean is at least $\tent(0)$. The first bullet point statement follows, and the second is proved in the same fashion. This proves Lemma~\ref{l.belowtent}(2).
\qed

\medskip 


For $H,H' \subseteq \R$, write 
$$
\goodk^{H,H'} = \Big\{ J(i,0) - \tent(0)  \in H \, , \, J(i,\ipdval) - \tent(\ipdval)  \in H' \Big\} \, ,
$$
and abbreviate $\goodk = \goodk^{H,H'}$ when the choice 
$$
H = \big( - \sigma_1 \const \big( \log \e^{-1} \big)^{1/2} , \infty\big) \, \, \, \textrm{and} \, \, \, 
H' = \big( - \sigma_2 \const \big( \log \e^{-1} \big)^{1/2} , \infty\big) 
$$
is made.
By Lemma~\ref{l.belowtent}(1),
\begin{equation}\label{e.goodkest}
 \PP_\mc{F} \big( \goodk^c \big) \cdot {{\bf 1}}_{\mpg \cap \{ \pole \cap [0,\ipdval] \not= \emptyset \}} \, \leq 2 \, (2\pi)^{-1/2} \const^{-1} \big( \log \e^{-1} \big)^{-1/2} \exp \big\{ - \tfrac{1}{2} \const^2 \log \e^{-1} \big\} \leq  \e^{\const^2/2} \, ,
\end{equation}
the latter bound since $\const \geq 1$ and $\e < e^{-1}$.

Note that
\begin{eqnarray}
\PP_\mc{F} \Big( J^{[0,\ipdval]} \in A  \Big) {{\bf 1}}_{\mpg} & \leq &
 \PP_\mc{F} \Big( J^{[0,\ipdval]} \in A  \, , \, \goodk  \Big) {{\bf 1}}_{\mpg \cap  \{ \pole \cap [0,\ipdval] \not= \emptyset \}} \label{e.jkk} \\
  & +  & 
\PP_\mc{F} \Big( J^{[0,\ipdval]} \in A   \Big) {{\bf 1}}_{\mpg \cap  \{ \pole \cap [0,\ipdval] = \emptyset \} }
\, + \, 
\PP_\mc{F} \big( \goodk^c \big) {{\bf 1}}_{\mpg \cap \{ \pole \cap [0,\ipdval] \not= \emptyset \}} \, . \nonumber 
\end{eqnarray}
It is the first term on the right-hand side
that is the most subtle to analyse. 
In order to do so, we introduce another augmentation of the $\sigma$-algebra~$\mc{F}$. Let $\fik$ denote the $\sigma$-algebra generated by $\mc{F}$ and the curve $J(i,\cdot)$
on $[\mfl,0] \cup [\ipdval,\mfr]$. (If either of these intervals is badly specified because the endpoints are out of order, we may treat the interval as the empty-set. However, the occurrence of $\mpg$ 
 will avoid this formal difficulty.)

Using this new device, the first term on the right-hand side 
of~(\ref{e.jkk})
equals
$$
\EE_\mc{F} \Big[ \PP_{\fik} \big(  J^{[0,\ipdval]}(i,\cdot) \in A \big)  {{\bf 1}}_{\goodk} \Big] \cdot {{\bf 1}}_{\mpg \cap  \{ \pole \cap [0,\ipdval] \not= \emptyset \}} \, ;
$$
in adopting this point of view, we are led to ask: `what is the conditional distribution of $J^{[0,\ipdval]}$ under $\PP_{\fik}$?'
For an instance of $\mc{F}$-measurable data that dictates the occurrence of $\mpg \cap \{ \pole \cap [0,\ipdval] \not= \emptyset \}$,  the process $J(i,\cdot)$ on $[0,\ipdval]$ under $\PP_{\fik}$ is distributed as Brownian bridge $B$ on $[0,\ipdval]$ under the law $\mc{B}_{1;J(i,0),J(i,\ipdval)}^{[0,\ipdval]}$ conditionally on $B(\pp) \geq \mc{L}_n(k+1,\pp)$. Thus,
\begin{eqnarray*}
 & & \PP_{\fik} \Big(  J^{[0,\ipdval]}(i,\cdot) \in A \Big) \cdot {\bf 1}_{\mpg \cap \{ \pole \cap [0,\ipdval] \not= \emptyset \}} \\
 & = & \mc{B}_{1;J(i,0),J(i,\ipdval)}^{[0,\ipdval]} \Big(  B^{[0,\ipdval]} \in A \, \Big\vert \,  B(\pp) \geq \mc{L}_n(k+1,\pp)  \Big) \\
 & \leq &  \mc{B}_{1;0,0}^{[0,\ipdval]} \Big(  B \in A \Big) \cdot  \mc{B}_{1;J(i,0),J(i,\ipdval)}^{[0,\ipdval]} \Big(  B(\pp) \geq \mc{L}_n(k+1,\pp)  \Big)^{-1} \, .
\end{eqnarray*}

Note that $B(\pp)$ under $\mc{B}_{1;J(i,0),J(i,\ipdval)}^{[0,\ipdval]}$
is normally distributed with mean 
$\big(1 - \pp \ipdval^{-1}\big)J(i,0) + \pp \ipdval^{-1} J(i,\ipdval)$ and variance $\sigma_3^2$, where recall that $\sigma_3^2 = \pp (\ipdval-\pp) \ipdval^{-1}$. The tent map is affine on $[0,\pp]$ and $[\pp,\ipdval]$, with a slope in each section of absolute value at most $4\eln$: thus,
$$
\tent(0) \geq \mc{L}_n(k+1,\pp) - 4\eln \pp \, \, \, \textrm{and} \, \, \, 
\tent(\ipdval) \geq \mc{L}_n(k+1,\pp) - 4\eln (\ipdval-\pp) \, .
$$

If we set $x_1 =  J(i,0) - \tent(0)$ and $x_2 =  J(i,\ipdval) - \tent(\ipdval)$, 
we see then that this mean is at least
$$
\mc{L}_n(k+1,\pp) - 8 \sigma_3^2 \eln +  \big( 1 -   \pp \ipdval^{-1} \big) x_1 + \pp \ipdval^{-1} x_2 \, . 
$$

Set 
$H_1 = \big( - \sigma_1 \const \big( \log \e^{-1} \big)^{1/2} , 0 \big]$ and 
$I_1 = \big( - \sigma_2 \const \big( \log \e^{-1} \big)^{1/2} ,  0 \big]$ as well as $H_2 = I_2 = [0,\infty)$.

Note then that
\begin{eqnarray*}
 & & 
\EE_\mc{F} \bigg[  \mc{B}_{1;J(i,0),J(i,\ipdval)}^{[0,\ipdval]} \Big(  B(\pp) \geq \mc{L}_n(k+1,\pp)  \Big)^{-1} {{\bf 1}}_{\goodk} \bigg] {{\bf 1}}_{\mpg \cap \{ \pole \cap [0,\ipdval] \not= \emptyset \}} \\
 & = &  
 \sum \EE_\mc{F} \bigg[  \mc{B}_{1;J(i,0),J(i,\ipdval)}^{[0,\ipdval]} \Big(  B(\pp) \geq \mc{L}_n(k+1,\pp)  \Big)^{-1} {{\bf 1}}_{\goodk^{H,I}} \bigg] {{\bf 1}}_{\mpg \cap \{ \pole \cap [0,\ipdval] \not= \emptyset \}} \, ,
\end{eqnarray*}  
where the sum is performed over $(u,v) \in \{1,2\}^2$, and the $(u,v)$-indexed summand, which we call $A_{uv}$, is specified by setting $H = H_u$ and $I = I_v$.

\begin{figure}[ht]
\begin{center}
\includegraphics[height=12cm]{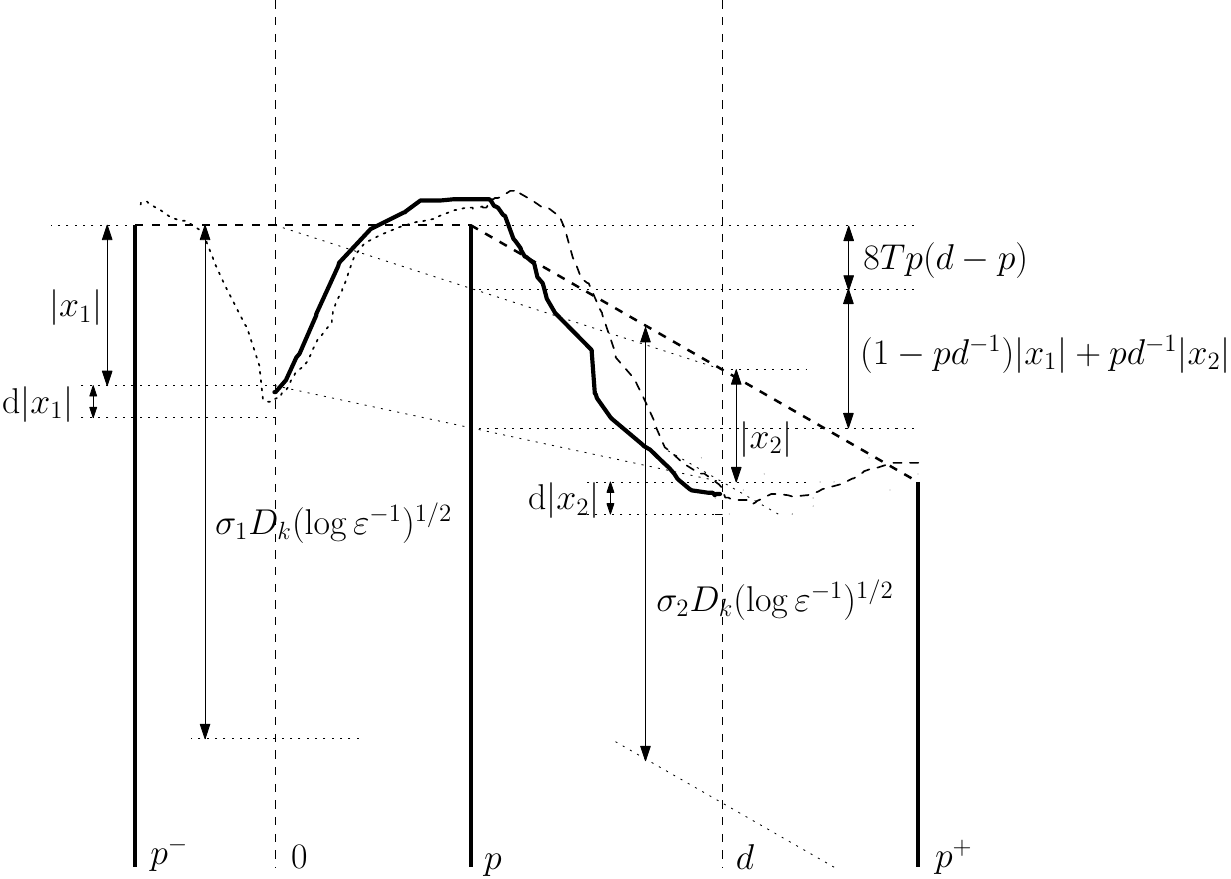}
\caption{This figure discusses a key aspect of the proof of Proposition~\ref{p.vaultbrownian}: the form of the integrand in~(\ref{e.a11ub}). Three vertical solid poles support the dashed tent map. The Gaussian factors  $g_{0,\sigma_1^2}(x_1)$  and $g_{0,\sigma_2^2}(x_2)$  arise from Lemma~\ref{l.belowtent}; the reason for their presence is illustrated by the dotted and dashed jump ensemble curves on $[p^-,p]$ and $[p,p^+]$, which must pay these kinetic costs to reach $\vert x_1 \vert$ and $\vert x_2 \vert$ below tent at $0$ and $\ipdval$. The pole at $p$ rises distance $ 8 \sigma_3^2 \eln  +  \big( 1 - \pp \ipdval^{-1} \big) \vert x_1 \vert + \pp \ipdval^{-1} \vert x_2 \vert$ above its intersection with the line segment joining $(0,x_1)$ and $(\ipdval,x_2)$ (where, in identifying these points, the vertical coordinate is measured relative to the tent map). 
The remaining factor in the integrand is the reciprocal of the probability that Brownian bridge between these two points vaults over the pole at $p$. The factor is at least as large as $\exp \{ O(x_1^2 + x_2^2)  \}$ (if $\ipdval$ is of unit order and $p$ is close to $\ipdval/2$, say), a quantity that is intolerably, polynomially-in-$\e^{-1}$, high when $\vert x_1 \vert \vee \vert x_2 \vert = \Theta\big( (\log \e^{-1})^{1/2} \big)$. However, the factor's dominant term is cancelled by the earlier two Gaussian factors. The picture illustrates this: the thick solid Brownian bridge that vaults over the pole at $p$ has the same first order kinetic costs associated to fluctuation during $[0,p]$ and $[p,\ipdval]$ as the dotted and dashed motions attached to the first two Gaussian factors.}  
\label{f.tentp}
\end{center}
\end{figure}

Using the notation and the statement of Lemma~\ref{l.belowtent}, we may note that 
\begin{eqnarray}
 A_{11} & \leq & \int_{H_1 \times I_1}  h^{\mc{F}}_{2;0,\ipdval}(x_1,x_2)   \nu_{0,\sigma_3^2} \Big( 8  \sigma_3^2 \eln  -  \big( 1 - \pp\ipdval^{-1} \big)x_1 -  \pp \ipdval^{-1} x_2 , \infty \Big)^{-1} \dd x_1 \dd x_2 \nonumber \\
& \leq & \int_{H_1 \times I_1}   g_{0,\sigma_1^2}(x_1) \cdot g_{0,\sigma_2^2}(x_2) \cdot \nu_{0,\sigma_3^2} \Big( 8 \sigma_3^2 \eln  -  \big(1 - \pp \ipdval^{-1} \big)x_1 -  \pp \ipdval^{-1}  x_2 , \infty \Big)^{-1} \, \dd x_1 \dd x_2 \, . \label{e.a11ub}
\end{eqnarray}
The latter integrand here has three factors, two small allies and one large opponent. The third term is unmanageably large in isolation and we will depend on a cancellation of first order kinetic costs between the third term and the first two. See Figure~\ref{f.tentp}. 

Let $N_1$ and $N_2$ denote independent normal random variables of  zero mean and respective variance $\sigma_1^2$ and $\sigma_2^2$. Set $N = \big(1 - \pp \ipdval^{-1} \big)N_1 + \pp \ipdval^{-1} N_2$. Let $V$ equal $N$ in the event that 
$$
N_1 \in \Big( - \sigma_1  \const \big( \log \e^{-1} \big)^{1/2} , 0 \Big) \, \, \,  \textrm{and} \, \, \,  
N_2 \in \Big( - \sigma_2  \const \big( \log \e^{-1} \big)^{1/2} , 0 \Big) \, ; 
$$
in the other case, we may take $V = \infty$. Writing $\EE$ for the expectation associated to these random variables, the expression in line~(\ref{e.a11ub}) equals 
$$
\EE \, \Big[ \, \nu_{0,\sigma_3^2} \big( 8  \sigma_3^2  \eln  - V \, , \, \infty \big)^{-1} \cdot {{\bf 1}}_{(N_1,N_2) \in H_1 \times I_1} \, \Big] \, = \,
\EE \, \Big[ \, \nu_{0,\sigma_3^2} \big( 8  \sigma_3^2  \eln  - V \, , \, \infty \big)^{-1} \cdot {{\bf 1}}_{V < \infty} \, \Big] \, .
$$ 
Further set $W$ equal to $N$ if 
$$
N \in \Big( - \big( \sigma_1 (1 - \pp\ipdval^{-1}) + \sigma_2 \pp\ipdval^{-1} \big) \const \big( \log \e^{-1} \big)^{1/2} , 0 \Big) \, ; 
$$
take $W = \infty$ otherwise. Note then that when $V$ is finite, $W$ equals $V$. Thus,
\begin{eqnarray} 
 & & \EE \Big[ \nu_{0,\sigma_3^2} \big( 8  \sigma_3^2  \eln  - V \, , \, \infty \big)^{-1} \cdot {{\bf 1}}_{V < \infty} \Big] 
  \leq  \EE \Big[ \nu_{0,\sigma_3^2} \big( 8 \sigma_3^2  \eln  - W \, , \, \infty \big)^{-1}  \cdot {{\bf 1}}_{W < \infty} \Big] \nonumber \\ 
 & \leq & \int_{\R}  g_{0,\sigma_4^2}(x)  \nu_{0,\sigma_3^2} \Big( 8\eln  \sigma_3^2 -  x , \infty \Big)^{-1} {{\bf 1 }}_{\big[- \big( \sigma_1 (1 - \pp \ipdval^{-1} ) + \sigma_2 \pp \ipdval^{-1} \big) \const ( \log \e^{-1} )^{1/2} , 0 \big]} \dd x  \, , \label{e.foursigma}
\end{eqnarray} 
where $\sigma_4^2$ equals the variance of $N$. We pause to collect some bounds satisfied by the four $\sigma^2$ from~(\ref{e.sigma12}) and~(\ref{e.sigma34}) that we are using. 
\begin{lemma}\label{l.sigmabounds}
We have that
\begin{enumerate}
 \item $\sigma^2_3/8 \leq \sigma^2_4 \leq \sigma^2_3$,
 \item  and
$\sigma_1 \big( 1 - \pp \ipdval^{-1} \big) + \sigma_2 \pp \ipdval^{-1}  \leq \sigma_3 \Big( \big( 1 - \pp \ipdval^{-1} \big)^{1/2} + \big( \pp \ipdval^{-1}  \big)^{1/2} \Big)$.
 \end{enumerate}
\end{lemma}
\noindent{\bf Proof: (1).} If $\pp \leq  \ipdval/2$, then $\tfrac{-\pp^-}{\pp - \pp^-} \geq 1/2$ and thus $\sigma_1^2 \geq \pp/2$ so that $\sigma_4^2 \geq \pp/8$. If $\pp \geq  \ipdval/2$, then $\sigma_4^2 \geq (\ipdval - \pp)/8$. Thus, $\sigma_4^2 \geq \sigma_3^2/8$. 
 Since $\sigma_1^2 \leq \pp$ and $\sigma_2^2 \leq \ipdval - \pp$,
we have that $\sigma_4^2 \leq \sigma_3^2$. 

\noindent{\bf (2)}. This follows from  $0 \in [\pp^-,\pp]$ and $\ipdval \in [\pp,\pp^+]$. \qed

\medskip

Setting
$$
 I = \bigg[ \, 0 \, , \,
\sigma_3 \Big( \big(1 - \pp\ipdval^{-1} \big)^{1/2} + \big(\pp\ipdval^{-1} \big)^{1/2} \Big) \,
 \const ( \log \e^{-1} )^{1/2} \, \bigg] \, ,
$$
and using the notation $-I = \{ -x: x \in I \}$,
the quantity in~(\ref{e.foursigma}) is seen by means of Lemma~\ref{l.sigmabounds}(2) to be at most
\begin{eqnarray*} 
 &  & 
\int_{\R}  g_{0,\sigma_4^2}(x)  \nu_{0,\sigma_3^2} \Big( 8\eln \sigma_3^2 -  x , \infty \Big)^{-1} {{\bf 1 }}_{-I}(x) \dd x = 
\int_I  g_{0,\sigma_4^2}(x)  \nu_{0,\sigma_3^2} \Big( 8\eln \sigma_3^2 + x , \infty \Big)^{-1}  \dd x  \\
 & \leq & 
\int_I  g_{0,\sigma_4^2}(x)  \nu_{0,\sigma_3^2} \Big( 8\eln \sigma_3^2 +     \sigma_3 + x , \infty \Big)^{-1}  \dd x 
\end{eqnarray*}
The addition of the $+ \sigma_3$ term in the last expression enables us to apply the Gaussian tail lower bound Lemma~\ref{l.gaussiantail} with $t = (8 \eln \sigma_3^2 + \sigma_3 + x) \sigma_3^{-1}$ safe in the knowledge that $t \geq 1$; the last integral is thus seen to be at most
\begin{eqnarray*}  
 &  & \int_I  g_{0,\sigma_4^2}(x)  g_{0,\sigma_3^2}(x + 8\eln \sigma_3^2 + \sigma_3)^{-1} 2 \big( x + 8 \eln \sigma_3^2 + \sigma_3 \big) \sigma_3^{-1}  \dd x \\
 & \leq & \int_I  g_{0,\sigma_4^2}(x)  g_{0,\sigma_3^2}(x + 9  \ipdval^{1/2} \eln \sigma_3)^{-1} 2 \big( x + 9 \ipdval^{1/2} \eln \sigma_3  \big) \sigma_3^{-1}  \dd x \\
 & \leq & \int_I  g_{0,\sigma_4^2}(x)  g_{0,\sigma_3^2}(x)^{-1} \exp \big\{ 9 \ipdval^{1/2} x \eln  \sigma_3^{-1}  + 9^2 \ipdval \eln^2/2 \big\} \cdot 2 \big( x + 9\ipdval^{1/2} \eln \sigma_3 \big) \sigma_3^{-1}  \dd x \\
 & \leq & \int_{\R}  \exp \big\{ 9 \ipdval^{1/2} \eln  x \sigma_3^{-1} + 41 \ipdval \eln^2  \big\} \cdot  2 \big( x + 9\ipdval^{1/2} \eln \sigma_3 \big) \sigma_4^{-1} {{\bf 1}}_{\big[0, \sqrt{2} \sigma_3  \, \const ( \log \e^{-1} )^{1/2} \big]}(x) \dd x \\
 & \leq &    \sqrt{2} \sigma_3  \, \const ( \log \e^{-1} )^{1/2} \exp \big\{ 9  \ipdval^{1/2} \eln   \sqrt{2}  \, \const ( \log \e^{-1} )^{1/2} + 41 \ipdval \eln^2  \big\} \\
  & & \qquad \qquad \qquad \qquad \qquad \qquad \qquad \qquad \times \, \, \,  2 \big( \sqrt{2} \sigma_3  \, \const ( \log \e^{-1} )^{1/2} + 9\ipdval^{1/2} \eln \sigma_3 \big) \sigma_4^{-1} \\
 & \leq &   4  \, \const ( \log \e^{-1} )^{1/2} \exp \big\{ 13  \ipdval^{1/2} \eln     \, \const ( \log \e^{-1} )^{1/2} + 41 \ipdval \eln^2  \big\} \cdot \big( 3  \ipdval^{1/2} \const ( \log \e^{-1} )^{1/2} + 18 \ipdval^{1/2} \eln \big)  \\
 & \leq &   84  \ipdval^{1/2} \const^2  \log \e^{-1}  \cdot \exp \big\{ 54  \ipdval     \const^2 ( \log \e^{-1} )^{5/6}    \big\} \, .  
\end{eqnarray*}
The first displayed inequality depended on $\sigma_3 \leq \ipdval^{1/2}$ and $\ipdval \wedge \eln \geq 1$ while the third made use of
$g_{0,\sigma_4^2}(x) g_{0,\sigma_3^2}(x)^{-1} \leq \sigma_3 \sigma_4^{-1}$, which is a consequence of $\sigma_4 \leq \sigma_3$. The fifth used $\sigma_3 \leq \ipdval^{1/2}$ and $\sigma_4 \geq 8^{-1/2 }\sigma_3$. That $\eln = \const \big( \log \e^{-1} \big)^{1/3}$, $\ipdval \geq 1$ and $\e < e^{-1}$ were applied to obtain the final inequality.

We find then that 
\begin{equation}\label{e.a11ubfinal}
 A_{11} \leq  84 \ipdval^{1/2} \const^2  \log \e^{-1}  \cdot \exp \big\{ 54   \ipdval    \const^2 ( \log \e^{-1} )^{5/6}    \big\} \, . 
\end{equation}

Note that 
\begin{eqnarray*}
A_{12} & = & 
 \EE_\mc{F} \bigg[  \mc{B}_{1;J(i,0),J(i,\ipdval)}^{[0,\ipdval]} \Big(  B(\pp) \geq \mc{L}_n(k+1,\pp)  \Big)^{-1} {{\bf 1}}_{\goodk^{H_1,I_2}} \bigg] {{\bf 1}}_{\mpg \cap \{ \pole \cap [0,\ipdval] \not= \emptyset \}} \\
  & \leq & 
 \EE_\mc{F} \bigg[  \mc{B}_{1;J(i,0),\tent(\ipdval)}^{[0,\ipdval]} \Big(  B(\pp) \geq \mc{L}_n(k+1,\pp)  \Big)^{-1} {{\bf 1}}_{\goodk^{H_1,[0,\infty)}} \bigg] {{\bf 1}}_{\mpg \cap \{ \pole \cap [0,\ipdval] \not= \emptyset \}} \\
 & \leq & \int_{H_1}  h^{\mc{F}}_{1;0}(x)   \nu_{0,\sigma_3^2} \Big( 8\sigma_3^2 \eln  -  \big( 1 - \pp\ipdval^{-1} \big) x  , \infty \Big)^{-1} \dd x \\
 & \leq & \int_{H_1}  g_{0,\sigma_1^2}(x)   \nu_{0,\sigma_3^2} \Big( 8  \sigma_3^2 \eln  -  \big(1 - \pp\ipdval^{-1}\big)x  , \infty \Big)^{-1} \dd x \\
 & \leq & Z[2,I_1]^{-1} \int_{H_1 \times I_1}  g_{0,\sigma_1^2}(x_1)  g_{0,\sigma_2^2}(x_2)   \nu_{0,\sigma_3^2} \Big( 8 \sigma_3^2 \eln  -  \big(1 - \pp\ipdval^{-1} \big)x_1 -  \pp \ipdval^{-1}  x_2  , \infty \Big)^{-1} \dd x \, .\\
\end{eqnarray*}
Lemma~\ref{l.belowtent}(1) was used in the penultimate inequality. In the final line, the quantity~$Z[2,I_1]$ denotes  $\int_{I_1}  g_{0,\sigma_2^2}(x)  \dd x$, which is the probability that $N_2$ assumes a value in $I_1$; $\const \geq 1$ and $\e < e^{-1}$ imply that 
$$
Z[2,I_1] \, \geq \, \frac{1}{2} \,  - \, (2\pi)^{-1/2} \const^{-1} \big( \log \e^{-1} \big)^{-1/2} \exp \big\{ - \tfrac{1}{2} \const^2 \log \e^{-1} \big\} \, \geq \, \frac{1}{4} \, \, . 
$$
Thus, $A_{12}$ is at most four times the upper bound on $A_{11}$ appearing in~(\ref{e.a11ub}). Thus, $A_{12}$ satisfies the bound~(\ref{e.a11ubfinal}) when the right-hand side is quadrupled. By a similar analysis, $A_{21}$ may be shown to satisfy the same bound. 

Finally, note that
\begin{eqnarray*} 
A_{22} & = & 
 \EE_\mc{F} \bigg[  \mc{B}_{1;J(i,0),J(i,\ipdval)}^{[0,\ipdval]} \Big(  B(\pp) \geq \mc{L}_n(k+1,\pp)  \Big)^{-1} {{\bf 1}}_{\goodk^{H_2,I_2}} \bigg] {{\bf 1}}_{\mpg \cap \{ \pole \cap [0,\ipdval] \not= \emptyset \}} \\
  & \leq & 
 \EE_\mc{F} \bigg[  \mc{B}_{1;\tent(0),\tent(\ipdval)}^{[0,\ipdval]} \Big(  B(\pp) \geq \mc{L}_n(k+1,\pp)  \Big)^{-1} {{\bf 1}}_{\goodk^{[0,\infty),[0,\infty)}} \bigg] {{\bf 1}}_{\mpg \cap \{ \pole \cap [0,\ipdval] \not= \emptyset \}} \\
 & \leq &   \nu_{0,\sigma_3^2} \Big( 8 \sigma_3^2 \eln  \,  , \, \infty \Big)^{-1} \, ;
\end{eqnarray*} 
were $\sigma_3$ known to be bounded away from zero, we would find this last quantity to be bounded and thus that it would easily satisfy the bound~(\ref{e.a11ubfinal}). Lacking this knowledge, we instead bound the quantity above by 
$$  
Z[1,H_1]^{-1}  Z[2,I_1]^{-1} \int_{H_1 \times I_1}  g_{0,\sigma_1^2}(x_1)  g_{0,\sigma_2^2}(x_2)   \nu_{0,\sigma_3^2} \Big( 8 \sigma_3^2 \eln  -  \big(1 - \pp \ipdval^{-1} \big) x_1 - \pp \ipdval^{-1} x_2 \, , \, \infty \Big)^{-1} \dd x \, ,
$$
where $Z[1,H_1] = \int_{H_1}  g_{0,\sigma_1^2}(x)  \dd x \geq 1/4$. Thus, $A_{22}$ is at most a multiple of sixteen of the expression in~(\ref{e.a11ub}) and thus also of the right-hand side of~(\ref{e.a11ubfinal}).

We find then that
\begin{eqnarray*}
 & & 
 \PP_\mc{F} \Big( J^{[0,\ipdval]} \in A  \, , \, \goodk \Big) {{\bf 1}}_{\mpg  \cap  \{ \pole \cap [0,\ipdval] \not= \emptyset \} } \\
  & \leq  &  \big( 1 + 4 + 4 + 16 \big) \cdot  84  \ipdval^{1/2} \const^2  \log \e^{-1}  \cdot \exp \big\{ 54   \ipdval    \const^2 ( \log \e^{-1} )^{5/6}    \big\}  \cdot \mc{B}_{1;0,0}^{[0,\ipdval]} \Big(  B \in A \Big) \, .
\end{eqnarray*}

The second term on (\ref{e.jkk})'s right-hand side equals
$$
\PP_\mc{F} \Big( J^{[0,\ipdval]} \in A   \Big) {{\bf 1}}_{\mpg \cap  \{ \pole \cap [0,\ipdval] = \emptyset \} } =  \mc{B}_{1;0,0}^{[0,\ipdval]} \Big(  B \in A \Big) {{\bf 1}}_{\mpg \cap  \{ \pole \cap [0,\ipdval] = \emptyset \} } \, .
$$

Applying the last two estimates along with~(\ref{e.goodkest}) to~(\ref{e.jkk}),
$$
\PP_\mc{F} \Big( J^{[0,\ipdval]} \in A  \Big) {{\bf 1}}_{\mpg} \leq 
 2101 \ipdval^{1/2} \const^2 \log \e^{-1} \cdot    \exp \big\{ 54 \ipdval \const^2 \big( \log \e^{-1} \big)^{5/6} \big\} \cdot \mc{B}_{1;0,0}^{[0,\ipdval]} \Big(  B \in A \Big) \, + \, \e^{\const^2/2} \, .
$$

This completes the proof of Proposition~\ref{p.vaultbrownian}. \qed

\begin{proposition}\label{p.lnbrownian}
Let $\ipdval \geq 1$, and suppose that $\e \in \big(0, e^{-2^3 \ipdval^3 \const^{-3}} \big)$. 
For any $i \in \intint{k}$ and any measurable standard bridge collection $A \subseteq \mc{C}_{0,0}\big([0,\ipdval],\R \big)$ that satisfies $\mc{B}_{1;0,0}^{[0,\ipdval]} \big(  B \in A \big) \geq \e^{\const^2/2}$,
$$
\PP \Big( \mc{L}_n^{[0,\ipdval]} \big( i, \cdot  \big) \in A  \Big)  \leq 
  2103 \, \ipdval^{1/2}    \exp \Big\{ 4028 \, \ipdval^2 k^{7/2}  \const^2 \big( \log \e^{-1} \big)^{5/6} \Big\} \cdot \mc{B}_{1;0,0}^{[0,\ipdval]} \Big(  B \in A \Big)  \, .
$$ 
\end{proposition}
\noindent{\bf Proof.}
Note that 
\begin{eqnarray*}
& & \PP_\mc{F} \Big( \mc{L}_n^{[0,\ipdval]} \big( i, \cdot  \big) \in A  \Big) \cdot {{\bf 1}}_{\mpg} 
  = 
\PP_\mc{F} \Big( J^{[0,\ipdval]}\big( i, \cdot  \big) \in A \, \Big\vert \,   \test_4(J) = 1 \Big) \cdot {{\bf 1}}_{\mpg} \\
 & \leq & 
 \frac{\PP_\mc{F} \Big( J^{[0,\ipdval]}\big( i, \cdot  \big) \in A \Big)}{\PP_\mc{F} \big( \,   \test_4(J) = 1 \big)} \cdot {{\bf 1}}_{\mpg} \, .
 \end{eqnarray*}

Recall that the inter-pole distance parameter $\ipd$ has been set equal to $\ipdval$. Applying Proposition~\ref{p.vaultsuccess} and Corollary~\ref{c.vaultbrownian}, we find that, 
provided that $\mc{B}_{1;0,0}^{[0,\ipdval]} \big(  B \in A \big) \geq \e^{\const^2/2}$, the last quantity is at most
$$
  \exp \Big\{  3973 k^{7/2} \ipdval^2 \const^2 \big( \log \e^{-1} \big)^{2/3} \Big\} 
  \cdot 2102 \ipdval^{1/2} \const^2 \log \e^{-1} \cdot    \exp \Big\{ 54 \ipdval \const^2 \big( \log \e^{-1} \big)^{5/6} \Big\} \cdot \mc{B}_{1;0,0}^{[0,1]} \Big(  B \in A \Big) \, .
$$
(The hypothesis $\ipdval \leq \eln/2$ of Corollary~\ref{c.vaultbrownian} holds since $\e \leq e^{-2^3 \ipdval^3 \const^{-3}}$.)

By Lemma~\ref{l.glub.new}, we have $\PP \big( \mpg^c  \big) \leq \e^{\const^2/2}$ provided that $\const \geq 16 c_k^{-1}$. Using this bound, we find that
\begin{eqnarray*}
& & \PP \Big( \mc{L}_n^{[0,\ipdval]} \big( i, \cdot  \big) \in A  \Big)  
 \leq  \PP \Big( \mc{L}_n^{[0,\ipdval]} \big( i, \cdot  \big) \in A  \, , \, \mpg \Big) \, + \, \PP \big( \mpg^c \big) \\
& \leq & \EE \bigg[ \, \PP_\mc{F} \Big( \mc{L}_n^{[0,\ipdval]} \big( i, \cdot  \big) \in A  \Big) \cdot {{\bf 1}}_\mpg \, \bigg] \, + \,  \e^{\const^2/2} \\
 & \leq &    2102 \ipdval^{1/2} \const^2 \, \log \e^{-1} \cdot    \exp \Big\{ 4027 \ipdval^2 k^{7/2}  \const^2 \big( \log \e^{-1} \big)^{5/6} \Big\} \cdot \mc{B}_{1;0,0}^{[0,\ipdval]} \Big(  B \in A \Big)\, + \,  \e^{\const^2/2} \, , 
\end{eqnarray*}
where in the latter inequality, we applied the bound just derived, alongside $\ipd \geq 1$, $\const \geq 1$ and $\e < e^{-1}$.
Noting that $\e < e^{-2^{3/2} \const^{-3}}$ implies that $\const^2 \log \e^{-1} \leq \exp \big\{  \const^2 ( \log \e^{-1} )^{5/6} \big\}$ (since $\const^4 (\log \e^{-1})^{5/3}/2$ lies in the interval between these two quantities), 
we obtain Proposition~\ref{p.lnbrownian}. \qed

\medskip

Equipped with Proposition~\ref{p.lnbrownian},
we are now able to conclude this section with the proofs of our principal results on Brownian regularity: 
in order, Theorem~\ref{t.airytail.ln}; Theorem~\ref{t.airytail}; and Theorem~\ref{t.rnbound}.

\noindent{\bf Proof of Theorem~\ref{t.airytail.ln}.}
First it may be useful to recall that the general 
hypotheses that govern the jump ensemble method are in force. That is, the strengthening~(\ref{e.constvalue.br}) of~(\ref{e.constvalue}) holds; and 
(\ref{e.epsilonupperbound}), 
(\ref{e.nlbone}) and~(\ref{e.nlbtwo}) also hold. 
Note that the hypotheses of Theorem~\ref{t.airytail.ln} enforce these requirements when we take, as we will, $\e = a$. In particular, the condition~(\ref{e.nlbone}) is validated by the hypothesis~(\ref{e.alowerbound}) that stipulates an $n$-dependent lower bound on $a$ of the form $a \geq \exp \big\{ - O(1) n^{\phimac_1/4 \, \wedge \, \phimac_2/4 \, \wedge \, \phimac_3/2} \big\}$, where the $O(1)$ carries a $k$-dependence.

By Lemma~\ref{l.parabolicinv}, if we prove Theorem~\ref{t.airytail.ln} with  $K$ taken equal to zero and the appearance of $\const$ replaced by $2^{-2/3} \const$, then we will succeed in proving the theorem itself, with $\const$ assuming its presently specified value given by~(\ref{e.constvalue.br}). This is
because the thus specified $\const$ regarded as a function of the regular ensemble parameter $\rsc > 0$ satisfies $\const(\rsc/2) \leq 2^{1/3} \const(\rsc)$. Since $\const^2 \geq 2$, we may apply Proposition~\ref{p.lnbrownian} with $\e = a$ and  $i = k$; since $4028 \cdot 2^{2/3} \leq 6395$, we then arrive at the desired result. \qed

\noindent{\bf Proof of Theorem~\ref{t.airytail}(1).}
Take $n = \infty$ and $\mc{L}_n = \mc{L}$ in Theorem~\ref{t.airytail.ln}. The sought result emerges when we set its parameters according to 
$\beta_k = 6395  k^{7/2}  \const^2$ and $\gamma_k =  e^{-1} \wedge (17)^{-1/k} C_k^{-1/k} \const^{-1}$; since  $\const \geq 2^{5/3}$ implies that  $\exp \big\{ -2^5 \ipdval^3 \const^{-3} \big\} \geq \exp \big\{ - \ipdval^3 \big\}$, and recalling~(\ref{e.formere}) and~(\ref{e.littlec}), we obtain the stated bounds on the growth and decay rates of the two just specified quantities. \qed

\noindent{\bf Proof of Theorem~\ref{t.airytail}(2).} Harmlessly take $K =0$ in view of the stationarity of the Airy line ensemble. For $s > 0$, take $A \subset  \mc{C}_{0,0}\big( [0,\ipdval] , \R \big)$ equal to the collection  of standard bridges $f$ on $[0,\ipdval]$ such that $\sup_{x \in [0,\ipdval]} \vert f(x) \vert > s \ipdval^{1/2}$. By Lemma~\ref{l.maxfluc},
$\mc{B}_{1;0,0}^{[0,\ipdval]} ( A ) \in  e^{-2s^2} \cdot [1,2]$.
We may thus apply Theorem~\ref{t.airytail}(1) with $a$ having the value $r e^{-2s^2}$ for $r$ lying somewhere in the interval $[1,2]$. The hypothesis $r e^{-2s^2} \leq \gamma_k \wedge e^{-\ipdval^3}$ being valid in view of $s \geq 2^{-1/2} \big( (\log \gamma_k^{-1}) \vee \ipdval^3 \,  + \,  \log 2 \big)^{1/2}$, we find that
$$
 \PP \Big(   \mc{L}^{[0,\ipdval]} ( k, \cdot ) \in A  \Big) \leq r e^{-2s^2} \cdot 2103 \, \ipdval^{1/2} \exp  \bigg\{ \beta_k \Big(   2s^2 - \log r  \Big)^{5/6} 
 \bigg\} \, .
$$
The right-hand side is at most $4206 \, e^{-2s^2} \exp \big\{ 2^{5/6} \beta_k s^{5/3} \big\} = 4206  \exp \big\{ - 2s^2 \big( 1 - 2^{-1/6} \beta_k s^{-1/3} \big) \big\}$. \qed

\noindent{\bf Proof of Theorem~\ref{t.rnbound}.}
The ensemble $\mc{L}:\N \times \R \to \R$ in the theorem's statement is formed by parabolically curving the Airy line ensemble -- see~(\ref{e.lairy}). 
We claim that there exist positive $\rsc$ and $\rsC$ such that, for any  $\bar\phimac \in (0,\infty)^3$,   $\mc{L}$ is a $\big(\bar\phimac,\rsc,\rsC\big)$-regular ensemble. 
Since the curve cardinality $n$ of $\mc{L}$ is infinite, this statement is technically an abuse of the meaning of Definition~\ref{d.regularsequence}.
Working however with that definition in the case $n = \infty$, we see that conditions involving $n$ disappear, whatever the value of $\bar\phimac$, making this vector parameter redundant.
Anyway, our task is to check that $\mc{L}$ is Brownian Gibbs, as well as regular. The Brownian Gibbs property is verified due to Proposition~\ref{p.airy}.
The condition $\rmreg(1)$ holds because $\xnmac$ equals $\infty$;
$\rmreg(2)$ and $\rmreg(3)$ with $\phimac_2, \phimac_3 > 0$ arbitrary are implied by \cite[(18)]{AiryLE}.
 That the value of the integral in Theorem~\ref{t.rnbound} is independent of $K \in \R$ follows from the stationarity of the Airy line ensemble. (In \cite{CorwinSun}, the stronger statement that the Airy line ensemble is ergodic with respect to horizontal shifts is proved.)

We will prove  Theorem~\ref{t.rnbound} with the setting
$\alpha_k =  \ctemp \const^{-2/5} \ipdval^{-12/5} k^{-21/5}$ where the constant $\ctemp$ equals $\tfrac{1}{42378}$.
Since we take the parameter $n$ equal to infinity, the jump ensemble method hypotheses~(\ref{e.nlbone})
and~(\ref{e.nlbtwo}) are vacuously satisfied. The parameter $\e > 0$ will be reset countably often, always consistently with the hypothesis~(\ref{e.epsilonupperbound}).

By the stationarity of the Airy line ensemble, it is enough to prove Theorem~\ref{t.rnbound} with  $K$ taken equal to zero. 
Write $\mu = \mc{B}_{1;0,0}^{[0,\ipdval]}$, and $\xi$ for  the distribution of  $\mc{L}^{[0,\ipdval]} \big( k, \cdot  \big)$ under~$\PP$; abbreviate $f = f_k$. For $m \in \N$, set $A_m \subseteq \mc{C}_{0,0}\big([0,\ipdval],\R \big)$ equal to $f^{-1}(m,m+1]$.  On $A_m$, $f \geq m$, and thus $\xi(A_m) > m \cdot \mu(A_m)$. We seek to specify $\e > 0$ so that  
$$
 2103 \ipdval^{1/2}    \exp \Big\{ 4028 \ipdval^2 k^{7/2}  \const^2 \big( \log \e^{-1} \big)^{5/6} \Big\} \, = \, m \, ;
$$
such a choice would entail that
$$
 \e =  \exp \bigg\{ - \Big( \tfrac{1}{4028} \const^{-2} \ipdval^{-2} k^{-7/2} \log \big(   m/2103 \cdot \ipdval^{-1/2} \big) \Big)^{6/5} \bigg\} \, .
$$
Since we demand of $\e > 0$ that it satisfy the upper bound~(\ref{e.epsilonupperbound}) as well as the upper bound hypothesised by Proposition~\ref{p.lnbrownian}, we specify $\e = \e_m > 0$ in this way only if $m 
 \geq \tau_k$ where $\tau_k = 2103  \ipdval^{1/2}   \exp \big\{ 4028 \ipdval^2 k^{7/2}  \const^2 ( \log q_k^{-1} )^{5/6} \big\}$; here, we set $q_k$ equal to  the minimum of the right-hand of~(\ref{e.epsilonupperbound}) and the quantity $\exp \big\{ - 2^3 \ipdval^3 \const^{-3} \big\}$.

When $m$ satisfies this bound, Proposition~\ref{p.lnbrownian} implies that $\mu(A_m) < \e^{\const^2/2}$. 
We find then that
\begin{eqnarray*}
 & & \int_{\mc{C}_{0,0}([0,\ipdval],\R)} \exp \Big\{ \ctemp  \const^{-2/5}  \ipdval^{-12/5}  k^{-21/5} \big( \log f \big)^{6/5} \Big\} \dd \mu \\
 & \leq & \sum_{m=0}^\infty \mu(A_m) \exp \Big\{  \ctemp  \const^{-2/5} \ipdval^{-12/5} k^{-21/5} \big( \log (m+1) \big)^{6/5} \Big\}  \\
  &  \leq &   \exp \Big\{  \ctemp \const^{-2/5} k^{-21/5} \big( \log (\tau_k+1) \big)^{6/5} \Big\}  \\
  & &   + \,\, 
  \sum_{m \geq \tau_k}  \exp \bigg\{ - \tfrac{1}{2} \const^2 \Big( \tfrac{1}{4028} \ipdval^{-2} \const^{-2} k^{-7/2} \log \big(   m/2103 \cdot \ipdval^{-1/2} \big) \Big)^{6/5} \\
& & \qquad \qquad \qquad \qquad \qquad \qquad 
  \, + \, \, \,  \ctemp  \const^{-2/5} \ipdval^{-12/5} k^{-21/5} \big( \log (m+1) \big)^{6/5} \bigg\}  \, .
\end{eqnarray*}
Recalling that $\ctemp$ equals $(42378)^{-1}$, the last sum is seen to be finite. \qed

\appendix

\chapter{Properties of regular Brownian Gibbs ensembles} 
  
 In this appendix, we establish some important properties of regular Brownian Gibbs ensembles.
 
  Proposition~\ref{p.lereg} is the result that makes the theory of such ensembles applicable to the Brownian last passage percolation setting, showing that scaled ensembles from that setting verify the regular ensemble axioms. 
In Section~\ref{s.scaledbrlpp}, we present the proof of this proposition.

   The higher curve index one-point lower tail bound Proposition~\ref{p.othercurves} is necessary to set up the jump ensemble method. Section~\ref{s.othercurves}  proves this proposition by employing a different Brownian Gibbs argument, which is a variant of one used in the KPZ line ensemble construction~\cite{KPZLE}.
 
Finally, in Section~\ref{s.collapsenearinfinity}, we furnish a proof of 
Proposition~\ref{p.collapsenearinfinity}, which asserts that regular ensemble curves decay rapidly very far from the origin. 
  
\section{Scaled Brownian LPP line ensembles are regular}\label{s.scaledbrlpp}

The aim of this section is to prove Proposition~\ref{p.lereg}.

Recall the Brownian LPP ensemble $L_n$ and its equality in law with Dyson's Brownian motion (in Proposition~\ref{p.brlppdbm}), 
  as well as the relation~(\ref{e.scl}) that specifies the scaled counterpart $\mc{L}^\scal_n$ of this ensemble.

\begin{lemma}\label{l.onepointupperlower}
There exist constants $C,c>0$ such that for every $n \in \N$ it is the case that
\begin{enumerate}
\item
for $s \in [0,2^{1/2}n^{1/3}]$,
$$
\PP \big( \mc{L}_n^\scal (1,0) \leq - s \big) \leq C \exp \big\{ - c s^{3/2} \big\} \, ,
$$
\item and, for $s \geq 0$,
$$
\PP \big( \mc{L}_n^\scal (1,0) \geq  s \big) \leq C \exp \big\{ - c s^{3/2} \big\} \, .
$$
\end{enumerate}
\end{lemma}
\noindent{\bf Proof. (1).}
By Proposition~\ref{p.brlppdbm} and Brownian scaling, $L_n(1,n)$ is equal in law to $2n L_n\big(1,(4n)^{-1} \big)$ and thus by Proposition~\ref{p.hbmgrabiner} to the {\rm GUE} top eigenvalue multiple $2n \lambda_n \big(1,(4n)^{-1}\big)$.
By (\ref{e.scl}), we see that
$$
\PP \big( \mc{L}^\scal_n(1,0) \leq - s \big)
=  
\PP \Big( L_n(1,n)  \leq 2n - 2^{1/2} n^{1/3} s \Big)
=
\PP \Big( \lambda_n\big(1,(4n)^{-1} \big) \leq 1 - 2^{-1/2} n^{-2/3}s \Big) \, .
$$
By (\ref{e.ledoux}), this quantity is at most $C' \exp \big\{ - 2^{3/2} s^{3/2} c' \big\}$
when $s \in [2^{1/2},2^{1/2}n^{1/3}]$. We set $C = C'$ and $c = 2^{3/2}c'$. Lemma~\ref{l.onepointupperlower}(1) follows by an increase if necessary in the value of $C$ in order to accomodate $s \in [0,2^{1/2})$.

\noindent{\bf (2).} 
We obtain this result
 with $C = \hat{C}$
and $c = 2^{-3/2} \hat{c}$ by applying Aubrun's bound~(\ref{e.aubrun}), using~(\ref{e.scl}) as in the first proof. \qed

\medskip
\medskip

\noindent{\bf Proof of Proposition~\ref{p.lereg}.}  We must show that each ensemble $\edgedbm_{n}$,  $n \in \N$, satisfies $\rmreg$ $(1)$, $(2)$ and $(3)$.

Condition $\rmreg(1)$ is satisfied since the left endpoint of $\edgedbm_n$ equals $- n^{1/3}/2$.

 Writing $X \eqdist Y$ to denote that the random variables $X$ and $Y$ are equal in law, note that
\begin{eqnarray}
 \dbm_{n} \big( 1, n + 2 n^{2/3} x \big) 
& \eqdist & \Big( 1 + 2n^{-1/3} x \Big)^{1/2}
 \dbm_{n} \big( 1, n \big) \label{e.disteq} \\
& \eqdist & \big( 1 + 2n^{-1/3} x \big)^{1/2} \Big( 2 n \,
 + \, 2^{1/2} n^{1/3}  \edgedbm_{n} \big(  1, 0 \big) \Big) \, . \nonumber
\end{eqnarray}

Introducing the variable $\phi  = n^{-1/3}x$, the latter expression may be written
\begin{eqnarray*}
& & 
  \Big( 1 + \phi - \phi^2/2 + O(\phi^3) \Big)  \Big( 2 n \,
 + \, 2^{1/2} n^{1/3}  \edgedbm_{n} ( 1, 0 ) \Big) \\
& = &  2n \, + \, 2n^{2/3}x \, - \, x^2 n^{1/3}   \, + \, n \cdot O \big( \phi^3 \big) \, +  \, 2^{1/2} n^{1/3}  \edgedbm_{n} ( 1, 0 )  \Big( 1 + O(\phi) \Big) \, ,
\end{eqnarray*}
where the big-$O$ notation implies a bounded factor associated to the term in question on any given compact interval of $\phi$-values that contains zero.

In light of~(\ref{e.scl}), we see that, if  $\vert x \vert \leq \smallc n^{1/3}$ for a small constant $\smallc > 0$, then $\edgedbm_{n} \big( 1 , x \big)$
is equal in law to a random variable that satisfies
$$
 - 2^{-1/2} x^2 + n^{2/3}  O \big(\phi^3 \big) + \edgedbm_{n}\big(1,0\big) \big( 1 + O(\phi) \big) \, .
$$

 When $\vert x \vert \leq n^{1/9}$, we have that  $\vert \phi \vert \leq n^{-2/9}$, in which case
the displayed random variable is $- 2^{-1/2} x^2 +  \edgedbm_{n}\big(1,0\big) \big( 1 + O(n^{-2/9}) \big) + O(1)$.  We obtain the one-point upper tail $\rmreg(3)$ with $\phimac_2 = 1/9$ from Lemma~\ref{l.onepointupperlower}(2), and one-point lower tail $\rmreg(2)$ with $\phimac_3 = 1/3$ from   Lemma~\ref{l.onepointupperlower}(1).  \qed

\section{The lower tail of the lower curves}\label{s.othercurves}

Here we prove Proposition \ref{p.othercurves}.
Recall that $\para: \R \to \R$ denotes the parabola $\para(x) =  2^{-1/2} x^2$ that appears in the definition of a regular ensemble.

Set $\rcon = 5 (3 - 2^{3/2})^{-1}$, $r_1 = 2^{3/2}$, and $r_k = \max \{ 5^3 , \rcon r_{k-1} \big\}$ for $k \geq 2$.

\begin{proposition}\label{p.strongothercurves}
For $\bar\phimac \in (0,\infty)^3$, $C,c > 0$ and $n \in \N$, let 
$$
\mc{L}_n:\intint{n} \times \big[-\xnmac,\infty\big) \to \R  
$$ 
be a    $\big(\bar\phimac,\rsc,\rsC\big)$-regular ensemble defined under the law~$\PP$. 
For $k \in \N$, let  $\Cstrong_k = 20^{k-1} 2^{k(k-1)/2} \Cstrong_1$ and $c_k =   \big( (3 - 2^{3/2})^{3/2} 2^{-1} 5^{-3/2} \big)^{k-1} c_1$ where $\Cstrong_1 = 10C$ and $c_1 = 2^{-5/2} c \wedge 1/8$. 
Set $\delta = \phimac_1/2 \wedge \phimac_2/2 \wedge \phimac_3$. Whenever  $(n,k) \in \N^2$ satisfies $n \geq k 
  \vee  (c/3)^{-2(\phimac_1 \wedge \phimac_2)^{-1}} \vee  6^{2/\delta}$, 
  $t \in \big[ 5 \vee (3 - 2^{3/2})^{-1/2} r_{k-1}^{1/2} ,n^\delta \big]$, $r \in \big[ r_k \, , \, 2n^\delta \big]$
  and $y \in \rsc/2 \cdot [-  n^\delta,  n^\delta]$,
$$
\PP \Big( \inf_{x \in [y-t,y+t]} \big( \mc{L}_n(k,x) + \para(x) \big) \leq - r \Big) \, \leq \, t^k \cdot \Cstrong_k \exp \big\{ - c_k r^{3/2} \big\} \, .
$$
\end{proposition}

\noindent{\bf Proof of Proposition~\ref{p.othercurves}.}
Choosing $t =  5 \vee (3 - 2^{3/2})^{-1/2} r_{k-1}^{1/2}$
in Proposition~\ref{p.strongothercurves}, we see that
$$
\PP \Big(  \mc{L}_n(k,y) + \para(y)  \leq - r \Big) \, \leq \,   \Big(  5 \vee (3 - 2^{3/2})^{-1/2} r_{k-1}^{1/2} \Big)^k \cdot \Cstrong_k \exp \big\{ - c_k r^{3/2} \big\} \, .
$$
for  $r \in \big[ r_k \, , \, 2n^\delta \big]$
  and $y \in c/2 \cdot [-n^\delta,n^\delta]$. Recalling~(\ref{e.formere}),  it is easily verified that, for $k \geq 1$,
 $$
   \max \bigg\{   \Big(  5 \vee (3 - 2^{3/2})^{-1/2} r_{k-1}^{1/2} \Big)^k  \Cstrong_k \, \, , \,  \exp \big\{ c_k r_k^{3/2} \big\} \bigg\} \leq C_k  \, ;
$$
thus, we find that 
$$
\PP \big(  \mc{L}_n(k,y) + \para(y)  \leq - r \big) \, \leq \,  \formerE_k \exp \big\{ - c_k r^{3/2} \big\}
$$
for $r \in [ 0 \, , \, 2n^\delta ]$
  and $y \in c/2 \cdot [-n^\delta,n^\delta]$. \qed

\medskip

\noindent{\bf Proof of Proposition~\ref{p.strongothercurves}.}
Lemma~\ref{l.parabolicinv} promptly permits us to reduce to verifying this assertion when $y$ equals zero, except for the detail that we must prove the version of the $y=0$ statement in which appearances of the quantity $c$ are replaced by $2c$. In an abuse of notation adopted to cope with this detail, we take $c_1$ equal to $2^{-3/2} c \wedge 1/8$ henceforth in this proof.

We will prove the $y=0$ assertion by induction on $k \in \N$.
Explicitly, our inductive hypothesis states that
if  
$$
(n,k,t,r) \in \N \times \N  \times \big[ 5 \vee (3 - 2^{3/2})^{-1/2} r_{k-1}^{1/2} ,n^\delta \big] \times \big[ r_k \, , \, 2n^\delta \big] 
$$
 satisfies $n \geq k 
  \vee  (2c/3)^{-2(\phimac_1 \wedge\phimac_2)^{-1}} \vee  6^{2/\delta}$,
then
\begin{equation}\label{e.expindhyp}
\PP \Big( \inf_{x \in [-t,t]} \big( \mc{L}_n(k,x)  + \para(x) \big) \leq - r \Big) \, \leq \, t^k \cdot \Cstrong_k \exp \big\{ - c_k r^{3/2} \big\} \, .
\end{equation}

We explain first why this hypothesis is valid in the base case $k=1$. Note that $n^\delta \cdot [-1,1] \subseteq \rsc n^{\phimac_2} \cdot[-1,1]$ because $\delta \leq \phimac_2/2$ and  $n \geq c^{-2\phimac_2^{-1}}$. Thus we may apply the {\rm one-point lower tail bound}~$\rmreg(2)$ with~$k=1$ 
at points in the set $[-t,t] \cap \Z \, \cup \{-t\} \cup \{t\}$, because all such points have absolute value at most $n^{\delta}$;
then a union bound yields
\begin{equation}\label{e.infzt}
\PP \Big( \inf_{x \in [-t,t] \cap \Z \, \cup \{-t\} \cup \{t\}} \big( \mc{L}_n(1,x)  + \para(x) \big) \leq - r/2 \Big) \, \leq \, (2t + 3) \cdot C \exp \big\{ - c (r/2)^{3/2} \big\} 
\end{equation}
when  $1 \leq r/2 \leq n^{\phimac_3}$ (which upper bound is ensured by $r \leq 2n^{\delta}$). 
In a temporary notation, write~$\PP^*$ for the law~$\PP$ conditioned on the occurrence of the complement of the event in the last display. Our first claim regarding~$\PP^*$ is that, when $r \geq 2^{3/2}$,
\begin{equation}\label{e.pstarclaim}
 \PP^* \Big(  \inf_{x \in [0,1]} \big( \mc{L}_n(1,x)  + \para(x) \big) \leq - r \Big) \, \leq \,   \exp \big\{ -   r^2/8 \big\} \, .
\end{equation}
To verify this, note that,
by the monotonicity Lemmas~\ref{l.monotoneone} and~\ref{l.monotonetwo},
the conditional distribution of $\mc{L}_n(1,\cdot):[0,1] \to \R$  under $\PP^*$
stochastically dominates $\mc{B}_{1;-r/2,-\para(1) - r/2}^{[0,1]}$. Since $\para(1) = 2^{-1/2}$, Lemma~\ref{l.maxfluc}  implies that, for $s > 0$,
$$
\mc{B}_{1;-r/2,-\para(1) - r/2}^{[0,1]} \Big( \inf_{x \in [0,1]} B(x) \leq - s - r/2 - 2^{-1/2} \Big) \leq \exp \big\{ - 2s^2 \big\} \, .
$$
Taking $s = r/4$ and using $r \geq 2^{3/2}$,
we learn that
$$
\PP^* \Big( \, \inf_{x \in [0,1]}  \mc{L}_n(1,x)  \leq - r  \,\Big) \leq \exp \big\{ -  r^2/8 \big\} \, . 
$$
Since $\para$ is non-negative, we obtain~(\ref{e.pstarclaim}).

The argument just given shows that the claim~(\ref{e.pstarclaim}) is equally true if in the definition of~$\PP^*$ we further condition on more information concerning the curve $\mc{L}_n(1,\cdot)$ outside $[0,1]$. Since  $\delta$ equals $\phimac_1/2  \wedge  \phimac_2/2  \wedge  \phimac_3$, this quantity is bounded above by the quantity $\deltapi = \phimac_1 \wedge \phimac_2$ appearing in the parabolic invariance Lemma~\ref{l.parabolicinv}, and thus the claim is also valid if we replace the interval $[0,1]$ by any other of the form $[m,m+1]$ with $[m,m+1] \subset [-n^\delta,n^\delta]$. Finally, the intervals $[ \lfloor t \rfloor, t ]$ and  $[ -t, \lceil -t \rceil  ]$ may also be accomodated, because the bound obtained using Lemma~\ref{l.maxfluc} is valid when an interval of less than unit length is used instead.
Our conclusion then is that
$$
 \PP^* \Big( \inf_{x \in [-t,t]} \big( \mc{L}_n(1,x)  + \para(x) \big) \leq - r \Big) \, \leq \, (2t + 3)  \exp \big\{ -   r^2/8 \big\} 
$$
when $r \geq 2^{3/2}$. By the definition of $\PP^*$, we may replace $\PP^*$ by $\PP$ here at the expense of intersecting the event in the left-hand side with the complement of the event in~(\ref{e.infzt}).
Combining then with~(\ref{e.infzt}), we find that  
$$
\PP \Big( \inf_{x \in [-t,t]} \big( \mc{L}_n(1,x)  + \para(x) \big) \leq - r \Big) \, \leq \, (2t + 3) \cdot C \exp \big\{ - c (r/2)^{3/2} \big\} \, + \, (2t + 3)  \exp \big\{ -   r^2/8 \big\}  \, .
$$
when $r \in [2^{3/2},2n^\delta]$.
Since $t\geq 1$ and $r \geq 1$, the quantity $10 t C \exp \big\{ - 1/8 \wedge c  2^{-3/2} \cdot r^{3/2} \big\}$ is an upper bound on the right-hand side. Thus, we confirm the inductive hypothesis at $k=1$ with the choice $\Cstrong_1 = 10 C$ and $c_1 = 2^{-3/2}c \wedge 1/8$.

We now verify the inductive hypothesis at general index $k > 1$, assuming the validity of this assertion for index~$k-1$. Our proof follows the same approach as that by which Lemmas~$7.2$ and~$7.3$ are proved in~\cite{KPZLE}. 

We begin with a fact about the basic parabola $-\para$.  Consider any real interval of length $t$, and the affine function whose values at the interval's endpoints coincide with the parabola's. Then the difference between the value of the parabola and the value of the affine function at the interval's midpoint will be independent of the interval in question and equal to $L t^2$, where here we introduce the quantity $L = 2^{-5/2}$. It is moreover the case that the maximal value of this difference evaluated at any point in the interval is achieved at the interval's midpoint.

In the first instance, we will verify the inductive hypothesis under the additional assumption that the parameter $t$ is at most  $3 n^{\delta/2}$. Suppose then that $t$ satisfies
\begin{equation}\label{e.trange}
t \in \big[ 5 \vee (3 - 2^{3/2})^{-1/2} r_{k-1}^{1/2} \, , \, 3 n^{\delta/2} \big] \, .
\end{equation} 
Define the event 
$$
\low_k^{[t,2t]} = \Big\{ \sup_{x \in [t,2t]} \big( \mc{L}_n(k,x) + \para(x) \big) \leq -  \tfrac{1}{2} L t^2 - t^{1/2} \Big\} 
$$
and its counterpart 
$\low_k^{[-2t,-t]}$ associated to the interval $[-2t,-t]$. 

\begin{figure}[ht]
\begin{center}
\includegraphics[height=12cm]{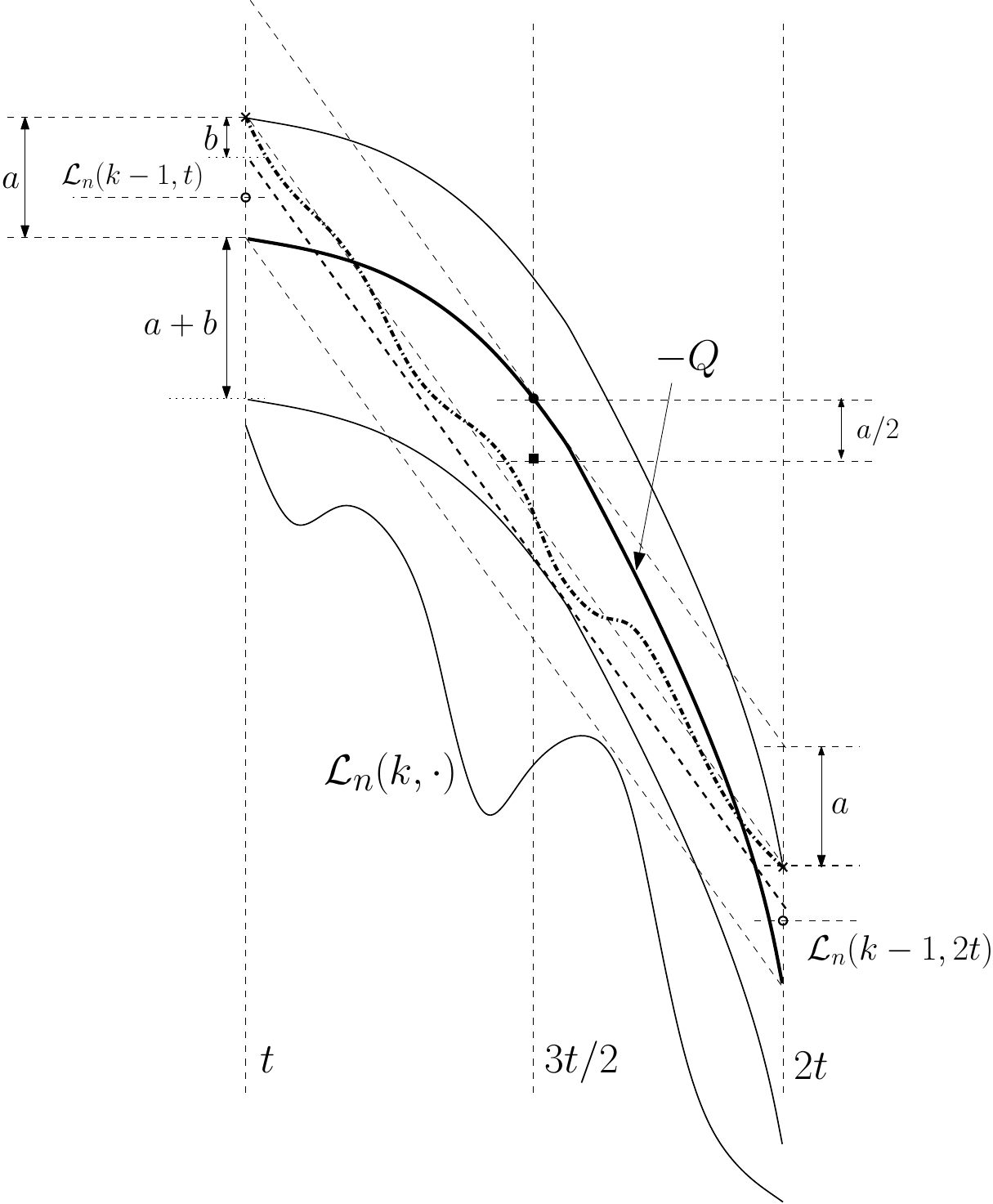}
\caption{Illustrating the proof of~(\ref{e.twolows}).
In the labels in the sketch, we write $a = \tfrac{1}{2} L t^2$ and $b = t^{1/2}$. The central concave curve is $-\para$; the other two are obtained by vertical translation by $a$ and $-a-b$.
The highest of the four parallel dashed sloping line segments has a point of tangency with $-\para$ at $3t/2$. When $\low_k^{[t,2t]}$ occurs, the curve $\mc{L}_n(k,\cdot):[t,2t] \to \R$ lies below the lowest of the three parabolas. In this location, it can offer no support to the $(k-1)\textsuperscript{st}$ curve. Indeed, in the typical circumstance that 
 $\mathsf{G}_{k-1}(s)$ occurs for $s \in \{ t, 2t  \}$, the curve $\mc{L}_n(k-1,\cdot)$ at $t$ and $2t$ lies below the endpoints of the highest parabola. This curve on $[t,2t]$ only rises stochastically if instead it begins and ends at this parabola's endpoints and the lower indexed curves disappear.  Brownian bridge with such endpoints easily bypasses $\mc{L}_n(k,\cdot)$ because it may avoid this curve by remaining above the third highest of the sloping dashed line segments, as the dashed-dotted curve illustrates. Thus, $\mc{L}_n(k-1,\cdot)$
 tends to adopt a linear, rather than a parabolic, route during $[t,2t]$, giving it a uniformly positive probability to lie below the point marked by a square at the midpoint time~$3t/2$. However, this point is very low judged in parabolically curved coordinates, so this eventuality is known to be rare by the one-point lower tail for the $(k-1)\textsuperscript{st}$ curve. Thus, $\low_k^{[t,2t]}$ is unlikely.}
\label{f.knosupport}
\end{center}
\end{figure}

We next establish 
a key part of the proof of the inductive step, namely the bound
\begin{equation}\label{e.twolows}
  \PP \big( \low_k^{[t,2t]} \big) \vee  \PP \big( \low_k^{[-2t,-t]} \big) \leq 
4C \exp \big\{  - c \cdot 2^{-27/4} t^3  \big\} \, .
\end{equation}
We will prove this upper bound for $\PP \big( \low_k^{[t,2t]} \big)$; the other proof is identical. The argument is illustrated by Figure~\ref{f.knosupport}.

For $s \geq - \xnmac$, define 
$$
 \mathsf{G}_{k-1}(s) = \Big\{ \mc{L}_n(k-1,s) + \para(s) \leq \tfrac{1}{2} L t^2 \Big\} \, .
$$

We set
$$
 \mathsf{A}_k^{[t,2t]} =   \mathsf{G}_{k-1}(t) \cap  \mathsf{G}_{k-1}(2t) \cap 
\low_k^{[t,2t]} \, .
$$

The monotonicity Lemmas~\ref{l.monotoneone} and~\ref{l.monotonetwo} imply that, 
under $\PP$ given $\mathsf{A}_k^{[t,2t]}$,
 the conditional distribution of $\mc{L}_n \big( k-1, \cdot \big):[t,2t] \to \R$ is stochastically dominated by~$\mc{B}_{1;u,v}^{[t,2t]}\big( \cdot \big\vert \notouch_f \big)$, where 
 $$
 u = -\para(t) + \tfrac{1}{2} L t^2 \, \, \, , \, \, \, v = - \para(2t) + \tfrac{1}{2} L t^2
 $$ 
 and $f:[t,2t] \to \R$ equals $f(x) = - \para(x) -  \tfrac{1}{2} L t^2 - t^{1/2}$.

Consider now the affine function $\ell$ whose graph contains the points $(t,u)$ and $(2t,v)$. Suppose the graph of $\ell$ is translated downwards until it makes contact with the graph of $f$. Our parabolic fact implies that first contact will occur at $x$-coordinate $3t/2$, when the distance of translation equals~$t^{1/2}$. In this light, we find that Lemma~\ref{l.maxfluc}
has the consequence that 
$$
\mc{B}_{1;u,v}^{[t,2t]} \big( \notouch_f \big) \geq 1 - e^{-2} \, .
$$
The location $-\para(3t/2) - \tfrac{1}{4} L t^2$ is the midpoint of the interval $\big[ \ell(3t/2) , - \para(3t/2) \big]$. Thus,
\begin{eqnarray*}
  &  &   \mc{B}_{1;u,v}^{[t,2t]} \Big( B(3t/2) \geq - \para(3t/2) - \tfrac{1}{4} L t^2 \Big) \, = \,  \mc{B}_{1;u,v}^{[t,2t]} \Big( B(3t/2) \geq  \ell (3t/2) + \tfrac{1}{4} L t^2 \Big) \\
  & = &  \nu_{0,t/4} \big(   \tfrac{1}{4} L t^2  , \infty \big) 
=    \nu_{0,1} \big(   \tfrac{1}{2} L t^{3/2}  , \infty \big) \\
  &  \leq & (2\pi)^{-1/2} \big( 2L^{-1} t^{-3/2} \big)^{1/2}  \exp \big\{ - 8^{-1} L^2 t^3 \big\} \leq  \exp \big\{ - 2^{-8} t^3 \big\}  \, ,
\end{eqnarray*}
where the last inequality depended on the certainly valid $t \geq 2 \pi^{-2/3}$.

The quantity $L/4$ equals $2^{-9/2}$.
We find then that, for such values of $t$,
\begin{eqnarray*}
  &  &  \PP \Big( \mc{L}_n\big( k-1,3t/2 \big) \geq - \para(3t/2) - 2^{-9/2} t^2 \, \Big\vert \, \mathsf{A}_k^{[t,2t]} \Big) \\
   & \leq &   \mc{B}_{1;u,v}^{[t,2t]}   \Big( B \big( 3t/2 \big) \geq - \para(3t/2) - 2^{-9/2} t^2 \, \Big\vert \, \notouch_f \Big) \\
   & \leq & \big( 1 - e^{-2} \big)^{-1}   \mc{B}_{1;u,v}^{[t,2t]}   \Big( B \big( 3t/2 \big) \geq - \para(3t/2) - 2^{-9/2} t^2  \Big) \leq \big( 1 - e^{-2} \big)^{-1}  \exp \big\{ - 2^{-8} t^3 \big\} \, . 
\end{eqnarray*} 
If $a \in (0,1)$ and two events $E$ and $F$ satisfy $\PP(E\vert F) \leq a$, then $\PP(F) \leq (1-a)^{-1} \PP(E^c)$. Expressing the last derived inequality in these terms, we infer that
$$
 \PP \big( \mathsf{A}_k^{[t,2t]} \big) \leq 2 \, \PP \Big( \mc{L}_n\big( k-1,3t/2 \big) \geq - \para(3t/2) - 2^{-9/2} t^2 \Big) \, ,
$$
since $a =  \big( 1 - e^{-2} \big)^{-1}  \exp \big\{ - 2^{-8} t^3 \big\}$ is at most one-half provided that $t^3 \geq 2^{8} \log \big( 2e^2 (e^2 -1)^{-1} \big)$ (which is valid since $t \geq 5$).
Using the one-point lower tail bound $\rmreg(2)$, we learn that, since $\vert 3t/2 \vert \leq c n^{\phimac_2}$ and $2^{-9/2} t^2 \in [1,n^{\phimac_3}]$,
$$
 \PP \big( \mathsf{A}_k^{[t,2t]} \big) \leq 2 C  \exp \big\{ - c \cdot 2^{-27/4} t^3 \big\} \, .
$$
(These bounds on $t$ are ensured by our choice of this parameter, with $t \leq 2  c/3 \cdot n^{\phimac_2}$ following from $t \leq n^\delta$, $\delta \leq \phimac_2/2$ and $n \geq (2c/3)^{-2\phimac_2^{-1}}$.)
Hence,
\begin{eqnarray*}
  \PP \big( \low_k^{[t,2t]} \big) & \leq &  \PP \big( \mathsf{A}_k^{[t,2t]} \big) + \PP \big( \neg \, \mathsf{G}_1(t) \big) + \PP \big( \neg \, \mathsf{G}_1(2t) \big) \\
   & \leq & 2C \exp \big\{  - c \cdot 2^{-27/4} t^3  \big\} \, + \,  2C  \exp \big\{  - c \cdot 2^{-21/4} t^3  \big\} \leq 
4C \exp \big\{  - c \cdot 2^{-27/4} t^3  \big\} \, .
\end{eqnarray*}
This completes the derivation of~(\ref{e.twolows}). 

We now present a further argument that leads from the bound~(\ref{e.twolows}) to the end of the inductive step. To do so, we now introduce the event
$$
\up_{k-1}^{[-2t,2t]}(s) = \Big\{ \inf_{x \in [-2t,2t]} \mc{L}_n(k-1,x) \geq - \para(2t) - s \Big\}
$$
where $s \geq 0$. Since $\para(2t) = \sup_{x \in [-2t,2t]} \para(x)$, we may 
apply the inductive hypothesis~(\ref{e.expindhyp}) at index~$k-1$ to bound above the failure probability of this event. 
The parameter $t$ appearing in the hypothesis will not be its presently assigned value, which is associated to this stage~$k$, but rather twice that value. We also take $r$ equal to $s$. These choices are permissible
if we insist that $s \in \big[2^{3/2} \vee r_{k-1} \, , \, 2n^\delta \big]$,
  because the other requirement, that $t \in 1/2 \cdot \big[ 5 \vee (3 - 2^{3/2})^{-1/2} r_{k-1}^{1/2} ,n^\delta \big]$,  is due to assumption and $n \geq 6^{2/\delta}$. For such $s$, the inductive hypothesis tells us that
\begin{equation}\label{e.noup}
\PP \Big(  \neg \, 
\up_{k-1}^{[-2t,2t]}(s) \Big) \leq (2 t)^{k-1} \cdot \Cstrong_{k-1} \exp \big\{ - c_{k-1} s^{3/2}  \big\} \, .
\end{equation}

Our plan for completing the inductive step is to consider the event
$$
 \mathsf{E} =   \nolow_k^{[-2t,-t]} \cap \nolow_k^{[t,2t]} \cap
\up_{k-1}^{[-2t,2t]}(s) \, ,
$$
where for now the parameter $s \geq 0$ remains unspecified. From what we have already learnt, and as we will record shortly in~(\ref{e.probebound}), we know that $\mathsf{E}^c$ is unlikely if $s$ is high. We now seek to argue that under $\PP \big( \cdot \big\vert E \big)$, $\mc{L}_n(k,\cdot)$ is unlikely to drop low {\em anywhere} in~$[-t,t]$.

\begin{figure}[ht]
\begin{center}
\includegraphics[height=12cm]{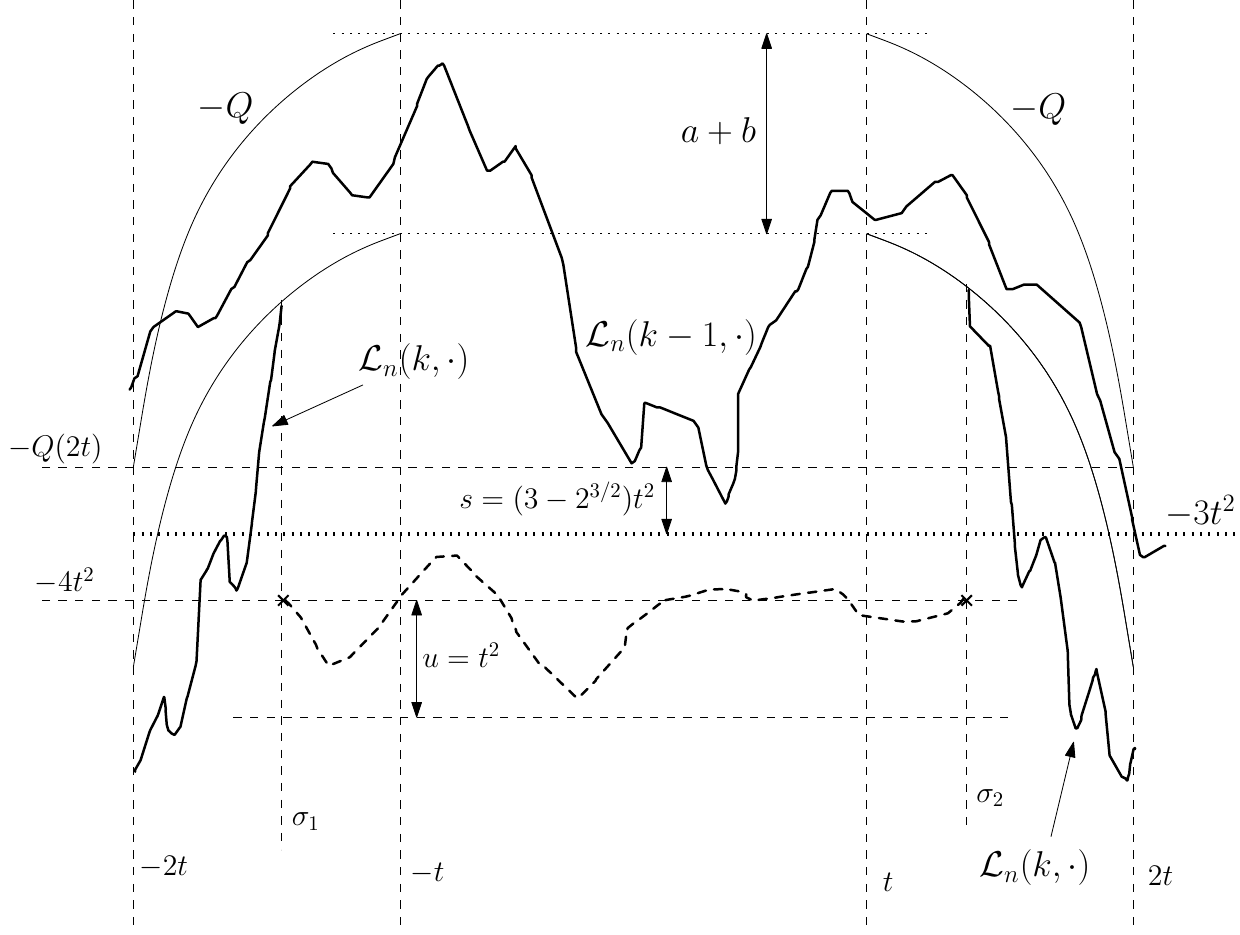}
\caption{The study of the conditional law $\PP \big( \cdot \big\vert E \big)$  leads to the derivation of~(\ref{e.econseq}).
As in Figure~\ref{f.knosupport}, we write $a = \tfrac{1}{2} L t^2$ and $b = t^{1/2}$. 
The sketch depicts an occurrence of~$\mathsf{E}$. The conditional distribution of $\mc{L}_n(k,\cdot)$ on $[\sigma_1,\sigma_2]$ is stochastically decreased by using endpoint values $-4t^2$, marked with crosses, and by removing curves of higher index.  The Brownian bridge dashed curve on $[\sigma_1,\sigma_2]$
is unlikely to drop by a further $u = t^2$ units; thus, nor is $\mathcal{L}_n(k,\cdot)$, because the curve $\mc{L}_n(k-1,\cdot)$ is at a comfortable distance upwards, staying above the coordinate $-\para(2t)- s = -3t^2$ that is indicated with a bold dotted line.}
\label{f.uptwosides}
\end{center}
\end{figure}

To argue that this is so, define random variables $\sigma_1 \in [-2t,-t]$ and  $\sigma_2 \in [t,2t]$ so that
$$
\sigma_1 = \inf \Big\{ x \in [-2t,-t]: \mc{L}_n(k,x) + \para(x) \geq - 2^{-7/2} t^2 - t^{1/2 }\Big\}
$$
and
$$
\sigma_2 = \sup \Big\{ x \in [t,2t]: \mc{L}_n(k,x) + \para(x) \geq - 2^{-7/2} t^2 - t^{1/2 }\Big\} \, ;
$$
it is immaterial how these random variables are specified when the infimums are taken over the empty-set, because this circumstance does not arise when $\nolow_k^{[-2t,-t]}\cap \nolow_k^{[t,2t]}$ takes place, and we will now use the random variables only in a situation where this event occurs. Note that $\intint{k} \times (\sigma_1,\sigma_2)$ is a stopping domain.

Consider the law $\PP \big( \cdot \big\vert E \big)$ (and consult Figure~\ref{f.uptwosides} for a visual explanation of the argument concerning this law). By the strong Gibbs and monotonicity   Lemmas~\ref{l.stronggibbslemma},~\ref{l.monotoneone} and~\ref{l.monotonetwo}, the resulting conditional distribution of $\mc{L}_n \big( k,\cdot \big): [\sigma_1,\sigma_2] \to \R$ stochastically dominates 
$$
 \mc{B}_{1;v_1,v_2}^{[\sigma_1,\sigma_2]} \big( \cdot \big\vert \notouch_g \big) \, ,
$$
where $v_1 = - \para(\sigma_1) - 2^{-7/2} t^2 - t^{1/2}$,
$v_2 = - \para(\sigma_2) - 2^{-7/2} t^2 - t^{1/2}$ and the function $g$ is identically equal to $-\para(2t) - s$.
Note that, since $\vert \sigma_1 \vert \vee \vert \sigma_2 \vert \leq 2t$ and $t \geq 1$, 
$$
 v_1 \wedge v_2 \geq - \big( 2^{3/2} + 2^{-3/2} \big) t^2 \geq -4t^2 \, .
$$
We now specify the parameter $s$ to be equal to $(3-2^{3/2})t^2$, so that $g$ equals the constant function~$-3t^2$.

 (In order that $s$ lie in the interval $[2^{3/2} \vee r_{k-1},2n^\delta]$, 
 as we have demanded, it is sufficient that 
 $t \in \big[5 \vee  (3 - 2^{3/2})^{-1/2} r_{k-1}^{1/2}  , 3 n^{\delta/2}\big]$,
  a requirement that we have imposed in~(\ref{e.trange}).)
 
  This choice of~$s$ yields
\begin{equation}\label{e.probebound}
\PP \big( \mathsf{E}^c \big) \leq 
2 \cdot 4C \exp \big\{  - c \cdot 2^{-27/4} t^3  \big\} \, + \,  (2 t)^{k-1}  \cdot \Cstrong_{k-1} \exp \big\{ - c_{k-1} (3 - 2^{3/2})^{3/2} t^3  \big\}
\end{equation}
by means of~(\ref{e.twolows}) and (\ref{e.noup}).

We have seen that, under $\PP \big( \cdot \big\vert E \big)$, the conditional distribution of the $k\textsuperscript{th}$ $\mc{
L}_n$-curve on~$[\sigma_1,\sigma_2]$ stochastically dominates
$\mc{B}_{1;-4t^2,-4t^2}^{[\sigma_1,\sigma_2]} \big( \cdot \big\vert \notouch_{-3t^2} \big)$.

Note that, since $\sigma_2 - \sigma_1 \leq 4t$,
$$
\mc{B}_{1;-4t^2,-4t^2}^{[\sigma_1,\sigma_2]} \big( \notouch_{-3t^2} \big) \geq 
\mc{B}_{1;0,0}^{[0,4t]} \big( \notouch_{t^2} \big) = 1 - 
\exp \big\{ - \tfrac{1}{2} t^3 \big\} \geq 1/2 \, ,
$$
where the equality is due to Lemma~\ref{l.maxfluc}
and the latter inequality to $t \geq \big(2 \log 2 \big)^{1/3}$.

By  Lemmas~\ref{l.stronggibbslemma} and~\ref{l.monotonetwo}, it follows that, for $u \geq 0$,
\begin{eqnarray*}
 & & \PP \Big( \inf_{x \in [\sigma_1,\sigma_2]} \mc{L}_n \big( k, x \big) \leq - 4t^2 - u  \, \Big\vert \, E \Big) \\
 & \leq & 
\mc{B}_{1;-4t^2,-4t^2}^{[\sigma_1,\sigma_2]} \Big(  \inf_{x \in [\sigma_1,\sigma_2]} B (  x ) \leq - 4t^2 - u \, \Big\vert \, \notouch_{-3t^2} \Big) \\
 & \leq & 2 \, 
\mc{B}_{1;-4t^2,-4t^2}^{[\sigma_1,\sigma_2]} \Big(  \inf_{x \in [\sigma_1,\sigma_2]} B (  x ) \leq - 4t^2 - u \, \Big) \leq 2 \exp \big\{ - \tfrac{1}{2} u^2 t^{-1} \big\} \\
\end{eqnarray*}
where the final inequality depended on $\sigma_2 - \sigma_1 \leq 4t$ and Lemma~\ref{l.maxfluc}.
Since $[-t,t] \subseteq [\sigma_1,\sigma_2]$,
we find that
$$
 \PP \Big( \inf_{x \in [-t,t]} \mc{L}_n \big( k, x \big) \leq - 4t^2 - u   \Big) \leq  2 \exp \big\{ - \tfrac{1}{2} u^2 t^{-1} \big\} \PP \big( E \big) \, + \, \PP \big( E^c \big) \, .  
$$
We now set $u = t^2$.
Using~(\ref{e.probebound}), we see that 
\begin{eqnarray}
 & & 
 \PP \Big( \inf_{x \in [-t,t]} \mc{L}_n \big( k, x \big) \leq - 5 t^2     \Big) \label{e.econseq}  \\
& \leq & 2 \exp \big\{ - \tfrac{1}{2} t^3 \big\} \, + \,
8 C \exp \big\{  - c \cdot 2^{-27/4} t^3  \big\} \, + \,  (2 t)^{k-1}  \cdot \Cstrong_{k-1} \exp \big\{ - c_{k-1} (3 - 2^{3/2})^{3/2} t^3  \big\} \nonumber \, .
\end{eqnarray}
Now setting $h = 5t^2$ and using $[-1,1] \subset [-t,t]$, 
we find that
$$
 \PP \Big( \inf_{x \in [-1,1]} \big( \mc{L}_n \big( k, x \big) + Q(x) \big) \leq - h     \Big) \leq   5 (2 t)^{k-1}  \cdot \Cstrong_{k-1} \exp \big\{ - c_{k-1} (3 - 2^{3/2})^{3/2} 5^{-3/2} \cdot h^{3/2}  \big\}  \, ,
$$
where $Q \geq 0$ has permitted the introduction of the $+ Q(x)$ term.
Note that the choice $h = 5t^2$  allows  $h$  to take any value in 
$\big[ 5^3 \vee 5(3 - 2^{3/2})^{-1} r_{k-1} , 45n^\delta\big]$
 by specifying $t$ suitably within its permitted range~(\ref{e.trange}). Note that the lower bound on $h$ here equals $r_k$.

By Lemma~\ref{l.parabolicinv}, the replacement of $c_{k-1}$ by $c_{k-1}/2$ in the exponent on the right-hand side of the preceding display produces an inequality which is equally valid when the interval~$[-1,1]$ over which the infimum is taken is replaced by any interval of length two that lies in $[-n^\delta,n^{\delta}]$ (since $n^\delta \leq \rsc/2 \cdot n^{\phimac_1 \wedge \phimac_2}$). By taking a union bound of the resulting inequalities, we learn that, if $t \leq n^\delta$ and 
$h \in 
\big[ r_k , 45n^\delta\big]$,
$$
 \PP \Big( \inf_{x \in [-t,t]} \big( \mc{L}_n ( k, x ) + \para(x) \big) \leq - h     \Big) \leq   (2t+1) \cdot 5 (2 t)^{k-1}  \cdot \Cstrong_{k-1} \exp \big\{ - c_{k-1} 2^{-1} (3 - 2^{3/2})^{3/2} 5^{-3/2} \cdot h^{3/2}  \big\}  \, .
$$
Setting $h=r$, and noting that $t \geq 1$, we verify the inductive hypothesis~(\ref{e.expindhyp}) at index~$k$ if we note that $\Cstrong_k = 20 \cdot 2^{k-1} \Cstrong_{k-1}$ and $c_k =  c_{k-1} (3 - 2^{3/2})^{3/2} 2^{-1} 5^{-3/2}$. 
The inductive step completed, we have obtained Proposition~\ref{p.strongothercurves}. \qed

\section{Regular ensemble curves collapse near infinity
}\label{s.collapsenearinfinity}

\noindent{\bf Proof of Proposition~\ref{p.collapsenearinfinity}.}
Recall that $Q(x) = - 2^{-1/2} x^2$ denotes the parabola that appears in the definiton of a regular ensemble.
Set $h = 2^{-5/2} \eta^2 n^{2\phimac_2}$.
Note that $\ell(x) = Q(x) - h$ when $x = \eta/2 \cdot n^{\phimac_2}$
and $\ell(x) = Q(x) + h$  when $x = \eta n^{\phimac_2}$.
Recall that the restriction of $\ell$ to $[0,\infty)$ is affine.

For $y \in [ \eta n^{\phimac_2} , \infty ]$, write $\mathsf{RightHigh}(y)$
for the event that  
$$
\mc{L}_n\big(1,z\big) \geq \ell(z) \, \, \textrm{for some} \, \,  z \in  [ \eta n^{\phimac_2} , y)  \, ,
$$
and $\mathsf{LeftHigh}(-z_n)$
for the event that  
$$
\mc{L}_n\big(1,z\big) \geq \ell(z) \, \, \textrm{for some} \, \,  z \in  [ -z_n, - \eta n^{\phimac_2})  \, .
$$

The proposition is implied by
\begin{equation}\label{e.righthighbound} 
\PP \big( \mathsf{RightHigh}(\infty) \big) \vee 
\PP \big( \mathsf{LeftHigh}(-z_n) \big)   \leq 
3C \exp \Big\{ - c \eta^3  2^{-15/4}   n^{3\phimac_2 \wedge 3\phimac_3/2} \Big\} \, .
\end{equation}
We will prove this bound for the first probability. 
The argument for the second is the same.
In regard to the first,
it is enough to argue that the bound holds on
$\PP \big( \mathsf{RightHigh}(y) \big)$ for any  $y \in [ \eta n^{\phimac_2} , \infty )$, and this we now do.

Thus we let  $y \in [ \eta n^{\phimac_2} , \infty )$ be given.
We set $\chi = \sup \big\{ x \in [ \eta n^{\phimac_2}, y ] :  \mc{L}_n\big(1,z\big) \geq \ell(z) \big\}$, with $\chi = -\infty$ when the infimum is over the empty-set. 
The event that $\chi > -\infty$ equals $\mathsf{RightHigh}(y)$, and thus 
we aim to bound above $\PP \big( \chi > - \infty \big)$. The random set $\{1 \} \times \big[ \eta/2 \cdot n^{\phimac_2}  ,  \chi \vee   \eta/2 \cdot n^{\phimac_2} \big]$ is a stopping domain. Applying Lemma~\ref{l.stronggibbslemma} and then dispensing with the lower boundary condition $\mc{L}_n(2,\cdot): \big[
 \eta/2 \cdot n^{\phimac_2}  ,  \chi \vee   \eta/2 \cdot n^{\phimac_2}\big] \to \R$
via the strong Gibbs Lemma~\ref{l.stronggibbslemma}, we find that
$$
\PP \Big( \, \mc{L}_n \big( 1, \eta n^{\phimac_2} \big) \geq  \ell \big(  \eta n^{\phimac_2} \big) \, \Big\vert \, \mathsf{RightHigh}(y) \, , \,   \mc{L}_n \big( 1, \eta/2 \cdot  n^{\phimac_2} \big) \geq  \ell \big(  \eta/2 \cdot n^{\phimac_2} \big)  \Big) \, \geq \, 1/2 \, .
$$
Thus,
\begin{equation}\label{e.lefthigh}
 \PP \big( \mathsf{RightHigh}(y)  \big) \leq 2 \, \PP \Big( \mc{L}_n \big( 1, \eta n^{\phimac_2} \big) \geq   \ell \big(  \eta n^{\phimac_2} \big) \Big) +  \PP \Big(  \mc{L}_n\big(1,  \eta/2 \cdot n^{\phimac_2}   \big) <  \ell \big(  \eta/2 \cdot n^{\phimac_2} \big) \Big) \, .
\end{equation}
Recall that $\ell \big(  \eta/2 \cdot n^{\phimac_2} \big) = Q \big(  \eta/2 \cdot n^{\phimac_2} \big) - h$
and $\ell \big(  \eta  n^{\phimac_2} \big) = Q \big(  \eta n^{\phimac_2} \big) + h$ with  $h = 2^{-5/2} \eta^2 n^{2\phimac_2}$.

Applying one-point lower tail ${\rm Reg}(2)$ with ${\bf z} =  \eta/2 \cdot n^{\phimac_2}$ and ${\bf s}
= h \wedge n^{\phimac_3}$,
we find that 
$$
 \PP \Big(  \mc{L}_n\big(1,  \eta/2 \cdot n^{\phimac_2}   \big) <  \ell \big(  \eta/2 \cdot n^{\phimac_2} \big)
  \Big) \leq C \exp \Big\{ - c \eta^3  2^{-15/4}   n^{3\phimac_2 \wedge 3\phimac_3/2} \Big\} \, ,
$$
where we used $\eta \leq \rsc \leq 2^5$ in the form ${\bf s} \geq 2^{-5/2} \eta^2 n^{2\phimac_2 \wedge \phimac_3}$.
(We are using a boldface notation here to indicate the parameters in the ${\rm Reg}$ condition.)
Note that
$$
  \PP \Big( \mc{L}_n \big( 1, \eta n^{\phimac_2} \big) \geq   \ell \big(  \eta n^{\phimac_2} \big) \Big) 
  =  \PP \Big(  \mc{L}_n \big( 1, \eta n^{\phimac_2} \big) \geq  Q \big(  \eta n^{\phimac_2} \big) + h \Big)
$$
whose right-hand side is seen by the one-point upper tail ${\rm Reg}(3)$ with ${\bf s} = h$ to be at most
$$
C \exp \Big\{ - \rsc \eta^3 2^{-15/4}  n^{3\phimac_2} \Big\} \, ,
$$
where we used Proposition~\ref{p.collapsenearinfinity}'s hypothesis $n \geq \big( 2^{5/4} \rsc^{-1} \big)^{\phimac_2^{-1}}$, and $\eta \leq \rsc$, to ensure that $h \geq 1$.
Returning to~(\ref{e.lefthigh}), we find then that
$$
 \PP \big( \mathsf{RightHigh}(y)  \big) \leq 
3C \exp \Big\{ - c \eta^3  2^{-15/4}   n^{3\phimac_2 \wedge 3\phimac_3/2} \Big\} \, .
$$ 
As we have noted, this suffices to prove~(\ref{e.righthighbound}). 
This completes the proof of Proposition~\ref{p.collapsenearinfinity}. \qed

\bibliographystyle{plain}

\bibliography{airy}

\end{document}